\documentclass{article}
\usepackage{amsthm}
\usepackage{amsmath}
\usepackage{amssymb}
\usepackage{epstopdf}
\usepackage{graphicx}
\usepackage{hyperref}
\usepackage{float}
\usepackage{pgfplots}
\newtheorem{theorem}{Theorem}

\newcommand{\sumend}{}

\begin{document}
\title%Optional Short Title
{New results for Euler sums}
\author{Ross C. McPhedran,\\
School of Physics, University of Sydney,\\
and David H. Bailey,\\
Lawrence Berkeley National Laboratory}
\maketitle

\begin{abstract}We present a large number of analytic evaluations of Euler sums, namely sums such as
\begin{align}
M(m,n_0,n_1,n_2, \ldots, n_t) &= \sum_{k=1}^\infty \frac{H(k)^m}{k^{n_0} (k+1)^{n_1} (k+2)^{n_2} \cdots (k+t)^{n_t}}, \nonumber
\end{align}
for nonnegative integers $m$ and $(n_i)$, with $m \geq 1$ and $n_0 + n_1 + \cdots + n_t \geq 2$, where $H(k) = \sum_{j=1}^k 1/j$ is the harmonic function. These results were obtained either by algebraic manipulations, or else by very high-precision numerical evaluations combined with an integer relation algorithm to obtain the analytic formulas. We show how many of these results can be derived from a few basic facts, and that these techniques are applicable to Euler sums of even more general forms than the above cases. We then show that these results permit the calculation of constants for Euler sums resembling the Stieltjes $\gamma$ constants arising in the theory of the Riemann zeta function, and we also present some preliminary results on the asymptotic behavior of these constants.\footnote{The authors dedicate this paper to the memory of Jonathan and Peter Borwein, two giants of mathematical research who recently passed away. Jonathan in particular investigated Euler sums in some earlier studies that we reference.}
\end{abstract}

\section{Introduction}

The investigations reported here had their origins in work \cite{sigsam} on the Keiper-Li criterion for the Riemann hypothesis \cite{keiper,li}. The Keiper-Li criterion involves positive valued coefficients $a_n$ arising in expansions of the Riemann zeta function. The new representation of the $a_n$ reported in \cite{sigsam} involved a combination of two sets of coefficients $C_{n,p}$ and $\Sigma^{\xi}_p$, again positive valued. This representation enabled the accurate calculation of the first 4000 coefficients $a_n$. The coefficients $C_{n,p}$ obeyed a recurrence relation, and had representations involving the classical Euler sums. A deeper understanding of the asymptotic behaviour of the $C_{n,p}$  as the two integer parameters $n$ and $p$ tended to infinity was sought, and naturally involved results from the extensive literature on Euler sums \cite{ChoiSri}-\cite{Zheng}.

In this paper we address Euler sums of the form
\begin{align}
M(m,n_0,n_1,n_2, \ldots, n_t) &= \sum_{k=1}^\infty \frac{H(k)^m}{k^{n_0} (k+1)^{n_1} (k+2)^{n_2} \cdots (k+t)^{n_t}},
\end{align}
for nonnegative integers $m$ and $(n_i)$, with $m \geq 1$ and $n_0 + n_1 + \cdots + n_t \geq 2$, where $H(k) = H_k = \sum_{j=1}^k 1/j$ is the harmonic function (we use both notations interchangeably below). However, the techniques presented below are applicable to Euler sums of even more general forms. We focus on Euler sums having a common order $r$, where $r = m + n_0 + n_1 + \cdots + n_t$. We combine results from the literature with many new results, in an effort to say as much as possible about systems of order $r$ ranging from 3 to 12. The complexity of these analyses increases rapidly with $r$.

Among the most striking results of this paper are the linkages between Euler sums and Stieltjes constants. The 
latter can be  defined by:
\begin{equation}
\lim_{n\rightarrow \infty}\left[\sum_{k=1}^n \frac{(\log k)^p}{k}-\frac{(\log n)^{p+1}}{p+1}\right]=\gamma_p.
\label{intro1}
\end{equation}
The Stieltjes constants $\gamma_p$ have a sign which varies in a complicated way as $p$ increases and their modulus increases. We define their equivalent for Euler sums as
\begin{equation}
\lim_{n\rightarrow \infty}\left[\sum_{k=1}^n \frac{H_k^{p-1}}{k}-\frac{1}{p} H_n^p\right]=\gamma_p^H.
\label{intro2}
\end{equation}
Here the harmonic Stieltjes constants $\gamma_p^H$ are all positive, and again increase rapidly as $p$ increases. The harmonic Stieltjes constants are shown to be given by a sum over a set of primitive sums which form a basis for the Euler sums of order $p$.

\section{Previous results on Euler Sums}\label{sec:Previous-Results}

We now present selected results from the literature on Euler sums, including  relatively recent results due to the late Jonathan Borwein and collaborators \cite{BBG1,BBG2}, which were obtained using the techniques of experimental mathematics to complement analysis. Two of the sets of sums they consider are:
\begin{align}
s_h(m,n) &= \sum_{k=1}^\infty \frac{H_k^m}{(k+1)^n} \label{EBBG1} \\
\sigma_h (m,n) &= \sum_{k=1}^\infty \frac{H_k^{(m)}}{(k+1)^n}, \label{EBBG1a}
\end{align}
where $H_k^{(m)} = \sum_{j=1}^k 1/j^m$. We define slight modifications of these:
\begin{align}
{\cal I}_h(m,n) &= \sum_{k=1}^\infty \frac{H_k^m}{k^n}, \hspace{2em} {\cal J}_h(m,n) = \sum_{k=1}^\infty \frac{H_k^{(m)}}{k^n}.
\label{mydefs1}
\end{align}
Then \cite{BBG1}:
\begin{align}
{\cal J}_h(m,n) &= \sum_{k=1}^\infty \frac{H_k^{(m)}}{k^n}=\sigma_h (m,n)+\zeta (m+n),
\label{EBBG3}
\end{align}
where $\zeta(p) = \sum_{k \geq 1} 1/k^p$ is the Riemann zeta function. For $m=1$, this is
\begin{align}
{\cal I}_h(1,n)-s_h(1,n) &= \zeta (n+1).
\label{myeqn}
\end{align}
For the special case $m=2$ under particular investigation in  \cite{BBG1,BBG2},
\begin{align}
{\cal I}_h(2,n)-s_h(2,n) &= 2 {\cal I}_h(1,n+1)-\zeta (n+2)=2 s_h(1,n+1)+\zeta (n+2).
\label{EBBG4}
\end{align}
For both these cases, the right-hand side tends down to unity as $n\rightarrow \infty$.

The sums ${\cal I}_h(n,p)$ can be represented in terms of the  $s_h(m,n)$  as follows:
\begin{align}
{\cal I}_h(n,p) &= \sum_{k=1}^\infty \frac{H_k^n}{k^p}=\zeta (n+p)+\sum_{m=1}^{n} \binom{n}{m} s_h(m,p+n-m).
\label{EBBG4a}
\end{align}
We can define the order of this expression to be $n+p$, i.e., the sum of the powers of $H_k$ and $k$ on the left-hand side. The sum of the arguments of $s_h$ on the right-hand side
is also $n+p$. The dual expression to equation (\ref{EBBG4a}) is
\begin{align}
s_h(n,p) &= \sum_{k=1}^\infty \frac{H_k^n}{(k+1)^p}=(-1)^n\zeta (n+p)+\sum_{m=1}^{n} \binom{n}{m} (-1)^{m-n} {\cal I}_h (m,p+n-m).
\label{EBBG4aa}
\end{align}
Euler provided the solution for $\sigma_h (1,m)=s_h(1,m)$ and thus for ${\cal J}_h(1,m)={\cal I}_h(1,m)$ for all $m\ge 2$:
\begin{align}
\sigma_h (1,m) &= \frac{m}{2} \zeta (m+1)-\frac{1}{2}  \sum_{k=1}^{m-2} \zeta (m-k) \zeta (k+1).
\label{EBBG5}
\end{align}
Another useful relationship is the reflection formula, valid for $m,n\ge 2$:
\begin{align}
\sigma_h (m,n)+\sigma_h (n,m) &= \zeta(m)\zeta(n)-\zeta(m+n),
\label{EBBG6}
\end{align}
or, written with different notation,
\begin{align}
{\cal J}_h(m,n)+{\cal J}_h(n,m) &= \zeta(m)\zeta(n)+\zeta(m+n),
\label{EBBG7}
\end{align}
so that $2{\cal J}_h(m,m)=\zeta(m)^2+\zeta(2m)$ for $m\ge 2$.

Euler \cite{BBG2} was able to derive the following expansions in terms of zeta functions, for the particular case where the sum of parameters $s+t$ is odd, and $t>1$. The first is for
$s$ odd, $t$ even:
\begin{align}
\sigma_h(s,t) &= \frac{1}{2} \left[ \binom{s+t}{s}-1\right]\zeta(s+t)+\zeta(s) \zeta(t) \nonumber \\
 & \hspace{1em} -\sum_{j=2}^{(s+t-1)/2} \left[\binom{2 j-2}{s-1}+ \binom{2 j-2}{t-1}\right] \zeta(2 j-1) \zeta(s+t-2 j+1). \label{EBBG8}
\end{align}
For $s$ even, $t$ odd:
\begin{align}
\sigma_h(s,t) &= -\frac{1}{2} \left[ \binom{s+t}{s}+1\right]\zeta(s+t)  \nonumber \\
 & \hspace{1em} +\sum_{j=2}^{(s+t-1)/2} \left[\binom{2 j-2}{s-1}+ \binom{2 j-2}{t-1}\right] \zeta(2 j-1) \zeta(s+t-2 j+1). \label{EBBG9}
\end{align}
A valuable result derived in \cite{BBG2} is:
\begin{align}
s_h(2,2 n-1) &= \frac{1}{6} (2 n^2-7 n-3) \zeta (2 n+1)+\zeta (2) \zeta(2 n-1)-\frac{1}{2} \sum_{k=1}^{n-2} (2 k-1)  \zeta (2 n-1-2 k) \zeta (2 k+2) \nonumber \\
 & \hspace{1em} +\frac{1}{3}\sum_{k=1}^{n-2} \zeta (2 k+1) \sum_{j=1}^{n-2-k} \zeta(2 j+1)  \zeta ( k+1-j) \zeta (2n -1-2 k-2 j) .
\label{sh2n}
\end{align}

Table 5 of \cite{BBG1} gives a list of sums $s_h$ for which the authors were unable to find representations in terms of zeta functions or zeta functions complemented by powers of logarithms of integers and polylogarithms of argument $1/2$, using various search algorithms. These results highlight the difficulty of  finding closed form representations of all the Euler-type sums arising in treatments of the sums ${\cal C}_{n,p}$ for $p$ large.

The literature on Euler sums \cite{Chu}-\cite{Zheng} concentrates on the sums $s_h(m,n)$, $\sigma_h(m,n)$, ${\cal I}_h(m,n)$ and ${\cal J}_h(m,n)$. Below we will analyze the mixed sums $M(m,n,p,q)$, which include the existing results for $s_h(m,n)$ (setting $n = q = 0$) and ${\cal I}_h(m,n)$ (setting $p = q = 0$) as special cases. The results we present in Appendix 2 include literature results up to order 7, extend them to orders 8--11 and also include selected results for order 12, as well as adding many extra results of all orders ({\em inter alia} for $p, q \neq 0$).

Note that the extension to orders larger than 7 is not straightforward. Order 8 was painstakingly investigated by Bailey, Borwein and Girgensohn \cite{BBG1}, and only the single analytic result for ${\cal I}_h(1,7)$ was found. By refining and extending the numerical methods used, and regarding ${\cal I}_h(2,6)$ as a known quantity, we have been able to obtain all the other sums for order eight in closed form, as reported in Appendix 2. Similar methods have been applied for orders 9, 10, 11 and 12, with the addition of the sums ${\cal I}_h(2,6), \, {\cal I}_h(2,8), \, {\cal I}_h(3,8), \, {\cal I}_h (2,10), \, {\cal I}_h (4,8)$ to the set of assumed constants. In Apppendix 1, the numerical values of these five assumed constants are given to 400 figures accuracy.

Some evaluations of ${\cal I}(m,n)$ are now presented \cite{Chu,BBG1,BBG2,ChoiSri,Zheng}, arranged according to the order $m+n$ (in the remainder of this section we will drop the $h$ subscript on ${\cal I}_h$). For order 3, there is only one:
\begin{align}
{\cal I}(1,2) &= \sum_{k=1}^\infty \frac{H(k)}{k^{2}} =  2\zeta(3). \label{pow3}
\end{align}
For order 4, there are two:
\begin{align}
{\cal I}(1,3) &= \sum_{k=1}^\infty \frac{H(k)}{k^{3}} = \frac{ 1}{4}\left(  5\zeta(4)\right), \hspace{2em} {\cal I}(2,2) = \sum_{k=1}^\infty \frac{H(k)^{2}}{k^{2}} = \frac{ 1}{4}\left(  17\zeta(4)\right). \label{pow4}
\end{align}
For order 5, there are three:
\begin{align}
{\cal I}(1,4) &= \sum_{k=1}^\infty \frac{H(k)}{k^{4}} =  3\zeta(5) -\zeta(2)\zeta(3), \hspace{2em} {\cal I}(2,3) = \sum_{k=1}^\infty \frac{H(k)^{2}}{k^{3}} = \frac{ 1}{2}\left(  7\zeta(5) -2\zeta(2)\zeta(3)\right), \nonumber \\
{\cal I}(3,2) &= \sum_{k=1}^\infty \frac{H(k)^{3}}{k^{2}} =  10\zeta(5) +\zeta(2)\zeta(3) . \label{pow5}
\end{align}
For order 6, there are four:
\begin{align}
\sum_{k=1}^\infty \frac{H(k)}{k^{5}} &= \frac{1}{4}\left(  7\zeta(6) -2\zeta(3)^2\right), \hspace{2em} \sum_{k=1}^\infty \frac{H(k)^{2}}{k^{4}} = \frac{1}{24}\left(  97\zeta(6) -48\zeta(3)^2\right), \nonumber \\  
\sum_{k=1}^\infty \frac{H(k)^{3}}{k^{3}} &= \frac{ 1}{16}\left(  93\zeta(6) -40\zeta(3)^2\right), \hspace{2em} \sum_{k=1}^\infty \frac{H(k)^{4}}{k^{2}} = \frac{1}{24}\left(  979\zeta(6) +72\zeta(3)^2\right). \label{pow6}
\end{align}
For order 7, there are five:
\begin{align}
{\cal I}(1,6) &= \sum_{k=1}^\infty \frac{H(k)}{k^6} =-\zeta (4) \zeta (3)-\zeta(2) \zeta (5) +4 \zeta (7), \nonumber \\ 
{\cal I}(2,5) &= \sum_{k=1}^\infty \frac{H(k)^2}{k^5} =-\frac{5}{2}\zeta (4) \zeta (3)-\zeta(2) \zeta (5) +6 \zeta (7), \nonumber \\ 
{\cal I}(3,4) &= \sum_{k=1}^\infty \frac{H(k)^3}{k^4} = \frac{693}{48} \zeta (7)+2 \zeta (5) \zeta (2)-\frac{51}{4} \zeta (4) \zeta (3), \nonumber \\ 
{\cal I}(4,3) &= \sum_{k=1}^\infty \frac{H(k)^4}{k^3} =\frac{185}{8} \zeta (7)+5 \zeta (5) \zeta (2)-\frac{43}{2} \zeta (4) \zeta (3), \nonumber \\ 
{\cal I}(5,2) &= \sum_{k=1}^\infty \frac{H(k)^5}{k^2} =\frac{2051}{16} \zeta (7)+\frac{57}{2} \zeta (5) \zeta (2)+ 33\zeta (4) \zeta (3). \label{pow7}
\end{align}

For order eight, there is a paucity of results in the literature. A careful study of this case was given by Bailey, Borwein and Girgensohn \cite{BBG1}. It employed an Euler-Maclaurin scheme for the high-precision evaluation of these sums, an enhanced version of which we describe in Section \ref{sec:Computational} below. The only analytic formula we can give in complete form is one studied by Euler:
\begin{equation}
{\cal I}(1,7)=\sum_{k=1}^\infty \frac{H(k)}{k^7} =\frac{9}{2} \zeta (8)-\zeta (6) \zeta (2)-\zeta (5) \zeta (3) -\frac{1}{2} \zeta (4)^2 ,
\label{pow8a}
\end{equation}
which can be simplified using the analytic expressions for $\zeta (2 n)$ to
\begin{align}
{\cal I}(1,7) &= \frac{1}{4}\left(  9\zeta(8) -4\zeta(3)\zeta(5)\right). \label{pow8b}
\end{align}

We have been able to establish solutions for four additional ${\cal I}$ constants if we express them in terms of the set of constants $\zeta (8)$, $\zeta(3) \zeta (5)$ and $\zeta (2) \zeta (3)^2$, together with ${\cal I}(2,6)$:
\begin{align}
{\cal I}(3,5) &= \frac{1}{96}\left(  595\zeta(8) +120\zeta(2)\zeta(3)^2 -576\zeta(3)\zeta(5) -264 {\cal I}(2,6)\right) \label{pow8part1} \\
{\cal I}(4,4) &= \frac{ 1}{144}\left( -14833\zeta(8) -4032\zeta(2)\zeta(3)^2 +16704\zeta(3)\zeta(5) +3744 {\cal I}(2,6)\right)  \label{pow8part2} \\
{\cal I}(5,3) &= \frac{1}{288}\left(  67811\zeta(8) +19080\zeta(2)\zeta(3)^2 -78768\zeta(3)\zeta(5) -16920 { \cal I}(2,6)\right) \label{pow8part3} \\
{\cal I}(6,2) &= \frac{ 1}{8}\left(  5843\zeta(8) -328\zeta(2)\zeta(3)^2 +3896\zeta(3)\zeta(5) +456 {\cal I}(2,6)\right) \label{pow8part4}
\end{align}

An earlier study \cite{BBG1} gives the expansions for all $s_h(m,n)$ with $m+n=9$, apart from those coming from (\ref{EBBG5}) and (\ref{sh2n}). The basis of function values needed is $\zeta (9)$, $\zeta(2) \zeta (7)$, $\zeta(3) \zeta (6)$, $\zeta (4) \zeta (5)$ and $\zeta (3)^3$, the last coming from the double sum in (\ref{sh2n}). These may be used with (\ref{EBBG4a}) to produce the following evaluations of ${\cal I}(m,n)$ for order $m+n = 9$:
\begin{align}
{\cal I}(1,8) &= 5\zeta(9) -\zeta(3)\zeta(6) -\zeta(4)\zeta(5) -\zeta(2)\zeta(7) \\
{\cal I}(2,7) &= \frac{ 1}{6}\left(  55\zeta(9) -21\zeta(3)\zeta(6) -15\zeta(4)\zeta(5) -6\zeta(2)\zeta(7) +2\zeta(3)^3\right) \\
{\cal I}(3,6) &= \frac{1}{24}\left(  521\zeta(9) -291\zeta(3)\zeta(6) -306\zeta(4)\zeta(5) +72\zeta(2)\zeta(7) +48\zeta(3)^3\right) \\
{\cal I}(4,5) &= \frac{1}{12}\left(  436\zeta(9) -279\zeta(3)\zeta(6) -258\zeta(4)\zeta(5) +84\zeta(2)\zeta(7) +40\zeta(3)^3\right) \\
{\cal I}(5,4) &= \frac{1}{72}\left(  9442\zeta(9) -14685\zeta(3)\zeta(6) +4752\zeta(4)\zeta(5) +2385\zeta(2)\zeta(7) -360\zeta(3)^3\right) \\
{\cal I}(6,3) &= \frac{ 1}{24}\left(  7474\zeta(9) -13122\zeta(3)\zeta(6) +6048\zeta(4)\zeta(5) +1953\zeta(2)\zeta(7) -544\zeta(3)^3\right)\\
{\cal I}(7,2) &= \frac{ 1}{72}\left(  276341\zeta(9) +88665\zeta(3)\zeta(6) +143163\zeta(4)\zeta(5) +59166\zeta(2)\zeta(7) +4032\zeta(3)^3\right)
\end{align}
The approximate numerical value of  ${\cal I}(7,2)$ is $9043.54574728044$; its integral estimate is $8976.6033415307$.

We next present the first results we know of for order $m+n = 10$. We originally obtained these results using the method described in Section \ref{sec:Computational}. The eight basic sums ${\cal I}$ are obtained with two sums ${\cal I}(2,6)$ and ${\cal I}(2,8)$ assumed known:
\begin{align}
{\cal I}(1,9) &= \frac{1}{4}\left(  11\zeta(10) -4\zeta(3)\zeta(7) -2\zeta(5)^2\right) \\
{\cal I}(3,7) &= \frac{ 1}{160}\left( -1661\zeta(10) +1280\zeta(3)\zeta(7) +80\zeta(3)^2\zeta(4) -560\zeta(2)\zeta(3)\zeta(5) +720\zeta(5)^2 \right. \nonumber \\ &\left. \hspace{1em}
+520 {\cal I}(2,8)\right)
\end{align}

\begin{align}
{\cal I}(4,6) &= \frac{ 1}{640}\left( -68823\zeta(10) +60000\zeta(3)\zeta(7) +1000\zeta(3)^2\zeta(4) -21680\zeta(2)\zeta(3)\zeta(5)  \right. \nonumber \\ &\left. \hspace{1em}
+23560\zeta(5)^2 +12120 {\cal I}(2,8) +1280\zeta(2) {\cal I}(2,6)\right) \\
{\cal I}(5,5) &= \frac{1}{256}\left(  64433\zeta(10) -57760\zeta(3)\zeta(7) +360\zeta(3)^2\zeta(4)  +20560\zeta(2)\zeta(3)\zeta(5) \right. \nonumber \\ &\left. \hspace{1em}
-22648\zeta(5)^2 -10920 {\cal I}(2,8) -1280\zeta(2) {\cal I}(2,6)\right) \\
{\cal I}(6,4) &= \frac{ 1}{128}\left( -271367\zeta(10) +176560\zeta(3)\zeta(7) -84648\zeta(3)^2\zeta(4) -400\zeta(2)\zeta(3)\zeta(5)  \right. \nonumber \\ &\left. \hspace{1em}
+121688\zeta(5)^2 +34376 {\cal I}(2,8) +15040\zeta(2) {\cal I}(2,6)\right) \\
{\cal I}(7,3) &= \frac{1}{2560}\left(  16614991\zeta(10) -10315520\zeta(3)\zeta(7) +5879160\zeta(3)^2\zeta(4) -705040\zeta(2)\zeta(3)\zeta(5) \right. \nonumber \\ &\left. \hspace{1em}
-7710760\zeta(5)^2 -2021880 {\cal I}(2,8) -1008000\zeta(2) {\cal I}(2,6)\right) \\
{\cal I}(8,2) &= \frac{ 1}{480}\left(  18741581\zeta(10) +6689520\zeta(3)\zeta(7) -524640\zeta(3)^2\zeta(4) +1452480\zeta(2)\zeta(3)\zeta(5) \right. \nonumber \\ &\left. \hspace{1em}
+4247040\zeta(5)^2 +485280 {\cal I}(2,8) +299520\zeta(2) {\cal I}(2,6)\right)
\end{align}

We now present results for order $m+n = 11$, which again are new in this study, and which again were originally found by us using the methods described below in Section \ref{sec:Computational}. These formulas involve the two sums ${\cal I}(2,6)$ and ${\cal I}(3,8)$.
\begin{align}
{\cal I}(1,10) &= 6\zeta(11) -\zeta(2)\zeta(9) -\zeta(3)\zeta(8) -\zeta(4)\zeta(7) -\zeta(5)\zeta(6) \\
{\cal I}(2,9) &= \frac{1}{2}\left(  26\zeta(11) -2\zeta(2)\zeta(9) -9\zeta(3)\zeta(8) -5\zeta(4)\zeta(7) -7\zeta(5)\zeta(6) +2\zeta(3)^2\zeta(5)\right)
 \\
{\cal I}(4,7) &= \frac{ 1}{48}\left( -2877\zeta(11) -272\zeta(2)\zeta(9) +1190\zeta(3)\zeta(8) +1212\zeta(4)\zeta(7) +1018\zeta(5)\zeta(6)  \right. \nonumber \\ &\left. \hspace{1em}
+80\zeta(2)\zeta(3)^3 -576\zeta(3)^2\zeta(5) +176 {\cal I}(3,8)\right) \\
{\cal I}(5,6) &= \frac{ 1}{576}\left( -781671\zeta(11) -88016\zeta(2)\zeta(9) +296660\zeta(3)\zeta(8) +411984\zeta(4)\zeta(7)  \right. \nonumber \\ &\left. \hspace{1em}
+220080\zeta(5)\zeta(6) +21120\zeta(2)\zeta(3)^3 -141120\zeta(3)^2\zeta(5) +8640\zeta(3) {\cal I}(2,6)  \right. \nonumber \\ &\left. \hspace{1em}
+27840 {\cal I}(3,8)\right) \\
{\cal I}(6,5) &= \frac{1}{192}\left(  734643\zeta(11)  83472\zeta(2)\zeta(9) -271244\zeta(3)\zeta(8) -395088\zeta(4)\zeta(7)  \right. \nonumber \\ &\left. \hspace{1em}
-205424\zeta(5)\zeta(6) -19360\zeta(2)\zeta(3)^3 +130176\zeta(3)^2\zeta(5) -9120\zeta(3) {\cal I}(2,6)  \right. \nonumber \\ &\left. \hspace{1em}
-25600 {\cal I}(3,8)\right) \\
{\cal I}(7,4) &= \frac{1}{1152}\left(  16370805\zeta(11)  1684144\zeta(2)\zeta(9) +5889744\zeta(3)\zeta(8) -10724760\zeta(4)\zeta(7) \right. \nonumber \\ &\left. \hspace{1em}
-10480104\zeta(5)\zeta(6) +844032\zeta(2)\zeta(3)^3 -2330496\zeta(3)^2\zeta(5) -1431360\zeta(3) {\cal I}(2,6) \right. \nonumber \\ &\left. \hspace{1em}
-630336 {\cal I}(3,8)\right) \\
{\cal I}(8,3) &= \frac{1}{72}\left(  2824380\zeta(11)  277304\zeta(2)\zeta(9) +1926401\zeta(3)\zeta(8) -1998972\zeta(4)\zeta(7) \right. \nonumber \\ &\left. \hspace{1em}
-2270310\zeta(5)\zeta(6) +243648\zeta(2)\zeta(3)^3 -803808\zeta(3)^2\zeta(5) -341280\zeta(3) {\cal I}(2,6) \right. \nonumber \\ &\left. \hspace{1em}
-113760 {\cal I}(3,8)\right) \\
{\cal I}(9,2) &= \frac{ 1}{64}\left(  7739347\zeta(11) +2048432\zeta(2)\zeta(9) +5357920\zeta(3)\zeta(8)  \right. \nonumber \\ &\left. \hspace{1em}
+8811792\zeta(4)\zeta(7) +10526056\zeta(5)\zeta(6) -294208\zeta(2)\zeta(3)^3 +2064192\zeta(3)^2\zeta(5)  \right. \nonumber \\ &\left. \hspace{1em}
+540096\zeta(3) {\cal I}(2,6) +199936 {\cal I}(3,8)\right)
\end{align}

Finally, we present results for order $m+n = 12$, which as before are new to this study, having been originally obtained by us using the methods described in Section \ref{sec:Computational}. These results involve the two sums ${\cal I}(2,10)$ and ${\cal I}(4,8)$:
\begin{align}
\sum_{k=1}^\infty \frac{H(k)}{k^{11}} &= \frac{1}{4}\left(  13\zeta(12) -4\zeta(3)\zeta(9) -4\zeta(5)\zeta(7)\right) \\
\sum_{k=1}^\infty \frac{H(k)^{3}}{k^{9}} &= \frac{1}{22112}\left(  355355\zeta(12) -221120\zeta(3)\zeta(9) -265344\zeta(5)\zeta(7)  \right. \nonumber \\ &\left. \hspace{1em}
-33168\zeta(3)^2\zeta(6) +5528\zeta(3)^4 +49752\zeta(2)\zeta(5)^2 +99504\zeta(2)\zeta(3)\zeta(7)  \right. \nonumber \\ &\left. \hspace{1em}
-82920 {\cal I}(2,10)\right) \\
\sum_{k=1}^\infty \frac{H(k)^{5}}{k^{7}} &= \frac{ 1}{265344}\left(  3612841\zeta(12) -884480\zeta(3)\zeta(9) -597024\zeta(5)\zeta(7)  \right. \nonumber \\ &\left. \hspace{1em}
+364848\zeta(3)^2\zeta(6) +221120\zeta(3)^4 +364848\zeta(2)\zeta(5)^2 +729696\zeta(2)\zeta(3)\zeta(7)  \right. \nonumber \\ &\left. \hspace{1em}
-3250464\zeta(3)\zeta(4)\zeta(5) -1028208 {\cal I}(2,10) +663360 {\cal I}(4,8)\right) \\
\sum_{k=1}^\infty \frac{H(k)^{6}}{k^{6}} &= \frac{ 1}{530688}\left( -4262917573\zeta(12) +2820739392\zeta(3)\zeta(9) +2446737024\zeta(5)\zeta(7)  \right. \nonumber \\ &\left. \hspace{1em}
+112663404\zeta(3)^2\zeta(6) -41128320\zeta(3)^4 -402626352\zeta(2)\zeta(5)^2  \right. \nonumber \\ &\left. \hspace{1em}
-741769152\zeta(2)\zeta(3)\zeta(7) -205077744\zeta(3)\zeta(4)\zeta(5) +52538112\zeta(4) {\cal I}(2,6)  \right. \nonumber \\ &\left. \hspace{1em}
+84213552\zeta(2) {\cal I}(2,8) +519676224 {\cal I}(2,10) -22554240 {\cal I}(4,8)\right) \\
\sum_{k=1}^\infty \frac{H(k)^{7}}{k^{5}} &= \frac{ 1}{1061376}\left( -29991036967\zeta(12) +19798731008\zeta(3)\zeta(9) +17219233536\zeta(5)\zeta(7)  \right. \nonumber \\ &\left. \hspace{1em}
+722473668\zeta(3)^2\zeta(6) -292232192\zeta(3)^4 -2832315024\zeta(2)\zeta(5)^2  \right. \nonumber \\ &\left. \hspace{1em}
-5220245184\zeta(2)\zeta(3)\zeta(7) -1329671952\zeta(3)\zeta(4)\zeta(5) +381697344\zeta(4) {\cal I}(2,6)  \right. \nonumber \\ &\left. \hspace{1em}
+589494864\zeta(2) {\cal I}(2,8) +3662808576 {\cal I}(2,10) -167166720 {\cal I}(4,8)\right) \\ 
\sum_{k=1}^\infty \frac{H(k)^{8}}{k^{4}} &= \frac{ 1}{199008}\left( -6469168763\zeta(12) -4417645920\zeta(3)\zeta(9) +2316436536\zeta(5)\zeta(7)  \right. \nonumber \\ &\left. \hspace{1em}
-7185432600\zeta(3)^2\zeta(6) +210815808\zeta(3)^4 +2292190728\zeta(2)\zeta(5)^2  \right. \nonumber \\ &\left. \hspace{1em}
+3705761136\zeta(2)\zeta(3)\zeta(7) +4396086720\zeta(3)\zeta(4)\zeta(5) +2171077776\zeta(4) {\cal I}(2,6)  \right. \nonumber \\ &\left. \hspace{1em}
+241230864\zeta(2) {\cal I}(2,8) -842782296 {\cal I}(2,10) +116552352 {\cal I}(4,8)\right) \\
\sum_{k=1}^\infty \frac{H(k)^{9}}{k^{3}} &= \frac{ 1}{176896}\left(  4340755723\zeta(12) -37498812096\zeta(3)\zeta(9) -8003239392\zeta(5)\zeta(7)  \right. \nonumber \\ &\left. \hspace{1em}
-29417337684\zeta(3)^2\zeta(6) +1136645248\zeta(3)^4 +12010630320\zeta(2)\zeta(5)^2  \right. \nonumber \\ &\left. \hspace{1em}
+20062394496\zeta(2)\zeta(3)\zeta(7) +18880585488\zeta(3)\zeta(4)\zeta(5) +8292663360\zeta(4) {\cal I}(2,6)  \right. \nonumber \\ &\left. \hspace{1em}
+375428592\zeta(2) {\cal I}(2,8) -7045878240 {\cal I}(2,10) +635233536 {\cal I}(4,8)\right) \\
\sum_{k=1}^\infty \frac{H(k)^{10}}{k^{2}} &= \frac{ 1}{176896}\left(  702828643635\zeta(12) +39514453568\zeta(3)\zeta(9) +93510608736\zeta(5)\zeta(7)  \right. \nonumber \\ &\left. \hspace{1em}
-23538514220\zeta(3)^2\zeta(6) +2706951040\zeta(3)^4 +35094519056\zeta(2)\zeta(5)^2  \right. \nonumber \\ &\left. \hspace{1em}
+62104868800\zeta(2)\zeta(3)\zeta(7) +96955381936\zeta(3)\zeta(4)\zeta(5) +16400028160\zeta(4) {\cal I}(2,6)  \right. \nonumber \\ &\left. \hspace{1em}
+954077520\zeta(2) {\cal I}(2,8) -12973442080 {\cal I}(2,10) +1115329280 {\cal I}(4,8)\right)
\end{align}

\subsection{Asymptotics of the sum ${\cal I} (m,2)$}

The summands of the Euler sums ${\cal I} (m,n)$ are always positive, and increase as $m$ increases, while decreasing as $n$ increases. The values of  sums depending on the ${\cal I} (m,n)$  discussed in this paper tend to be dominated  by the lowest sum ${\cal I} (m,2)$ for large values of $m$, and so it is valuable to have asymptotic approximations for it. The value of the sum can be well estimated by an integral, given that the maximum of the truncated summand $(\log k+\gamma)^m/k^2$ occurs for $k=\exp (m/2)$, large enough for the discrete sum to be well approximated by the corresponding integral. For a general positive integer $q$, the result follows from the recursion
\begin{align}
{\cal N}_{q+1} &= \int_1^\infty \frac{(\log k+\gamma)^{q+1} \, {\rm d}k}{k^2}=\gamma^{q+1}+(q+1) {\cal N}_q, \hspace{2em} {\cal N}_1=1+\gamma ,
\label{intrecur}
\end{align}
where $\gamma = 0.5772156649\ldots$ is Euler's constant. This recurrence can be solved exactly, giving 
\begin{align}
\int_1^\infty \frac{(\log k+\gamma)^m}{k^2} {\rm d}k &= m! \left[e^\gamma\right]_m . \label{form:pow8c}
\end{align}
Here we have introduced the notation for the truncated exponential:
\begin{equation}
\left[e^\gamma\right]_m =1+\sum_{q=1}^m \frac{\gamma^q}{q!} . \label{pow8cc}
\end{equation}
Note that for large $m$, ${\cal N}_m/m! \rightarrow \exp(\gamma)$.
 
Although the integral in equation (\ref{form:pow8c}) is exactly evaluated, its use in approximating the sums ${\cal I} (m,2)$ for $m$ large depends on two approximations: the sum is well approximated by an integral, and the two-term  asymptotic series of the harmonic number function gives a sufficiently accurate representation for the integrand. These approximations are tested in Table \ref{tab:tabIm2}, which shows that the integral approximation gains relative accuracy rapidly as $m$ increases, until at $m=9$ it is accurate to two parts in 1000.
  
\begin{table}
\begin{center}
\begin{tabular}{|c|c|c|c|}\hline
$m$&${\cal I} (m,2)$& Formula \eqref{form:pow8c} &ratio \\
\hline
 1 & 2.4041138063 & 1.5772156649 & 0.656048 \\
 2 & 4.5998737432 & 3.4876092536 & 0.758196 \\
 3 & 12.346581901 & 10.655143277 & 0.863003 \\
 4 & 45.833941465 & 42.731580639 & 0.932313 \\
 5 & 220.80305576 & 213.72197848 & 0.967930 \\
 6 & 1302.2827194 & 1282.3688561 & 0.984708 \\
 7 & 9043.5457472 & 8976.6033415 & 0.992597 \\
 8 & 72074.045293 & 71812.839054 & 0.996375 \\
 9 & 647472.79308 & 646315.55860 & 0.998212 \\
\hline
\end{tabular}
\end{center}
\caption{The sums ${\cal I} (m,2)$ are compared with their integral approximation (\ref{form:pow8c}), together with ratios.}\label{tab:tabIm2}
\end{table}

\section{Mixed Euler sums}

In the previous sections, we have focused on the ${\cal I}$ sums, whose denominators are powers of $k$, and on the $s_h$ sums, which have powers of $k+1$. But one is immediately led to consider more general denominators, which have not been previously studied in the literature in any detail. To that end we now consider ``mixed Euler sums,'' namely sums such as 
\begin{align}
M(m,n_0,n_1,n_2, \ldots, n_t) &= \sum_{k=1}^\infty \frac{H(k)^m}{k^{n_0} (k+1)^{n_1} (k+2)^{n_2} \cdots (k+t)^{n_t}}, \label{form:Mixed}
\end{align}
for nonnegative integers $m$ and $(n_i)$, with $m \geq 1$ and $n_0 + n_1 + \cdots + n_t \geq 2$, where $H(k) = H_k = \sum_{j=1}^k 1/j$ is the harmonic function as before, and where $r = m + n_0 + n_1 + \cdots + n_t$ is the order. It is clear that the $s_h$ and ${\cal I}$ sums are merely special cases: $s_h(m,n) = M(m,0,n)$ and ${\cal I}(m,n) = M(m,n)$, so hereafter we will use the $M$ notation. We first demonstrate, by means of examples, why Euler sums with more complicated denominators can be reduced to the basic $M(m,n)$ cases.

\vspace{1ex}
\begin{theorem}\label{thm:Thm1}
If the order of a mixed Euler sum of the form \eqref{form:Mixed} is $12$ or less, then it is expressible as a rational linear sum of terms chosen from the following list, depending on the order as shown (constants for a given order include all those of smaller orders, plus the listed ``additional constants''):

\vspace{1ex} \noindent
Constants for order 3: $1, \, \zeta(2), \, \zeta(3)$ \\
Additional constant for order 4: $\zeta(4)$ \\
Additional constants for order 5: $\zeta(5), \, \zeta(2)\zeta(3)$ \\
Additional constants for order 6: $\zeta(6), \, \zeta(3)^2$ \\
Additional constants for order 7: $\zeta(7), \, \zeta(2)\zeta(5), \, \zeta(3)\zeta(4)$ \\
Additional constants for order 8: $\zeta(8), \, \zeta(2)\zeta(3)^2, \, \zeta(3)\zeta(5), \, M(2,6)$ \\
Additional constants for order 9: $\zeta(9), \, \zeta(2)\zeta(7), \, \zeta(3)\zeta(6), \, \zeta(4)\zeta(5), \, \zeta(3)^3$ \\
Additional constants for order 10: $\zeta(10), \, \zeta(3)\zeta(7), \, \zeta(3)^2\zeta(4), \, \zeta(2)\zeta(3)\zeta(5), \, \zeta(5)^2, \, \zeta(2)M(2,6), \, M(2,8)$ \\
Additional constants for order 11: $\zeta(11), \, \zeta(2)\zeta(9), \, \zeta(3)\zeta(8), \, \zeta(4)\zeta(7), \, \zeta(5)\zeta(6), \, \zeta(2)\zeta(3)^3, \, \zeta(5)\zeta(3)^2$, \\
\hspace*{2em} $\zeta(3)M(2,6), \, M(3,8)$ \\
Additional constants for order 12: $\zeta(12), \, \zeta(3)\zeta(9), \, \zeta(5)\zeta(7), \, 
\zeta(2)\zeta(5)^2, \, \zeta(2)\zeta(3)\zeta(7), \, \zeta(3)\zeta(4)\zeta(5),$ \\
\hspace*{2em} $\zeta(3)^2\zeta(6), \, \zeta(3)^4, \, \zeta(4)M(2,6), \, \zeta(2)M(2,8), \, M(2,10), \, M(4,8)$
\end{theorem}

\noindent
{\em Note:} We conjecture that the representation of a order-12 or less mixed Euler sum as a rational linear combination of the constants in the list in Theorem \ref{thm:Thm1} above is unique, since integer relation computations on this set rule out any relations with reasonable-sized coefficients (see next paragraph for details), but we have no proof of this. We also conjecture that a result similar to Theorem 1 applies for all higher orders: most likely it only remains to identify the appropriate ``atoms,'' akin to the list in Theorem \ref{thm:Thm1}.

We should also clarify that Theorem 1 relies in part on some results in the previous section that were obtained using the computational techniques described below in Section \ref{sec:Computational}.

Note that the above list includes the constants $M(2,6), \, M(2,8), \, M(3,8), \, M(2,10)$ and $M(4,8)$. These constants appear to be linearly independent from the rest of the set, as indicated by the fact that a multipair PSLQ computer run (see Section \ref{sec:Computational}) with the full set of order 8 constants shown above finds no integer relation with Euclidean norm less than $5.88 \cdot 10^{22}$; the full set of order 10 constants above produces no integer relation with Euclidean norm less than $1.28 \cdot 10^{13}$; and the full set of order 12 constants produces no integer relation with Euclidean norm less than $2.13 \cdot 10^6$. Nevertheless, the question of whether $M(2,6), \, M(2,8), \, M(3,8), \, M(2,10)$ and $M(4,8)$, singly or collectively, can be expressed analytically in terms of zetas or other well-known mathematical constants remains open. As an aid to further research, we include 400-digit values of these constants in Appendix 1 (Section \ref{sec:Digits}).

\vspace{1ex}\noindent
{\em Sketch of proof:} We first observe (see Section \ref{sec:Previous-Results} above) that \emph{each} of the basic Euler sums $M(m,n) = {\cal I}(m,n)$ with order $m+n \leq 12$ is reducible to a rational linear sum of the above-listed ``atomic'' constants. We now argue that any general mixed Euler sum \eqref{form:Mixed} of order 12 or less can be reduced to a rational linear combination of the basic Euler sums $M(m,n)$ of the same order or less, and thus to a rational linear combination of the constants in Theorem \ref{thm:Thm1}, by the application (possibly repeated) of these two algebraic techniques:
\begin{enumerate}
\item Changing sums with expressions involving $(k+1), \, (k+2)$ or $(k+w)$ for any integer $w > 0$ to sums involving only $k$, by means of a process akin to ``completing the square'' of elementary algebra.
\item Applying a partial fraction decomposition: Recall that any rational function can be written uniquely as the sum of terms based on the factorization of the denominator polynomial, as in the example
\begin{align}
\frac{1}{(k + 1)(k + 2)^2} &= \frac{1}{k+1} - \frac{1}{k+2} - \frac{1}{(k+2)^2}.
\end{align}
This can be produced in \emph{Wolfram Mathematica} by the command: \tt{Apart[1/((k+1)*(k+2)\^{}2)]}.
\end{enumerate}

\noindent
To illustrate these techniques, note that one can write $M(2,0,2) = \sum_{k=1}^\infty H(k)^2 / (k+1)^2$ as
\begin{align}
M(2,0,2) &= \sum_{k=1}^\infty \frac{H(k)^2}{(k+1)^2} \nonumber \\
  &= \frac{1}{2^2} + \frac{(1+1/2)^2}{3^2} + \frac{(1+1/2+1/3)^2}{4^2} + \frac{(1+1/2+1/3+1/4)^2}{5^2} + \cdots \nonumber \\
  &= \left(\frac{(1+1/2)^2}{2^2} - \frac{2/2}{2^2} - \frac{1/4}{2^2}\right) + \left(\frac{(1+1/2+1/3)^2}{3^2} - \frac{2/3(1+1/2)}{3^2} - \frac{1/9}{3^2}\right) \nonumber \\
  &\hspace*{1em} + \left(\frac{(1+1/2+1/3+1/4)^2}{4^2} - \frac{2/4(1+1/2+1/3)}{4^2} - \frac{1/16}{4^2}\right) + \cdots \nonumber \\
  &= \left(\frac{(1+1/2)^2}{2^2} + \frac{(1+1/2+1/3)^2}{3^2} + \cdots\right) - 2 \left(\frac{1}{2^3} + \frac{(1+1/2)}{3^3} \cdots\right) - \left(\frac{1}{2^4} + \frac{1}{3^4} + \cdots\right) \nonumber \\
  &= \left(\sum_{k=1}^\infty \frac{H(k)^2}{k^2} - 1\right) - 2 \sum_{k=1}^\infty \frac{H(k)}{(k+1)^3} - (\zeta(4) - 1) \nonumber \\
  &= M(2,2) - 2 \, M(1,0,3) - \zeta(4). \label{form:M202}
\end{align}
Note, crucially, that this manipulation rewrites the mixed Euler sum $M(2,0,2)$ (of order 4) to an expression involving $M(2,2)$ (of order 4), the mixed sum $M(1,0,3) = \sum_{k=1}^\infty H(k)/(k+1)^3$, (also of order 4), and the constant $\zeta(4)$ (again of order 4). A similar manipulation can now be performed on $M(1,0,3)$:
\begin{align}
M(1,0,3) &= \sum_{k=1}^\infty \frac{H(k)}{(k+1)^3} \; = \; \frac{1}{2^3} + \frac{(1+1/2)}{3^3} + \frac{(1+1/2+1/3)}{4^3} + \cdots \nonumber \\
 &= \left(\frac{(1+1/2)}{2^3} - \frac{1/2}{2^3}\right) + \left(\frac{(1+1/2+1/3)}{3^3} - \frac{1/3}{3^3}\right) + \cdots \nonumber \\
 &= \left(\sum_{k=1}^\infty \frac{H(k)}{k^3} - 1\right) - (\zeta(4) - 1) \nonumber \\
 &= M(1,3) - \zeta(4) \; = \; 5/4 \, \zeta(4) - \zeta(4) \; = \; 1/4 \, \zeta(4), \label{form:M103}
\end{align}
so that $M(2,0,2) = M(2,2) - 2 \, M(1,0,3) - \zeta(4) = 17/4 \, \zeta(4) - 1/2 \, \zeta(4) - \zeta(4) = 11/4 \, \zeta(4)$, which is of order 4. Note that none of these algebraic manipulations increased the order.

A second example of this technique is $M(2,1,1) = \sum_{k \geq 1} H(k)^2/(k(k+1))$. Note that by employing a manipulation similar to that used above in \eqref{form:M202} and \eqref{form:M103}, combined with the partial fraction decomposition
\begin{align}
\frac{1}{k(k+1)} &= \frac{1}{k} - \frac{1}{k+1},
\end{align}
this can be written
\begin{align}
M(2,1,1) &= \sum_{k=1}^\infty \frac{H(k)^2}{k(k+1)} \; = \; \sum_{k=1}^\infty \left(\frac{H(k)^2}{k} - \frac{H(k)^2}{k+1}\right) \nonumber \\
&= \left(\frac{1}{1} + \frac{(1+1/2)^2}{2} + \frac{(1+1/2+1/3)^2}{3} + \cdots\right) - \left(\frac{1}{2} + \frac{(1+1/2)^2}{3} + \frac{(1+1/2+1/3)^2}{4} + \cdots\right) \nonumber \\
&= \left(\frac{1}{1} + \frac{(1+1/2)^2}{2} + \frac{(1+1/2+1/3)^2}{3} + \cdots\right) \nonumber \\
& \hspace*{1em} - \left[\left(\frac{(1+1/2)^2}{2} - \frac{2/2}{2} - \frac{1/4}{2}\right) + \left(\frac{(1+1/2+1/3)^2}{3} - \frac{2/3 \,(1+1/2)}{3} - \frac{1/9}{3}\right) \right. \nonumber \\
& \hspace*{1em} \left. + \left(\frac{(1+1/2+1/3+1/4)^2}{4} - \frac{2/4 \, (1+1/2+1/3)}{4} - \frac{1/16}{4}\right) + \cdots\right] \nonumber \\
&= 1 + 2 \sum_{k=1}^\infty \frac{H(k)}{(k+1)^2} + (\zeta(3) - 1) \; = \; 2 M(1,2) + \zeta(3) \; = \; 3 \, \zeta(3). \label{form:M111}
\end{align}
Note again that none of these operations increased the order; in fact, in this case the order of the final result, namely 3, is less than the order of the original problem, namely 4.

Consider now a more complicated sum such as $M(2,2,2) = \sum_{k \geq 1} H(k)^2 / (k^2 (k+1)^2$. Sums like this can be readily reduced by means of a partial fraction decomposition, which in this case is:
\begin{align}
\frac{1}{k^2(k+1)^2} &= -2\left(\frac{1}{k} - \frac{1}{k+1}\right) + \frac{1}{k^2} + \frac{1}{(k+1)^2},
\end{align}
so that
\begin{align}
M(2,2,2) &= \sum_{k=1}^\infty \frac{H(k)^2}{k^2 (k+1)^2} \nonumber \\
 &= -2 \sum_{k=1}^\infty H(k)^2 \left(\frac{1}{k} - \frac{1}{k+1}\right) + \sum_{k=1}^\infty \frac{H(k)^2}{k^2} + \sum_{k=1}^\infty \frac{H(k)^2}{(k+1)^2} \nonumber \\
 &= -6 \zeta(3) + 17/4 \, \zeta(4) + 11/4 \, \zeta(4) \nonumber \\
 &= 7 \zeta(4) - 6 \zeta(3),
\end{align}
where we have employed results from \eqref{form:M202}, \eqref{form:M103} and \eqref{form:M111} above.

One example involving $(k+2)$ is $M(2,0,0,2) = \sum_{k \geq 1} H(k)^2/(k+2)^2$. This can be reduced as follows (omitting details of some intermediate evaluations using the above techniques):
\begin{align}
M(2,0,0,2) &= \sum_{k=1}^\infty \frac{H(k)^2}{(k+2)^2} \nonumber \\
  &= \frac{1^2}{3^2} + \frac{(1+1/2)^2}{4^2} + \frac{(1+1/2+1/3)^2}{5^2} + \frac{(1+1/2+1/3+1/4)^2}{6^2} + \cdots \nonumber \\
  &= \left(\frac{(1+1/2+1/3)^2}{3^2} - \frac{2(1/2+1/3)}{3^2} - \frac{(1/2+1/3)^2}{3^2} \right) \nonumber \\
  & \hspace*{1em} + \left(\frac{(1/2+1/3+1/4)^2}{4^2} - \frac{2(1+1/2)(1/3+1/4)}{4^2} - \frac{(1/3+1/4)^2}{4^2} \right) \nonumber \\
  & \hspace*{1em} + \left(\frac{1/2+1/3+1/4+1/5)^2}{5^2} - \frac{2(1+1/2+1/3)(1/4+1/5)}{5^2} - \frac{(1/4+1/5)^2}{5^2}\right) + \cdots \nonumber
\end{align}

\begin{align}
  &= \sum_{k=3}^\infty \frac{H(k)^2}{k^2} - 2 \sum_{k=1}^\infty \frac{H(k)}{(k+1)(k+2)^2} - 2 \sum_{k=1}^\infty \frac{H(k)}{(k+2)^3} - \sum_{k=1}^\infty \frac{1}{(k+1)^2 (k+2)^2} \nonumber \\
  & \hspace*{1em} -  2 \sum_{k=1}^ \infty \frac{1}{(k+1)(k+2)^3} - \sum_{k=1}^\infty \frac{1}{(k+2)^4} \nonumber \\
  &= \left(-\frac{25}{16} + \frac{17}{4} \, \zeta(4)\right) - 2\left(3 - \zeta(2) - \zeta(3)\right) - 2 \left(-3 + \zeta(2) + \zeta(3) + \frac{1}{4} \zeta(4) \right) \nonumber \\
  & \hspace*{1em} - \left(-\frac{13}{4} + 2 \zeta(2)\right) - 2 \left(\frac{23}{8} - \zeta(2) - \zeta(3)\right) - \left(-\frac{17}{16} + \zeta(4)\right) \nonumber \\
  &= -3 + 2 \zeta(3) + \frac{11}{4} \zeta(4) \label{M2002}
\end{align}

The same techniques work for denominators involving $(k+w)$ for any integer $w > 2$. For example, after rather laborious effort one can deduce that
\begin{align}
M(2,0,0,0,0,2) &= \sum_{k=1}^\infty \frac{H(k)^2}{(k+4)^2} \; = \; \frac{1}{3456}\left(-1045 + 
288 \, \zeta(2) + 648 \, \zeta(3) \right).
\end{align}

We should note, however, that the algebraic manipulations required in these evaluations grow very sharply in complexity with increasing powers of $H(k)$ in the numerator and increasing terms in the denominator. Thus we have found, quite frankly, that in most cases these analytic formulas are more easily obtained by the computational methods we describe below in Section \ref{sec:Computational}.

\section{Euler sum analogues to the Stieltjes constants}

As an application of these techniques, we address a finite Euler sum result from Choi and Srivastava\cite{ChoiSri}:
\begin{equation}
 \sum_{k=1}^{n-1} \frac{H_k}{k}=\frac{1}{2} H_{n-1}^2+\frac{1}{2} H_{n-1}^{(2)}.
\label{ChoiSri}
\end{equation}
(Note that in this section it is more convenient to use the subscript notation for the harmonic function: $H_k = H(k)$.) We note also  the following result from \cite{mathstack2}:
\begin{equation}
\sum_{k=1}^n\frac{H_k^{(2)}}{k}=H_n^{(2)} H_n-\sum_{k=1}^n \frac{H_k}{k^2}+H_n^{(3)}.
\label{harm9}
\end{equation}
This is easily generalised to an arbitrary harmonic number of order $p$:
\begin{equation}
\sum_{k=1}^n\frac{H_k^{(p)}}{k}=H_n^{(p)} H_n-\sum_{k=1}^n \frac{H_k}{k^p}+H_n^{(p+1)}.
\label{harm10}
\end{equation}

We continue with the sums over a finite range, commencing with a result  given in \cite{mathstack1}:
\begin{equation}
\sum_{k=1}^n \frac{H_k^2}{k}=\frac{1}{3} H_n^3-\frac{1}{3} H_n^{(3)}+\sum_{k=1}^n \frac{H_k}{k^2}.
\label{form:1onk1}
\end{equation}
The method used to derive \eqref{form:1onk1} employs Abel's summation formula and may easily be generalised to higher values of $p$. When this was done, a pattern emerged  for all the sums:
\begin{equation}
\sum_{k=1}^n \frac{H_k^{p-1}}{k}=\frac{1}{p} H_n^p+ {\cal D}_{p,1} H_n^{(p)}+\sum_{q=2}^{p-1} {\cal D}_{p,q} \sum_{k=1}^n \frac{H_k^{q-1}}{k^{p-q+1}}.
\label{form:1onk2}
\end{equation}
Employing this pattern, it is easy to evaluate the coefficients ${\cal D}_{p,q}$ by choosing the same number of values of $n$ as the number of unknowns, and solving linear equations for the $p-1$ unknowns. The values obtained can easily be checked for other values of $n$. Note that the coefficients in the linear equations are exactly known, and the values for the ${\cal D}_{p,q}$ are also exact. Some values are given in Table \ref{tabDpq}.

\begin{table}
\begin{tabular}{||c|c||}\hline
$p$ &  ${\cal D}_{p,q}$ \\ \hline
2 & 1/2 \\ \hline
3 & -1/3,  1 \\ \hline
4 & 1/4, -1, 3/2 \\ \hline
5 & -1/5, 1, -2, 2 \\ \hline
6 & 1/6, -1, 5/2, -10/3, 5/2 \\ \hline
7 & -1/7, 1, -3, 5, -5, 3 \\ \hline
8 & 1/8, -1, 7/2, -7, 35/4, -7, 7/2 \\ \hline
9 & -1/9, 1, -4, 28/3, -14, 14, -28/3, 4 \\ \hline
10 & 1/10, -1, 9/2, -12, 21, -126/5, 21, -12, 9/2 \\ \hline
11 & -1/11, 1, -5, 15, -30, 42, -42, 30, -15, 5 \\ \hline
12 & 1/12, -1, 11/2, -55/3, 165/4, -66, 77, -66, 165/4, -55/3, 11/2 \\ \hline
13 & -1/13, 1, -6, 22, -55, 99, -132, 132, -99, 55, -22, 6 \\ \hline
14 & 1/14, -1, 13/2, -26, 143/2, -143, 429/2, -1716/7, 429/2, -143, 143/2, -26, 13/2 \\ \hline
15 & -1/15, 1, -7, 91/3, -91, 1001/5, -1001/3, 429, -429, 1001/3, -1001/5, 91, -91/3, 7 \\ \hline
16 & 1/16, -1, 15/2, -35, 455/4, -273, 1001/2, -715, 6435/8, -715, 1001/2, -273, 455/4, -35, 15/2 \\ \hline
17 & -1/17, 1, -8, 40, -140, 364, -728, 1144, -1430, 1430, -1144, 728, -364, 140, -40, 8 \\ \hline
18 & 1/18, -1, 17/2, -136/3, 170, -476, 3094/3, -1768, 2431, -1768, 3094/3, -476, 170, -136/3, 17/2 \\ \hline
19&1/19, 1, -9, 51, -204, 612, -1428, 2652, -3978, 4862, -4862, 3978, -2652, 1428, -612, 204, -51, 9 \\ \hline
20&1/20, -1, 19/2, -57, 969/4, -3876/5, 1938, -3876, 12597/2, -8398, 9237, \\
 & -8398, 12597/2, -3876, 1938, -3876/5, 969/4, -57, 19/2 \\ \hline
\end{tabular}
\label{tabDpq}
\caption{The quantities ${\cal D}_{p,q}$ in equation \eqref{form:1onk2} and \eqref{form:1onk4} for various values of $p$.}
\end{table}

The  lists of coefficients in Table \ref{tabDpq} have some evident properties. The sum of the ${\cal D}_{p,q}$ over $q$ when combined with $1/p$ from the first term on the right-hand side in  equation (\ref{form:1onk2}) is required to be unity, so that the results for $n=1$  on both sides of equation (\ref{form:1onk2}) match. For $p$ even,  ${\cal D}_{p,1}=1/p$ and later coefficients show an even symmetry. For $p$ odd, later coefficients show an odd symmetry. For $p$ odd, the second coefficient in Table \ref{tabDpq} is 1, while for $p$ even it is -1. There are $p-1$ coefficients ${\cal D}$, with the first two being $1/p,\pm 1$. The rest of the  ${\cal D}$'s  fall into $(p-3)/2$ pairs which combine subtractively for $p$ odd, or $(p-4)/2$ additive pairs and a central element for $p$ even.

We take the limit as $n\rightarrow \infty$ in equation (\ref{form:1onk2}) to define a set of harmonic sum Stieltjes constants $\gamma_p^H$, where:
\begin{equation}
\lim_{n\rightarrow \infty}\left[\sum_{k=1}^n \frac{H_k^{p-1}}{k}-\frac{1}{p} H_n^p\right]=\gamma_p^H\label{form:1onk3},
\end{equation}
and 
\begin{eqnarray}
\gamma_p^H &=&  {\cal D}_{p,1} \zeta (p)+\sum_{q=2}^{p-1} {\cal D}_{p,q} \sum_{k=1}^\infty \frac{H_k^{q-1}}{k^{p-q+1}}\nonumber \\
 &=& {\cal D}_{p,1} \zeta (p)+\sum_{q=2}^{p-1} {\cal D}_{p,q} M_{q-1,p-q+1}.
\label{form:1onk4}
\end{eqnarray}
The most slowly convergent by direct summation of the terms in \eqref{form:1onk4} occurs for $q=p-1$, where the summand goes to zero as $\log(k)^{p-2}/k^2$.

By applying formula \eqref{form:1onk4}, together with the computational techniques described below in Section \ref{sec:Computational}, we were able to obtain these results:
\begin{align}
\gamma_2^H &= \frac{1}{2}\left(  \zeta(2)\right) \\
\gamma_3^H &= \frac{ 1}{3}\left(  5\zeta(3)\right) \\
\gamma_4^H &= \frac{ 1}{8}\left(  43\zeta(4)\right) \\
\gamma_5^H &= \frac{ 1}{5}\left(  79\zeta(5) +15\zeta(2)\zeta(3)\right) \\
\gamma_6^H &= \frac{ 1}{24}\left(  2187\zeta(6) +272\zeta(3)^2\right) \\
\gamma_7^H &= \frac{1}{56}\left(  18311\zeta(7) +4060\zeta(2)\zeta(5) +8358\zeta(3)\zeta(4)\right) \\
\gamma_8^H &= \frac{1}{576}\left(  1926401\zeta(8) +48384\zeta(2)\zeta(3)^2 +440064\zeta(3)\zeta(5)\right) \\
\gamma_9^H &= \frac{1}{36}\left(  501978\zeta(9) +266355\zeta(3)\zeta(6) +241794\zeta(4)\zeta(5) +105273\zeta(2)\zeta(7) +12104\zeta(3)^3\right) \\
\gamma_{10}^H &= \frac{1}{80}\left(  17061619\zeta(10) +3161210\zeta(3)\zeta(7) +705180\zeta(3)^2\zeta(4) +928080\zeta(2)\zeta(3)\zeta(5)  \right. \nonumber \\ &\left. \hspace{1em}
+1770112\zeta(5)^2 +37320\zeta(2) M(2,6)\right) \\
\gamma_{11}^H &= \frac{1}{264}\left(  230253219\zeta(11) +49094276\zeta(2)\zeta(9) +165822855\zeta(3)\zeta(8) +130449891\zeta(4)\zeta(7)  \right. \nonumber \\ &\left. \hspace{1em}
+156493260\zeta(5)\zeta(6) +805200\zeta(2)\zeta(3)^3 +19281504\zeta(3)^2\zeta(5) +1849320\zeta(3) M(2,6) \right. \nonumber \\ &\left. \hspace{1em}
+1232880 M(3,8)\right)
\end{align}

Using the integral estimate ${\cal N}_n\rightarrow [ \exp(\gamma)]_n n!$ and replacing $n$ by $p-2$, we  multiply this by $(p-1)/2$, the final ${\cal D}_{p,q}$ in each line of Table 2.
This gives the estimate for $\gamma_p^H$:
\begin{equation}
\gamma_p^H\approx [\exp(\gamma)]_{p-2} \frac{(p-1)! }{2} .
\label{form:gammapest} 
\end{equation}

The approximation \eqref{form:gammapest} is compared with the $\gamma_p^H$ for $p$ ranging from 3 to 11 in Table \ref{tabgapH}. The trend is clearly for the relative accuracy to improve as $p$ increases.

\begin{table}
\begin{center}
\begin{tabular}{|c|c|c|c|}\hline
$p$ & $\gamma_p^H$ & Formula \eqref{form:gammapest} & Ratio \\
\hline
 3 & 2.0034281719 & 1.5772156649 & 0.787258 \\
 4 & 5.8174873811 & 5.2314138804 & 0.899256 \\
 5 & 22.315371582 & 21.310286555 & 0.954959 \\
 6 & 109.08138223 & 106.82895159 & 0.979350 \\
 7 & 647.55020378 & 641.16593544 & 0.990140 \\
 8 & 4510.0214667 & 4488.2909965 & 0.995181 \\
 9 & 35992.013221 & 35906.413366 & 0.997621 \\
10 & 323539.34424 & 323157.77574 & 0.998820 \\
11 & 3233473.9305 & 3231577.7930 & 0.999413 \\
\hline
\end{tabular}
\end{center}
\label{tabgapH}
\caption{The values $\gamma_p^H$ are compared with their integral approximations (\ref{form:gammapest}), together with ratios.}
\end{table}

We turn now to equivalent expressions for the case when the denominator in the basic sum is $k+1$ rather than $k$. From \cite{ChoiSri},
\begin{equation}
\sum_{k=1}^{n-1} \frac{H_k}{k+1}=\frac{1}{2} H_{n}^2-\frac{1}{2} H_{n}^{(2)}.
\label{fstkp1}
\end{equation}
The equivalent of equation (\ref{form:1onk1}) is
\begin{equation}
\sum_{k=1}^{n-1} \frac{H_k^2}{k+1}=\frac{1}{3} H_n^3-\frac{1}{3} H_n^{(3)}-\sum_{k=1}^n \frac{H_k}{(k+1)^2}.
\label{1onk1p}
\end{equation}
Once  again, this may be extended to higher powers of $H_k$ in the numerator, giving results of the following form:
\begin{equation}
\sum_{k=1}^{n-1} \frac{H_k^{p-1}}{k+1}=\frac{1}{p} H_n^p + {\cal E}_{p,1} H_n^{(p)}+\sum_{q=2}^{p-1} {\cal E}_{p,q} \sum_{k=1}^n \frac{H_k^{q-1}}{(k+1)^{p-q+1}}.
\label{1onk2p}
\end{equation}
The sums on the right-hand side are the quantities $s_h$ studied inter alia in \cite{BBG1, BBG2}, and connected with the ${\cal I}$ by (\ref {EBBG4aa}). The coefficients ${\cal E}_{p,q} $ are  given for $p$ up to 21 in Table \ref{tabEpq}. Note that all the ${\cal E}_{p,q} $ are negative, unlike the alternating sign behaviour of the ${\cal D}_{p,q} $.

We can define an alternate set of Stieltjes-like constants from equation (\ref{1onk2p}):
\begin{equation}
\lim_{n\rightarrow \infty} \left[\sum_{k=1}^{n-1} \frac{H_k^{p-1}}{k+1}-\frac{1}{p} H_n^p \right]=\gamma^h_p,
\label{hgas}
\end{equation}
where
\begin{equation}
\gamma^h_p= {\cal E}_{p,1} \zeta (p)+\sum_{q=2}^{p-1} {\cal E}_{p,q} s_h (q-1,p-q+1).
\label{hgasa}
\end{equation}
The $\gamma^h_p$ would then all be negative, and by virtue of (\ref{EBBG4}) with a modulus  well approximated asymptotically by (\ref{form:gammapest}).

\begin{table}
\begin{tabular}{||c|c||}\hline
$p$ &  ${\cal E}_{p,q}$ \\ \hline
2& -1/2  \\ \hline
3 & -1/3, -1 \\ \hline
4 & -1/4, -1, -3/2 \\ \hline
5 & -1/5, -1, -2, -2 \\ \hline
6 & -1/6, -1, -5/2, -10/3, -5/2 \\ \hline
7 & -1/7, -1, -3, -5, -5, -3\\ \hline
8 & -1/8, -1, -7/2, -7, -35/4, -7, -7/2 \\ \hline
9 & -1/9, -1, -4, -28/3, -14, -14, -28/3, -4 \\ \hline
10 & -1/10, -1, -9/2, -12, -21, -126/5, -21, -12, -9/2 \\ \hline
11 & -1/11, -1, -5, -15, -30, -42, -42, -30, -15, -5 \\ \hline
12 & -1/12,-1, -(11/2), -(55/3), -(165/4), -66, -77, -66, -(165/4), -(55/3), -(11/2) \\ \hline
13 & -1/13,-1, -6, -22, -55, -99, -132, -132, -99, -55, -22, -6 \\ \hline
14 & -1/14,-1, -(13/2), -26, -(143/2), -143, -(429/2), -(1716/7), -(429/2), -143, -(143/2), -26, -(13/2)\\ \hline
15 & -1/15,-1, -7, -(91/3), -91, -(1001/5), -(1001/3), -429, -429, -(1001/3), -(1001/5), -91, -(91/3), -7 \\ \hline
16 & -1/16,-1, -(15/2), -35, -(455/4), -273, -(1001/2), -715, -(6435/8),\\
&  -715, -(1001/2), -273, -(455/4), -35, -(15/2) \\ \hline
17 & -1/17,-1, -8, -40, -140, -364, -728, -1144, -1430, -1430, -1144, -728, -364, -140, -40, -8 \\ \hline
18 & -1/18,-1, -(17/2), -(136/3), -170, -476, -(3094/3), -1768, -2431, -(24310/9),\\
&  -2431, -1768, -(3094/3), -476, -170, -(136/3), -(17/2)\\ \hline
19 & -1/19, -1, -9, -51, -204, -612, -1428, -2652, -3978, -4862, -4862, -3978, -2652, -1428, -612, -204, -51, -9 \\ \hline
20& -1/20, -1, -(19/2), -57, -(969/4), -(3876/5), -1938, -3876, -(12597/2), -8398, -(46189/5), \\
& -8398, -(12597/2), -3876, -1938, -(3876/5), -(969/4), -57, -(19/2)\\ \hline
21& -1/21,-1, -10, -(190/3), -285, -969, -2584, -(38760/7), -9690, -(41990/3), -16796,\\
 &  -16796, -(41990/3), -9690, -(38760/7), -2584, -969, -285, -(190/3), -10\\ \hline
\end{tabular}
\caption{The coefficients  ${\cal E}_{p,q}$ in equation (\ref{form:1onk2}) for various values of $p$.}
\label{tabEpq}
\end{table}

\section{Computational techniques} \label{sec:Computational}

We initially obtained many of the formulas presented above and in Appendix 2 (Section \ref{sec:Formulas}) by a computational procedure that utilizes advanced techniques to produce a very high-precision numerical value of the sum, then employs an integer relation algorithm to identify the numerical value as a rational linear sum of constants from the list in Theorem \ref{thm:Thm1}. We present here a brief summary of these techniques, which are based in part on schemes described in \cite{BBG1,Bailey2025}.

In this way, we have found formulas for the much more numerous set of mixed Euler cases
\begin{align}
M(m,n,p,q) &= \sum_{k=1}^\infty \frac{H_k^m}{k^n (k+1)^p (k+2)^q} \label{form:Mnpq}
\end{align}
for orders $r = m + n + p + q = 3$ through $11$, and also a selection of results of order 12. As noted above, this class includes the cases $s_h(m,n)$ and ${\cal I}(m,n)$ as subsets.  We present the full collection of these formulas in Appendix 2 (Section \ref{sec:Formulas}).

\subsection{Computing Euler sums to high precision}

One key tool for these computations is the Euler-Maclaurin summation formula \cite[pg.~285]{Atkinson1989}, which approximates a summation as an integral with high-order corrections (here $f(t)$ is assumed to have $(2s+2)$-th order derivatives on $[a,b]$):
\begin{align}
\sum_{j=a}^b f(j) &= \int_a^b f(t) \, {\rm d}t \, + \, \frac{1}{2} \left(f(a) + f(b)\right) \, + \,\sum_{j=1}^s \frac{B_{2j} \left(D^{2j-1} f(b) - D^{2j-1} f(a)\right)}{(2j)!} \, + \, R_s(a,b), \label{form:eulermac}
\end{align}
where $B_k$ is the $k$-th Bernoulli number \cite{NIST}, $D^k f(a)$ is the $k$-th derivative of $f(t)$ evaluated at $t = a$, and
\begin{align}
R_s(a,b) &= \frac{-1}{(2s+2)!} \int_a^b B_{2s+2} (t - [t]) D^{2s+2} f(t) \, {\rm d}t,\end{align}
where $[\cdot]$ denotes greatest integer and $B_{k}(\cdot)$ is the $k$-th Bernoulli polynomial \cite{NIST} (note $B_k = B_k(0)$).

Applying the Euler-Maclaurin summation formula to the harmonic function $H(t) = \sum_{j=1}^t 1/j$ yields
\begin{align}
H(t) &= \gamma + \log(t) + \frac{1}{2t} - \sum_{j=1}^s \frac{B_{2j}}{2j t^{2j}} + R_s(t),
\end{align}
where $\gamma = 0.5772156649\ldots$ is Euler's constant and $|R_s(t)| \leq |B_{2s+2}|/((2s+2) t^{2s+2})$; see \cite{BBG1} for full details. In the computations for the present study, we set $s = 21$, so that $H(t)$ is approximated by
\begin{align}
\hat{H}(t) &= \gamma + \log(t) + \frac{1}{2 t} - \frac{1}{12 t^2} + \frac{1}{120 t^4} - \frac{1}{252 t^6} + \frac{1}{240 t^8} - \frac{1}{132 t^{10}} + \frac{691}{32760 t^{12}} - \frac{1}{12 t^{14}} \nonumber \\
 & + \frac{3617}{8160 t^{16}} - \frac{43867}{14364 t^{18}} + \frac{174611}{6600 t^{20}} - \frac{77683}{276 t^{22}} + \frac{236364091}{65520 t^{24}} - \frac{657931}{12 t^{26}} + \frac{3392780147}{3480 t^{28}} \nonumber \\
 & - \frac{1723168255201}{85932 t^{30}} + \frac{7709321041217}{16320 t^{32}} - \frac{151628697551}{12 t^{34}} + \frac{26315271553053477373}{69090840 t^{36}} \nonumber \\
 & - \frac{154210205991661}{12 t^{38}} + \frac{261082718496449122051}{541200 t^{40}} - \frac{1520097643918070802691}{75852 t^{42}}, \label{form:happrox}
\end{align}
which approximates $H(t)$ to within roughly $t^{-44}$ for large $t$. The expression \eqref{form:happrox} can be obtained using \emph{Wolfram Mathematica} with the command  {\tt Series[HarmonicNumber[t],\{t,Infinity,42\}]}.

Given $M(m,n,p,q)$, denote $\hat{G}(t) = \hat{H}(t)^m /(t^n (t+1)^p (t+2)^q)$. Using the Euler-Maclaurin summation formula \eqref{form:eulermac} once again, one can write
\begin{align}
M(m,n,p,q) &= \sum_{j=1}^k \frac{H(j)^m}{j^n (j+1)^p (j+2)^q} + \sum_{j=k+1}^\infty \frac{H(j)^m}{j^n (j+1)^p (j+2)^q} \; \approx \; \sum_{j=1}^k \frac{H(j)^m}{j^n (j+1)^p (j+2)^q} + \sum_{j=k+1}^\infty \hat{G}(j) \nonumber \\ 
 & \approx \sum_{j=1}^k \frac{H(j)^m}{j^n (j+1)^p (j+2)^q} \, + \int_{k+1}^\infty \hat{G}(t) \, {\rm d}t \, + \frac{1}{2} \hat{G}(k + 1) \, - \sum_{j=1}^{s} \frac{B_{2j} D^{2j-1} \hat{G} (k+1)}{(2j)!}, \label{form:Mapprox}
\end{align}
where $s = 21$, which is accurate to within roughly $k^{-44}$. Initially we set $k = 10^8 = 100,000,000$, so the approximation in the second line of \eqref{form:Mapprox} is correct to within roughly $10^{-354}$, which was sufficient for our early investigations.  For larger cases, and for all runs listed in Appendix 2, we set $k = 10^9 = 1,000,000,000$, so this approximation is correct to within roughly $10^{-396}$.

We evaluated the first term of \eqref{form:Mapprox} (the explicit summation) using an arbitrary precision package \cite{Bailey2024}. Using $k = 10^8$ and a working precision of 360 digits (producing roughly 350 good digits) required 5--9 minutes CPU time per case on a 2024 Apple Mac Studio system with an M4 processor; using $k = 10^9$ and a working precision of 420 digits (producing roughly 400 good digits) required 50--90 minutes per case. For our 400-digit computations, we evaluated the second term (the integral) using the exp-sinh quadrature algorithm \cite{Bailey2024,Bailey2005}, with the arbitrary precision software set to 400 digits; this required only 3--4 seconds per case (we first tried to evaluate these integrals using \emph{Wolfram Mathematica} version 14.2, but this failed for larger $m$). The third term is straightforward. The fourth term, which involves the symbolic expansion and numerical evaluation to 400-digit accuracy of high-order derivatives of the approximation function $\hat{G}(t) = \hat{H}(t)^m /(t^n (t+1)^p (t+2)^q)$, where $\hat{H}(t)$ is given by the expression \eqref{form:happrox}, was computed using \emph{Wolfram Mathematica}; this required up to 400 seconds CPU time per case for larger $m$.

\subsection{Using an integer relation algorithm to find formulas}

Once a 400-digit value for a given mixed Euler constant was obtained, we employed the multipair PSLQ algorithm to search for integer relations with known constants \cite{Bailey2025,Bailey2000,Ferguson1999}. Given an $v$-long vector $x = (x_0, x_1, \cdots, x_{v-1})$ of high-precision floating-point reals, the multipair PSLQ algorithm searches for integers $(a_0, a_1, \cdots, a_{v-1})$ such that
$a_0 x_0 + a_1 x_1 + \cdots + a_{v-1} x_{v-1} = 0$ to within available precision, or else establishes that there is no such integer relation within a given bound. The algorithm operates by generating an iterative sequence of $v \times v$ integer matrices $B$, so that the entries of the vector $y = B \cdot x$ become progressively smaller, until one entry of $y$ is numerically zero, at which iteration the algorithm halts, with the relation given by the row of $B$ corresponding to the zero entry of $y$. In the application here, we set $x_0$ to the 400-digit value of $M(m,n,p,q)$. For the other entries of the input $x$ vector, we specified 400-digit values of constants listed in Theorem \ref{thm:Thm1}, depending on the order $r$.

Integer relation detection by any algorithm requires very high precision (at least $v \cdot \max_i \log_{10} |a_i|$ digits) to produce numerically reliable results, since otherwise the real relation, if any, will be lost in a sea of numerical artifacts. An effective check of numerical reliability with the multipair PSLQ algorithm is to note the dynamic range of the entries of the $y$ vector at the iteration of detection. In the computer runs for results presented above and in Appendix \ref{sec:Formulas}, this dynamic range always exceeded $10^{63}$, and in most cases exceeded $10^{300}$. In other words, each of these relations holds to at least 63 digits (and in most cases to more than 300 digits) beyond the level required to discover the relation. However, these results should not be regarded as formally proven by these computations.

Figure \ref{fig:pslqplot} illustrates the process of finding a relation using the multipair PSLQ algorithm and assessing the numerical reliability of the result. This shows the base-10 logarithm of the minimum absolute value of the $y$ vector (vertical axis), plotted against the iteration number (horizontal axis), in the multipair PSLQ computer run that the present authors employed to discover the order-10 formula
\begin{align}
M(1,3,6,0) &= \sum_{k=1}^\infty \frac{H(k)}{k^{3}(k+1)^{6}} = \frac{ 1}{4}\left(  84\zeta(2) -108\zeta(3) -5\zeta(4) -48\zeta(5) +24\zeta(2)\zeta(3) -9\zeta(6) +6\zeta(3)^2  \right. \nonumber \\ &\left. \hspace{1em}
-12\zeta(7) +4\zeta(2)\zeta(5) +4\zeta(3)\zeta(4)\right). \label{form:m1360}
\end{align}
Note that as the algorithm proceeds, the minimum absolute value of the $y$ vector slowly decreases, from approximately $10^{-3}$ to approximately $10^{-65}$, but at iteration 311 abruptly drops to approximately $10^{-405}$, a drop of 340 orders of magnitude. Note that since we are using 400-digit precision, $10^{-405}$ is effectively zero, so the algorithm terminates here with the relation $(4, -84, 108, 5, 48, -24, 9, -6, 12, -4, -4)$. In other words, formula \eqref{form:m1360} holds to roughly 340 digits beyond the precision level required to discover it. This dynamic range at the iteration of detection can thus be considered a ``confidence level'' of the result's numerical reliability.

\begin{figure}[h]
\begin{tikzpicture}
\begin{axis}[
 width=15cm,
 height=10cm,
 xmajorgrids=true,
 ymajorgrids=true,
 grid style=dashed,
 line width=0.5pt,
 mark size=0.1pt,
 xlabel={Iteration},
 xmin=0, xmax=325,
 xtick={50,100,150,200,250,300},
 ymin=-425, ymax=0,
 ytick={-400,-350,-300,-250,-200,-150,-100,-50,0},
 ylabel={$\log_{10} \min_i |y_i|$}]
 \addplot[color=black,solid,mark=*,mark options={solid}]table{NES-pslq-plot.dat};
\end{axis}
\end{tikzpicture}
\caption{Plot of $\log_{10} \min_i |y_i|$ in the multipair PSLQ computer run for $M(1,3,6,0)$, showing the detection of the relation at iteration 311.}\label{fig:pslqplot}
\end{figure}

\subsection{Computational results}

The process described above succeeded in finding relations for \emph{each} of the cases of form \eqref{form:Mnpq} with orders between 3 and 11, a total of 960 cases, plus an additional 54 selected cases of order 12. See Appendix 2 (Section \ref{sec:Formulas}) for a complete listing of these formulas. As noted below, in order to minimize the possibility of transcription errors, \emph{the LaTeX code for each section of results was generated automatically by a computer program from the output computer files, and this LaTeX code is included here without any alteration}.

The techniques described in this section are applicable to more general classes of Euler sums, including Euler sums of orders higher than 12 and Euler sums with more complicated polynomial denominators. However, the computational cost increases with the precision required and the number of selected right-hand-side constants. The principal challenges here are the first and fourth term of \eqref{form:Mapprox}, namely ($a$) the cost of explicitly computing and summing to high precision a large number of terms of the mixed Euler sum series, and ($b$) the symbolic expansion and numerical evaluation of high-order derivatives of the function $\hat{G}(t)$. Perhaps further investigation into the underlying theory of Euler sums will yield computational schemes that are more efficient for large problems.

\section{Conclusions}
We have presented techniques, both algebraic and computational, for finding analytic evaluations of a significantly larger class of Euler sums than studied previously. We believe that most of these formulas are new to the literature. Along this line, we have found that \emph{Wolfram Mathematica} (version 14.2) can evaluate many of the basic cases, but a large majority are not evaluated by this software.

These methods appear to be applicable to even more general Euler sums. For example, by applying the  methods described above, we have obtained these intriguing computational results, among others:
\begin{align}
\sum_{k=1}^\infty \frac{H(k)}{k(2k+1)} &=  2\log(2)^2 \\
\sum_{k=1}^\infty \frac{H(k)}{k^{2}(2k+1)} &= 2\zeta(3) -4\log(2)^2 \\
\sum_{k=1}^\infty \frac{H(k)}{(2k+1)^{2}} &= \frac{ 1}{4}\left(  7\zeta(3) -6\log(2)\zeta(2)\right) \\
\sum_{k=1}^\infty \frac{H(k)}{k^{2}(2k+1)^{2}} &=  9\zeta(3) -6\log(2)\zeta(2) -8\log(2)^2 \\
%\sum_{k=1}^\infty \frac{H(k)}{(2k+1)^{3}} &= \frac{1}{32}\left(  45\zeta(4) -56\log(2)\zeta(3)\right) \\
\sum_{k=1}^\infty \frac{H(k)}{(2k+1)^{4}} &= \frac{ 1}{16}\left(  62\zeta(5) -21\zeta(2)\zeta(3) -30\log(2)\zeta(4)\right)
\end{align}
Note the appearance of $\log(2)$ in these formulas. Each of these formulas holds to nearly 400-digit accuracy (approximately 350 digits beyond the level required to discover them), but at present we do not yet know how they can be rigorously proven.

The new results presented in this study also highlight the benefits of attempting to solve for all Euler sums of a given order. The results given for sums of mixed type may be of use in indicating the zeta function values likely to arise in attempts to numerically solve for recalcitrant sums like those for order eight and higher. It is hoped that the asymptotic form inferred for the constants $\gamma_p^H$ can be deduced rigorously, as it may well prove useful in other applications of high-order Euler sums.

\newpage
\section{Appendix 1: 400-digit values of three key constants}\label{sec:Digits}

We present here 400-digit approximations of $M(2,6), \, M(2,8), \, M(3,8), \, M(2,10), \, M(4,8)$:

\vspace{2ex}
\noindent
1.041413395855265060833934370636480151499859280096830090748511645153773087302971 \\
78483751544719684852509976852158376374740737268847953695380222383595172532123654 \\
63963612795034976112760332996361625685218808108323018034356756036322549570832977 \\
08604139265652530043836463078378465035583569011375448218307043216126923803712749 \\
23988797094981204968396475470138806138535478550733612009250215922048841374239723 \\
645442685850 \\

\noindent
1.009386471889869832518544227219279156409372942639652641202049549364385367847079 \\
49180863769095027121905627225975982985135460410529740749826141104503536876835470 \\
18469301862442802589875242849768895787689895958104331283788277223328457927340866 \\
40158920385626435450329285165922784555461987108701748322359094830741802548831985 \\
88668354450902612335818964472409228594433865464246509588184931824643738162119198 \\
662316661058 \\

\noindent
1.014305290895216264339827024366251554326370696089068947073583456867667637693523 \\
94671117513355550815252027023563642862142136424802329417381850641867359561102369 \\
07708608852232885420834448581394420559852401108798464519014241848466439384418357 \\
64122407964525143823389069592803884034573487533288088530610292952331243167418134 \\
67198428033583320677784291408584463648948157254030603048103031772735772545074505 \\
256977622009 \\

\noindent
1.002258993186511461546882204200782204716716526446955625961726703382258341612187 \\
43078937691748374313286427718378216976785972900132909275418926271525587620211343 \\
23828486038109593077991669275749305259201098402321866062725818826804223344328667 \\
24311160745124744110382924162704634065128094036087399151400598689180216783658166 \\
11238941907596547797317455957463173904843228565114429853394147887119147441919740 \\
167418480233 \\

\noindent
1.021889991239632409955119439812528407213142628943801388819377166083386890422403 \\
37811800890914656945480609918216413770587578233999208809687116705691954898386671 \\
55265235208787528310809835281206739252549491701207237871226480178257316430518840 \\
84028017686753324378706573579491021902617968550914371718501262794713100873318573 \\
13123823952255996488552702615920587505801228078818063210192756692513776767867839 \\
561290224768   

\section{Appendix 2: Formulas for orders 3 through 12}\label{sec:Formulas}

We present here the full set of results for $M(m,n,p,q)$ for orders 3 through 11, plus some additional selected cases of order 12. Each of these formulas holds to at least 380-digit precision, which is at least 63 digits (and in most cases more than 300 digits) beyond the level required to discover the relation. However, these formulas should not be regarded as formally proven solely by these computations.

\emph{To minimize the possibility of transcription errors, in each section below the formulas were produced by a computer program that parses the computer run output files, extracts the formulas, sorts them lexiographically and then generates LaTeX code (including all spacing, line breaks and page breaks). We then applied a separate program, which parses this LaTeX code, numerically evaluates each of the left- and right-hand sides, and verifies equality to 200-digit precision. No errors were found.}

\emph{We have included this LaTeX code below without any alteration. }

\vspace{2ex}
Formulas for order $r = m + n + p + q = 3$:
\begin{align}
\sum_{k=1}^\infty \frac{H(k)}{k^{2}} \sumend &= -\left( -2\zeta(3)\right) \label{eq03001} \\
\sum_{k=1}^\infty \frac{H(k)}{k(k+1)} \sumend &= \left(  \zeta(2)\right) \label{eq03002} \\
\sum_{k=1}^\infty \frac{H(k)}{(k+1)^{2}} \sumend &= -\left( -\zeta(3)\right) \label{eq03003} \\
\sum_{k=1}^\infty \frac{H(k)}{k(k+2)} \sumend &= \frac{1}{2}\left(  1 +\zeta(2)\right) \label{eq03004} \\
\sum_{k=1}^\infty \frac{H(k)}{(k+1)(k+2)} \sumend &= -\left( -1\right) \label{eq03005} \\
\sum_{k=1}^\infty \frac{H(k)}{(k+2)^{2}} \sumend &= \left( -2 +\zeta(2) +\zeta(3)\right) \label{eq03006}
\end{align}

Formulas for order $r = m + n + p + q = 4$:
\begin{align}
\sum_{k=1}^\infty \frac{H(k)}{k^{3}} \sumend &= \frac{-1}{4}\left( -5\zeta(4)\right) \label{eq04001} \\
\sum_{k=1}^\infty \frac{H(k)}{k^{2}(k+1)} \sumend &= -\left(  \zeta(2) -2\zeta(3)\right) \label{eq04002} \\
\sum_{k=1}^\infty \frac{H(k)}{k(k+1)^{2}} \sumend &= \left(  \zeta(2) -\zeta(3)\right) \label{eq04003} \\
\sum_{k=1}^\infty \frac{H(k)}{(k+1)^{3}} \sumend &= \frac{1}{4}\left(  \zeta(4)\right) \label{eq04004} \\
\sum_{k=1}^\infty \frac{H(k)}{k^{2}(k+2)} \sumend &= \frac{-1}{4}\left(  1 +\zeta(2) -4\zeta(3)\right) \label{eq04005} \\
\sum_{k=1}^\infty \frac{H(k)}{k(k+1)(k+2)} \sumend &= \frac{1}{2}\left( -1 +\zeta(2)\right) \label{eq04006} \\
\sum_{k=1}^\infty \frac{H(k)}{(k+1)^{2}(k+2)} \sumend &= -\left(  1 -\zeta(3)\right) \label{eq04007} \\
\sum_{k=1}^\infty \frac{H(k)}{k(k+2)^{2}} \sumend &= \frac{1}{4}\left(  5 -\zeta(2) -2\zeta(3)\right) \label{eq04008} \\
\sum_{k=1}^\infty \frac{H(k)}{(k+1)(k+2)^{2}} \sumend &= -\left( -3 +\zeta(2) +\zeta(3)\right) \label{eq04009} \\
\sum_{k=1}^\infty \frac{H(k)}{(k+2)^{3}} \sumend &= \frac{-1}{4}\left(  12 -4\zeta(2) -4\zeta(3) -\zeta(4)\right) \label{eq04010} \\
\sum_{k=1}^\infty \frac{H(k)^{2}}{k^{2}} \sumend &= \frac{-1}{4}\left( -17\zeta(4)\right) \label{eq04011} \\
\sum_{k=1}^\infty \frac{H(k)^{2}}{k(k+1)} \sumend &= \left(  3\zeta(3)\right) \label{eq04012} \\
\sum_{k=1}^\infty \frac{H(k)^{2}}{(k+1)^{2}} \sumend &= \frac{-1}{4}\left( -11\zeta(4)\right) \label{eq04013} \\
\sum_{k=1}^\infty \frac{H(k)^{2}}{k(k+2)} \sumend &= \frac{1}{2}\left(  1 +\zeta(2) +3\zeta(3)\right) \label{eq04014} \\
\sum_{k=1}^\infty \frac{H(k)^{2}}{(k+1)(k+2)} \sumend &= -\left( -1 -\zeta(2)\right) \label{eq04015} \\
\sum_{k=1}^\infty \frac{H(k)^{2}}{(k+2)^{2}} \sumend &= \frac{-1}{4}\left(  12 -8\zeta(3) -11\zeta(4)\right) \label{eq04016}
\end{align}

\newpage
Formulas for order $r = m + n + p + q = 5$:
\begin{align}
\sum_{k=1}^\infty \frac{H(k)}{k^{4}} \sumend &= -\left( -3\zeta(5) +\zeta(2)\zeta(3)\right) \label{eq05001} \\
\sum_{k=1}^\infty \frac{H(k)}{k^{3}(k+1)} \sumend &= \frac{-1}{4}\left( -4\zeta(2) +8\zeta(3) -5\zeta(4)\right) \label{eq05002} \\
\sum_{k=1}^\infty \frac{H(k)}{k^{2}(k+1)^{2}} \sumend &= \left( -2\zeta(2) +3\zeta(3)\right) \label{eq05003} \\
\sum_{k=1}^\infty \frac{H(k)}{k(k+1)^{3}} \sumend &= \frac{1}{4}\left(  4\zeta(2) -4\zeta(3) -\zeta(4)\right) \label{eq05004} \\
\sum_{k=1}^\infty \frac{H(k)}{(k+1)^{4}} \sumend &= \left(  2\zeta(5) -\zeta(2)\zeta(3)\right) \label{eq05005} \\
\sum_{k=1}^\infty \frac{H(k)}{k^{3}(k+2)} \sumend &= \frac{-1}{8}\left( -1 -\zeta(2) +4\zeta(3) -5\zeta(4)\right) \label{eq05006} \\
\sum_{k=1}^\infty \frac{H(k)}{k^{2}(k+1)(k+2)} \sumend &= \frac{-1}{4}\left( -1 +3\zeta(2) -4\zeta(3)\right) \label{eq05007} \\
\sum_{k=1}^\infty \frac{H(k)}{k(k+1)^{2}(k+2)} \sumend &= \frac{1}{2}\left(  1 +\zeta(2) -2\zeta(3)\right) \label{eq05008} \\
\sum_{k=1}^\infty \frac{H(k)}{(k+1)^{3}(k+2)} \sumend &= \frac{1}{4}\left(  4 -4\zeta(3) +\zeta(4)\right) \label{eq05009} \\
\sum_{k=1}^\infty \frac{H(k)}{k^{2}(k+2)^{2}} \sumend &= \frac{1}{4}\left( -3 +3\zeta(3)\right) \label{eq05010} \\
\sum_{k=1}^\infty \frac{H(k)}{k(k+1)(k+2)^{2}} \sumend &= \frac{1}{4}\left( -7 +3\zeta(2) +2\zeta(3)\right) \label{eq05011} \\
\sum_{k=1}^\infty \frac{H(k)}{(k+1)^{2}(k+2)^{2}} \sumend &= \left( -4 +\zeta(2) +2\zeta(3)\right) \label{eq05012} \\
\sum_{k=1}^\infty \frac{H(k)}{k(k+2)^{3}} \sumend &= \frac{1}{8}\left(  17 -5\zeta(2) -6\zeta(3) -\zeta(4)\right) \label{eq05013} \\
\sum_{k=1}^\infty \frac{H(k)}{(k+1)(k+2)^{3}} \sumend &= \frac{1}{4}\left(  24 -8\zeta(2) -8\zeta(3) -\zeta(4)\right) \label{eq05014} \\
\sum_{k=1}^\infty \frac{H(k)}{(k+2)^{4}} \sumend &= -\left(  4 -\zeta(2) -\zeta(3) -\zeta(4) -2\zeta(5) +\zeta(2)\zeta(3)\right) \label{eq05015} \\
\sum_{k=1}^\infty \frac{H(k)^{2}}{k^{3}} \sumend &= \frac{-1}{2}\left( -7\zeta(5) +2\zeta(2)\zeta(3)\right) \label{eq05016}
\end{align}
 
\begin{align}
\sum_{k=1}^\infty \frac{H(k)^{2}}{k^{2}(k+1)} \sumend &= \frac{1}{4}\left( -12\zeta(3) +17\zeta(4)\right) \label{eq05017} \\
\sum_{k=1}^\infty \frac{H(k)^{2}}{k(k+1)^{2}} \sumend &= \frac{-1}{4}\left( -12\zeta(3) +11\zeta(4)\right) \label{eq05018} \\
\sum_{k=1}^\infty \frac{H(k)^{2}}{(k+1)^{3}} \sumend &= \frac{-1}{2}\left(  3\zeta(5) -2\zeta(2)\zeta(3)\right) \label{eq05019} \\
\sum_{k=1}^\infty \frac{H(k)^{2}}{k^{2}(k+2)} \sumend &= \frac{-1}{8}\left(  2 +2\zeta(2) +6\zeta(3) -17\zeta(4)\right) \label{eq05020} \\
\sum_{k=1}^\infty \frac{H(k)^{2}}{k(k+1)(k+2)} \sumend &= \frac{1}{2}\left( -1 -\zeta(2) +3\zeta(3)\right) \label{eq05021} \\
\sum_{k=1}^\infty \frac{H(k)^{2}}{(k+1)^{2}(k+2)} \sumend &= \frac{-1}{4}\left(  4 +4\zeta(2) -11\zeta(4)\right) \label{eq05022} \\
\sum_{k=1}^\infty \frac{H(k)^{2}}{k(k+2)^{2}} \sumend &= \frac{-1}{8}\left( -14 -2\zeta(2) +2\zeta(3) +11\zeta(4)\right) \label{eq05023} \\
\sum_{k=1}^\infty \frac{H(k)^{2}}{(k+1)(k+2)^{2}} \sumend &= \frac{-1}{4}\left( -16 -4\zeta(2) +8\zeta(3) +11\zeta(4)\right) \label{eq05024} \\
\sum_{k=1}^\infty \frac{H(k)^{2}}{(k+2)^{3}} \sumend &= \frac{1}{2}\left( -12 +2\zeta(2) +6\zeta(3) +\zeta(4) -3\zeta(5) +2\zeta(2)\zeta(3)\right) \label{eq05025} \\
\sum_{k=1}^\infty \frac{H(k)^{3}}{k^{2}} \sumend &= -\left( -10\zeta(5) -\zeta(2)\zeta(3)\right) \label{eq05026} \\
\sum_{k=1}^\infty \frac{H(k)^{3}}{k(k+1)} \sumend &= -\left( -10\zeta(4)\right) \label{eq05027} \\
\sum_{k=1}^\infty \frac{H(k)^{3}}{(k+1)^{2}} \sumend &= \frac{-1}{2}\left( -15\zeta(5) -2\zeta(2)\zeta(3)\right) \label{eq05028} \\
\sum_{k=1}^\infty \frac{H(k)^{3}}{k(k+2)} \sumend &= \frac{-1}{2}\left( -1 -2\zeta(2) -4\zeta(3) -10\zeta(4)\right) \label{eq05029} \\
\sum_{k=1}^\infty \frac{H(k)^{3}}{(k+1)(k+2)} \sumend &= \left(  1 +2\zeta(2) +4\zeta(3)\right) \label{eq05030} \\
\sum_{k=1}^\infty \frac{H(k)^{3}}{(k+2)^{2}} \sumend &= \frac{1}{4}\left( -16 -12\zeta(2) -4\zeta(3) +33\zeta(4) +30\zeta(5) +4\zeta(2)\zeta(3)\right) \label{eq05031}
\end{align}

\newpage
Formulas for order $r = m + n + p + q = 6$:
\begin{align}
\sum_{k=1}^\infty \frac{H(k)}{k^{5}} \sumend &= \frac{1}{4}\left(  7\zeta(6) -2\zeta(3)^2\right) \label{eq06001} \\
\sum_{k=1}^\infty \frac{H(k)}{k^{4}(k+1)} \sumend &= \frac{1}{4}\left( -4\zeta(2) +8\zeta(3) -5\zeta(4) +12\zeta(5) -4\zeta(2)\zeta(3)\right) \label{eq06002} \\
\sum_{k=1}^\infty \frac{H(k)}{k^{3}(k+1)^{2}} \sumend &= \frac{1}{4}\left(  12\zeta(2) -20\zeta(3) +5\zeta(4)\right) \label{eq06003} \\
\sum_{k=1}^\infty \frac{H(k)}{k^{2}(k+1)^{3}} \sumend &= \frac{-1}{4}\left(  12\zeta(2) -16\zeta(3) -\zeta(4)\right) \label{eq06004} \\
\sum_{k=1}^\infty \frac{H(k)}{k(k+1)^{4}} \sumend &= \frac{1}{4}\left(  4\zeta(2) -4\zeta(3) -\zeta(4) -8\zeta(5) +4\zeta(2)\zeta(3)\right) \label{eq06005} \\
\sum_{k=1}^\infty \frac{H(k)}{(k+1)^{5}} \sumend &= \frac{-1}{4}\left( -3\zeta(6) +2\zeta(3)^2\right) \label{eq06006} \\
\sum_{k=1}^\infty \frac{H(k)}{k^{4}(k+2)} \sumend &= \frac{1}{16}\left( -1 -\zeta(2) +4\zeta(3) -5\zeta(4) +24\zeta(5) -8\zeta(2)\zeta(3)\right) \label{eq06007} \\
\sum_{k=1}^\infty \frac{H(k)}{k^{3}(k+1)(k+2)} \sumend &= \frac{-1}{8}\left(  1 -7\zeta(2) +12\zeta(3) -5\zeta(4)\right) \label{eq06008} \\
\sum_{k=1}^\infty \frac{H(k)}{k^{2}(k+1)^{2}(k+2)} \sumend &= \frac{-1}{4}\left(  1 +5\zeta(2) -8\zeta(3)\right) \label{eq06009} \\
\sum_{k=1}^\infty \frac{H(k)}{k(k+1)^{3}(k+2)} \sumend &= \frac{-1}{4}\left(  2 -2\zeta(2) +\zeta(4)\right) \label{eq06010} \\
\sum_{k=1}^\infty \frac{H(k)}{(k+1)^{4}(k+2)} \sumend &= \frac{1}{4}\left( -4 +4\zeta(3) -\zeta(4) +8\zeta(5) -4\zeta(2)\zeta(3)\right) \label{eq06011} \\
\sum_{k=1}^\infty \frac{H(k)}{k^{3}(k+2)^{2}} \sumend &= \frac{-1}{16}\left( -7 -\zeta(2) +10\zeta(3) -5\zeta(4)\right) \label{eq06012} \\
\sum_{k=1}^\infty \frac{H(k)}{k^{2}(k+1)(k+2)^{2}} \sumend &= \frac{-1}{4}\left( -4 +3\zeta(2) -\zeta(3)\right) \label{eq06013} \\
\sum_{k=1}^\infty \frac{H(k)}{k(k+1)^{2}(k+2)^{2}} \sumend &= \frac{-1}{4}\left( -9 +\zeta(2) +6\zeta(3)\right) \label{eq06014} \\
\sum_{k=1}^\infty \frac{H(k)}{(k+1)^{3}(k+2)^{2}} \sumend &= \frac{1}{4}\left(  20 -4\zeta(2) -12\zeta(3) +\zeta(4)\right) \label{eq06015} \\
\sum_{k=1}^\infty \frac{H(k)}{k^{2}(k+2)^{3}} \sumend &= \frac{-1}{16}\left(  23 -5\zeta(2) -12\zeta(3) -\zeta(4)\right) \label{eq06016}
\end{align}
 
\begin{align}
\sum_{k=1}^\infty \frac{H(k)}{k(k+1)(k+2)^{3}} \sumend &= \frac{-1}{8}\left(  31 -11\zeta(2) -10\zeta(3) -\zeta(4)\right) \label{eq06017} \\
\sum_{k=1}^\infty \frac{H(k)}{(k+1)^{2}(k+2)^{3}} \sumend &= \frac{-1}{4}\left(  40 -12\zeta(2) -16\zeta(3) -\zeta(4)\right) \label{eq06018} \\
\sum_{k=1}^\infty \frac{H(k)}{k(k+2)^{4}} \sumend &= \frac{1}{16}\left(  49 -13\zeta(2) -14\zeta(3) -9\zeta(4) -16\zeta(5) +8\zeta(2)\zeta(3)\right) \label{eq06019} \\
\sum_{k=1}^\infty \frac{H(k)}{(k+1)(k+2)^{4}} \sumend &= \frac{1}{4}\left(  40 -12\zeta(2) -12\zeta(3) -5\zeta(4) -8\zeta(5) +4\zeta(2)\zeta(3)\right) \label{eq06020} \\
\sum_{k=1}^\infty \frac{H(k)}{(k+2)^{5}} \sumend &= \frac{-1}{4}\left(  20 -4\zeta(2) -4\zeta(3) -4\zeta(4) -4\zeta(5) -3\zeta(6) +2\zeta(3)^2\right) \label{eq06021} \\
\sum_{k=1}^\infty \frac{H(k)^{2}}{k^{4}} \sumend &= \frac{1}{24}\left(  97\zeta(6) -48\zeta(3)^2\right) \label{eq06022} \\
\sum_{k=1}^\infty \frac{H(k)^{2}}{k^{3}(k+1)} \sumend &= \frac{1}{4}\left(  12\zeta(3) -17\zeta(4) +14\zeta(5) -4\zeta(2)\zeta(3)\right) \label{eq06023} \\
\sum_{k=1}^\infty \frac{H(k)^{2}}{k^{2}(k+1)^{2}} \sumend &= \left( -6\zeta(3) +7\zeta(4)\right) \label{eq06024} \\
\sum_{k=1}^\infty \frac{H(k)^{2}}{k(k+1)^{3}} \sumend &= \frac{-1}{4}\left( -12\zeta(3) +11\zeta(4) -6\zeta(5) +4\zeta(2)\zeta(3)\right) \label{eq06025} \\
\sum_{k=1}^\infty \frac{H(k)^{2}}{(k+1)^{4}} \sumend &= \frac{-1}{24}\left( -37\zeta(6) +24\zeta(3)^2\right) \label{eq06026} \\
\sum_{k=1}^\infty \frac{H(k)^{2}}{k^{3}(k+2)} \sumend &= \frac{-1}{16}\left( -2 -2\zeta(2) -6\zeta(3) +17\zeta(4) -28\zeta(5) +8\zeta(2)\zeta(3)\right) \label{eq06027} \\
\sum_{k=1}^\infty \frac{H(k)^{2}}{k^{2}(k+1)(k+2)} \sumend &= \frac{1}{8}\left(  2 +2\zeta(2) -18\zeta(3) +17\zeta(4)\right) \label{eq06028} \\
\sum_{k=1}^\infty \frac{H(k)^{2}}{k(k+1)^{2}(k+2)} \sumend &= \frac{-1}{4}\left( -2 -2\zeta(2) -6\zeta(3) +11\zeta(4)\right) \label{eq06029} \\
\sum_{k=1}^\infty \frac{H(k)^{2}}{(k+1)^{3}(k+2)} \sumend &= \frac{1}{4}\left(  4 +4\zeta(2) -11\zeta(4) -6\zeta(5) +4\zeta(2)\zeta(3)\right) \label{eq06030} \\
\sum_{k=1}^\infty \frac{H(k)^{2}}{k^{2}(k+2)^{2}} \sumend &= \frac{1}{4}\left( -4 -\zeta(2) -\zeta(3) +7\zeta(4)\right) \label{eq06031} \\
\sum_{k=1}^\infty \frac{H(k)^{2}}{k(k+1)(k+2)^{2}} \sumend &= \frac{1}{8}\left( -18 -6\zeta(2) +14\zeta(3) +11\zeta(4)\right) \label{eq06032}
\end{align}
 
\begin{align}
\sum_{k=1}^\infty \frac{H(k)^{2}}{(k+1)^{2}(k+2)^{2}} \sumend &= \frac{-1}{2}\left(  10 +4\zeta(2) -4\zeta(3) -11\zeta(4)\right) \label{eq06033} \\
\sum_{k=1}^\infty \frac{H(k)^{2}}{k(k+2)^{3}} \sumend &= \frac{-1}{16}\left( -62 +6\zeta(2) +26\zeta(3) +15\zeta(4) -12\zeta(5)  \right. \nonumber \\ &\left. \hspace{1em}
+8\zeta(2)\zeta(3)\right) \label{eq06034} \\
\sum_{k=1}^\infty \frac{H(k)^{2}}{(k+1)(k+2)^{3}} \sumend &= \frac{-1}{4}\left( -40 +20\zeta(3) +13\zeta(4) -6\zeta(5) +4\zeta(2)\zeta(3)\right) \label{eq06035} \\
\sum_{k=1}^\infty \frac{H(k)^{2}}{(k+2)^{4}} \sumend &= \frac{1}{24}\left( -240 +48\zeta(2) +96\zeta(3) +36\zeta(4) +96\zeta(5) -48\zeta(2)\zeta(3)  \right. \nonumber \\ &\left. \hspace{1em}
+37\zeta(6) -24\zeta(3)^2\right) \label{eq06036} \\
\sum_{k=1}^\infty \frac{H(k)^{3}}{k^{3}} \sumend &= \frac{-1}{16}\left( -93\zeta(6) +40\zeta(3)^2\right) \label{eq06037} \\
\sum_{k=1}^\infty \frac{H(k)^{3}}{k^{2}(k+1)} \sumend &= \left( -10\zeta(4) +10\zeta(5) +\zeta(2)\zeta(3)\right) \label{eq06038} \\
\sum_{k=1}^\infty \frac{H(k)^{3}}{k(k+1)^{2}} \sumend &= \frac{1}{2}\left(  20\zeta(4) -15\zeta(5) -2\zeta(2)\zeta(3)\right) \label{eq06039} \\
\sum_{k=1}^\infty \frac{H(k)^{3}}{(k+1)^{3}} \sumend &= \frac{-1}{16}\left(  33\zeta(6) -32\zeta(3)^2\right) \label{eq06040} \\
\sum_{k=1}^\infty \frac{H(k)^{3}}{k^{2}(k+2)} \sumend &= \frac{1}{4}\left( -1 -2\zeta(2) -4\zeta(3) -10\zeta(4) +20\zeta(5) +2\zeta(2)\zeta(3)\right) \label{eq06041} \\
\sum_{k=1}^\infty \frac{H(k)^{3}}{k(k+1)(k+2)} \sumend &= \frac{1}{2}\left( -1 -2\zeta(2) -4\zeta(3) +10\zeta(4)\right) \label{eq06042} \\
\sum_{k=1}^\infty \frac{H(k)^{3}}{(k+1)^{2}(k+2)} \sumend &= \frac{1}{2}\left( -2 -4\zeta(2) -8\zeta(3) +15\zeta(5) +2\zeta(2)\zeta(3)\right) \label{eq06043} \\
\sum_{k=1}^\infty \frac{H(k)^{3}}{k(k+2)^{2}} \sumend &= \frac{1}{8}\left(  18 +16\zeta(2) +12\zeta(3) -13\zeta(4) -30\zeta(5) -4\zeta(2)\zeta(3)\right) \label{eq06044} \\
\sum_{k=1}^\infty \frac{H(k)^{3}}{(k+1)(k+2)^{2}} \sumend &= \frac{1}{4}\left(  20 +20\zeta(2) +20\zeta(3) -33\zeta(4) -30\zeta(5)  \right. \nonumber \\ &\left. \hspace{1em}
-4\zeta(2)\zeta(3)\right) \label{eq06045} \\
\sum_{k=1}^\infty \frac{H(k)^{3}}{(k+2)^{3}} \sumend &= \frac{1}{16}\left( -160 -48\zeta(2) +48\zeta(3) +144\zeta(4) -72\zeta(5) +48\zeta(2)\zeta(3)  \right. \nonumber \\ &\left. \hspace{1em}
-33\zeta(6) +32\zeta(3)^2\right) \label{eq06046}
\end{align}
 
\begin{align}
\sum_{k=1}^\infty \frac{H(k)^{4}}{k^{2}} \sumend &= \frac{1}{24}\left(  979\zeta(6) +72\zeta(3)^2\right) \label{eq06047} \\
\sum_{k=1}^\infty \frac{H(k)^{4}}{k(k+1)} \sumend &= \left(  30\zeta(5) +6\zeta(2)\zeta(3)\right) \label{eq06048} \\
\sum_{k=1}^\infty \frac{H(k)^{4}}{(k+1)^{2}} \sumend &= \frac{-1}{24}\left( -859\zeta(6) -72\zeta(3)^2\right) \label{eq06049} \\
\sum_{k=1}^\infty \frac{H(k)^{4}}{k(k+2)} \sumend &= \frac{1}{4}\left(  2 +6\zeta(2) +22\zeta(3) +37\zeta(4) +60\zeta(5) +12\zeta(2)\zeta(3)\right) \label{eq06050} \\
\sum_{k=1}^\infty \frac{H(k)^{4}}{(k+1)(k+2)} \sumend &= \frac{-1}{2}\left( -2 -6\zeta(2) -22\zeta(3) -37\zeta(4)\right) \label{eq06051} \\
\sum_{k=1}^\infty \frac{H(k)^{4}}{(k+2)^{2}} \sumend &= \frac{-1}{24}\left(  120 +192\zeta(2) +432\zeta(3) +48\zeta(4) -720\zeta(5) -96\zeta(2)\zeta(3)  \right. \nonumber \\ &\left. \hspace{1em}
-859\zeta(6) -72\zeta(3)^2\right) \label{eq06052}
\end{align}

\newpage
Formulas for order $r = m + n + p + q = 7$:
\begin{align}
\sum_{k=1}^\infty \frac{H(k)}{k^{6}} \sumend &= \left(  4\zeta(7) -\zeta(2)\zeta(5) -\zeta(3)\zeta(4)\right) \label{eq07001} \\
\sum_{k=1}^\infty \frac{H(k)}{k^{5}(k+1)} \sumend &= \frac{-1}{4}\left( -4\zeta(2) +8\zeta(3) -5\zeta(4) +12\zeta(5) -4\zeta(2)\zeta(3) -7\zeta(6)  \right. \nonumber \\ &\left. \hspace{1em}
+2\zeta(3)^2\right) \label{eq07002} \\
\sum_{k=1}^\infty \frac{H(k)}{k^{4}(k+1)^{2}} \sumend &= \frac{-1}{2}\left(  8\zeta(2) -14\zeta(3) +5\zeta(4) -6\zeta(5) +2\zeta(2)\zeta(3)\right) \label{eq07003} \\
\sum_{k=1}^\infty \frac{H(k)}{k^{3}(k+1)^{3}} \sumend &= -\left( -6\zeta(2) +9\zeta(3) -\zeta(4)\right) \label{eq07004} \\
\sum_{k=1}^\infty \frac{H(k)}{k^{2}(k+1)^{4}} \sumend &= \frac{1}{2}\left( -8\zeta(2) +10\zeta(3) +\zeta(4) +4\zeta(5) -2\zeta(2)\zeta(3)\right) \label{eq07005} \\
\sum_{k=1}^\infty \frac{H(k)}{k(k+1)^{5}} \sumend &= \frac{-1}{4}\left( -4\zeta(2) +4\zeta(3) +\zeta(4) +8\zeta(5) -4\zeta(2)\zeta(3) +3\zeta(6)  \right. \nonumber \\ &\left. \hspace{1em}
-2\zeta(3)^2\right) \label{eq07006} \\
\sum_{k=1}^\infty \frac{H(k)}{(k+1)^{6}} \sumend &= -\left( -3\zeta(7) +\zeta(2)\zeta(5) +\zeta(3)\zeta(4)\right) \label{eq07007} \\
\sum_{k=1}^\infty \frac{H(k)}{k^{5}(k+2)} \sumend &= \frac{1}{32}\left(  1 +\zeta(2) -4\zeta(3) +5\zeta(4) -24\zeta(5) +8\zeta(2)\zeta(3) +28\zeta(6)  \right. \nonumber \\ &\left. \hspace{1em}
-8\zeta(3)^2\right) \label{eq07008} \\
\sum_{k=1}^\infty \frac{H(k)}{k^{4}(k+1)(k+2)} \sumend &= \frac{1}{16}\left(  1 -15\zeta(2) +28\zeta(3) -15\zeta(4) +24\zeta(5)  \right. \nonumber \\ &\left. \hspace{1em}
-8\zeta(2)\zeta(3)\right) \label{eq07009} \\
\sum_{k=1}^\infty \frac{H(k)}{k^{3}(k+1)^{2}(k+2)} \sumend &= \frac{1}{8}\left(  1 +17\zeta(2) -28\zeta(3) +5\zeta(4)\right) \label{eq07010} \\
\sum_{k=1}^\infty \frac{H(k)}{k^{2}(k+1)^{3}(k+2)} \sumend &= \frac{-1}{4}\left( -1 +7\zeta(2) -8\zeta(3) -\zeta(4)\right) \label{eq07011} \\
\sum_{k=1}^\infty \frac{H(k)}{k(k+1)^{4}(k+2)} \sumend &= \frac{1}{2}\left(  1 +\zeta(2) -2\zeta(3) -4\zeta(5) +2\zeta(2)\zeta(3)\right) \label{eq07012} \\
\sum_{k=1}^\infty \frac{H(k)}{(k+1)^{5}(k+2)} \sumend &= \frac{-1}{4}\left( -4 +4\zeta(3) -\zeta(4) +8\zeta(5) -4\zeta(2)\zeta(3) -3\zeta(6)  \right. \nonumber \\ &\left. \hspace{1em}
+2\zeta(3)^2\right) \label{eq07013} \\
\sum_{k=1}^\infty \frac{H(k)}{k^{4}(k+2)^{2}} \sumend &= \frac{-1}{16}\left(  4 +\zeta(2) -7\zeta(3) +5\zeta(4) -12\zeta(5) +4\zeta(2)\zeta(3)\right) \label{eq07014}
\end{align}
 
\begin{align}
\sum_{k=1}^\infty \frac{H(k)}{k^{3}(k+1)(k+2)^{2}} \sumend &= \frac{1}{16}\left( -9 +13\zeta(2) -14\zeta(3) +5\zeta(4)\right) \label{eq07015} \\
\sum_{k=1}^\infty \frac{H(k)}{k^{2}(k+1)^{2}(k+2)^{2}} \sumend &= \frac{-1}{4}\left(  5 +2\zeta(2) -7\zeta(3)\right) \label{eq07016} \\
\sum_{k=1}^\infty \frac{H(k)}{k(k+1)^{3}(k+2)^{2}} \sumend &= \frac{1}{4}\left( -11 +3\zeta(2) +6\zeta(3) -\zeta(4)\right) \label{eq07017} \\
\sum_{k=1}^\infty \frac{H(k)}{(k+1)^{4}(k+2)^{2}} \sumend &= \frac{1}{2}\left( -12 +2\zeta(2) +8\zeta(3) -\zeta(4) +4\zeta(5) -2\zeta(2)\zeta(3)\right) \label{eq07018} \\
\sum_{k=1}^\infty \frac{H(k)}{k^{3}(k+2)^{3}} \sumend &= \frac{-1}{16}\left( -15 +2\zeta(2) +11\zeta(3) -2\zeta(4)\right) \label{eq07019} \\
\sum_{k=1}^\infty \frac{H(k)}{k^{2}(k+1)(k+2)^{3}} \sumend &= \frac{-1}{16}\left( -39 +17\zeta(2) +8\zeta(3) +\zeta(4)\right) \label{eq07020} \\
\sum_{k=1}^\infty \frac{H(k)}{k(k+1)^{2}(k+2)^{3}} \sumend &= \frac{-1}{8}\left( -49 +13\zeta(2) +22\zeta(3) +\zeta(4)\right) \label{eq07021} \\
\sum_{k=1}^\infty \frac{H(k)}{(k+1)^{3}(k+2)^{3}} \sumend &= \left(  15 -4\zeta(2) -7\zeta(3)\right) \label{eq07022} \\
\sum_{k=1}^\infty \frac{H(k)}{k^{2}(k+2)^{4}} \sumend &= \frac{1}{16}\left( -36 +9\zeta(2) +13\zeta(3) +5\zeta(4) +8\zeta(5) -4\zeta(2)\zeta(3)\right) \label{eq07023} \\
\sum_{k=1}^\infty \frac{H(k)}{k(k+1)(k+2)^{4}} \sumend &= \frac{1}{16}\left( -111 +35\zeta(2) +34\zeta(3) +11\zeta(4) +16\zeta(5)  \right. \nonumber \\ &\left. \hspace{1em}
-8\zeta(2)\zeta(3)\right) \label{eq07024} \\
\sum_{k=1}^\infty \frac{H(k)}{(k+1)^{2}(k+2)^{4}} \sumend &= \frac{-1}{2}\left(  40 -12\zeta(2) -14\zeta(3) -3\zeta(4) -4\zeta(5)  \right. \nonumber \\ &\left. \hspace{1em}
+2\zeta(2)\zeta(3)\right) \label{eq07025} \\
\sum_{k=1}^\infty \frac{H(k)}{k(k+2)^{5}} \sumend &= \frac{-1}{32}\left( -129 +29\zeta(2) +30\zeta(3) +25\zeta(4) +32\zeta(5) -8\zeta(2)\zeta(3)  \right. \nonumber \\ &\left. \hspace{1em}
+12\zeta(6) -8\zeta(3)^2\right) \label{eq07026} \\
\sum_{k=1}^\infty \frac{H(k)}{(k+1)(k+2)^{5}} \sumend &= \frac{-1}{4}\left( -60 +16\zeta(2) +16\zeta(3) +9\zeta(4) +12\zeta(5) -4\zeta(2)\zeta(3)  \right. \nonumber \\ &\left. \hspace{1em}
+3\zeta(6) -2\zeta(3)^2\right) \label{eq07027} \\
\sum_{k=1}^\infty \frac{H(k)}{(k+2)^{6}} \sumend &= \left( -6 +\zeta(2) +\zeta(3) +\zeta(4) +\zeta(5) +\zeta(6) +3\zeta(7) -\zeta(2)\zeta(5)  \right. \nonumber \\ &\left. \hspace{1em}
-\zeta(3)\zeta(4)\right) \label{eq07028}
\end{align}
 
\begin{align}
\sum_{k=1}^\infty \frac{H(k)^{2}}{k^{5}} \sumend &= \frac{1}{2}\left(  12\zeta(7) -2\zeta(2)\zeta(5) -5\zeta(3)\zeta(4)\right) \label{eq07029} \\
\sum_{k=1}^\infty \frac{H(k)^{2}}{k^{4}(k+1)} \sumend &= \frac{-1}{24}\left(  72\zeta(3) -102\zeta(4) +84\zeta(5) -24\zeta(2)\zeta(3) -97\zeta(6)  \right. \nonumber \\ &\left. \hspace{1em}
+48\zeta(3)^2\right) \label{eq07030} \\
\sum_{k=1}^\infty \frac{H(k)^{2}}{k^{3}(k+1)^{2}} \sumend &= \frac{1}{4}\left(  36\zeta(3) -45\zeta(4) +14\zeta(5) -4\zeta(2)\zeta(3)\right) \label{eq07031} \\
\sum_{k=1}^\infty \frac{H(k)^{2}}{k^{2}(k+1)^{3}} \sumend &= \frac{-1}{4}\left(  36\zeta(3) -39\zeta(4) +6\zeta(5) -4\zeta(2)\zeta(3)\right) \label{eq07032} \\
\sum_{k=1}^\infty \frac{H(k)^{2}}{k(k+1)^{4}} \sumend &= \frac{-1}{24}\left( -72\zeta(3) +66\zeta(4) -36\zeta(5) +24\zeta(2)\zeta(3) +37\zeta(6)  \right. \nonumber \\ &\left. \hspace{1em}
-24\zeta(3)^2\right) \label{eq07033} \\
\sum_{k=1}^\infty \frac{H(k)^{2}}{(k+1)^{5}} \sumend &= \frac{1}{2}\left( -2\zeta(7) +2\zeta(2)\zeta(5) -\zeta(3)\zeta(4)\right) \label{eq07034} \\
\sum_{k=1}^\infty \frac{H(k)^{2}}{k^{4}(k+2)} \sumend &= \frac{1}{96}\left( -6 -6\zeta(2) -18\zeta(3) +51\zeta(4) -84\zeta(5) +24\zeta(2)\zeta(3)  \right. \nonumber \\ &\left. \hspace{1em}
+194\zeta(6) -96\zeta(3)^2\right) \label{eq07035} \\
\sum_{k=1}^\infty \frac{H(k)^{2}}{k^{3}(k+1)(k+2)} \sumend &= \frac{-1}{16}\left(  2 +2\zeta(2) -42\zeta(3) +51\zeta(4) -28\zeta(5)  \right. \nonumber \\ &\left. \hspace{1em}
+8\zeta(2)\zeta(3)\right) \label{eq07036} \\
\sum_{k=1}^\infty \frac{H(k)^{2}}{k^{2}(k+1)^{2}(k+2)} \sumend &= \frac{-1}{8}\left(  2 +2\zeta(2) +30\zeta(3) -39\zeta(4)\right) \label{eq07037} \\
\sum_{k=1}^\infty \frac{H(k)^{2}}{k(k+1)^{3}(k+2)} \sumend &= \frac{-1}{2}\left(  1 +\zeta(2) -3\zeta(3) -3\zeta(5) +2\zeta(2)\zeta(3)\right) \label{eq07038} \\
\sum_{k=1}^\infty \frac{H(k)^{2}}{(k+1)^{4}(k+2)} \sumend &= \frac{-1}{24}\left(  24 +24\zeta(2) -66\zeta(4) -36\zeta(5) +24\zeta(2)\zeta(3) -37\zeta(6)  \right. \nonumber \\ &\left. \hspace{1em}
+24\zeta(3)^2\right) \label{eq07039} \\
\sum_{k=1}^\infty \frac{H(k)^{2}}{k^{3}(k+2)^{2}} \sumend &= \frac{-1}{32}\left( -18 -6\zeta(2) -10\zeta(3) +45\zeta(4) -28\zeta(5)  \right. \nonumber \\ &\left. \hspace{1em}
+8\zeta(2)\zeta(3)\right) \label{eq07040} \\
\sum_{k=1}^\infty \frac{H(k)^{2}}{k^{2}(k+1)(k+2)^{2}} \sumend &= \frac{-1}{8}\left( -10 -4\zeta(2) +16\zeta(3) -3\zeta(4)\right) \label{eq07041}
\end{align}
 
\begin{align}
\sum_{k=1}^\infty \frac{H(k)^{2}}{k(k+1)^{2}(k+2)^{2}} \sumend &= \frac{1}{8}\left(  22 +10\zeta(2) -2\zeta(3) -33\zeta(4)\right) \label{eq07042} \\
\sum_{k=1}^\infty \frac{H(k)^{2}}{(k+1)^{3}(k+2)^{2}} \sumend &= \frac{1}{4}\left(  24 +12\zeta(2) -8\zeta(3) -33\zeta(4) -6\zeta(5)  \right. \nonumber \\ &\left. \hspace{1em}
+4\zeta(2)\zeta(3)\right) \label{eq07043} \\
\sum_{k=1}^\infty \frac{H(k)^{2}}{k^{2}(k+2)^{3}} \sumend &= \frac{-1}{32}\left(  78 -2\zeta(2) -22\zeta(3) -43\zeta(4) +12\zeta(5)  \right. \nonumber \\ &\left. \hspace{1em}
-8\zeta(2)\zeta(3)\right) \label{eq07044} \\
\sum_{k=1}^\infty \frac{H(k)^{2}}{k(k+1)(k+2)^{3}} \sumend &= \frac{1}{16}\left( -98 -6\zeta(2) +54\zeta(3) +37\zeta(4) -12\zeta(5)  \right. \nonumber \\ &\left. \hspace{1em}
+8\zeta(2)\zeta(3)\right) \label{eq07045} \\
\sum_{k=1}^\infty \frac{H(k)^{2}}{(k+1)^{2}(k+2)^{3}} \sumend &= \frac{1}{4}\left( -60 -8\zeta(2) +28\zeta(3) +35\zeta(4) -6\zeta(5)  \right. \nonumber \\ &\left. \hspace{1em}
+4\zeta(2)\zeta(3)\right) \label{eq07046} \\
\sum_{k=1}^\infty \frac{H(k)^{2}}{k(k+2)^{4}} \sumend &= \frac{-1}{96}\left( -666 +114\zeta(2) +270\zeta(3) +117\zeta(4) +156\zeta(5) -72\zeta(2)\zeta(3)  \right. \nonumber \\ &\left. \hspace{1em}
+74\zeta(6) -48\zeta(3)^2\right) \label{eq07047} \\
\sum_{k=1}^\infty \frac{H(k)^{2}}{(k+1)(k+2)^{4}} \sumend &= \frac{-1}{24}\left( -480 +48\zeta(2) +216\zeta(3) +114\zeta(4) +60\zeta(5)  \right. \nonumber \\ &\left. \hspace{1em}
-24\zeta(2)\zeta(3) +37\zeta(6) -24\zeta(3)^2\right) \label{eq07048} \\
\sum_{k=1}^\infty \frac{H(k)^{2}}{(k+2)^{5}} \sumend &= \frac{1}{2}\left( -30 +6\zeta(2) +10\zeta(3) +5\zeta(4) +10\zeta(5) -4\zeta(2)\zeta(3) +3\zeta(6)  \right. \nonumber \\ &\left. \hspace{1em}
-2\zeta(3)^2 -2\zeta(7) +2\zeta(2)\zeta(5) -\zeta(3)\zeta(4)\right) \label{eq07049} \\
\sum_{k=1}^\infty \frac{H(k)^{3}}{k^{4}} \sumend &= \frac{1}{16}\left(  231\zeta(7) +32\zeta(2)\zeta(5) -204\zeta(3)\zeta(4)\right) \label{eq07050} \\
\sum_{k=1}^\infty \frac{H(k)^{3}}{k^{3}(k+1)} \sumend &= \frac{-1}{16}\left( -160\zeta(4) +160\zeta(5) +16\zeta(2)\zeta(3) -93\zeta(6)  \right. \nonumber \\ &\left. \hspace{1em}
+40\zeta(3)^2\right) \label{eq07051} \\
\sum_{k=1}^\infty \frac{H(k)^{3}}{k^{2}(k+1)^{2}} \sumend &= \frac{-1}{2}\left(  40\zeta(4) -35\zeta(5) -4\zeta(2)\zeta(3)\right) \label{eq07052} \\
\sum_{k=1}^\infty \frac{H(k)^{3}}{k(k+1)^{3}} \sumend &= \frac{-1}{16}\left( -160\zeta(4) +120\zeta(5) +16\zeta(2)\zeta(3) -33\zeta(6)  \right. \nonumber \\ &\left. \hspace{1em}
+32\zeta(3)^2\right) \label{eq07053}
\end{align}
 
\begin{align}
\sum_{k=1}^\infty \frac{H(k)^{3}}{(k+1)^{4}} \sumend &= \frac{1}{16}\left(  119\zeta(7) +32\zeta(2)\zeta(5) -132\zeta(3)\zeta(4)\right) \label{eq07054} \\
\sum_{k=1}^\infty \frac{H(k)^{3}}{k^{3}(k+2)} \sumend &= \frac{-1}{32}\left( -4 -8\zeta(2) -16\zeta(3) -40\zeta(4) +80\zeta(5) +8\zeta(2)\zeta(3)  \right. \nonumber \\ &\left. \hspace{1em}
-93\zeta(6) +40\zeta(3)^2\right) \label{eq07055} \\
\sum_{k=1}^\infty \frac{H(k)^{3}}{k^{2}(k+1)(k+2)} \sumend &= \frac{-1}{4}\left( -1 -2\zeta(2) -4\zeta(3) +30\zeta(4) -20\zeta(5)  \right. \nonumber \\ &\left. \hspace{1em}
-2\zeta(2)\zeta(3)\right) \label{eq07056} \\
\sum_{k=1}^\infty \frac{H(k)^{3}}{k(k+1)^{2}(k+2)} \sumend &= \frac{-1}{2}\left( -1 -2\zeta(2) -4\zeta(3) -10\zeta(4) +15\zeta(5)  \right. \nonumber \\ &\left. \hspace{1em}
+2\zeta(2)\zeta(3)\right) \label{eq07057} \\
\sum_{k=1}^\infty \frac{H(k)^{3}}{(k+1)^{3}(k+2)} \sumend &= \frac{1}{16}\left(  16 +32\zeta(2) +64\zeta(3) -120\zeta(5) -16\zeta(2)\zeta(3) -33\zeta(6)  \right. \nonumber \\ &\left. \hspace{1em}
+32\zeta(3)^2\right) \label{eq07058} \\
\sum_{k=1}^\infty \frac{H(k)^{3}}{k^{2}(k+2)^{2}} \sumend &= \frac{-1}{16}\left(  20 +20\zeta(2) +20\zeta(3) +7\zeta(4) -70\zeta(5)  \right. \nonumber \\ &\left. \hspace{1em}
-8\zeta(2)\zeta(3)\right) \label{eq07059} \\
\sum_{k=1}^\infty \frac{H(k)^{3}}{k(k+1)(k+2)^{2}} \sumend &= \frac{-1}{8}\left(  22 +24\zeta(2) +28\zeta(3) -53\zeta(4) -30\zeta(5)  \right. \nonumber \\ &\left. \hspace{1em}
-4\zeta(2)\zeta(3)\right) \label{eq07060} \\
\sum_{k=1}^\infty \frac{H(k)^{3}}{(k+1)^{2}(k+2)^{2}} \sumend &= \frac{1}{4}\left( -24 -28\zeta(2) -36\zeta(3) +33\zeta(4) +60\zeta(5)  \right. \nonumber \\ &\left. \hspace{1em}
+8\zeta(2)\zeta(3)\right) \label{eq07061} \\
\sum_{k=1}^\infty \frac{H(k)^{3}}{k(k+2)^{3}} \sumend &= \frac{-1}{32}\left( -196 -80\zeta(2) +24\zeta(3) +170\zeta(4) -12\zeta(5) +56\zeta(2)\zeta(3)  \right. \nonumber \\ &\left. \hspace{1em}
-33\zeta(6) +32\zeta(3)^2\right) \label{eq07062} \\
\sum_{k=1}^\infty \frac{H(k)^{3}}{(k+1)(k+2)^{3}} \sumend &= \frac{1}{16}\left(  240 +128\zeta(2) +32\zeta(3) -276\zeta(4) -48\zeta(5)  \right. \nonumber \\ &\left. \hspace{1em}
-64\zeta(2)\zeta(3) +33\zeta(6) -32\zeta(3)^2\right) \label{eq07063} \\
\sum_{k=1}^\infty \frac{H(k)^{3}}{(k+2)^{4}} \sumend &= \frac{1}{16}\left( -320 -32\zeta(2) +128\zeta(3) +172\zeta(4) +24\zeta(5) +74\zeta(6)  \right. \nonumber \\ &\left. \hspace{1em}
-48\zeta(3)^2 +119\zeta(7) +32\zeta(2)\zeta(5) -132\zeta(3)\zeta(4)\right) \label{eq07064}
\end{align}
 
\begin{align}
\sum_{k=1}^\infty \frac{H(k)^{4}}{k^{3}} \sumend &= \frac{-1}{8}\left( -185\zeta(7) -40\zeta(2)\zeta(5) +172\zeta(3)\zeta(4)\right) \label{eq07065} \\
\sum_{k=1}^\infty \frac{H(k)^{4}}{k^{2}(k+1)} \sumend &= \frac{1}{24}\left( -720\zeta(5) -144\zeta(2)\zeta(3) +979\zeta(6) +72\zeta(3)^2\right) \label{eq07066} \\
\sum_{k=1}^\infty \frac{H(k)^{4}}{k(k+1)^{2}} \sumend &= \frac{1}{24}\left(  720\zeta(5) +144\zeta(2)\zeta(3) -859\zeta(6) -72\zeta(3)^2\right) \label{eq07067} \\
\sum_{k=1}^\infty \frac{H(k)^{4}}{(k+1)^{3}} \sumend &= \frac{1}{8}\left( -109\zeta(7) -40\zeta(2)\zeta(5) +148\zeta(3)\zeta(4)\right) \label{eq07068} \\
\sum_{k=1}^\infty \frac{H(k)^{4}}{k^{2}(k+2)} \sumend &= \frac{-1}{48}\left(  12 +36\zeta(2) +132\zeta(3) +222\zeta(4) +360\zeta(5) +72\zeta(2)\zeta(3)  \right. \nonumber \\ &\left. \hspace{1em}
-979\zeta(6) -72\zeta(3)^2\right) \label{eq07069} \\
\sum_{k=1}^\infty \frac{H(k)^{4}}{k(k+1)(k+2)} \sumend &= \frac{-1}{4}\left(  2 +6\zeta(2) +22\zeta(3) +37\zeta(4) -60\zeta(5)  \right. \nonumber \\ &\left. \hspace{1em}
-12\zeta(2)\zeta(3)\right) \label{eq07070} \\
\sum_{k=1}^\infty \frac{H(k)^{4}}{(k+1)^{2}(k+2)} \sumend &= \frac{1}{24}\left( -24 -72\zeta(2) -264\zeta(3) -444\zeta(4) +859\zeta(6)  \right. \nonumber \\ &\left. \hspace{1em}
+72\zeta(3)^2\right) \label{eq07071} \\
\sum_{k=1}^\infty \frac{H(k)^{4}}{k(k+2)^{2}} \sumend &= \frac{1}{48}\left(  132 +228\zeta(2) +564\zeta(3) +270\zeta(4) -360\zeta(5) -24\zeta(2)\zeta(3)  \right. \nonumber \\ &\left. \hspace{1em}
-859\zeta(6) -72\zeta(3)^2\right) \label{eq07072} \\
\sum_{k=1}^\infty \frac{H(k)^{4}}{(k+1)(k+2)^{2}} \sumend &= \frac{1}{24}\left(  144 +264\zeta(2) +696\zeta(3) +492\zeta(4) -720\zeta(5)  \right. \nonumber \\ &\left. \hspace{1em}
-96\zeta(2)\zeta(3) -859\zeta(6) -72\zeta(3)^2\right) \label{eq07073} \\
\sum_{k=1}^\infty \frac{H(k)^{4}}{(k+2)^{3}} \sumend &= \frac{-1}{8}\left(  120 +112\zeta(2) +160\zeta(3) -124\zeta(4) -168\zeta(5) -80\zeta(2)\zeta(3)  \right. \nonumber \\ &\left. \hspace{1em}
+66\zeta(6) -64\zeta(3)^2 +109\zeta(7) +40\zeta(2)\zeta(5) -148\zeta(3)\zeta(4)\right) \label{eq07074} \\
\sum_{k=1}^\infty \frac{H(k)^{5}}{k^{2}} \sumend &= \frac{1}{16}\left(  2051\zeta(7) +456\zeta(2)\zeta(5) +528\zeta(3)\zeta(4)\right) \label{eq07075} \\
\sum_{k=1}^\infty \frac{H(k)^{5}}{k(k+1)} \sumend &= \frac{-1}{2}\left( -357\zeta(6) -45\zeta(3)^2\right) \label{eq07076} \\
\sum_{k=1}^\infty \frac{H(k)^{5}}{(k+1)^{2}} \sumend &= \frac{-1}{16}\left( -1855\zeta(7) -456\zeta(2)\zeta(5) -528\zeta(3)\zeta(4)\right) \label{eq07077}
\end{align}
 
\begin{align}
\sum_{k=1}^\infty \frac{H(k)^{5}}{k(k+2)} \sumend &= \frac{-1}{8}\left( -4 -16\zeta(2) -84\zeta(3) -251\zeta(4) -284\zeta(5) -60\zeta(2)\zeta(3)  \right. \nonumber \\ &\left. \hspace{1em}
-714\zeta(6) -90\zeta(3)^2\right) \label{eq07078} \\
\sum_{k=1}^\infty \frac{H(k)^{5}}{(k+1)(k+2)} \sumend &= \frac{1}{4}\left(  4 +16\zeta(2) +84\zeta(3) +251\zeta(4) +284\zeta(5)  \right. \nonumber \\ &\left. \hspace{1em}
+60\zeta(2)\zeta(3)\right) \label{eq07079} \\
\sum_{k=1}^\infty \frac{H(k)^{5}}{(k+2)^{2}} \sumend &= \frac{-1}{48}\left(  288 +720\zeta(2) +2784\zeta(3) +4704\zeta(4) -192\zeta(5)  \right. \nonumber \\ &\left. \hspace{1em}
+240\zeta(2)\zeta(3) -8590\zeta(6) -720\zeta(3)^2 -5565\zeta(7) -1368\zeta(2)\zeta(5)  \right. \nonumber \\ &\left. \hspace{1em}
-1584\zeta(3)\zeta(4)\right) \label{eq07080}
\end{align}

\newpage
Formulas for order $r = m + n + p + q = 8$:
\begin{align}
\sum_{k=1}^\infty \frac{H(k)}{k^{7}} \sumend &= \frac{1}{4}\left(  9\zeta(8) -4\zeta(3)\zeta(5)\right) \label{eq08001} \\
\sum_{k=1}^\infty \frac{H(k)}{k^{6}(k+1)} \sumend &= \frac{1}{4}\left( -4\zeta(2) +8\zeta(3) -5\zeta(4) +12\zeta(5) -4\zeta(2)\zeta(3) -7\zeta(6)  \right. \nonumber \\ &\left. \hspace{1em}
+2\zeta(3)^2 +16\zeta(7) -4\zeta(2)\zeta(5) -4\zeta(3)\zeta(4)\right) \label{eq08002} \\
\sum_{k=1}^\infty \frac{H(k)}{k^{5}(k+1)^{2}} \sumend &= \frac{1}{4}\left(  20\zeta(2) -36\zeta(3) +15\zeta(4) -24\zeta(5) +8\zeta(2)\zeta(3) +7\zeta(6)  \right. \nonumber \\ &\left. \hspace{1em}
-2\zeta(3)^2\right) \label{eq08003} \\
\sum_{k=1}^\infty \frac{H(k)}{k^{4}(k+1)^{3}} \sumend &= \frac{1}{2}\left( -20\zeta(2) +32\zeta(3) -7\zeta(4) +6\zeta(5) -2\zeta(2)\zeta(3)\right) \label{eq08004} \\
\sum_{k=1}^\infty \frac{H(k)}{k^{3}(k+1)^{4}} \sumend &= \frac{1}{2}\left(  20\zeta(2) -28\zeta(3) +\zeta(4) -4\zeta(5) +2\zeta(2)\zeta(3)\right) \label{eq08005} \\
\sum_{k=1}^\infty \frac{H(k)}{k^{2}(k+1)^{5}} \sumend &= \frac{1}{4}\left( -20\zeta(2) +24\zeta(3) +3\zeta(4) +16\zeta(5) -8\zeta(2)\zeta(3) +3\zeta(6)  \right. \nonumber \\ &\left. \hspace{1em}
-2\zeta(3)^2\right) \label{eq08006} \\
\sum_{k=1}^\infty \frac{H(k)}{k(k+1)^{6}} \sumend &= \frac{-1}{4}\left( -4\zeta(2) +4\zeta(3) +\zeta(4) +8\zeta(5) -4\zeta(2)\zeta(3) +3\zeta(6)  \right. \nonumber \\ &\left. \hspace{1em}
-2\zeta(3)^2 +12\zeta(7) -4\zeta(2)\zeta(5) -4\zeta(3)\zeta(4)\right) \label{eq08007} \\
\sum_{k=1}^\infty \frac{H(k)}{(k+1)^{7}} \sumend &= \frac{-1}{4}\left( -5\zeta(8) +4\zeta(3)\zeta(5)\right) \label{eq08008} \\
\sum_{k=1}^\infty \frac{H(k)}{k^{6}(k+2)} \sumend &= \frac{-1}{64}\left(  1 +\zeta(2) -4\zeta(3) +5\zeta(4) -24\zeta(5) +8\zeta(2)\zeta(3) +28\zeta(6)  \right. \nonumber \\ &\left. \hspace{1em}
-8\zeta(3)^2 -128\zeta(7) +32\zeta(2)\zeta(5) +32\zeta(3)\zeta(4)\right) \label{eq08009} \\
\sum_{k=1}^\infty \frac{H(k)}{k^{5}(k+1)(k+2)} \sumend &= \frac{-1}{32}\left(  1 -31\zeta(2) +60\zeta(3) -35\zeta(4) +72\zeta(5) -24\zeta(2)\zeta(3)  \right. \nonumber \\ &\left. \hspace{1em}
-28\zeta(6) +8\zeta(3)^2\right) \label{eq08010} \\
\sum_{k=1}^\infty \frac{H(k)}{k^{4}(k+1)^{2}(k+2)} \sumend &= \frac{-1}{16}\left(  1 +49\zeta(2) -84\zeta(3) +25\zeta(4) -24\zeta(5)  \right. \nonumber \\ &\left. \hspace{1em}
+8\zeta(2)\zeta(3)\right) \label{eq08011} \\
\sum_{k=1}^\infty \frac{H(k)}{k^{3}(k+1)^{3}(k+2)} \sumend &= \frac{-1}{8}\left(  1 -31\zeta(2) +44\zeta(3) -3\zeta(4)\right) \label{eq08012} \\
\sum_{k=1}^\infty \frac{H(k)}{k^{2}(k+1)^{4}(k+2)} \sumend &= \frac{-1}{4}\left(  1 +9\zeta(2) -12\zeta(3) -\zeta(4) -8\zeta(5)  \right. \nonumber \\ &\left. \hspace{1em}
+4\zeta(2)\zeta(3)\right) \label{eq08013}
\end{align}
 
\begin{align}
\sum_{k=1}^\infty \frac{H(k)}{k(k+1)^{5}(k+2)} \sumend &= \frac{1}{4}\left( -2 +2\zeta(2) -\zeta(4) -3\zeta(6) +2\zeta(3)^2\right) \label{eq08014} \\
\sum_{k=1}^\infty \frac{H(k)}{(k+1)^{6}(k+2)} \sumend &= \frac{1}{4}\left( -4 +4\zeta(3) -\zeta(4) +8\zeta(5) -4\zeta(2)\zeta(3) -3\zeta(6) +2\zeta(3)^2  \right. \nonumber \\ &\left. \hspace{1em}
+12\zeta(7) -4\zeta(2)\zeta(5) -4\zeta(3)\zeta(4)\right) \label{eq08015} \\
\sum_{k=1}^\infty \frac{H(k)}{k^{5}(k+2)^{2}} \sumend &= \frac{-1}{64}\left( -9 -3\zeta(2) +18\zeta(3) -15\zeta(4) +48\zeta(5) -16\zeta(2)\zeta(3)  \right. \nonumber \\ &\left. \hspace{1em}
-28\zeta(6) +8\zeta(3)^2\right) \label{eq08016} \\
\sum_{k=1}^\infty \frac{H(k)}{k^{4}(k+1)(k+2)^{2}} \sumend &= \frac{-1}{16}\left( -5 +14\zeta(2) -21\zeta(3) +10\zeta(4) -12\zeta(5)  \right. \nonumber \\ &\left. \hspace{1em}
+4\zeta(2)\zeta(3)\right) \label{eq08017} \\
\sum_{k=1}^\infty \frac{H(k)}{k^{3}(k+1)^{2}(k+2)^{2}} \sumend &= \frac{1}{16}\left(  11 +21\zeta(2) -42\zeta(3) +5\zeta(4)\right) \label{eq08018} \\
\sum_{k=1}^\infty \frac{H(k)}{k^{2}(k+1)^{3}(k+2)^{2}} \sumend &= \frac{1}{4}\left(  6 -5\zeta(2) +\zeta(3) +\zeta(4)\right) \label{eq08019} \\
\sum_{k=1}^\infty \frac{H(k)}{k(k+1)^{4}(k+2)^{2}} \sumend &= \frac{1}{4}\left(  13 -\zeta(2) -10\zeta(3) +\zeta(4) -8\zeta(5) +4\zeta(2)\zeta(3)\right) \label{eq08020} \\
\sum_{k=1}^\infty \frac{H(k)}{(k+1)^{5}(k+2)^{2}} \sumend &= \frac{-1}{4}\left( -28 +4\zeta(2) +20\zeta(3) -3\zeta(4) +16\zeta(5) -8\zeta(2)\zeta(3)  \right. \nonumber \\ &\left. \hspace{1em}
-3\zeta(6) +2\zeta(3)^2\right) \label{eq08021} \\
\sum_{k=1}^\infty \frac{H(k)}{k^{4}(k+2)^{3}} \sumend &= \frac{1}{32}\left( -19 +\zeta(2) +18\zeta(3) -7\zeta(4) +12\zeta(5) -4\zeta(2)\zeta(3)\right) \label{eq08022} \\
\sum_{k=1}^\infty \frac{H(k)}{k^{3}(k+1)(k+2)^{3}} \sumend &= \frac{-1}{16}\left(  24 -15\zeta(2) +3\zeta(3) -3\zeta(4)\right) \label{eq08023} \\
\sum_{k=1}^\infty \frac{H(k)}{k^{2}(k+1)^{2}(k+2)^{3}} \sumend &= \frac{-1}{16}\left(  59 -9\zeta(2) -36\zeta(3) -\zeta(4)\right) \label{eq08024} \\
\sum_{k=1}^\infty \frac{H(k)}{k(k+1)^{3}(k+2)^{3}} \sumend &= \frac{-1}{8}\left(  71 -19\zeta(2) -34\zeta(3) +\zeta(4)\right) \label{eq08025} \\
\sum_{k=1}^\infty \frac{H(k)}{(k+1)^{4}(k+2)^{3}} \sumend &= \frac{1}{2}\left( -42 +10\zeta(2) +22\zeta(3) -\zeta(4) +4\zeta(5)  \right. \nonumber \\ &\left. \hspace{1em}
-2\zeta(2)\zeta(3)\right) \label{eq08026} \\
\sum_{k=1}^\infty \frac{H(k)}{k^{3}(k+2)^{4}} \sumend &= \frac{1}{32}\left(  51 -11\zeta(2) -24\zeta(3) -3\zeta(4) -8\zeta(5) +4\zeta(2)\zeta(3)\right) \label{eq08027}
\end{align}
 
\begin{align}
\sum_{k=1}^\infty \frac{H(k)}{k^{2}(k+1)(k+2)^{4}} \sumend &= \frac{-1}{16}\left( -75 +26\zeta(2) +21\zeta(3) +6\zeta(4) +8\zeta(5)  \right. \nonumber \\ &\left. \hspace{1em}
-4\zeta(2)\zeta(3)\right) \label{eq08028} \\
\sum_{k=1}^\infty \frac{H(k)}{k(k+1)^{2}(k+2)^{4}} \sumend &= \frac{1}{16}\left(  209 -61\zeta(2) -78\zeta(3) -13\zeta(4) -16\zeta(5)  \right. \nonumber \\ &\left. \hspace{1em}
+8\zeta(2)\zeta(3)\right) \label{eq08029} \\
\sum_{k=1}^\infty \frac{H(k)}{(k+1)^{3}(k+2)^{4}} \sumend &= \frac{1}{2}\left(  70 -20\zeta(2) -28\zeta(3) -3\zeta(4) -4\zeta(5)  \right. \nonumber \\ &\left. \hspace{1em}
+2\zeta(2)\zeta(3)\right) \label{eq08030} \\
\sum_{k=1}^\infty \frac{H(k)}{k^{2}(k+2)^{5}} \sumend &= \frac{-1}{64}\left(  201 -47\zeta(2) -56\zeta(3) -35\zeta(4) -48\zeta(5) +16\zeta(2)\zeta(3)  \right. \nonumber \\ &\left. \hspace{1em}
-12\zeta(6) +8\zeta(3)^2\right) \label{eq08031} \\
\sum_{k=1}^\infty \frac{H(k)}{k(k+1)(k+2)^{5}} \sumend &= \frac{1}{32}\left( -351 +99\zeta(2) +98\zeta(3) +47\zeta(4) +64\zeta(5) -24\zeta(2)\zeta(3)  \right. \nonumber \\ &\left. \hspace{1em}
+12\zeta(6) -8\zeta(3)^2\right) \label{eq08032} \\
\sum_{k=1}^\infty \frac{H(k)}{(k+1)^{2}(k+2)^{5}} \sumend &= \frac{1}{4}\left( -140 +40\zeta(2) +44\zeta(3) +15\zeta(4) +20\zeta(5) -8\zeta(2)\zeta(3)  \right. \nonumber \\ &\left. \hspace{1em}
+3\zeta(6) -2\zeta(3)^2\right) \label{eq08033} \\
\sum_{k=1}^\infty \frac{H(k)}{k(k+2)^{6}} \sumend &= \frac{1}{64}\left(  321 -61\zeta(2) -62\zeta(3) -57\zeta(4) -64\zeta(5) +8\zeta(2)\zeta(3)  \right. \nonumber \\ &\left. \hspace{1em}
-44\zeta(6) +8\zeta(3)^2 -96\zeta(7) +32\zeta(2)\zeta(5) +32\zeta(3)\zeta(4)\right) \label{eq08034} \\
\sum_{k=1}^\infty \frac{H(k)}{(k+1)(k+2)^{6}} \sumend &= \frac{1}{4}\left(  84 -20\zeta(2) -20\zeta(3) -13\zeta(4) -16\zeta(5) +4\zeta(2)\zeta(3)  \right. \nonumber \\ &\left. \hspace{1em}
-7\zeta(6) +2\zeta(3)^2 -12\zeta(7) +4\zeta(2)\zeta(5) +4\zeta(3)\zeta(4)\right) \label{eq08035} \\
\sum_{k=1}^\infty \frac{H(k)}{(k+2)^{7}} \sumend &= \frac{1}{4}\left( -28 +4\zeta(2) +4\zeta(3) +4\zeta(4) +4\zeta(5) +4\zeta(6) +4\zeta(7) +5\zeta(8)  \right. \nonumber \\ &\left. \hspace{1em}
-4\zeta(3)\zeta(5)\right) \label{eq08036} \\
\sum_{k=1}^\infty \frac{H(k)^{2}}{k^{6}} \sumend &= -\left( - M(2,6)\right) \label{eq08037} \\
\sum_{k=1}^\infty \frac{H(k)^{2}}{k^{5}(k+1)} \sumend &= \frac{-1}{24}\left( -72\zeta(3) +102\zeta(4) -84\zeta(5) +24\zeta(2)\zeta(3) +97\zeta(6)  \right. \nonumber \\ &\left. \hspace{1em}
-48\zeta(3)^2 -144\zeta(7) +24\zeta(2)\zeta(5) +60\zeta(3)\zeta(4)\right) \label{eq08038}
\end{align}
 
\begin{align}
\sum_{k=1}^\infty \frac{H(k)^{2}}{k^{4}(k+1)^{2}} \sumend &= \frac{-1}{24}\left(  288\zeta(3) -372\zeta(4) +168\zeta(5) -48\zeta(2)\zeta(3) -97\zeta(6)  \right. \nonumber \\ &\left. \hspace{1em}
+48\zeta(3)^2\right) \label{eq08039} \\
\sum_{k=1}^\infty \frac{H(k)^{2}}{k^{3}(k+1)^{3}} \sumend &= \left(  18\zeta(3) -21\zeta(4) +5\zeta(5) -2\zeta(2)\zeta(3)\right) \label{eq08040} \\
\sum_{k=1}^\infty \frac{H(k)^{2}}{k^{2}(k+1)^{4}} \sumend &= \frac{-1}{24}\left(  288\zeta(3) -300\zeta(4) +72\zeta(5) -48\zeta(2)\zeta(3) -37\zeta(6)  \right. \nonumber \\ &\left. \hspace{1em}
+24\zeta(3)^2\right) \label{eq08041} \\
\sum_{k=1}^\infty \frac{H(k)^{2}}{k(k+1)^{5}} \sumend &= \frac{1}{24}\left(  72\zeta(3) -66\zeta(4) +36\zeta(5) -24\zeta(2)\zeta(3) -37\zeta(6)  \right. \nonumber \\ &\left. \hspace{1em}
+24\zeta(3)^2 +24\zeta(7) -24\zeta(2)\zeta(5) +12\zeta(3)\zeta(4)\right) \label{eq08042} \\
\sum_{k=1}^\infty \frac{H(k)^{2}}{(k+1)^{6}} \sumend &= \frac{1}{2}\left( -7\zeta(8) +4\zeta(3)\zeta(5) +2 M(2,6)\right) \label{eq08043} \\
\sum_{k=1}^\infty \frac{H(k)^{2}}{k^{5}(k+2)} \sumend &= \frac{1}{192}\left(  6 +6\zeta(2) +18\zeta(3) -51\zeta(4) +84\zeta(5) -24\zeta(2)\zeta(3)  \right. \nonumber \\ &\left. \hspace{1em}
-194\zeta(6) +96\zeta(3)^2 +576\zeta(7) -96\zeta(2)\zeta(5) -240\zeta(3)\zeta(4)\right) \label{eq08044} \\
\sum_{k=1}^\infty \frac{H(k)^{2}}{k^{4}(k+1)(k+2)} \sumend &= \frac{-1}{96}\left( -6 -6\zeta(2) +270\zeta(3) -357\zeta(4) +252\zeta(5)  \right. \nonumber \\ &\left. \hspace{1em}
-72\zeta(2)\zeta(3) -194\zeta(6) +96\zeta(3)^2\right) \label{eq08045} \\
\sum_{k=1}^\infty \frac{H(k)^{2}}{k^{3}(k+1)^{2}(k+2)} \sumend &= \frac{-1}{16}\left( -2 -2\zeta(2) -102\zeta(3) +129\zeta(4) -28\zeta(5)  \right. \nonumber \\ &\left. \hspace{1em}
+8\zeta(2)\zeta(3)\right) \label{eq08046} \\
\sum_{k=1}^\infty \frac{H(k)^{2}}{k^{2}(k+1)^{3}(k+2)} \sumend &= \frac{1}{8}\left(  2 +2\zeta(2) -42\zeta(3) +39\zeta(4) -12\zeta(5)  \right. \nonumber \\ &\left. \hspace{1em}
+8\zeta(2)\zeta(3)\right) \label{eq08047} \\
\sum_{k=1}^\infty \frac{H(k)^{2}}{k(k+1)^{4}(k+2)} \sumend &= \frac{-1}{24}\left( -12 -12\zeta(2) -36\zeta(3) +66\zeta(4) +37\zeta(6)  \right. \nonumber \\ &\left. \hspace{1em}
-24\zeta(3)^2\right) \label{eq08048} \\
\sum_{k=1}^\infty \frac{H(k)^{2}}{(k+1)^{5}(k+2)} \sumend &= \frac{-1}{24}\left( -24 -24\zeta(2) +66\zeta(4) +36\zeta(5) -24\zeta(2)\zeta(3) +37\zeta(6)  \right. \nonumber \\ &\left. \hspace{1em}
-24\zeta(3)^2 +24\zeta(7) -24\zeta(2)\zeta(5) +12\zeta(3)\zeta(4)\right) \label{eq08049} \\
\sum_{k=1}^\infty \frac{H(k)^{2}}{k^{4}(k+2)^{2}} \sumend &= \frac{-1}{96}\left(  30 +12\zeta(2) +24\zeta(3) -93\zeta(4) +84\zeta(5) -24\zeta(2)\zeta(3)  \right. \nonumber \\ &\left. \hspace{1em}
-97\zeta(6) +48\zeta(3)^2\right) \label{eq08050}
\end{align}
 
\begin{align}
\sum_{k=1}^\infty \frac{H(k)^{2}}{k^{3}(k+1)(k+2)^{2}} \sumend &= \frac{-1}{32}\left(  22 +10\zeta(2) -74\zeta(3) +57\zeta(4) -28\zeta(5)  \right. \nonumber \\ &\left. \hspace{1em}
+8\zeta(2)\zeta(3)\right) \label{eq08051} \\
\sum_{k=1}^\infty \frac{H(k)^{2}}{k^{2}(k+1)^{2}(k+2)^{2}} \sumend &= \frac{-1}{4}\left(  6 +3\zeta(2) +7\zeta(3) -18\zeta(4)\right) \label{eq08052} \\
\sum_{k=1}^\infty \frac{H(k)^{2}}{k(k+1)^{3}(k+2)^{2}} \sumend &= \frac{1}{8}\left( -26 -14\zeta(2) +14\zeta(3) +33\zeta(4) +12\zeta(5)  \right. \nonumber \\ &\left. \hspace{1em}
-8\zeta(2)\zeta(3)\right) \label{eq08053} \\
\sum_{k=1}^\infty \frac{H(k)^{2}}{(k+1)^{4}(k+2)^{2}} \sumend &= \frac{-1}{24}\left(  168 +96\zeta(2) -48\zeta(3) -264\zeta(4) -72\zeta(5)  \right. \nonumber \\ &\left. \hspace{1em}
+48\zeta(2)\zeta(3) -37\zeta(6) +24\zeta(3)^2\right) \label{eq08054} \\
\sum_{k=1}^\infty \frac{H(k)^{2}}{k^{3}(k+2)^{3}} \sumend &= \frac{-1}{16}\left( -24 -\zeta(2) +3\zeta(3) +22\zeta(4) -10\zeta(5)  \right. \nonumber \\ &\left. \hspace{1em}
+4\zeta(2)\zeta(3)\right) \label{eq08055} \\
\sum_{k=1}^\infty \frac{H(k)^{2}}{k^{2}(k+1)(k+2)^{3}} \sumend &= \frac{-1}{32}\left( -118 -14\zeta(2) +86\zeta(3) +31\zeta(4) -12\zeta(5)  \right. \nonumber \\ &\left. \hspace{1em}
+8\zeta(2)\zeta(3)\right) \label{eq08056} \\
\sum_{k=1}^\infty \frac{H(k)^{2}}{k(k+1)^{2}(k+2)^{3}} \sumend &= \frac{1}{16}\left(  142 +26\zeta(2) -58\zeta(3) -103\zeta(4) +12\zeta(5)  \right. \nonumber \\ &\left. \hspace{1em}
-8\zeta(2)\zeta(3)\right) \label{eq08057} \\
\sum_{k=1}^\infty \frac{H(k)^{2}}{(k+1)^{3}(k+2)^{3}} \sumend &= \left(  21 +5\zeta(2) -9\zeta(3) -17\zeta(4)\right) \label{eq08058} \\
\sum_{k=1}^\infty \frac{H(k)^{2}}{k^{2}(k+2)^{4}} \sumend &= \frac{1}{96}\left( -450 +60\zeta(2) +168\zeta(3) +123\zeta(4) +60\zeta(5)  \right. \nonumber \\ &\left. \hspace{1em}
-24\zeta(2)\zeta(3) +37\zeta(6) -24\zeta(3)^2\right) \label{eq08059} \\
\sum_{k=1}^\infty \frac{H(k)^{2}}{k(k+1)(k+2)^{4}} \sumend &= \frac{1}{96}\left( -1254 +78\zeta(2) +594\zeta(3) +339\zeta(4) +84\zeta(5)  \right. \nonumber \\ &\left. \hspace{1em}
-24\zeta(2)\zeta(3) +74\zeta(6) -48\zeta(3)^2\right) \label{eq08060} \\
\sum_{k=1}^\infty \frac{H(k)^{2}}{(k+1)^{2}(k+2)^{4}} \sumend &= \frac{-1}{24}\left(  840 -384\zeta(3) -324\zeta(4) -24\zeta(5) -37\zeta(6)  \right. \nonumber \\ &\left. \hspace{1em}
+24\zeta(3)^2\right) \label{eq08061} \\
\sum_{k=1}^\infty \frac{H(k)^{2}}{k(k+2)^{5}} \sumend &= \frac{1}{192}\left(  2106 -402\zeta(2) -750\zeta(3) -357\zeta(4) -636\zeta(5)  \right. \nonumber \\ &\left. \hspace{1em}
+264\zeta(2)\zeta(3) -218\zeta(6) +144\zeta(3)^2 +96\zeta(7) -96\zeta(2)\zeta(5)  \right. \nonumber \\ &\left. \hspace{1em}
+48\zeta(3)\zeta(4)\right) \label{eq08062}
\end{align}
 
\begin{align}
\sum_{k=1}^\infty \frac{H(k)^{2}}{(k+1)(k+2)^{5}} \sumend &= \frac{1}{24}\left(  840 -120\zeta(2) -336\zeta(3) -174\zeta(4) -180\zeta(5)  \right. \nonumber \\ &\left. \hspace{1em}
+72\zeta(2)\zeta(3) -73\zeta(6) +48\zeta(3)^2 +24\zeta(7) -24\zeta(2)\zeta(5)  \right. \nonumber \\ &\left. \hspace{1em}
+12\zeta(3)\zeta(4)\right) \label{eq08063} \\
\sum_{k=1}^\infty \frac{H(k)^{2}}{(k+2)^{6}} \sumend &= \frac{1}{2}\left( -42 +8\zeta(2) +12\zeta(3) +7\zeta(4) +12\zeta(5) -4\zeta(2)\zeta(3) +5\zeta(6)  \right. \nonumber \\ &\left. \hspace{1em}
-2\zeta(3)^2 +12\zeta(7) -4\zeta(2)\zeta(5) -4\zeta(3)\zeta(4) -7\zeta(8) +4\zeta(3)\zeta(5)  \right. \nonumber \\ &\left. \hspace{1em}
+2 M(2,6)\right) \label{eq08064} \\
\sum_{k=1}^\infty \frac{H(k)^{3}}{k^{5}} \sumend &= \frac{1}{96}\left( -595\zeta(8) -120\zeta(2)\zeta(3)^2 +576\zeta(3)\zeta(5) +264 M(2,6)\right) \label{eq08065} \\
\sum_{k=1}^\infty \frac{H(k)^{3}}{k^{4}(k+1)} \sumend &= \frac{-1}{16}\left(  160\zeta(4) -160\zeta(5) -16\zeta(2)\zeta(3) +93\zeta(6) -40\zeta(3)^2  \right. \nonumber \\ &\left. \hspace{1em}
-231\zeta(7) -32\zeta(2)\zeta(5) +204\zeta(3)\zeta(4)\right) \label{eq08066} \\
\sum_{k=1}^\infty \frac{H(k)^{3}}{k^{3}(k+1)^{2}} \sumend &= \frac{1}{16}\left(  480\zeta(4) -440\zeta(5) -48\zeta(2)\zeta(3) +93\zeta(6)  \right. \nonumber \\ &\left. \hspace{1em}
-40\zeta(3)^2\right) \label{eq08067} \\
\sum_{k=1}^\infty \frac{H(k)^{3}}{k^{2}(k+1)^{3}} \sumend &= \frac{-1}{16}\left(  480\zeta(4) -400\zeta(5) -48\zeta(2)\zeta(3) +33\zeta(6)  \right. \nonumber \\ &\left. \hspace{1em}
-32\zeta(3)^2\right) \label{eq08068} \\
\sum_{k=1}^\infty \frac{H(k)^{3}}{k(k+1)^{4}} \sumend &= \frac{-1}{16}\left( -160\zeta(4) +120\zeta(5) +16\zeta(2)\zeta(3) -33\zeta(6) +32\zeta(3)^2  \right. \nonumber \\ &\left. \hspace{1em}
+119\zeta(7) +32\zeta(2)\zeta(5) -132\zeta(3)\zeta(4)\right) \label{eq08069} \\
\sum_{k=1}^\infty \frac{H(k)^{3}}{(k+1)^{5}} \sumend &= \frac{1}{96}\left( -43\zeta(8) -120\zeta(2)\zeta(3)^2 +288\zeta(3)\zeta(5) -24 M(2,6)\right) \label{eq08070} \\
\sum_{k=1}^\infty \frac{H(k)^{3}}{k^{4}(k+2)} \sumend &= \frac{-1}{64}\left(  4 +8\zeta(2) +16\zeta(3) +40\zeta(4) -80\zeta(5) -8\zeta(2)\zeta(3)  \right. \nonumber \\ &\left. \hspace{1em}
+93\zeta(6) -40\zeta(3)^2 -462\zeta(7) -64\zeta(2)\zeta(5) +408\zeta(3)\zeta(4)\right) \label{eq08071} \\
\sum_{k=1}^\infty \frac{H(k)^{3}}{k^{3}(k+1)(k+2)} \sumend &= \frac{-1}{32}\left(  4 +8\zeta(2) +16\zeta(3) -280\zeta(4) +240\zeta(5) +24\zeta(2)\zeta(3)  \right. \nonumber \\ &\left. \hspace{1em}
-93\zeta(6) +40\zeta(3)^2\right) \label{eq08072} \\
\sum_{k=1}^\infty \frac{H(k)^{3}}{k^{2}(k+1)^{2}(k+2)} \sumend &= \frac{1}{4}\left( -1 -2\zeta(2) -4\zeta(3) -50\zeta(4) +50\zeta(5)  \right. \nonumber \\ &\left. \hspace{1em}
+6\zeta(2)\zeta(3)\right) \label{eq08073}
\end{align}
 
\begin{align}
\sum_{k=1}^\infty \frac{H(k)^{3}}{k(k+1)^{3}(k+2)} \sumend &= \frac{1}{16}\left( -8 -16\zeta(2) -32\zeta(3) +80\zeta(4) +33\zeta(6) -32\zeta(3)^2\right) \label{eq08074} \\
\sum_{k=1}^\infty \frac{H(k)^{3}}{(k+1)^{4}(k+2)} \sumend &= \frac{-1}{16}\left(  16 +32\zeta(2) +64\zeta(3) -120\zeta(5) -16\zeta(2)\zeta(3) -33\zeta(6)  \right. \nonumber \\ &\left. \hspace{1em}
+32\zeta(3)^2 -119\zeta(7) -32\zeta(2)\zeta(5) +132\zeta(3)\zeta(4)\right) \label{eq08075} \\
\sum_{k=1}^\infty \frac{H(k)^{3}}{k^{3}(k+2)^{2}} \sumend &= \frac{-1}{64}\left( -44 -48\zeta(2) -56\zeta(3) -54\zeta(4) +220\zeta(5) +24\zeta(2)\zeta(3)  \right. \nonumber \\ &\left. \hspace{1em}
-93\zeta(6) +40\zeta(3)^2\right) \label{eq08076} \\
\sum_{k=1}^\infty \frac{H(k)^{3}}{k^{2}(k+1)(k+2)^{2}} \sumend &= \frac{-1}{16}\left( -24 -28\zeta(2) -36\zeta(3) +113\zeta(4) -10\zeta(5)\right) \label{eq08077} \\
\sum_{k=1}^\infty \frac{H(k)^{3}}{k(k+1)^{2}(k+2)^{2}} \sumend &= \frac{-1}{8}\left( -26 -32\zeta(2) -44\zeta(3) +13\zeta(4) +90\zeta(5)  \right. \nonumber \\ &\left. \hspace{1em}
+12\zeta(2)\zeta(3)\right) \label{eq08078} \\
\sum_{k=1}^\infty \frac{H(k)^{3}}{(k+1)^{3}(k+2)^{2}} \sumend &= \frac{1}{16}\left(  112 +144\zeta(2) +208\zeta(3) -132\zeta(4) -360\zeta(5)  \right. \nonumber \\ &\left. \hspace{1em}
-48\zeta(2)\zeta(3) -33\zeta(6) +32\zeta(3)^2\right) \label{eq08079} \\
\sum_{k=1}^\infty \frac{H(k)^{3}}{k^{2}(k+2)^{3}} \sumend &= \frac{1}{64}\left( -236 -120\zeta(2) -16\zeta(3) +156\zeta(4) +128\zeta(5)  \right. \nonumber \\ &\left. \hspace{1em}
+72\zeta(2)\zeta(3) -33\zeta(6) +32\zeta(3)^2\right) \label{eq08080} \\
\sum_{k=1}^\infty \frac{H(k)^{3}}{k(k+1)(k+2)^{3}} \sumend &= \frac{-1}{32}\left(  284 +176\zeta(2) +88\zeta(3) -382\zeta(4) -108\zeta(5)  \right. \nonumber \\ &\left. \hspace{1em}
-72\zeta(2)\zeta(3) +33\zeta(6) -32\zeta(3)^2\right) \label{eq08081} \\
\sum_{k=1}^\infty \frac{H(k)^{3}}{(k+1)^{2}(k+2)^{3}} \sumend &= \frac{-1}{16}\left(  336 +240\zeta(2) +176\zeta(3) -408\zeta(4) -288\zeta(5)  \right. \nonumber \\ &\left. \hspace{1em}
-96\zeta(2)\zeta(3) +33\zeta(6) -32\zeta(3)^2\right) \label{eq08082} \\
\sum_{k=1}^\infty \frac{H(k)^{3}}{k(k+2)^{4}} \sumend &= \frac{-1}{64}\left( -836 -144\zeta(2) +280\zeta(3) +514\zeta(4) +36\zeta(5) +56\zeta(2)\zeta(3)  \right. \nonumber \\ &\left. \hspace{1em}
+115\zeta(6) -64\zeta(3)^2 +238\zeta(7) +64\zeta(2)\zeta(5) -264\zeta(3)\zeta(4)\right) \label{eq08083} \\
\sum_{k=1}^\infty \frac{H(k)^{3}}{(k+1)(k+2)^{4}} \sumend &= \frac{-1}{16}\left( -560 -160\zeta(2) +96\zeta(3) +448\zeta(4) +72\zeta(5)  \right. \nonumber \\ &\left. \hspace{1em}
+64\zeta(2)\zeta(3) +41\zeta(6) -16\zeta(3)^2 +119\zeta(7) +32\zeta(2)\zeta(5)  \right. \nonumber \\ &\left. \hspace{1em}
-132\zeta(3)\zeta(4)\right) \label{eq08084}
\end{align}
 
\begin{align}
\sum_{k=1}^\infty \frac{H(k)^{3}}{(k+2)^{5}} \sumend &= \frac{-1}{96}\left(  3360 -1344\zeta(3) -1296\zeta(4) -816\zeta(5) +288\zeta(2)\zeta(3)  \right. \nonumber \\ &\left. \hspace{1em}
-660\zeta(6) +432\zeta(3)^2 +288\zeta(7) -288\zeta(2)\zeta(5) +144\zeta(3)\zeta(4) +43\zeta(8)  \right. \nonumber \\ &\left. \hspace{1em}
+120\zeta(2)\zeta(3)^2 -288\zeta(3)\zeta(5) +24 M(2,6)\right) \label{eq08085} \\
\sum_{k=1}^\infty \frac{H(k)^{4}}{k^{4}} \sumend &= \frac{-1}{144}\left(  14833\zeta(8) +4032\zeta(2)\zeta(3)^2 -16704\zeta(3)\zeta(5)  \right. \nonumber \\ &\left. \hspace{1em}
-3744 M(2,6)\right) \label{eq08086} \\
\sum_{k=1}^\infty \frac{H(k)^{4}}{k^{3}(k+1)} \sumend &= \frac{1}{24}\left(  720\zeta(5) +144\zeta(2)\zeta(3) -979\zeta(6) -72\zeta(3)^2 +555\zeta(7)  \right. \nonumber \\ &\left. \hspace{1em}
+120\zeta(2)\zeta(5) -516\zeta(3)\zeta(4)\right) \label{eq08087} \\
\sum_{k=1}^\infty \frac{H(k)^{4}}{k^{2}(k+1)^{2}} \sumend &= \frac{1}{12}\left( -720\zeta(5) -144\zeta(2)\zeta(3) +919\zeta(6) +72\zeta(3)^2\right) \label{eq08088} \\
\sum_{k=1}^\infty \frac{H(k)^{4}}{k(k+1)^{3}} \sumend &= \frac{1}{24}\left(  720\zeta(5) +144\zeta(2)\zeta(3) -859\zeta(6) -72\zeta(3)^2 +327\zeta(7)  \right. \nonumber \\ &\left. \hspace{1em}
+120\zeta(2)\zeta(5) -444\zeta(3)\zeta(4)\right) \label{eq08089} \\
\sum_{k=1}^\infty \frac{H(k)^{4}}{(k+1)^{4}} \sumend &= \frac{1}{144}\left( -12415\zeta(8) -3312\zeta(2)\zeta(3)^2 +13824\zeta(3)\zeta(5)  \right. \nonumber \\ &\left. \hspace{1em}
+3024 M(2,6)\right) \label{eq08090} \\
\sum_{k=1}^\infty \frac{H(k)^{4}}{k^{3}(k+2)} \sumend &= \frac{1}{96}\left(  12 +36\zeta(2) +132\zeta(3) +222\zeta(4) +360\zeta(5) +72\zeta(2)\zeta(3)  \right. \nonumber \\ &\left. \hspace{1em}
-979\zeta(6) -72\zeta(3)^2 +1110\zeta(7) +240\zeta(2)\zeta(5) -1032\zeta(3)\zeta(4)\right) \label{eq08091} \\
\sum_{k=1}^\infty \frac{H(k)^{4}}{k^{2}(k+1)(k+2)} \sumend &= \frac{1}{48}\left(  12 +36\zeta(2) +132\zeta(3) +222\zeta(4) -1080\zeta(5)  \right. \nonumber \\ &\left. \hspace{1em}
-216\zeta(2)\zeta(3) +979\zeta(6) +72\zeta(3)^2\right) \label{eq08092} \\
\sum_{k=1}^\infty \frac{H(k)^{4}}{k(k+1)^{2}(k+2)} \sumend &= \frac{-1}{24}\left( -12 -36\zeta(2) -132\zeta(3) -222\zeta(4) -360\zeta(5)  \right. \nonumber \\ &\left. \hspace{1em}
-72\zeta(2)\zeta(3) +859\zeta(6) +72\zeta(3)^2\right) \label{eq08093} \\
\sum_{k=1}^\infty \frac{H(k)^{4}}{(k+1)^{3}(k+2)} \sumend &= \frac{-1}{24}\left( -24 -72\zeta(2) -264\zeta(3) -444\zeta(4) +859\zeta(6) +72\zeta(3)^2  \right. \nonumber \\ &\left. \hspace{1em}
+327\zeta(7) +120\zeta(2)\zeta(5) -444\zeta(3)\zeta(4)\right) \label{eq08094} \\
\sum_{k=1}^\infty \frac{H(k)^{4}}{k^{2}(k+2)^{2}} \sumend &= \frac{1}{48}\left( -72 -132\zeta(2) -348\zeta(3) -246\zeta(4) -24\zeta(2)\zeta(3)  \right. \nonumber \\ &\left. \hspace{1em}
+919\zeta(6) +72\zeta(3)^2\right) \label{eq08095}
\end{align}
 
\begin{align}
\sum_{k=1}^\infty \frac{H(k)^{4}}{k(k+1)(k+2)^{2}} \sumend &= \frac{-1}{48}\left(  156 +300\zeta(2) +828\zeta(3) +714\zeta(4) -1080\zeta(5)  \right. \nonumber \\ &\left. \hspace{1em}
-168\zeta(2)\zeta(3) -859\zeta(6) -72\zeta(3)^2\right) \label{eq08096} \\
\sum_{k=1}^\infty \frac{H(k)^{4}}{(k+1)^{2}(k+2)^{2}} \sumend &= \frac{1}{12}\left( -84 -168\zeta(2) -480\zeta(3) -468\zeta(4) +360\zeta(5)  \right. \nonumber \\ &\left. \hspace{1em}
+48\zeta(2)\zeta(3) +859\zeta(6) +72\zeta(3)^2\right) \label{eq08097} \\
\sum_{k=1}^\infty \frac{H(k)^{4}}{k(k+2)^{3}} \sumend &= \frac{1}{96}\left(  852 +900\zeta(2) +1524\zeta(3) -474\zeta(4) -1368\zeta(5)  \right. \nonumber \\ &\left. \hspace{1em}
-504\zeta(2)\zeta(3) -463\zeta(6) -456\zeta(3)^2 +654\zeta(7) +240\zeta(2)\zeta(5)  \right. \nonumber \\ &\left. \hspace{1em}
-888\zeta(3)\zeta(4)\right) \label{eq08098} \\
\sum_{k=1}^\infty \frac{H(k)^{4}}{(k+1)(k+2)^{3}} \sumend &= \frac{-1}{24}\left( -504 -600\zeta(2) -1176\zeta(3) -120\zeta(4) +1224\zeta(5)  \right. \nonumber \\ &\left. \hspace{1em}
+336\zeta(2)\zeta(3) +661\zeta(6) +264\zeta(3)^2 -327\zeta(7) -120\zeta(2)\zeta(5)  \right. \nonumber \\ &\left. \hspace{1em}
+444\zeta(3)\zeta(4)\right) \label{eq08099} \\
\sum_{k=1}^\infty \frac{H(k)^{4}}{(k+2)^{4}} \sumend &= \frac{-1}{144}\left(  5040 +2880\zeta(2) +2304\zeta(3) -5040\zeta(4) -2880\zeta(5)  \right. \nonumber \\ &\left. \hspace{1em}
-1728\zeta(2)\zeta(3) -144\zeta(6) -288\zeta(3)^2 -4284\zeta(7) -1152\zeta(2)\zeta(5)  \right. \nonumber \\ &\left. \hspace{1em}
+4752\zeta(3)\zeta(4) +12415\zeta(8) +3312\zeta(2)\zeta(3)^2 -13824\zeta(3)\zeta(5)  \right. \nonumber \\ &\left. \hspace{1em}
-3024 M(2,6)\right) \label{eq08100} \\
\sum_{k=1}^\infty \frac{H(k)^{5}}{k^{3}} \sumend &= \frac{1}{288}\left( -67811\zeta(8) -19080\zeta(2)\zeta(3)^2 +78768\zeta(3)\zeta(5)  \right. \nonumber \\ &\left. \hspace{1em}
+16920 M(2,6)\right) \label{eq08101} \\
\sum_{k=1}^\infty \frac{H(k)^{5}}{k^{2}(k+1)} \sumend &= \frac{-1}{16}\left(  2856\zeta(6) +360\zeta(3)^2 -2051\zeta(7) -456\zeta(2)\zeta(5)  \right. \nonumber \\ &\left. \hspace{1em}
-528\zeta(3)\zeta(4)\right) \label{eq08102} \\
\sum_{k=1}^\infty \frac{H(k)^{5}}{k(k+1)^{2}} \sumend &= \frac{-1}{16}\left( -2856\zeta(6) -360\zeta(3)^2 +1855\zeta(7) +456\zeta(2)\zeta(5)  \right. \nonumber \\ &\left. \hspace{1em}
+528\zeta(3)\zeta(4)\right) \label{eq08103} \\
\sum_{k=1}^\infty \frac{H(k)^{5}}{(k+1)^{3}} \sumend &= \frac{-1}{288}\left( -65621\zeta(8) -17640\zeta(2)\zeta(3)^2 +72432\zeta(3)\zeta(5)  \right. \nonumber \\ &\left. \hspace{1em}
+15480 M(2,6)\right) \label{eq08104} \\
\sum_{k=1}^\infty \frac{H(k)^{5}}{k^{2}(k+2)} \sumend &= \frac{-1}{32}\left(  8 +32\zeta(2) +168\zeta(3) +502\zeta(4) +568\zeta(5) +120\zeta(2)\zeta(3)  \right. \nonumber \\ &\left. \hspace{1em}
+1428\zeta(6) +180\zeta(3)^2 -2051\zeta(7) -456\zeta(2)\zeta(5) -528\zeta(3)\zeta(4)\right) \label{eq08105}
\end{align}
 
\begin{align}
\sum_{k=1}^\infty \frac{H(k)^{5}}{k(k+1)(k+2)} \sumend &= \frac{1}{8}\left( -4 -16\zeta(2) -84\zeta(3) -251\zeta(4) -284\zeta(5) -60\zeta(2)\zeta(3)  \right. \nonumber \\ &\left. \hspace{1em}
+714\zeta(6) +90\zeta(3)^2\right) \label{eq08106} \\
\sum_{k=1}^\infty \frac{H(k)^{5}}{(k+1)^{2}(k+2)} \sumend &= \frac{1}{16}\left( -16 -64\zeta(2) -336\zeta(3) -1004\zeta(4) -1136\zeta(5)  \right. \nonumber \\ &\left. \hspace{1em}
-240\zeta(2)\zeta(3) +1855\zeta(7) +456\zeta(2)\zeta(5) +528\zeta(3)\zeta(4)\right) \label{eq08107} \\
\sum_{k=1}^\infty \frac{H(k)^{5}}{k(k+2)^{2}} \sumend &= \frac{1}{96}\left(  312 +816\zeta(2) +3288\zeta(3) +6210\zeta(4) +1512\zeta(5)  \right. \nonumber \\ &\left. \hspace{1em}
+600\zeta(2)\zeta(3) -4306\zeta(6) -180\zeta(3)^2 -5565\zeta(7) -1368\zeta(2)\zeta(5)  \right. \nonumber \\ &\left. \hspace{1em}
-1584\zeta(3)\zeta(4)\right) \label{eq08108} \\
\sum_{k=1}^\infty \frac{H(k)^{5}}{(k+1)(k+2)^{2}} \sumend &= \frac{-1}{48}\left( -336 -912\zeta(2) -3792\zeta(3) -7716\zeta(4) -3216\zeta(5)  \right. \nonumber \\ &\left. \hspace{1em}
-960\zeta(2)\zeta(3) +8590\zeta(6) +720\zeta(3)^2 +5565\zeta(7) +1368\zeta(2)\zeta(5)  \right. \nonumber \\ &\left. \hspace{1em}
+1584\zeta(3)\zeta(4)\right) \label{eq08109} \\
\sum_{k=1}^\infty \frac{H(k)^{5}}{(k+2)^{3}} \sumend &= \frac{-1}{288}\left(  6048 +10080\zeta(2) +30240\zeta(3) +30096\zeta(4) -18432\zeta(5)  \right. \nonumber \\ &\left. \hspace{1em}
-4320\zeta(2)\zeta(3) -45600\zeta(6) -10080\zeta(3)^2 +19620\zeta(7) +7200\zeta(2)\zeta(5)  \right. \nonumber \\ &\left. \hspace{1em}
-26640\zeta(3)\zeta(4) -65621\zeta(8) -17640\zeta(2)\zeta(3)^2 +72432\zeta(3)\zeta(5)  \right. \nonumber \\ &\left. \hspace{1em}
+15480 M(2,6)\right) \label{eq08110} \\
\sum_{k=1}^\infty \frac{H(k)^{6}}{k^{2}} \sumend &= \frac{-1}{8}\left( -5843\zeta(8) +328\zeta(2)\zeta(3)^2 -3896\zeta(3)\zeta(5) -456 M(2,6)\right) \label{eq08111} \\
\sum_{k=1}^\infty \frac{H(k)^{6}}{k(k+1)} \sumend &= -\left( -644\zeta(7) -145\zeta(2)\zeta(5) -297\zeta(3)\zeta(4)\right) \label{eq08112} \\
\sum_{k=1}^\infty \frac{H(k)^{6}}{(k+1)^{2}} \sumend &= \frac{-1}{24}\left( -17027\zeta(8) +924\zeta(2)\zeta(3)^2 -11328\zeta(3)\zeta(5)  \right. \nonumber \\ &\left. \hspace{1em}
-1308 M(2,6)\right) \label{eq08113} \\
\sum_{k=1}^\infty \frac{H(k)^{6}}{k(k+2)} \sumend &= \frac{1}{8}\left(  4 +20\zeta(2) +136\zeta(3) +571\zeta(4) +1142\zeta(5) +244\zeta(2)\zeta(3)  \right. \nonumber \\ &\left. \hspace{1em}
+2097\zeta(6) +268\zeta(3)^2 +2576\zeta(7) +580\zeta(2)\zeta(5) +1188\zeta(3)\zeta(4)\right) \label{eq08114} \\
\sum_{k=1}^\infty \frac{H(k)^{6}}{(k+1)(k+2)} \sumend &= \frac{1}{4}\left(  4 +20\zeta(2) +136\zeta(3) +571\zeta(4) +1142\zeta(5) +244\zeta(2)\zeta(3)  \right. \nonumber \\ &\left. \hspace{1em}
+2097\zeta(6) +268\zeta(3)^2\right) \label{eq08115}
\end{align}
 
\begin{align}
\sum_{k=1}^\infty \frac{H(k)^{6}}{(k+2)^{2}} \sumend &= \frac{-1}{24}\left(  168 +576\zeta(2) +3120\zeta(3) +9288\zeta(4) +10104\zeta(5)  \right. \nonumber \\ &\left. \hspace{1em}
+2448\zeta(2)\zeta(3) -303\zeta(6) +528\zeta(3)^2 -16695\zeta(7) -4104\zeta(2)\zeta(5)  \right. \nonumber \\ &\left. \hspace{1em}
-4752\zeta(3)\zeta(4) -17027\zeta(8) +924\zeta(2)\zeta(3)^2 -11328\zeta(3)\zeta(5) -1308 M(2,6)\right) \label{eq08116}
\end{align}

\newpage
Formulas for order $r = m + n + p + q = 9$:
\begin{align}
\sum_{k=1}^\infty \frac{H(k)}{k^{8}} \sumend &= -\left( -5\zeta(9) +\zeta(3)\zeta(6) +\zeta(4)\zeta(5) +\zeta(2)\zeta(7)\right) \label{eq09001} \\
\sum_{k=1}^\infty \frac{H(k)}{k^{7}(k+1)} \sumend &= \frac{-1}{4}\left( -4\zeta(2) +8\zeta(3) -5\zeta(4) +12\zeta(5) -4\zeta(2)\zeta(3) -7\zeta(6)  \right. \nonumber \\ &\left. \hspace{1em}
+2\zeta(3)^2 +16\zeta(7) -4\zeta(2)\zeta(5) -4\zeta(3)\zeta(4) -9\zeta(8) +4\zeta(3)\zeta(5)\right) \label{eq09002} \\
\sum_{k=1}^\infty \frac{H(k)}{k^{6}(k+1)^{2}} \sumend &= \frac{-1}{2}\left(  12\zeta(2) -22\zeta(3) +10\zeta(4) -18\zeta(5) +6\zeta(2)\zeta(3) +7\zeta(6)  \right. \nonumber \\ &\left. \hspace{1em}
-2\zeta(3)^2 -8\zeta(7) +2\zeta(2)\zeta(5) +2\zeta(3)\zeta(4)\right) \label{eq09003} \\
\sum_{k=1}^\infty \frac{H(k)}{k^{5}(k+1)^{3}} \sumend &= \frac{1}{4}\left(  60\zeta(2) -100\zeta(3) +29\zeta(4) -36\zeta(5) +12\zeta(2)\zeta(3)  \right. \nonumber \\ &\left. \hspace{1em}
+7\zeta(6) -2\zeta(3)^2\right) \label{eq09004} \\
\sum_{k=1}^\infty \frac{H(k)}{k^{4}(k+1)^{4}} \sumend &= \left( -20\zeta(2) +30\zeta(3) -4\zeta(4) +5\zeta(5) -2\zeta(2)\zeta(3)\right) \label{eq09005} \\
\sum_{k=1}^\infty \frac{H(k)}{k^{3}(k+1)^{5}} \sumend &= \frac{1}{4}\left(  60\zeta(2) -80\zeta(3) -\zeta(4) -24\zeta(5) +12\zeta(2)\zeta(3) -3\zeta(6)  \right. \nonumber \\ &\left. \hspace{1em}
+2\zeta(3)^2\right) \label{eq09006} \\
\sum_{k=1}^\infty \frac{H(k)}{k^{2}(k+1)^{6}} \sumend &= \frac{-1}{2}\left(  12\zeta(2) -14\zeta(3) -2\zeta(4) -12\zeta(5) +6\zeta(2)\zeta(3) -3\zeta(6)  \right. \nonumber \\ &\left. \hspace{1em}
+2\zeta(3)^2 -6\zeta(7) +2\zeta(2)\zeta(5) +2\zeta(3)\zeta(4)\right) \label{eq09007} \\
\sum_{k=1}^\infty \frac{H(k)}{k(k+1)^{7}} \sumend &= \frac{-1}{4}\left( -4\zeta(2) +4\zeta(3) +\zeta(4) +8\zeta(5) -4\zeta(2)\zeta(3) +3\zeta(6)  \right. \nonumber \\ &\left. \hspace{1em}
-2\zeta(3)^2 +12\zeta(7) -4\zeta(2)\zeta(5) -4\zeta(3)\zeta(4) +5\zeta(8) -4\zeta(3)\zeta(5)\right) \label{eq09008} \\
\sum_{k=1}^\infty \frac{H(k)}{(k+1)^{8}} \sumend &= \left(  4\zeta(9) -\zeta(3)\zeta(6) -\zeta(4)\zeta(5) -\zeta(2)\zeta(7)\right) \label{eq09009} \\
\sum_{k=1}^\infty \frac{H(k)}{k^{7}(k+2)} \sumend &= \frac{1}{128}\left(  1 +\zeta(2) -4\zeta(3) +5\zeta(4) -24\zeta(5) +8\zeta(2)\zeta(3) +28\zeta(6)  \right. \nonumber \\ &\left. \hspace{1em}
-8\zeta(3)^2 -128\zeta(7) +32\zeta(2)\zeta(5) +32\zeta(3)\zeta(4) +144\zeta(8)  \right. \nonumber \\ &\left. \hspace{1em}
-64\zeta(3)\zeta(5)\right) \label{eq09010} \\
\sum_{k=1}^\infty \frac{H(k)}{k^{6}(k+1)(k+2)} \sumend &= \frac{-1}{64}\left( -1 +63\zeta(2) -124\zeta(3) +75\zeta(4) -168\zeta(5) +56\zeta(2)\zeta(3)  \right. \nonumber \\ &\left. \hspace{1em}
+84\zeta(6) -24\zeta(3)^2 -128\zeta(7) +32\zeta(2)\zeta(5) +32\zeta(3)\zeta(4)\right) \label{eq09011} \\
\sum_{k=1}^\infty \frac{H(k)}{k^{5}(k+1)^{2}(k+2)} \sumend &= \frac{-1}{32}\left( -1 -129\zeta(2) +228\zeta(3) -85\zeta(4) +120\zeta(5)  \right. \nonumber \\ &\left. \hspace{1em}
-40\zeta(2)\zeta(3) -28\zeta(6) +8\zeta(3)^2\right) \label{eq09012}
\end{align}
 
\begin{align}
\sum_{k=1}^\infty \frac{H(k)}{k^{4}(k+1)^{3}(k+2)} \sumend &= \frac{1}{16}\left(  1 -111\zeta(2) +172\zeta(3) -31\zeta(4) +24\zeta(5)  \right. \nonumber \\ &\left. \hspace{1em}
-8\zeta(2)\zeta(3)\right) \label{eq09013} \\
\sum_{k=1}^\infty \frac{H(k)}{k^{3}(k+1)^{4}(k+2)} \sumend &= \frac{-1}{8}\left( -1 -49\zeta(2) +68\zeta(3) -\zeta(4) +16\zeta(5)  \right. \nonumber \\ &\left. \hspace{1em}
-8\zeta(2)\zeta(3)\right) \label{eq09014} \\
\sum_{k=1}^\infty \frac{H(k)}{k^{2}(k+1)^{5}(k+2)} \sumend &= \frac{1}{4}\left(  1 -11\zeta(2) +12\zeta(3) +2\zeta(4) +8\zeta(5) -4\zeta(2)\zeta(3)  \right. \nonumber \\ &\left. \hspace{1em}
+3\zeta(6) -2\zeta(3)^2\right) \label{eq09015} \\
\sum_{k=1}^\infty \frac{H(k)}{k(k+1)^{6}(k+2)} \sumend &= \frac{1}{2}\left(  1 +\zeta(2) -2\zeta(3) -4\zeta(5) +2\zeta(2)\zeta(3) -6\zeta(7)  \right. \nonumber \\ &\left. \hspace{1em}
+2\zeta(2)\zeta(5) +2\zeta(3)\zeta(4)\right) \label{eq09016} \\
\sum_{k=1}^\infty \frac{H(k)}{(k+1)^{7}(k+2)} \sumend &= \frac{-1}{4}\left( -4 +4\zeta(3) -\zeta(4) +8\zeta(5) -4\zeta(2)\zeta(3) -3\zeta(6) +2\zeta(3)^2  \right. \nonumber \\ &\left. \hspace{1em}
+12\zeta(7) -4\zeta(2)\zeta(5) -4\zeta(3)\zeta(4) -5\zeta(8) +4\zeta(3)\zeta(5)\right) \label{eq09017} \\
\sum_{k=1}^\infty \frac{H(k)}{k^{6}(k+2)^{2}} \sumend &= \frac{1}{64}\left( -5 -2\zeta(2) +11\zeta(3) -10\zeta(4) +36\zeta(5) -12\zeta(2)\zeta(3)  \right. \nonumber \\ &\left. \hspace{1em}
-28\zeta(6) +8\zeta(3)^2 +64\zeta(7) -16\zeta(2)\zeta(5) -16\zeta(3)\zeta(4)\right) \label{eq09018} \\
\sum_{k=1}^\infty \frac{H(k)}{k^{5}(k+1)(k+2)^{2}} \sumend &= \frac{1}{64}\left( -11 +59\zeta(2) -102\zeta(3) +55\zeta(4) -96\zeta(5) +32\zeta(2)\zeta(3)  \right. \nonumber \\ &\left. \hspace{1em}
+28\zeta(6) -8\zeta(3)^2\right) \label{eq09019} \\
\sum_{k=1}^\infty \frac{H(k)}{k^{4}(k+1)^{2}(k+2)^{2}} \sumend &= \frac{-1}{16}\left(  6 +35\zeta(2) -63\zeta(3) +15\zeta(4) -12\zeta(5)  \right. \nonumber \\ &\left. \hspace{1em}
+4\zeta(2)\zeta(3)\right) \label{eq09020} \\
\sum_{k=1}^\infty \frac{H(k)}{k^{3}(k+1)^{3}(k+2)^{2}} \sumend &= \frac{1}{16}\left( -13 +41\zeta(2) -46\zeta(3) +\zeta(4)\right) \label{eq09021} \\
\sum_{k=1}^\infty \frac{H(k)}{k^{2}(k+1)^{4}(k+2)^{2}} \sumend &= \frac{-1}{4}\left(  7 +4\zeta(2) -11\zeta(3) -8\zeta(5) +4\zeta(2)\zeta(3)\right) \label{eq09022} \\
\sum_{k=1}^\infty \frac{H(k)}{k(k+1)^{5}(k+2)^{2}} \sumend &= \frac{1}{4}\left( -15 +3\zeta(2) +10\zeta(3) -2\zeta(4) +8\zeta(5) -4\zeta(2)\zeta(3)  \right. \nonumber \\ &\left. \hspace{1em}
-3\zeta(6) +2\zeta(3)^2\right) \label{eq09023} \\
\sum_{k=1}^\infty \frac{H(k)}{(k+1)^{6}(k+2)^{2}} \sumend &= \frac{1}{2}\left( -16 +2\zeta(2) +12\zeta(3) -2\zeta(4) +12\zeta(5) -6\zeta(2)\zeta(3)  \right. \nonumber \\ &\left. \hspace{1em}
-3\zeta(6) +2\zeta(3)^2 +6\zeta(7) -2\zeta(2)\zeta(5) -2\zeta(3)\zeta(4)\right) \label{eq09024}
\end{align}
 
\begin{align}
\sum_{k=1}^\infty \frac{H(k)}{k^{5}(k+2)^{3}} \sumend &= \frac{-1}{128}\left( -47 -\zeta(2) +54\zeta(3) -29\zeta(4) +72\zeta(5) -24\zeta(2)\zeta(3)  \right. \nonumber \\ &\left. \hspace{1em}
-28\zeta(6) +8\zeta(3)^2\right) \label{eq09025} \\
\sum_{k=1}^\infty \frac{H(k)}{k^{4}(k+1)(k+2)^{3}} \sumend &= \frac{-1}{32}\left( -29 +29\zeta(2) -24\zeta(3) +13\zeta(4) -12\zeta(5)  \right. \nonumber \\ &\left. \hspace{1em}
+4\zeta(2)\zeta(3)\right) \label{eq09026} \\
\sum_{k=1}^\infty \frac{H(k)}{k^{3}(k+1)^{2}(k+2)^{3}} \sumend &= \frac{-1}{16}\left( -35 -6\zeta(2) +39\zeta(3) -2\zeta(4)\right) \label{eq09027} \\
\sum_{k=1}^\infty \frac{H(k)}{k^{2}(k+1)^{3}(k+2)^{3}} \sumend &= \frac{-1}{16}\left( -83 +29\zeta(2) +32\zeta(3) -3\zeta(4)\right) \label{eq09028} \\
\sum_{k=1}^\infty \frac{H(k)}{k(k+1)^{4}(k+2)^{3}} \sumend &= \frac{-1}{8}\left( -97 +21\zeta(2) +54\zeta(3) -3\zeta(4) +16\zeta(5)  \right. \nonumber \\ &\left. \hspace{1em}
-8\zeta(2)\zeta(3)\right) \label{eq09029} \\
\sum_{k=1}^\infty \frac{H(k)}{(k+1)^{5}(k+2)^{3}} \sumend &= \frac{-1}{4}\left( -112 +24\zeta(2) +64\zeta(3) -5\zeta(4) +24\zeta(5) -12\zeta(2)\zeta(3)  \right. \nonumber \\ &\left. \hspace{1em}
-3\zeta(6) +2\zeta(3)^2\right) \label{eq09030} \\
\sum_{k=1}^\infty \frac{H(k)}{k^{4}(k+2)^{4}} \sumend &= \frac{-1}{32}\left(  35 -6\zeta(2) -21\zeta(3) +2\zeta(4) -10\zeta(5) +4\zeta(2)\zeta(3)\right) \label{eq09031} \\
\sum_{k=1}^\infty \frac{H(k)}{k^{3}(k+1)(k+2)^{4}} \sumend &= \frac{-1}{32}\left(  99 -41\zeta(2) -18\zeta(3) -9\zeta(4) -8\zeta(5)  \right. \nonumber \\ &\left. \hspace{1em}
+4\zeta(2)\zeta(3)\right) \label{eq09032} \\
\sum_{k=1}^\infty \frac{H(k)}{k^{2}(k+1)^{2}(k+2)^{4}} \sumend &= \frac{-1}{16}\left(  134 -35\zeta(2) -57\zeta(3) -7\zeta(4) -8\zeta(5)  \right. \nonumber \\ &\left. \hspace{1em}
+4\zeta(2)\zeta(3)\right) \label{eq09033} \\
\sum_{k=1}^\infty \frac{H(k)}{k(k+1)^{3}(k+2)^{4}} \sumend &= \frac{-1}{16}\left(  351 -99\zeta(2) -146\zeta(3) -11\zeta(4) -16\zeta(5)  \right. \nonumber \\ &\left. \hspace{1em}
+8\zeta(2)\zeta(3)\right) \label{eq09034} \\
\sum_{k=1}^\infty \frac{H(k)}{(k+1)^{4}(k+2)^{4}} \sumend &= -\left(  56 -15\zeta(2) -25\zeta(3) -\zeta(4) -4\zeta(5) +2\zeta(2)\zeta(3)\right) \label{eq09035} \\
\sum_{k=1}^\infty \frac{H(k)}{k^{3}(k+2)^{5}} \sumend &= \frac{-1}{128}\left( -303 +69\zeta(2) +104\zeta(3) +41\zeta(4) +64\zeta(5) -24\zeta(2)\zeta(3)  \right. \nonumber \\ &\left. \hspace{1em}
+12\zeta(6) -8\zeta(3)^2\right) \label{eq09036}
\end{align}
 
\begin{align}
\sum_{k=1}^\infty \frac{H(k)}{k^{2}(k+1)(k+2)^{5}} \sumend &= \frac{-1}{64}\left( -501 +151\zeta(2) +140\zeta(3) +59\zeta(4) +80\zeta(5)  \right. \nonumber \\ &\left. \hspace{1em}
-32\zeta(2)\zeta(3) +12\zeta(6) -8\zeta(3)^2\right) \label{eq09037} \\
\sum_{k=1}^\infty \frac{H(k)}{k(k+1)^{2}(k+2)^{5}} \sumend &= \frac{1}{32}\left(  769 -221\zeta(2) -254\zeta(3) -73\zeta(4) -96\zeta(5)  \right. \nonumber \\ &\left. \hspace{1em}
+40\zeta(2)\zeta(3) -12\zeta(6) +8\zeta(3)^2\right) \label{eq09038} \\
\sum_{k=1}^\infty \frac{H(k)}{(k+1)^{3}(k+2)^{5}} \sumend &= \frac{1}{4}\left(  280 -80\zeta(2) -100\zeta(3) -21\zeta(4) -28\zeta(5) +12\zeta(2)\zeta(3)  \right. \nonumber \\ &\left. \hspace{1em}
-3\zeta(6) +2\zeta(3)^2\right) \label{eq09039} \\
\sum_{k=1}^\infty \frac{H(k)}{k^{2}(k+2)^{6}} \sumend &= \frac{1}{64}\left( -261 +54\zeta(2) +59\zeta(3) +46\zeta(4) +56\zeta(5) -12\zeta(2)\zeta(3)  \right. \nonumber \\ &\left. \hspace{1em}
+28\zeta(6) -8\zeta(3)^2 +48\zeta(7) -16\zeta(2)\zeta(5) -16\zeta(3)\zeta(4)\right) \label{eq09040} \\
\sum_{k=1}^\infty \frac{H(k)}{k(k+1)(k+2)^{6}} \sumend &= \frac{-1}{64}\left(  1023 -259\zeta(2) -258\zeta(3) -151\zeta(4) -192\zeta(5)  \right. \nonumber \\ &\left. \hspace{1em}
+56\zeta(2)\zeta(3) -68\zeta(6) +24\zeta(3)^2 -96\zeta(7) +32\zeta(2)\zeta(5)  \right. \nonumber \\ &\left. \hspace{1em}
+32\zeta(3)\zeta(4)\right) \label{eq09041} \\
\sum_{k=1}^\infty \frac{H(k)}{(k+1)^{2}(k+2)^{6}} \sumend &= \frac{1}{2}\left( -112 +30\zeta(2) +32\zeta(3) +14\zeta(4) +18\zeta(5) -6\zeta(2)\zeta(3)  \right. \nonumber \\ &\left. \hspace{1em}
+5\zeta(6) -2\zeta(3)^2 +6\zeta(7) -2\zeta(2)\zeta(5) -2\zeta(3)\zeta(4)\right) \label{eq09042} \\
\sum_{k=1}^\infty \frac{H(k)}{k(k+2)^{7}} \sumend &= \frac{1}{128}\left(  769 -125\zeta(2) -126\zeta(3) -121\zeta(4) -128\zeta(5) +8\zeta(2)\zeta(3)  \right. \nonumber \\ &\left. \hspace{1em}
-108\zeta(6) +8\zeta(3)^2 -160\zeta(7) +32\zeta(2)\zeta(5) +32\zeta(3)\zeta(4) -80\zeta(8)  \right. \nonumber \\ &\left. \hspace{1em}
+64\zeta(3)\zeta(5)\right) \label{eq09043} \\
\sum_{k=1}^\infty \frac{H(k)}{(k+1)(k+2)^{7}} \sumend &= \frac{-1}{4}\left( -112 +24\zeta(2) +24\zeta(3) +17\zeta(4) +20\zeta(5) -4\zeta(2)\zeta(3)  \right. \nonumber \\ &\left. \hspace{1em}
+11\zeta(6) -2\zeta(3)^2 +16\zeta(7) -4\zeta(2)\zeta(5) -4\zeta(3)\zeta(4) +5\zeta(8)  \right. \nonumber \\ &\left. \hspace{1em}
-4\zeta(3)\zeta(5)\right) \label{eq09044} \\
\sum_{k=1}^\infty \frac{H(k)}{(k+2)^{8}} \sumend &= -\left(  8 -\zeta(2) -\zeta(3) -\zeta(4) -\zeta(5) -\zeta(6) -\zeta(7) -\zeta(8) -4\zeta(9)  \right. \nonumber \\ &\left. \hspace{1em}
+\zeta(3)\zeta(6) +\zeta(4)\zeta(5) +\zeta(2)\zeta(7)\right) \label{eq09045} \\
\sum_{k=1}^\infty \frac{H(k)^{2}}{k^{7}} \sumend &= \frac{-1}{6}\left( -55\zeta(9) +21\zeta(3)\zeta(6) +15\zeta(4)\zeta(5) +6\zeta(2)\zeta(7)  \right. \nonumber \\ &\left. \hspace{1em}
-2\zeta(3)^3\right) \label{eq09046}
\end{align}
 
\begin{align}
\sum_{k=1}^\infty \frac{H(k)^{2}}{k^{6}(k+1)} \sumend &= \frac{1}{24}\left( -72\zeta(3) +102\zeta(4) -84\zeta(5) +24\zeta(2)\zeta(3) +97\zeta(6)  \right. \nonumber \\ &\left. \hspace{1em}
-48\zeta(3)^2 -144\zeta(7) +24\zeta(2)\zeta(5) +60\zeta(3)\zeta(4) +24 M(2,6)\right) \label{eq09047} \\
\sum_{k=1}^\infty \frac{H(k)^{2}}{k^{5}(k+1)^{2}} \sumend &= \frac{-1}{12}\left( -180\zeta(3) +237\zeta(4) -126\zeta(5) +36\zeta(2)\zeta(3) +97\zeta(6)  \right. \nonumber \\ &\left. \hspace{1em}
-48\zeta(3)^2 -72\zeta(7) +12\zeta(2)\zeta(5) +30\zeta(3)\zeta(4)\right) \label{eq09048} \\
\sum_{k=1}^\infty \frac{H(k)^{2}}{k^{4}(k+1)^{3}} \sumend &= \frac{1}{24}\left( -720\zeta(3) +876\zeta(4) -288\zeta(5) +96\zeta(2)\zeta(3) +97\zeta(6)  \right. \nonumber \\ &\left. \hspace{1em}
-48\zeta(3)^2\right) \label{eq09049} \\
\sum_{k=1}^\infty \frac{H(k)^{2}}{k^{3}(k+1)^{4}} \sumend &= \frac{-1}{24}\left( -720\zeta(3) +804\zeta(4) -192\zeta(5) +96\zeta(2)\zeta(3) +37\zeta(6)  \right. \nonumber \\ &\left. \hspace{1em}
-24\zeta(3)^2\right) \label{eq09050} \\
\sum_{k=1}^\infty \frac{H(k)^{2}}{k^{2}(k+1)^{5}} \sumend &= \frac{1}{12}\left( -180\zeta(3) +183\zeta(4) -54\zeta(5) +36\zeta(2)\zeta(3) +37\zeta(6)  \right. \nonumber \\ &\left. \hspace{1em}
-24\zeta(3)^2 -12\zeta(7) +12\zeta(2)\zeta(5) -6\zeta(3)\zeta(4)\right) \label{eq09051} \\
\sum_{k=1}^\infty \frac{H(k)^{2}}{k(k+1)^{6}} \sumend &= \frac{1}{24}\left(  72\zeta(3) -66\zeta(4) +36\zeta(5) -24\zeta(2)\zeta(3) -37\zeta(6)  \right. \nonumber \\ &\left. \hspace{1em}
+24\zeta(3)^2 +24\zeta(7) -24\zeta(2)\zeta(5) +12\zeta(3)\zeta(4) +84\zeta(8) -48\zeta(3)\zeta(5)  \right. \nonumber \\ &\left. \hspace{1em}
-24 M(2,6)\right) \label{eq09052} \\
\sum_{k=1}^\infty \frac{H(k)^{2}}{(k+1)^{7}} \sumend &= \frac{-1}{6}\left( -\zeta(9) +9\zeta(3)\zeta(6) +3\zeta(4)\zeta(5) -6\zeta(2)\zeta(7)  \right. \nonumber \\ &\left. \hspace{1em}
-2\zeta(3)^3\right) \label{eq09053} \\
\sum_{k=1}^\infty \frac{H(k)^{2}}{k^{6}(k+2)} \sumend &= \frac{1}{384}\left( -6 -6\zeta(2) -18\zeta(3) +51\zeta(4) -84\zeta(5) +24\zeta(2)\zeta(3)  \right. \nonumber \\ &\left. \hspace{1em}
+194\zeta(6) -96\zeta(3)^2 -576\zeta(7) +96\zeta(2)\zeta(5) +240\zeta(3)\zeta(4) +192 M(2,6)\right) \label{eq09054} \\
\sum_{k=1}^\infty \frac{H(k)^{2}}{k^{5}(k+1)(k+2)} \sumend &= \frac{-1}{64}\left(  2 +2\zeta(2) -186\zeta(3) +255\zeta(4) -196\zeta(5)  \right. \nonumber \\ &\left. \hspace{1em}
+56\zeta(2)\zeta(3) +194\zeta(6) -96\zeta(3)^2 -192\zeta(7) +32\zeta(2)\zeta(5)  \right. \nonumber \\ &\left. \hspace{1em}
+80\zeta(3)\zeta(4)\right) \label{eq09055} \\
\sum_{k=1}^\infty \frac{H(k)^{2}}{k^{4}(k+1)^{2}(k+2)} \sumend &= \frac{1}{96}\left( -6 -6\zeta(2) -882\zeta(3) +1131\zeta(4) -420\zeta(5)  \right. \nonumber \\ &\left. \hspace{1em}
+120\zeta(2)\zeta(3) +194\zeta(6) -96\zeta(3)^2\right) \label{eq09056}
\end{align}
 
\begin{align}
\sum_{k=1}^\infty \frac{H(k)^{2}}{k^{3}(k+1)^{3}(k+2)} \sumend &= \frac{-1}{16}\left(  2 +2\zeta(2) -186\zeta(3) +207\zeta(4) -52\zeta(5)  \right. \nonumber \\ &\left. \hspace{1em}
+24\zeta(2)\zeta(3)\right) \label{eq09057} \\
\sum_{k=1}^\infty \frac{H(k)^{2}}{k^{2}(k+1)^{4}(k+2)} \sumend &= \frac{-1}{24}\left(  6 +6\zeta(2) +162\zeta(3) -183\zeta(4) +36\zeta(5)  \right. \nonumber \\ &\left. \hspace{1em}
-24\zeta(2)\zeta(3) -37\zeta(6) +24\zeta(3)^2\right) \label{eq09058} \\
\sum_{k=1}^\infty \frac{H(k)^{2}}{k(k+1)^{5}(k+2)} \sumend &= \frac{-1}{2}\left(  1 +\zeta(2) -3\zeta(3) -3\zeta(5) +2\zeta(2)\zeta(3) -2\zeta(7)  \right. \nonumber \\ &\left. \hspace{1em}
+2\zeta(2)\zeta(5) -\zeta(3)\zeta(4)\right) \label{eq09059} \\
\sum_{k=1}^\infty \frac{H(k)^{2}}{(k+1)^{6}(k+2)} \sumend &= \frac{1}{24}\left( -24 -24\zeta(2) +66\zeta(4) +36\zeta(5) -24\zeta(2)\zeta(3) +37\zeta(6)  \right. \nonumber \\ &\left. \hspace{1em}
-24\zeta(3)^2 +24\zeta(7) -24\zeta(2)\zeta(5) +12\zeta(3)\zeta(4) -84\zeta(8) +48\zeta(3)\zeta(5)  \right. \nonumber \\ &\left. \hspace{1em}
+24 M(2,6)\right) \label{eq09060} \\
\sum_{k=1}^\infty \frac{H(k)^{2}}{k^{5}(k+2)^{2}} \sumend &= \frac{1}{384}\left(  66 +30\zeta(2) +66\zeta(3) -237\zeta(4) +252\zeta(5)  \right. \nonumber \\ &\left. \hspace{1em}
-72\zeta(2)\zeta(3) -388\zeta(6) +192\zeta(3)^2 +576\zeta(7) -96\zeta(2)\zeta(5)  \right. \nonumber \\ &\left. \hspace{1em}
-240\zeta(3)\zeta(4)\right) \label{eq09061} \\
\sum_{k=1}^\infty \frac{H(k)^{2}}{k^{4}(k+1)(k+2)^{2}} \sumend &= \frac{1}{96}\left(  36 +18\zeta(2) -246\zeta(3) +264\zeta(4) -168\zeta(5)  \right. \nonumber \\ &\left. \hspace{1em}
+48\zeta(2)\zeta(3) +97\zeta(6) -48\zeta(3)^2\right) \label{eq09062} \\
\sum_{k=1}^\infty \frac{H(k)^{2}}{k^{3}(k+1)^{2}(k+2)^{2}} \sumend &= \frac{1}{32}\left(  26 +14\zeta(2) +130\zeta(3) -201\zeta(4) +28\zeta(5)  \right. \nonumber \\ &\left. \hspace{1em}
-8\zeta(2)\zeta(3)\right) \label{eq09063} \\
\sum_{k=1}^\infty \frac{H(k)^{2}}{k^{2}(k+1)^{3}(k+2)^{2}} \sumend &= \frac{1}{8}\left(  14 +8\zeta(2) -28\zeta(3) +3\zeta(4) -12\zeta(5)  \right. \nonumber \\ &\left. \hspace{1em}
+8\zeta(2)\zeta(3)\right) \label{eq09064} \\
\sum_{k=1}^\infty \frac{H(k)^{2}}{k(k+1)^{4}(k+2)^{2}} \sumend &= \frac{1}{24}\left(  90 +54\zeta(2) -6\zeta(3) -165\zeta(4) -36\zeta(5)  \right. \nonumber \\ &\left. \hspace{1em}
+24\zeta(2)\zeta(3) -37\zeta(6) +24\zeta(3)^2\right) \label{eq09065} \\
\sum_{k=1}^\infty \frac{H(k)^{2}}{(k+1)^{5}(k+2)^{2}} \sumend &= \frac{-1}{12}\left( -96 -60\zeta(2) +24\zeta(3) +165\zeta(4) +54\zeta(5)  \right. \nonumber \\ &\left. \hspace{1em}
-36\zeta(2)\zeta(3) +37\zeta(6) -24\zeta(3)^2 +12\zeta(7) -12\zeta(2)\zeta(5)  \right. \nonumber \\ &\left. \hspace{1em}
+6\zeta(3)\zeta(4)\right) \label{eq09066}
\end{align}
 
\begin{align}
\sum_{k=1}^\infty \frac{H(k)^{2}}{k^{4}(k+2)^{3}} \sumend &= \frac{-1}{192}\left(  174 +18\zeta(2) +6\zeta(3) -225\zeta(4) +144\zeta(5)  \right. \nonumber \\ &\left. \hspace{1em}
-48\zeta(2)\zeta(3) -97\zeta(6) +48\zeta(3)^2\right) \label{eq09067} \\
\sum_{k=1}^\infty \frac{H(k)^{2}}{k^{3}(k+1)(k+2)^{3}} \sumend &= \frac{-1}{32}\left(  70 +12\zeta(2) -80\zeta(3) +13\zeta(4) -8\zeta(5)\right) \label{eq09068} \\
\sum_{k=1}^\infty \frac{H(k)^{2}}{k^{2}(k+1)^{2}(k+2)^{3}} \sumend &= \frac{-1}{32}\left(  166 +38\zeta(2) -30\zeta(3) -175\zeta(4) +12\zeta(5)  \right. \nonumber \\ &\left. \hspace{1em}
-8\zeta(2)\zeta(3)\right) \label{eq09069} \\
\sum_{k=1}^\infty \frac{H(k)^{2}}{k(k+1)^{3}(k+2)^{3}} \sumend &= \frac{1}{16}\left( -194 -54\zeta(2) +86\zeta(3) +169\zeta(4) +12\zeta(5)  \right. \nonumber \\ &\left. \hspace{1em}
-8\zeta(2)\zeta(3)\right) \label{eq09070} \\
\sum_{k=1}^\infty \frac{H(k)^{2}}{(k+1)^{4}(k+2)^{3}} \sumend &= \frac{-1}{24}\left(  672 +216\zeta(2) -264\zeta(3) -672\zeta(4) -72\zeta(5)  \right. \nonumber \\ &\left. \hspace{1em}
+48\zeta(2)\zeta(3) -37\zeta(6) +24\zeta(3)^2\right) \label{eq09071} \\
\sum_{k=1}^\infty \frac{H(k)^{2}}{k^{3}(k+2)^{4}} \sumend &= \frac{1}{192}\left(  594 -54\zeta(2) -186\zeta(3) -255\zeta(4) -37\zeta(6)  \right. \nonumber \\ &\left. \hspace{1em}
+24\zeta(3)^2\right) \label{eq09072} \\
\sum_{k=1}^\infty \frac{H(k)^{2}}{k^{2}(k+1)(k+2)^{4}} \sumend &= \frac{-1}{96}\left( -804 +18\zeta(2) +426\zeta(3) +216\zeta(4) +24\zeta(5) +37\zeta(6)  \right. \nonumber \\ &\left. \hspace{1em}
-24\zeta(3)^2\right) \label{eq09073} \\
\sum_{k=1}^\infty \frac{H(k)^{2}}{k(k+1)^{2}(k+2)^{4}} \sumend &= \frac{-1}{96}\left( -2106 -78\zeta(2) +942\zeta(3) +957\zeta(4) +12\zeta(5)  \right. \nonumber \\ &\left. \hspace{1em}
+24\zeta(2)\zeta(3) +74\zeta(6) -48\zeta(3)^2\right) \label{eq09074} \\
\sum_{k=1}^\infty \frac{H(k)^{2}}{(k+1)^{3}(k+2)^{4}} \sumend &= \frac{-1}{24}\left( -1344 -120\zeta(2) +600\zeta(3) +732\zeta(4) +24\zeta(5) +37\zeta(6)  \right. \nonumber \\ &\left. \hspace{1em}
-24\zeta(3)^2\right) \label{eq09075} \\
\sum_{k=1}^\infty \frac{H(k)^{2}}{k^{2}(k+2)^{5}} \sumend &= \frac{-1}{384}\left(  3006 -522\zeta(2) -1086\zeta(3) -603\zeta(4) -756\zeta(5)  \right. \nonumber \\ &\left. \hspace{1em}
+312\zeta(2)\zeta(3) -292\zeta(6) +192\zeta(3)^2 +96\zeta(7) -96\zeta(2)\zeta(5)  \right. \nonumber \\ &\left. \hspace{1em}
+48\zeta(3)\zeta(4)\right) \label{eq09076} \\
\sum_{k=1}^\infty \frac{H(k)^{2}}{k(k+1)(k+2)^{5}} \sumend &= \frac{1}{64}\left( -1538 +186\zeta(2) +646\zeta(3) +345\zeta(4) +268\zeta(5)  \right. \nonumber \\ &\left. \hspace{1em}
-104\zeta(2)\zeta(3) +122\zeta(6) -80\zeta(3)^2 -32\zeta(7) +32\zeta(2)\zeta(5)  \right. \nonumber \\ &\left. \hspace{1em}
-16\zeta(3)\zeta(4)\right) \label{eq09077}
\end{align}
 
\begin{align}
\sum_{k=1}^\infty \frac{H(k)^{2}}{(k+1)^{2}(k+2)^{5}} \sumend &= \frac{1}{12}\left( -840 +60\zeta(2) +360\zeta(3) +249\zeta(4) +102\zeta(5)  \right. \nonumber \\ &\left. \hspace{1em}
-36\zeta(2)\zeta(3) +55\zeta(6) -36\zeta(3)^2 -12\zeta(7) +12\zeta(2)\zeta(5)  \right. \nonumber \\ &\left. \hspace{1em}
-6\zeta(3)\zeta(4)\right) \label{eq09078} \\
\sum_{k=1}^\infty \frac{H(k)^{2}}{k(k+2)^{6}} \sumend &= \frac{1}{384}\left(  6138 -1170\zeta(2) -1902\zeta(3) -1029\zeta(4) -1788\zeta(5)  \right. \nonumber \\ &\left. \hspace{1em}
+648\zeta(2)\zeta(3) -698\zeta(6) +336\zeta(3)^2 -1056\zeta(7) +288\zeta(2)\zeta(5)  \right. \nonumber \\ &\left. \hspace{1em}
+432\zeta(3)\zeta(4) +672\zeta(8) -384\zeta(3)\zeta(5) -192 M(2,6)\right) \label{eq09079} \\
\sum_{k=1}^\infty \frac{H(k)^{2}}{(k+1)(k+2)^{6}} \sumend &= \frac{-1}{24}\left( -1344 +216\zeta(2) +480\zeta(3) +258\zeta(4) +324\zeta(5)  \right. \nonumber \\ &\left. \hspace{1em}
-120\zeta(2)\zeta(3) +133\zeta(6) -72\zeta(3)^2 +120\zeta(7) -24\zeta(2)\zeta(5) -60\zeta(3)\zeta(4)  \right. \nonumber \\ &\left. \hspace{1em}
-84\zeta(8) +48\zeta(3)\zeta(5) +24 M(2,6)\right) \label{eq09080} \\
\sum_{k=1}^\infty \frac{H(k)^{2}}{(k+2)^{7}} \sumend &= \frac{1}{6}\left( -168 +30\zeta(2) +42\zeta(3) +27\zeta(4) +42\zeta(5) -12\zeta(2)\zeta(3)  \right. \nonumber \\ &\left. \hspace{1em}
+21\zeta(6) -6\zeta(3)^2 +42\zeta(7) -12\zeta(2)\zeta(5) -12\zeta(3)\zeta(4) +15\zeta(8)  \right. \nonumber \\ &\left. \hspace{1em}
-12\zeta(3)\zeta(5) +\zeta(9) -9\zeta(3)\zeta(6) -3\zeta(4)\zeta(5) +6\zeta(2)\zeta(7)  \right. \nonumber \\ &\left. \hspace{1em}
+2\zeta(3)^3\right) \label{eq09081} \\
\sum_{k=1}^\infty \frac{H(k)^{3}}{k^{6}} \sumend &= \frac{1}{24}\left(  521\zeta(9) -291\zeta(3)\zeta(6) -306\zeta(4)\zeta(5) +72\zeta(2)\zeta(7)  \right. \nonumber \\ &\left. \hspace{1em}
+48\zeta(3)^3\right) \label{eq09082} \\
\sum_{k=1}^\infty \frac{H(k)^{3}}{k^{5}(k+1)} \sumend &= \frac{1}{96}\left(  960\zeta(4) -960\zeta(5) -96\zeta(2)\zeta(3) +558\zeta(6) -240\zeta(3)^2  \right. \nonumber \\ &\left. \hspace{1em}
-1386\zeta(7) -192\zeta(2)\zeta(5) +1224\zeta(3)\zeta(4) -595\zeta(8) -120\zeta(2)\zeta(3)^2  \right. \nonumber \\ &\left. \hspace{1em}
+576\zeta(3)\zeta(5) +264 M(2,6)\right) \label{eq09083} \\
\sum_{k=1}^\infty \frac{H(k)^{3}}{k^{4}(k+1)^{2}} \sumend &= \frac{1}{16}\left( -640\zeta(4) +600\zeta(5) +64\zeta(2)\zeta(3) -186\zeta(6) +80\zeta(3)^2  \right. \nonumber \\ &\left. \hspace{1em}
+231\zeta(7) +32\zeta(2)\zeta(5) -204\zeta(3)\zeta(4)\right) \label{eq09084} \\
\sum_{k=1}^\infty \frac{H(k)^{3}}{k^{3}(k+1)^{3}} \sumend &= \frac{1}{8}\left(  480\zeta(4) -420\zeta(5) -48\zeta(2)\zeta(3) +63\zeta(6)  \right. \nonumber \\ &\left. \hspace{1em}
-36\zeta(3)^2\right) \label{eq09085} \\
\sum_{k=1}^\infty \frac{H(k)^{3}}{k^{2}(k+1)^{4}} \sumend &= \frac{1}{16}\left( -640\zeta(4) +520\zeta(5) +64\zeta(2)\zeta(3) -66\zeta(6) +64\zeta(3)^2  \right. \nonumber \\ &\left. \hspace{1em}
+119\zeta(7) +32\zeta(2)\zeta(5) -132\zeta(3)\zeta(4)\right) \label{eq09086}
\end{align}
 
\begin{align}
\sum_{k=1}^\infty \frac{H(k)^{3}}{k(k+1)^{5}} \sumend &= \frac{1}{96}\left(  960\zeta(4) -720\zeta(5) -96\zeta(2)\zeta(3) +198\zeta(6) -192\zeta(3)^2  \right. \nonumber \\ &\left. \hspace{1em}
-714\zeta(7) -192\zeta(2)\zeta(5) +792\zeta(3)\zeta(4) +43\zeta(8) +120\zeta(2)\zeta(3)^2  \right. \nonumber \\ &\left. \hspace{1em}
-288\zeta(3)\zeta(5) +24 M(2,6)\right) \label{eq09087} \\
\sum_{k=1}^\infty \frac{H(k)^{3}}{(k+1)^{6}} \sumend &= \frac{-1}{24}\left( -197\zeta(9) +111\zeta(3)\zeta(6) +198\zeta(4)\zeta(5) -72\zeta(2)\zeta(7)  \right. \nonumber \\ &\left. \hspace{1em}
-24\zeta(3)^3\right) \label{eq09088} \\
\sum_{k=1}^\infty \frac{H(k)^{3}}{k^{5}(k+2)} \sumend &= \frac{-1}{384}\left( -12 -24\zeta(2) -48\zeta(3) -120\zeta(4) +240\zeta(5) +24\zeta(2)\zeta(3)  \right. \nonumber \\ &\left. \hspace{1em}
-279\zeta(6) +120\zeta(3)^2 +1386\zeta(7) +192\zeta(2)\zeta(5) -1224\zeta(3)\zeta(4) +1190\zeta(8)  \right. \nonumber \\ &\left. \hspace{1em}
+240\zeta(2)\zeta(3)^2 -1152\zeta(3)\zeta(5) -528 M(2,6)\right) \label{eq09089} \\
\sum_{k=1}^\infty \frac{H(k)^{3}}{k^{4}(k+1)(k+2)} \sumend &= \frac{1}{64}\left(  4 +8\zeta(2) +16\zeta(3) -600\zeta(4) +560\zeta(5) +56\zeta(2)\zeta(3)  \right. \nonumber \\ &\left. \hspace{1em}
-279\zeta(6) +120\zeta(3)^2 +462\zeta(7) +64\zeta(2)\zeta(5) -408\zeta(3)\zeta(4)\right) \label{eq09090} \\
\sum_{k=1}^\infty \frac{H(k)^{3}}{k^{3}(k+1)^{2}(k+2)} \sumend &= \frac{1}{32}\left(  4 +8\zeta(2) +16\zeta(3) +680\zeta(4) -640\zeta(5)  \right. \nonumber \\ &\left. \hspace{1em}
-72\zeta(2)\zeta(3) +93\zeta(6) -40\zeta(3)^2\right) \label{eq09091} \\
\sum_{k=1}^\infty \frac{H(k)^{3}}{k^{2}(k+1)^{3}(k+2)} \sumend &= \frac{1}{16}\left(  4 +8\zeta(2) +16\zeta(3) -280\zeta(4) +200\zeta(5)  \right. \nonumber \\ &\left. \hspace{1em}
+24\zeta(2)\zeta(3) -33\zeta(6) +32\zeta(3)^2\right) \label{eq09092} \\
\sum_{k=1}^\infty \frac{H(k)^{3}}{k(k+1)^{4}(k+2)} \sumend &= \frac{-1}{16}\left( -8 -16\zeta(2) -32\zeta(3) -80\zeta(4) +120\zeta(5) +16\zeta(2)\zeta(3)  \right. \nonumber \\ &\left. \hspace{1em}
+119\zeta(7) +32\zeta(2)\zeta(5) -132\zeta(3)\zeta(4)\right) \label{eq09093} \\
\sum_{k=1}^\infty \frac{H(k)^{3}}{(k+1)^{5}(k+2)} \sumend &= \frac{1}{96}\left(  96 +192\zeta(2) +384\zeta(3) -720\zeta(5) -96\zeta(2)\zeta(3)  \right. \nonumber \\ &\left. \hspace{1em}
-198\zeta(6) +192\zeta(3)^2 -714\zeta(7) -192\zeta(2)\zeta(5) +792\zeta(3)\zeta(4) -43\zeta(8)  \right. \nonumber \\ &\left. \hspace{1em}
-120\zeta(2)\zeta(3)^2 +288\zeta(3)\zeta(5) -24 M(2,6)\right) \label{eq09094} \\
\sum_{k=1}^\infty \frac{H(k)^{3}}{k^{4}(k+2)^{2}} \sumend &= \frac{-1}{64}\left(  24 +28\zeta(2) +36\zeta(3) +47\zeta(4) -150\zeta(5) -16\zeta(2)\zeta(3)  \right. \nonumber \\ &\left. \hspace{1em}
+93\zeta(6) -40\zeta(3)^2 -231\zeta(7) -32\zeta(2)\zeta(5) +204\zeta(3)\zeta(4)\right) \label{eq09095} \\
\sum_{k=1}^\infty \frac{H(k)^{3}}{k^{3}(k+1)(k+2)^{2}} \sumend &= \frac{-1}{64}\left(  52 +64\zeta(2) +88\zeta(3) -506\zeta(4) +260\zeta(5)  \right. \nonumber \\ &\left. \hspace{1em}
+24\zeta(2)\zeta(3) -93\zeta(6) +40\zeta(3)^2\right) \label{eq09096}
\end{align}
 
\begin{align}
\sum_{k=1}^\infty \frac{H(k)^{3}}{k^{2}(k+1)^{2}(k+2)^{2}} \sumend &= \frac{1}{16}\left( -28 -36\zeta(2) -52\zeta(3) -87\zeta(4) +190\zeta(5)  \right. \nonumber \\ &\left. \hspace{1em}
+24\zeta(2)\zeta(3)\right) \label{eq09097} \\
\sum_{k=1}^\infty \frac{H(k)^{3}}{k(k+1)^{3}(k+2)^{2}} \sumend &= \frac{1}{16}\left( -60 -80\zeta(2) -120\zeta(3) +106\zeta(4) +180\zeta(5)  \right. \nonumber \\ &\left. \hspace{1em}
+24\zeta(2)\zeta(3) +33\zeta(6) -32\zeta(3)^2\right) \label{eq09098} \\
\sum_{k=1}^\infty \frac{H(k)^{3}}{(k+1)^{4}(k+2)^{2}} \sumend &= \frac{-1}{16}\left(  128 +176\zeta(2) +272\zeta(3) -132\zeta(4) -480\zeta(5)  \right. \nonumber \\ &\left. \hspace{1em}
-64\zeta(2)\zeta(3) -66\zeta(6) +64\zeta(3)^2 -119\zeta(7) -32\zeta(2)\zeta(5)  \right. \nonumber \\ &\left. \hspace{1em}
+132\zeta(3)\zeta(4)\right) \label{eq09099} \\
\sum_{k=1}^\infty \frac{H(k)^{3}}{k^{3}(k+2)^{3}} \sumend &= \frac{-1}{64}\left( -140 -84\zeta(2) -36\zeta(3) +51\zeta(4) +174\zeta(5)  \right. \nonumber \\ &\left. \hspace{1em}
+48\zeta(2)\zeta(3) -63\zeta(6) +36\zeta(3)^2\right) \label{eq09100} \\
\sum_{k=1}^\infty \frac{H(k)^{3}}{k^{2}(k+1)(k+2)^{3}} \sumend &= \frac{-1}{64}\left( -332 -232\zeta(2) -160\zeta(3) +608\zeta(4) +88\zeta(5)  \right. \nonumber \\ &\left. \hspace{1em}
+72\zeta(2)\zeta(3) -33\zeta(6) +32\zeta(3)^2\right) \label{eq09101} \\
\sum_{k=1}^\infty \frac{H(k)^{3}}{k(k+1)^{2}(k+2)^{3}} \sumend &= \frac{1}{32}\left(  388 +304\zeta(2) +264\zeta(3) -434\zeta(4) -468\zeta(5)  \right. \nonumber \\ &\left. \hspace{1em}
-120\zeta(2)\zeta(3) +33\zeta(6) -32\zeta(3)^2\right) \label{eq09102} \\
\sum_{k=1}^\infty \frac{H(k)^{3}}{(k+1)^{3}(k+2)^{3}} \sumend &= \frac{1}{4}\left(  112 +96\zeta(2) +96\zeta(3) -135\zeta(4) -162\zeta(5)  \right. \nonumber \\ &\left. \hspace{1em}
-36\zeta(2)\zeta(3)\right) \label{eq09103} \\
\sum_{k=1}^\infty \frac{H(k)^{3}}{k^{2}(k+2)^{4}} \sumend &= \frac{-1}{64}\left(  536 +132\zeta(2) -132\zeta(3) -335\zeta(4) -82\zeta(5)  \right. \nonumber \\ &\left. \hspace{1em}
-64\zeta(2)\zeta(3) -41\zeta(6) +16\zeta(3)^2 -119\zeta(7) -32\zeta(2)\zeta(5)  \right. \nonumber \\ &\left. \hspace{1em}
+132\zeta(3)\zeta(4)\right) \label{eq09104} \\
\sum_{k=1}^\infty \frac{H(k)^{3}}{k(k+1)(k+2)^{4}} \sumend &= \frac{-1}{64}\left(  1404 +496\zeta(2) -104\zeta(3) -1278\zeta(4) -252\zeta(5)  \right. \nonumber \\ &\left. \hspace{1em}
-200\zeta(2)\zeta(3) -49\zeta(6) -238\zeta(7) -64\zeta(2)\zeta(5) +264\zeta(3)\zeta(4)\right) \label{eq09105} \\
\sum_{k=1}^\infty \frac{H(k)^{3}}{(k+1)^{2}(k+2)^{4}} \sumend &= \frac{1}{16}\left( -896 -400\zeta(2) -80\zeta(3) +856\zeta(4) +360\zeta(5)  \right. \nonumber \\ &\left. \hspace{1em}
+160\zeta(2)\zeta(3) +8\zeta(6) +16\zeta(3)^2 +119\zeta(7) +32\zeta(2)\zeta(5)  \right. \nonumber \\ &\left. \hspace{1em}
-132\zeta(3)\zeta(4)\right) \label{eq09106}
\end{align}
 
\begin{align}
\sum_{k=1}^\infty \frac{H(k)^{3}}{k(k+2)^{5}} \sumend &= \frac{1}{384}\left(  9228 +432\zeta(2) -3528\zeta(3) -4134\zeta(4) -1740\zeta(5)  \right. \nonumber \\ &\left. \hspace{1em}
+408\zeta(2)\zeta(3) -1665\zeta(6) +1056\zeta(3)^2 -138\zeta(7) -768\zeta(2)\zeta(5)  \right. \nonumber \\ &\left. \hspace{1em}
+1080\zeta(3)\zeta(4) +86\zeta(8) +240\zeta(2)\zeta(3)^2 -576\zeta(3)\zeta(5) +48 M(2,6)\right) \label{eq09107} \\
\sum_{k=1}^\infty \frac{H(k)^{3}}{(k+1)(k+2)^{5}} \sumend &= \frac{-1}{96}\left( -6720 -960\zeta(2) +1920\zeta(3) +3984\zeta(4) +1248\zeta(5)  \right. \nonumber \\ &\left. \hspace{1em}
+96\zeta(2)\zeta(3) +906\zeta(6) -528\zeta(3)^2 +426\zeta(7) +480\zeta(2)\zeta(5) -936\zeta(3)\zeta(4)  \right. \nonumber \\ &\left. \hspace{1em}
-43\zeta(8) -120\zeta(2)\zeta(3)^2 +288\zeta(3)\zeta(5) -24 M(2,6)\right) \label{eq09108} \\
\sum_{k=1}^\infty \frac{H(k)^{3}}{(k+2)^{6}} \sumend &= \frac{1}{24}\left( -1344 +72\zeta(2) +504\zeta(3) +414\zeta(4) +396\zeta(5) -144\zeta(2)\zeta(3)  \right. \nonumber \\ &\left. \hspace{1em}
+243\zeta(6) -144\zeta(3)^2 +144\zeta(7) -108\zeta(3)\zeta(4) -252\zeta(8) +144\zeta(3)\zeta(5)  \right. \nonumber \\ &\left. \hspace{1em}
+72 M(2,6) +197\zeta(9) -111\zeta(3)\zeta(6) -198\zeta(4)\zeta(5) +72\zeta(2)\zeta(7)  \right. \nonumber \\ &\left. \hspace{1em}
+24\zeta(3)^3\right) \label{eq09109} \\
\sum_{k=1}^\infty \frac{H(k)^{4}}{k^{5}} \sumend &= \frac{1}{12}\left(  436\zeta(9) -279\zeta(3)\zeta(6) -258\zeta(4)\zeta(5) +84\zeta(2)\zeta(7)  \right. \nonumber \\ &\left. \hspace{1em}
+40\zeta(3)^3\right) \label{eq09110} \\
\sum_{k=1}^\infty \frac{H(k)^{4}}{k^{4}(k+1)} \sumend &= \frac{-1}{144}\left(  4320\zeta(5) +864\zeta(2)\zeta(3) -5874\zeta(6) -432\zeta(3)^2  \right. \nonumber \\ &\left. \hspace{1em}
+3330\zeta(7) +720\zeta(2)\zeta(5) -3096\zeta(3)\zeta(4) +14833\zeta(8) +4032\zeta(2)\zeta(3)^2  \right. \nonumber \\ &\left. \hspace{1em}
-16704\zeta(3)\zeta(5) -3744 M(2,6)\right) \label{eq09111} \\
\sum_{k=1}^\infty \frac{H(k)^{4}}{k^{3}(k+1)^{2}} \sumend &= \frac{-1}{8}\left( -720\zeta(5) -144\zeta(2)\zeta(3) +939\zeta(6) +72\zeta(3)^2 -185\zeta(7)  \right. \nonumber \\ &\left. \hspace{1em}
-40\zeta(2)\zeta(5) +172\zeta(3)\zeta(4)\right) \label{eq09112} \\
\sum_{k=1}^\infty \frac{H(k)^{4}}{k^{2}(k+1)^{3}} \sumend &= \frac{-1}{8}\left(  720\zeta(5) +144\zeta(2)\zeta(3) -899\zeta(6) -72\zeta(3)^2 +109\zeta(7)  \right. \nonumber \\ &\left. \hspace{1em}
+40\zeta(2)\zeta(5) -148\zeta(3)\zeta(4)\right) \label{eq09113} \\
\sum_{k=1}^\infty \frac{H(k)^{4}}{k(k+1)^{4}} \sumend &= \frac{1}{144}\left(  4320\zeta(5) +864\zeta(2)\zeta(3) -5154\zeta(6) -432\zeta(3)^2  \right. \nonumber \\ &\left. \hspace{1em}
+1962\zeta(7) +720\zeta(2)\zeta(5) -2664\zeta(3)\zeta(4) +12415\zeta(8) +3312\zeta(2)\zeta(3)^2  \right. \nonumber \\ &\left. \hspace{1em}
-13824\zeta(3)\zeta(5) -3024 M(2,6)\right) \label{eq09114} \\
\sum_{k=1}^\infty \frac{H(k)^{4}}{(k+1)^{5}} \sumend &= \frac{-1}{12}\left(  174\zeta(9) -99\zeta(3)\zeta(6) -222\zeta(4)\zeta(5) +84\zeta(2)\zeta(7)  \right. \nonumber \\ &\left. \hspace{1em}
+32\zeta(3)^3\right) \label{eq09115}
\end{align}
 
\begin{align}
\sum_{k=1}^\infty \frac{H(k)^{4}}{k^{4}(k+2)} \sumend &= \frac{1}{576}\left( -36 -108\zeta(2) -396\zeta(3) -666\zeta(4) -1080\zeta(5)  \right. \nonumber \\ &\left. \hspace{1em}
-216\zeta(2)\zeta(3) +2937\zeta(6) +216\zeta(3)^2 -3330\zeta(7) -720\zeta(2)\zeta(5)  \right. \nonumber \\ &\left. \hspace{1em}
+3096\zeta(3)\zeta(4) -29666\zeta(8) -8064\zeta(2)\zeta(3)^2 +33408\zeta(3)\zeta(5)  \right. \nonumber \\ &\left. \hspace{1em}
+7488 M(2,6)\right) \label{eq09116} \\
\sum_{k=1}^\infty \frac{H(k)^{4}}{k^{3}(k+1)(k+2)} \sumend &= \frac{-1}{32}\left(  4 +12\zeta(2) +44\zeta(3) +74\zeta(4) -840\zeta(5)  \right. \nonumber \\ &\left. \hspace{1em}
-168\zeta(2)\zeta(3) +979\zeta(6) +72\zeta(3)^2 -370\zeta(7) -80\zeta(2)\zeta(5)  \right. \nonumber \\ &\left. \hspace{1em}
+344\zeta(3)\zeta(4)\right) \label{eq09117} \\
\sum_{k=1}^\infty \frac{H(k)^{4}}{k^{2}(k+1)^{2}(k+2)} \sumend &= \frac{-1}{16}\left(  4 +12\zeta(2) +44\zeta(3) +74\zeta(4) +600\zeta(5)  \right. \nonumber \\ &\left. \hspace{1em}
+120\zeta(2)\zeta(3) -899\zeta(6) -72\zeta(3)^2\right) \label{eq09118} \\
\sum_{k=1}^\infty \frac{H(k)^{4}}{k(k+1)^{3}(k+2)} \sumend &= \frac{1}{8}\left( -4 -12\zeta(2) -44\zeta(3) -74\zeta(4) +120\zeta(5) +24\zeta(2)\zeta(3)  \right. \nonumber \\ &\left. \hspace{1em}
+109\zeta(7) +40\zeta(2)\zeta(5) -148\zeta(3)\zeta(4)\right) \label{eq09119} \\
\sum_{k=1}^\infty \frac{H(k)^{4}}{(k+1)^{4}(k+2)} \sumend &= \frac{1}{144}\left( -144 -432\zeta(2) -1584\zeta(3) -2664\zeta(4) +5154\zeta(6)  \right. \nonumber \\ &\left. \hspace{1em}
+432\zeta(3)^2 +1962\zeta(7) +720\zeta(2)\zeta(5) -2664\zeta(3)\zeta(4) -12415\zeta(8)  \right. \nonumber \\ &\left. \hspace{1em}
-3312\zeta(2)\zeta(3)^2 +13824\zeta(3)\zeta(5) +3024 M(2,6)\right) \label{eq09120} \\
\sum_{k=1}^\infty \frac{H(k)^{4}}{k^{3}(k+2)^{2}} \sumend &= \frac{1}{64}\left(  52 +100\zeta(2) +276\zeta(3) +238\zeta(4) +120\zeta(5)  \right. \nonumber \\ &\left. \hspace{1em}
+40\zeta(2)\zeta(3) -939\zeta(6) -72\zeta(3)^2 +370\zeta(7) +80\zeta(2)\zeta(5)  \right. \nonumber \\ &\left. \hspace{1em}
-344\zeta(3)\zeta(4)\right) \label{eq09121} \\
\sum_{k=1}^\infty \frac{H(k)^{4}}{k^{2}(k+1)(k+2)^{2}} \sumend &= \frac{-1}{4}\left( -7 -14\zeta(2) -40\zeta(3) -39\zeta(4) +90\zeta(5)  \right. \nonumber \\ &\left. \hspace{1em}
+16\zeta(2)\zeta(3) -5\zeta(6)\right) \label{eq09122} \\
\sum_{k=1}^\infty \frac{H(k)^{4}}{k(k+1)^{2}(k+2)^{2}} \sumend &= \frac{-1}{16}\left( -60 -124\zeta(2) -364\zeta(3) -386\zeta(4) +120\zeta(5)  \right. \nonumber \\ &\left. \hspace{1em}
+8\zeta(2)\zeta(3) +859\zeta(6) +72\zeta(3)^2\right) \label{eq09123} \\
\sum_{k=1}^\infty \frac{H(k)^{4}}{(k+1)^{3}(k+2)^{2}} \sumend &= \frac{1}{8}\left(  64 +136\zeta(2) +408\zeta(3) +460\zeta(4) -240\zeta(5)  \right. \nonumber \\ &\left. \hspace{1em}
-32\zeta(2)\zeta(3) -859\zeta(6) -72\zeta(3)^2 -109\zeta(7) -40\zeta(2)\zeta(5)  \right. \nonumber \\ &\left. \hspace{1em}
+148\zeta(3)\zeta(4)\right) \label{eq09124}
\end{align}
 
\begin{align}
\sum_{k=1}^\infty \frac{H(k)^{4}}{k^{2}(k+2)^{3}} \sumend &= \frac{-1}{64}\left(  332 +388\zeta(2) +740\zeta(3) +6\zeta(4) -456\zeta(5)  \right. \nonumber \\ &\left. \hspace{1em}
-152\zeta(2)\zeta(3) -767\zeta(6) -200\zeta(3)^2 +218\zeta(7) +80\zeta(2)\zeta(5)  \right. \nonumber \\ &\left. \hspace{1em}
-296\zeta(3)\zeta(4)\right) \label{eq09125} \\
\sum_{k=1}^\infty \frac{H(k)^{4}}{k(k+1)(k+2)^{3}} \sumend &= \frac{-1}{32}\left(  388 +500\zeta(2) +1060\zeta(3) +318\zeta(4) -1176\zeta(5)  \right. \nonumber \\ &\left. \hspace{1em}
-280\zeta(2)\zeta(3) -727\zeta(6) -200\zeta(3)^2 +218\zeta(7) +80\zeta(2)\zeta(5)  \right. \nonumber \\ &\left. \hspace{1em}
-296\zeta(3)\zeta(4)\right) \label{eq09126} \\
\sum_{k=1}^\infty \frac{H(k)^{4}}{(k+1)^{2}(k+2)^{3}} \sumend &= \frac{1}{8}\left( -224 -312\zeta(2) -712\zeta(3) -352\zeta(4) +648\zeta(5)  \right. \nonumber \\ &\left. \hspace{1em}
+144\zeta(2)\zeta(3) +793\zeta(6) +136\zeta(3)^2 -109\zeta(7) -40\zeta(2)\zeta(5)  \right. \nonumber \\ &\left. \hspace{1em}
+148\zeta(3)\zeta(4)\right) \label{eq09127} \\
\sum_{k=1}^\infty \frac{H(k)^{4}}{k(k+2)^{4}} \sumend &= \frac{1}{576}\left(  12636 +8460\zeta(2) +9180\zeta(3) -11502\zeta(4) -9864\zeta(5)  \right. \nonumber \\ &\left. \hspace{1em}
-4968\zeta(2)\zeta(3) -1677\zeta(6) -1944\zeta(3)^2 -6606\zeta(7) -1584\zeta(2)\zeta(5)  \right. \nonumber \\ &\left. \hspace{1em}
+6840\zeta(3)\zeta(4) +24830\zeta(8) +6624\zeta(2)\zeta(3)^2 -27648\zeta(3)\zeta(5)  \right. \nonumber \\ &\left. \hspace{1em}
-6048 M(2,6)\right) \label{eq09128} \\
\sum_{k=1}^\infty \frac{H(k)^{4}}{(k+1)(k+2)^{4}} \sumend &= \frac{-1}{144}\left( -8064 -6480\zeta(2) -9360\zeta(3) +4320\zeta(4) +10224\zeta(5)  \right. \nonumber \\ &\left. \hspace{1em}
+3744\zeta(2)\zeta(3) +4110\zeta(6) +1872\zeta(3)^2 +2322\zeta(7) +432\zeta(2)\zeta(5)  \right. \nonumber \\ &\left. \hspace{1em}
-2088\zeta(3)\zeta(4) -12415\zeta(8) -3312\zeta(2)\zeta(3)^2 +13824\zeta(3)\zeta(5)  \right. \nonumber \\ &\left. \hspace{1em}
+3024 M(2,6)\right) \label{eq09129} \\
\sum_{k=1}^\infty \frac{H(k)^{4}}{(k+2)^{5}} \sumend &= \frac{-1}{24}\left(  1680 +600\zeta(2) +120\zeta(3) -1380\zeta(4) -672\zeta(5)  \right. \nonumber \\ &\left. \hspace{1em}
-240\zeta(2)\zeta(3) -318\zeta(6) +144\zeta(3)^2 -570\zeta(7) -336\zeta(2)\zeta(5) +864\zeta(3)\zeta(4)  \right. \nonumber \\ &\left. \hspace{1em}
+43\zeta(8) +120\zeta(2)\zeta(3)^2 -288\zeta(3)\zeta(5) +24 M(2,6) +348\zeta(9) -198\zeta(3)\zeta(6)  \right. \nonumber \\ &\left. \hspace{1em}
-444\zeta(4)\zeta(5) +168\zeta(2)\zeta(7) +64\zeta(3)^3\right) \label{eq09130} \\
\sum_{k=1}^\infty \frac{H(k)^{5}}{k^{4}} \sumend &= \frac{1}{72}\left(  9442\zeta(9) -14685\zeta(3)\zeta(6) +4752\zeta(4)\zeta(5) +2385\zeta(2)\zeta(7)  \right. \nonumber \\ &\left. \hspace{1em}
-360\zeta(3)^3\right) \label{eq09131} \\
\sum_{k=1}^\infty \frac{H(k)^{5}}{k^{3}(k+1)} \sumend &= \frac{1}{288}\left(  51408\zeta(6) +6480\zeta(3)^2 -36918\zeta(7) -8208\zeta(2)\zeta(5)  \right. \nonumber \\ &\left. \hspace{1em}
-9504\zeta(3)\zeta(4) -67811\zeta(8) -19080\zeta(2)\zeta(3)^2 +78768\zeta(3)\zeta(5)  \right. \nonumber \\ &\left. \hspace{1em}
+16920 M(2,6)\right) \label{eq09132}
\end{align}
 
\begin{align}
\sum_{k=1}^\infty \frac{H(k)^{5}}{k^{2}(k+1)^{2}} \sumend &= \frac{1}{8}\left( -2856\zeta(6) -360\zeta(3)^2 +1953\zeta(7) +456\zeta(2)\zeta(5)  \right. \nonumber \\ &\left. \hspace{1em}
+528\zeta(3)\zeta(4)\right) \label{eq09133} \\
\sum_{k=1}^\infty \frac{H(k)^{5}}{k(k+1)^{3}} \sumend &= \frac{1}{288}\left(  51408\zeta(6) +6480\zeta(3)^2 -33390\zeta(7) -8208\zeta(2)\zeta(5)  \right. \nonumber \\ &\left. \hspace{1em}
-9504\zeta(3)\zeta(4) -65621\zeta(8) -17640\zeta(2)\zeta(3)^2 +72432\zeta(3)\zeta(5)  \right. \nonumber \\ &\left. \hspace{1em}
+15480 M(2,6)\right) \label{eq09134} \\
\sum_{k=1}^\infty \frac{H(k)^{5}}{(k+1)^{4}} \sumend &= \frac{1}{72}\left(  7120\zeta(9) -12885\zeta(3)\zeta(6) +4752\zeta(4)\zeta(5)  \right. \nonumber \\ &\left. \hspace{1em}
+2385\zeta(2)\zeta(7) -360\zeta(3)^3\right) \label{eq09135} \\
\sum_{k=1}^\infty \frac{H(k)^{5}}{k^{3}(k+2)} \sumend &= \frac{1}{576}\left(  72 +288\zeta(2) +1512\zeta(3) +4518\zeta(4) +5112\zeta(5)  \right. \nonumber \\ &\left. \hspace{1em}
+1080\zeta(2)\zeta(3) +12852\zeta(6) +1620\zeta(3)^2 -18459\zeta(7) -4104\zeta(2)\zeta(5)  \right. \nonumber \\ &\left. \hspace{1em}
-4752\zeta(3)\zeta(4) -67811\zeta(8) -19080\zeta(2)\zeta(3)^2 +78768\zeta(3)\zeta(5)  \right. \nonumber \\ &\left. \hspace{1em}
+16920 M(2,6)\right) \label{eq09136} \\
\sum_{k=1}^\infty \frac{H(k)^{5}}{k^{2}(k+1)(k+2)} \sumend &= \frac{-1}{32}\left( -8 -32\zeta(2) -168\zeta(3) -502\zeta(4) -568\zeta(5)  \right. \nonumber \\ &\left. \hspace{1em}
-120\zeta(2)\zeta(3) +4284\zeta(6) +540\zeta(3)^2 -2051\zeta(7) -456\zeta(2)\zeta(5)  \right. \nonumber \\ &\left. \hspace{1em}
-528\zeta(3)\zeta(4)\right) \label{eq09137} \\
\sum_{k=1}^\infty \frac{H(k)^{5}}{k(k+1)^{2}(k+2)} \sumend &= \frac{-1}{16}\left( -8 -32\zeta(2) -168\zeta(3) -502\zeta(4) -568\zeta(5)  \right. \nonumber \\ &\left. \hspace{1em}
-120\zeta(2)\zeta(3) -1428\zeta(6) -180\zeta(3)^2 +1855\zeta(7) +456\zeta(2)\zeta(5)  \right. \nonumber \\ &\left. \hspace{1em}
+528\zeta(3)\zeta(4)\right) \label{eq09138} \\
\sum_{k=1}^\infty \frac{H(k)^{5}}{(k+1)^{3}(k+2)} \sumend &= \frac{1}{288}\left(  288 +1152\zeta(2) +6048\zeta(3) +18072\zeta(4) +20448\zeta(5)  \right. \nonumber \\ &\left. \hspace{1em}
+4320\zeta(2)\zeta(3) -33390\zeta(7) -8208\zeta(2)\zeta(5) -9504\zeta(3)\zeta(4) +65621\zeta(8)  \right. \nonumber \\ &\left. \hspace{1em}
+17640\zeta(2)\zeta(3)^2 -72432\zeta(3)\zeta(5) -15480 M(2,6)\right) \label{eq09139} \\
\sum_{k=1}^\infty \frac{H(k)^{5}}{k^{2}(k+2)^{2}} \sumend &= \frac{-1}{96}\left(  168 +456\zeta(2) +1896\zeta(3) +3858\zeta(4) +1608\zeta(5)  \right. \nonumber \\ &\left. \hspace{1em}
+480\zeta(2)\zeta(3) -11\zeta(6) +180\zeta(3)^2 -5859\zeta(7) -1368\zeta(2)\zeta(5)  \right. \nonumber \\ &\left. \hspace{1em}
-1584\zeta(3)\zeta(4)\right) \label{eq09140} \\
\sum_{k=1}^\infty \frac{H(k)^{5}}{k(k+1)(k+2)^{2}} \sumend &= \frac{-1}{96}\left(  360 +1008\zeta(2) +4296\zeta(3) +9222\zeta(4) +4920\zeta(5)  \right. \nonumber \\ &\left. \hspace{1em}
+1320\zeta(2)\zeta(3) -12874\zeta(6) -1260\zeta(3)^2 -5565\zeta(7) -1368\zeta(2)\zeta(5)  \right. \nonumber \\ &\left. \hspace{1em}
-1584\zeta(3)\zeta(4)\right) \label{eq09141}
\end{align}
 
\begin{align}
\sum_{k=1}^\infty \frac{H(k)^{5}}{(k+1)^{2}(k+2)^{2}} \sumend &= \frac{-1}{24}\left(  192 +552\zeta(2) +2400\zeta(3) +5364\zeta(4) +3312\zeta(5)  \right. \nonumber \\ &\left. \hspace{1em}
+840\zeta(2)\zeta(3) -4295\zeta(6) -360\zeta(3)^2 -5565\zeta(7) -1368\zeta(2)\zeta(5)  \right. \nonumber \\ &\left. \hspace{1em}
-1584\zeta(3)\zeta(4)\right) \label{eq09142} \\
\sum_{k=1}^\infty \frac{H(k)^{5}}{k(k+2)^{3}} \sumend &= \frac{1}{576}\left(  6984 +12528\zeta(2) +40104\zeta(3) +48726\zeta(4) -13896\zeta(5)  \right. \nonumber \\ &\left. \hspace{1em}
-2520\zeta(2)\zeta(3) -58518\zeta(6) -10620\zeta(3)^2 +2925\zeta(7) +3096\zeta(2)\zeta(5)  \right. \nonumber \\ &\left. \hspace{1em}
-31392\zeta(3)\zeta(4) -65621\zeta(8) -17640\zeta(2)\zeta(3)^2 +72432\zeta(3)\zeta(5)  \right. \nonumber \\ &\left. \hspace{1em}
+15480 M(2,6)\right) \label{eq09143} \\
\sum_{k=1}^\infty \frac{H(k)^{5}}{(k+1)(k+2)^{3}} \sumend &= \frac{-1}{288}\left( -8064 -15552\zeta(2) -52992\zeta(3) -76392\zeta(4) -864\zeta(5)  \right. \nonumber \\ &\left. \hspace{1em}
-1440\zeta(2)\zeta(3) +97140\zeta(6) +14400\zeta(3)^2 +13770\zeta(7) +1008\zeta(2)\zeta(5)  \right. \nonumber \\ &\left. \hspace{1em}
+36144\zeta(3)\zeta(4) +65621\zeta(8) +17640\zeta(2)\zeta(3)^2 -72432\zeta(3)\zeta(5)  \right. \nonumber \\ &\left. \hspace{1em}
-15480 M(2,6)\right) \label{eq09144} \\
\sum_{k=1}^\infty \frac{H(k)^{5}}{(k+2)^{4}} \sumend &= \frac{1}{144}\left( -8064 -9360\zeta(2) -22320\zeta(3) -11520\zeta(4) +17136\zeta(5)  \right. \nonumber \\ &\left. \hspace{1em}
+5760\zeta(2)\zeta(3) +22050\zeta(6) +6480\zeta(3)^2 +900\zeta(7) -720\zeta(2)\zeta(5)  \right. \nonumber \\ &\left. \hspace{1em}
+1440\zeta(3)\zeta(4) -62075\zeta(8) -16560\zeta(2)\zeta(3)^2 +69120\zeta(3)\zeta(5) +15120 M(2,6)  \right. \nonumber \\ &\left. \hspace{1em}
+14240\zeta(9) -25770\zeta(3)\zeta(6) +9504\zeta(4)\zeta(5) +4770\zeta(2)\zeta(7)  \right. \nonumber \\ &\left. \hspace{1em}
-720\zeta(3)^3\right) \label{eq09145} \\
\sum_{k=1}^\infty \frac{H(k)^{6}}{k^{3}} \sumend &= \frac{-1}{24}\left( -7474\zeta(9) +13122\zeta(3)\zeta(6) -6048\zeta(4)\zeta(5) -1953\zeta(2)\zeta(7)  \right. \nonumber \\ &\left. \hspace{1em}
+544\zeta(3)^3\right) \label{eq09146} \\
\sum_{k=1}^\infty \frac{H(k)^{6}}{k^{2}(k+1)} \sumend &= \frac{-1}{8}\left(  5152\zeta(7) +1160\zeta(2)\zeta(5) +2376\zeta(3)\zeta(4) -5843\zeta(8)  \right. \nonumber \\ &\left. \hspace{1em}
+328\zeta(2)\zeta(3)^2 -3896\zeta(3)\zeta(5) -456 M(2,6)\right) \label{eq09147} \\
\sum_{k=1}^\infty \frac{H(k)^{6}}{k(k+1)^{2}} \sumend &= \frac{1}{24}\left(  15456\zeta(7) +3480\zeta(2)\zeta(5) +7128\zeta(3)\zeta(4) -17027\zeta(8)  \right. \nonumber \\ &\left. \hspace{1em}
+924\zeta(2)\zeta(3)^2 -11328\zeta(3)\zeta(5) -1308 M(2,6)\right) \label{eq09148} \\
\sum_{k=1}^\infty \frac{H(k)^{6}}{(k+1)^{3}} \sumend &= \frac{-1}{24}\left(  6146\zeta(9) -12582\zeta(3)\zeta(6) +5832\zeta(4)\zeta(5)  \right. \nonumber \\ &\left. \hspace{1em}
+1953\zeta(2)\zeta(7) -536\zeta(3)^3\right) \label{eq09149}
\end{align}
 
\begin{align}
\sum_{k=1}^\infty \frac{H(k)^{6}}{k^{2}(k+2)} \sumend &= \frac{1}{16}\left( -4 -20\zeta(2) -136\zeta(3) -571\zeta(4) -1142\zeta(5) -244\zeta(2)\zeta(3)  \right. \nonumber \\ &\left. \hspace{1em}
-2097\zeta(6) -268\zeta(3)^2 -2576\zeta(7) -580\zeta(2)\zeta(5) -1188\zeta(3)\zeta(4) +5843\zeta(8)  \right. \nonumber \\ &\left. \hspace{1em}
-328\zeta(2)\zeta(3)^2 +3896\zeta(3)\zeta(5) +456 M(2,6)\right) \label{eq09150} \\
\sum_{k=1}^\infty \frac{H(k)^{6}}{k(k+1)(k+2)} \sumend &= \frac{-1}{8}\left(  4 +20\zeta(2) +136\zeta(3) +571\zeta(4) +1142\zeta(5) +244\zeta(2)\zeta(3)  \right. \nonumber \\ &\left. \hspace{1em}
+2097\zeta(6) +268\zeta(3)^2 -2576\zeta(7) -580\zeta(2)\zeta(5) -1188\zeta(3)\zeta(4)\right) \label{eq09151} \\
\sum_{k=1}^\infty \frac{H(k)^{6}}{(k+1)^{2}(k+2)} \sumend &= \frac{-1}{24}\left(  24 +120\zeta(2) +816\zeta(3) +3426\zeta(4) +6852\zeta(5)  \right. \nonumber \\ &\left. \hspace{1em}
+1464\zeta(2)\zeta(3) +12582\zeta(6) +1608\zeta(3)^2 -17027\zeta(8) +924\zeta(2)\zeta(3)^2  \right. \nonumber \\ &\left. \hspace{1em}
-11328\zeta(3)\zeta(5) -1308 M(2,6)\right) \label{eq09152} \\
\sum_{k=1}^\infty \frac{H(k)^{6}}{k(k+2)^{2}} \sumend &= \frac{-1}{48}\left( -180 -636\zeta(2) -3528\zeta(3) -11001\zeta(4) -13530\zeta(5)  \right. \nonumber \\ &\left. \hspace{1em}
-3180\zeta(2)\zeta(3) -5988\zeta(6) -1332\zeta(3)^2 +8967\zeta(7) +2364\zeta(2)\zeta(5)  \right. \nonumber \\ &\left. \hspace{1em}
+1188\zeta(3)\zeta(4) +17027\zeta(8) -924\zeta(2)\zeta(3)^2 +11328\zeta(3)\zeta(5) +1308 M(2,6)\right) \label{eq09153} \\
\sum_{k=1}^\infty \frac{H(k)^{6}}{(k+1)(k+2)^{2}} \sumend &= \frac{1}{24}\left(  192 +696\zeta(2) +3936\zeta(3) +12714\zeta(4) +16956\zeta(5)  \right. \nonumber \\ &\left. \hspace{1em}
+3912\zeta(2)\zeta(3) +12279\zeta(6) +2136\zeta(3)^2 -16695\zeta(7) -4104\zeta(2)\zeta(5)  \right. \nonumber \\ &\left. \hspace{1em}
-4752\zeta(3)\zeta(4) -17027\zeta(8) +924\zeta(2)\zeta(3)^2 -11328\zeta(3)\zeta(5) -1308 M(2,6)\right) \label{eq09154} \\
\sum_{k=1}^\infty \frac{H(k)^{6}}{(k+2)^{3}} \sumend &= \frac{1}{48}\left( -1344 -3312\zeta(2) -14832\zeta(3) -33120\zeta(4) -20592\zeta(5)  \right. \nonumber \\ &\left. \hspace{1em}
-5184\zeta(2)\zeta(3) +24396\zeta(6) +3024\zeta(3)^2 +23580\zeta(7) +4608\zeta(2)\zeta(5)  \right. \nonumber \\ &\left. \hspace{1em}
+22824\zeta(3)\zeta(4) +65621\zeta(8) +17640\zeta(2)\zeta(3)^2 -72432\zeta(3)\zeta(5) -15480 M(2,6)  \right. \nonumber \\ &\left. \hspace{1em}
-12292\zeta(9) +25164\zeta(3)\zeta(6) -11664\zeta(4)\zeta(5) -3906\zeta(2)\zeta(7)  \right. \nonumber \\ &\left. \hspace{1em}
+1072\zeta(3)^3\right) \label{eq09155} \\
\sum_{k=1}^\infty \frac{H(k)^{7}}{k^{2}} \sumend &= \frac{-1}{72}\left( -276341\zeta(9) -88665\zeta(3)\zeta(6) -143163\zeta(4)\zeta(5)  \right. \nonumber \\ &\left. \hspace{1em}
-59166\zeta(2)\zeta(7) -4032\zeta(3)^3\right) \label{eq09156} \\
\sum_{k=1}^\infty \frac{H(k)^{7}}{k(k+1)} \sumend &= \frac{-1}{18}\left( -119774\zeta(8) -3024\zeta(2)\zeta(3)^2 -27405\zeta(3)\zeta(5)\right) \label{eq09157} \\
\sum_{k=1}^\infty \frac{H(k)^{7}}{(k+1)^{2}} \sumend &= \frac{-1}{72}\left( -269402\zeta(9) -88665\zeta(3)\zeta(6) -141273\zeta(4)\zeta(5)  \right. \nonumber \\ &\left. \hspace{1em}
-59166\zeta(2)\zeta(7) -4032\zeta(3)^3\right) \label{eq09158}
\end{align}
 
\begin{align}
\sum_{k=1}^\infty \frac{H(k)^{7}}{k(k+2)} \sumend &= \frac{1}{288}\left(  144 +864\zeta(2) +7200\zeta(3) +38664\zeta(4) +108504\zeta(5)  \right. \nonumber \\ &\left. \hspace{1em}
+23184\zeta(2)\zeta(3) +352887\zeta(6) +45864\zeta(3)^2 +319554\zeta(7) +73080\zeta(2)\zeta(5)  \right. \nonumber \\ &\left. \hspace{1em}
+148932\zeta(3)\zeta(4) +958192\zeta(8) +24192\zeta(2)\zeta(3)^2 +219240\zeta(3)\zeta(5)\right) \label{eq09159} \\
\sum_{k=1}^\infty \frac{H(k)^{7}}{(k+1)(k+2)} \sumend &= \frac{-1}{48}\left( -48 -288\zeta(2) -2400\zeta(3) -12888\zeta(4) -36168\zeta(5)  \right. \nonumber \\ &\left. \hspace{1em}
-7728\zeta(2)\zeta(3) -117629\zeta(6) -15288\zeta(3)^2 -106518\zeta(7) -24360\zeta(2)\zeta(5)  \right. \nonumber \\ &\left. \hspace{1em}
-49644\zeta(3)\zeta(4)\right) \label{eq09160} \\
\sum_{k=1}^\infty \frac{H(k)^{7}}{(k+2)^{2}} \sumend &= \frac{-1}{144}\left(  1152 +5040\zeta(2) +34992\zeta(3) +146340\zeta(4) +287712\zeta(5)  \right. \nonumber \\ &\left. \hspace{1em}
+64512\zeta(2)\zeta(3) +525384\zeta(6) +76608\zeta(3)^2 -31041\zeta(7) -13104\zeta(2)\zeta(5)  \right. \nonumber \\ &\left. \hspace{1em}
+49140\zeta(3)\zeta(4) -715134\zeta(8) +38808\zeta(2)\zeta(3)^2 -475776\zeta(3)\zeta(5) -54936 M(2,6)  \right. \nonumber \\ &\left. \hspace{1em}
-538804\zeta(9) -177330\zeta(3)\zeta(6) -282546\zeta(4)\zeta(5) -118332\zeta(2)\zeta(7)  \right. \nonumber \\ &\left. \hspace{1em}
-8064\zeta(3)^3\right) \label{eq09161}
\end{align}

\newpage
Formulas for order $r = m + n + p + q = 10$:
\begin{align}
\sum_{k=1}^\infty \frac{H(k)}{k^{9}} \sumend &= \frac{1}{4}\left(  11\zeta(10) -4\zeta(3)\zeta(7) -2\zeta(5)^2\right) \label{eq10001} \\
\sum_{k=1}^\infty \frac{H(k)}{k^{8}(k+1)} \sumend &= \frac{1}{4}\left( -4\zeta(2) +8\zeta(3) -5\zeta(4) +12\zeta(5) -4\zeta(2)\zeta(3) -7\zeta(6)  \right. \nonumber \\ &\left. \hspace{1em}
+2\zeta(3)^2 +16\zeta(7) -4\zeta(2)\zeta(5) -4\zeta(3)\zeta(4) -9\zeta(8) +4\zeta(3)\zeta(5)  \right. \nonumber \\ &\left. \hspace{1em}
+20\zeta(9) -4\zeta(3)\zeta(6) -4\zeta(4)\zeta(5) -4\zeta(2)\zeta(7)\right) \label{eq10002} \\
\sum_{k=1}^\infty \frac{H(k)}{k^{7}(k+1)^{2}} \sumend &= \frac{1}{4}\left(  28\zeta(2) -52\zeta(3) +25\zeta(4) -48\zeta(5) +16\zeta(2)\zeta(3)  \right. \nonumber \\ &\left. \hspace{1em}
+21\zeta(6) -6\zeta(3)^2 -32\zeta(7) +8\zeta(2)\zeta(5) +8\zeta(3)\zeta(4) +9\zeta(8)  \right. \nonumber \\ &\left. \hspace{1em}
-4\zeta(3)\zeta(5)\right) \label{eq10003} \\
\sum_{k=1}^\infty \frac{H(k)}{k^{6}(k+1)^{3}} \sumend &= \frac{-1}{4}\left(  84\zeta(2) -144\zeta(3) +49\zeta(4) -72\zeta(5) +24\zeta(2)\zeta(3)  \right. \nonumber \\ &\left. \hspace{1em}
+21\zeta(6) -6\zeta(3)^2 -16\zeta(7) +4\zeta(2)\zeta(5) +4\zeta(3)\zeta(4)\right) \label{eq10004} \\
\sum_{k=1}^\infty \frac{H(k)}{k^{5}(k+1)^{4}} \sumend &= \frac{-1}{4}\left( -140\zeta(2) +220\zeta(3) -45\zeta(4) +56\zeta(5) -20\zeta(2)\zeta(3)  \right. \nonumber \\ &\left. \hspace{1em}
-7\zeta(6) +2\zeta(3)^2\right) \label{eq10005} \\
\sum_{k=1}^\infty \frac{H(k)}{k^{4}(k+1)^{5}} \sumend &= \frac{-1}{4}\left(  140\zeta(2) -200\zeta(3) +15\zeta(4) -44\zeta(5) +20\zeta(2)\zeta(3)  \right. \nonumber \\ &\left. \hspace{1em}
-3\zeta(6) +2\zeta(3)^2\right) \label{eq10006} \\
\sum_{k=1}^\infty \frac{H(k)}{k^{3}(k+1)^{6}} \sumend &= \frac{-1}{4}\left( -84\zeta(2) +108\zeta(3) +5\zeta(4) +48\zeta(5) -24\zeta(2)\zeta(3)  \right. \nonumber \\ &\left. \hspace{1em}
+9\zeta(6) -6\zeta(3)^2 +12\zeta(7) -4\zeta(2)\zeta(5) -4\zeta(3)\zeta(4)\right) \label{eq10007} \\
\sum_{k=1}^\infty \frac{H(k)}{k^{2}(k+1)^{7}} \sumend &= \frac{1}{4}\left( -28\zeta(2) +32\zeta(3) +5\zeta(4) +32\zeta(5) -16\zeta(2)\zeta(3) +9\zeta(6)  \right. \nonumber \\ &\left. \hspace{1em}
-6\zeta(3)^2 +24\zeta(7) -8\zeta(2)\zeta(5) -8\zeta(3)\zeta(4) +5\zeta(8) -4\zeta(3)\zeta(5)\right) \label{eq10008} \\
\sum_{k=1}^\infty \frac{H(k)}{k(k+1)^{8}} \sumend &= \frac{-1}{4}\left( -4\zeta(2) +4\zeta(3) +\zeta(4) +8\zeta(5) -4\zeta(2)\zeta(3) +3\zeta(6)  \right. \nonumber \\ &\left. \hspace{1em}
-2\zeta(3)^2 +12\zeta(7) -4\zeta(2)\zeta(5) -4\zeta(3)\zeta(4) +5\zeta(8) -4\zeta(3)\zeta(5)  \right. \nonumber \\ &\left. \hspace{1em}
+16\zeta(9) -4\zeta(3)\zeta(6) -4\zeta(4)\zeta(5) -4\zeta(2)\zeta(7)\right) \label{eq10009} \\
\sum_{k=1}^\infty \frac{H(k)}{(k+1)^{9}} \sumend &= \frac{-1}{4}\left( -7\zeta(10) +4\zeta(3)\zeta(7) +2\zeta(5)^2\right) \label{eq10010} \\
\sum_{k=1}^\infty \frac{H(k)}{k^{8}(k+2)} \sumend &= \frac{-1}{256}\left(  1 +\zeta(2) -4\zeta(3) +5\zeta(4) -24\zeta(5) +8\zeta(2)\zeta(3) +28\zeta(6)  \right. \nonumber \\ &\left. \hspace{1em}
-8\zeta(3)^2 -128\zeta(7) +32\zeta(2)\zeta(5) +32\zeta(3)\zeta(4) +144\zeta(8) -64\zeta(3)\zeta(5)  \right. \nonumber \\ &\left. \hspace{1em}
-640\zeta(9) +128\zeta(3)\zeta(6) +128\zeta(4)\zeta(5) +128\zeta(2)\zeta(7)\right) \label{eq10011}
\end{align}
 
\begin{align}
\sum_{k=1}^\infty \frac{H(k)}{k^{7}(k+1)(k+2)} \sumend &= \frac{1}{128}\left( -1 +127\zeta(2) -252\zeta(3) +155\zeta(4) -360\zeta(5) +120\zeta(2)\zeta(3)  \right. \nonumber \\ &\left. \hspace{1em}
+196\zeta(6) -56\zeta(3)^2 -384\zeta(7) +96\zeta(2)\zeta(5) +96\zeta(3)\zeta(4) +144\zeta(8)  \right. \nonumber \\ &\left. \hspace{1em}
-64\zeta(3)\zeta(5)\right) \label{eq10012} \\
\sum_{k=1}^\infty \frac{H(k)}{k^{6}(k+1)^{2}(k+2)} \sumend &= \frac{1}{64}\left( -1 -321\zeta(2) +580\zeta(3) -245\zeta(4) +408\zeta(5)  \right. \nonumber \\ &\left. \hspace{1em}
-136\zeta(2)\zeta(3) -140\zeta(6) +40\zeta(3)^2 +128\zeta(7) -32\zeta(2)\zeta(5)  \right. \nonumber \\ &\left. \hspace{1em}
-32\zeta(3)\zeta(4)\right) \label{eq10013} \\
\sum_{k=1}^\infty \frac{H(k)}{k^{5}(k+1)^{3}(k+2)} \sumend &= \frac{-1}{32}\left(  1 -351\zeta(2) +572\zeta(3) -147\zeta(4) +168\zeta(5)  \right. \nonumber \\ &\left. \hspace{1em}
-56\zeta(2)\zeta(3) -28\zeta(6) +8\zeta(3)^2\right) \label{eq10014} \\
\sum_{k=1}^\infty \frac{H(k)}{k^{4}(k+1)^{4}(k+2)} \sumend &= \frac{1}{16}\left( -1 -209\zeta(2) +308\zeta(3) -33\zeta(4) +56\zeta(5)  \right. \nonumber \\ &\left. \hspace{1em}
-24\zeta(2)\zeta(3)\right) \label{eq10015} \\
\sum_{k=1}^\infty \frac{H(k)}{k^{3}(k+1)^{5}(k+2)} \sumend &= \frac{-1}{8}\left(  1 -71\zeta(2) +92\zeta(3) +3\zeta(4) +32\zeta(5) -16\zeta(2)\zeta(3)  \right. \nonumber \\ &\left. \hspace{1em}
+6\zeta(6) -4\zeta(3)^2\right) \label{eq10016} \\
\sum_{k=1}^\infty \frac{H(k)}{k^{2}(k+1)^{6}(k+2)} \sumend &= \frac{-1}{4}\left(  1 +13\zeta(2) -16\zeta(3) -2\zeta(4) -16\zeta(5) +8\zeta(2)\zeta(3)  \right. \nonumber \\ &\left. \hspace{1em}
-3\zeta(6) +2\zeta(3)^2 -12\zeta(7) +4\zeta(2)\zeta(5) +4\zeta(3)\zeta(4)\right) \label{eq10017} \\
\sum_{k=1}^\infty \frac{H(k)}{k(k+1)^{7}(k+2)} \sumend &= \frac{-1}{4}\left(  2 -2\zeta(2) +\zeta(4) +3\zeta(6) -2\zeta(3)^2 +5\zeta(8)  \right. \nonumber \\ &\left. \hspace{1em}
-4\zeta(3)\zeta(5)\right) \label{eq10018} \\
\sum_{k=1}^\infty \frac{H(k)}{(k+1)^{8}(k+2)} \sumend &= \frac{1}{4}\left( -4 +4\zeta(3) -\zeta(4) +8\zeta(5) -4\zeta(2)\zeta(3) -3\zeta(6) +2\zeta(3)^2  \right. \nonumber \\ &\left. \hspace{1em}
+12\zeta(7) -4\zeta(2)\zeta(5) -4\zeta(3)\zeta(4) -5\zeta(8) +4\zeta(3)\zeta(5) +16\zeta(9)  \right. \nonumber \\ &\left. \hspace{1em}
-4\zeta(3)\zeta(6) -4\zeta(4)\zeta(5) -4\zeta(2)\zeta(7)\right) \label{eq10019} \\
\sum_{k=1}^\infty \frac{H(k)}{k^{7}(k+2)^{2}} \sumend &= \frac{1}{256}\left(  11 +5\zeta(2) -26\zeta(3) +25\zeta(4) -96\zeta(5) +32\zeta(2)\zeta(3)  \right. \nonumber \\ &\left. \hspace{1em}
+84\zeta(6) -24\zeta(3)^2 -256\zeta(7) +64\zeta(2)\zeta(5) +64\zeta(3)\zeta(4) +144\zeta(8)  \right. \nonumber \\ &\left. \hspace{1em}
-64\zeta(3)\zeta(5)\right) \label{eq10020} \\
\sum_{k=1}^\infty \frac{H(k)}{k^{6}(k+1)(k+2)^{2}} \sumend &= \frac{-1}{64}\left( -6 +61\zeta(2) -113\zeta(3) +65\zeta(4) -132\zeta(5)  \right. \nonumber \\ &\left. \hspace{1em}
+44\zeta(2)\zeta(3) +56\zeta(6) -16\zeta(3)^2 -64\zeta(7) +16\zeta(2)\zeta(5)  \right. \nonumber \\ &\left. \hspace{1em}
+16\zeta(3)\zeta(4)\right) \label{eq10021}
\end{align}
 
\begin{align}
\sum_{k=1}^\infty \frac{H(k)}{k^{5}(k+1)^{2}(k+2)^{2}} \sumend &= \frac{1}{64}\left(  13 +199\zeta(2) -354\zeta(3) +115\zeta(4) -144\zeta(5)  \right. \nonumber \\ &\left. \hspace{1em}
+48\zeta(2)\zeta(3) +28\zeta(6) -8\zeta(3)^2\right) \label{eq10022} \\
\sum_{k=1}^\infty \frac{H(k)}{k^{4}(k+1)^{3}(k+2)^{2}} \sumend &= \frac{-1}{16}\left( -7 +76\zeta(2) -109\zeta(3) +16\zeta(4) -12\zeta(5)  \right. \nonumber \\ &\left. \hspace{1em}
+4\zeta(2)\zeta(3)\right) \label{eq10023} \\
\sum_{k=1}^\infty \frac{H(k)}{k^{3}(k+1)^{4}(k+2)^{2}} \sumend &= \frac{-1}{16}\left( -15 -57\zeta(2) +90\zeta(3) -\zeta(4) +32\zeta(5)  \right. \nonumber \\ &\left. \hspace{1em}
-16\zeta(2)\zeta(3)\right) \label{eq10024} \\
\sum_{k=1}^\infty \frac{H(k)}{k^{2}(k+1)^{5}(k+2)^{2}} \sumend &= \frac{1}{4}\left(  8 -7\zeta(2) +\zeta(3) +2\zeta(4) +3\zeta(6) -2\zeta(3)^2\right) \label{eq10025} \\
\sum_{k=1}^\infty \frac{H(k)}{k(k+1)^{6}(k+2)^{2}} \sumend &= \frac{1}{4}\left(  17 -\zeta(2) -14\zeta(3) +2\zeta(4) -16\zeta(5) +8\zeta(2)\zeta(3)  \right. \nonumber \\ &\left. \hspace{1em}
+3\zeta(6) -2\zeta(3)^2 -12\zeta(7) +4\zeta(2)\zeta(5) +4\zeta(3)\zeta(4)\right) \label{eq10026} \\
\sum_{k=1}^\infty \frac{H(k)}{(k+1)^{7}(k+2)^{2}} \sumend &= \frac{1}{4}\left(  36 -4\zeta(2) -28\zeta(3) +5\zeta(4) -32\zeta(5) +16\zeta(2)\zeta(3)  \right. \nonumber \\ &\left. \hspace{1em}
+9\zeta(6) -6\zeta(3)^2 -24\zeta(7) +8\zeta(2)\zeta(5) +8\zeta(3)\zeta(4) +5\zeta(8)  \right. \nonumber \\ &\left. \hspace{1em}
-4\zeta(3)\zeta(5)\right) \label{eq10027} \\
\sum_{k=1}^\infty \frac{H(k)}{k^{6}(k+2)^{3}} \sumend &= \frac{1}{256}\left( -57 -5\zeta(2) +76\zeta(3) -49\zeta(4) +144\zeta(5) -48\zeta(2)\zeta(3)  \right. \nonumber \\ &\left. \hspace{1em}
-84\zeta(6) +24\zeta(3)^2 +128\zeta(7) -32\zeta(2)\zeta(5) -32\zeta(3)\zeta(4)\right) \label{eq10028} \\
\sum_{k=1}^\infty \frac{H(k)}{k^{5}(k+1)(k+2)^{3}} \sumend &= \frac{1}{128}\left( -69 +117\zeta(2) -150\zeta(3) +81\zeta(4) -120\zeta(5)  \right. \nonumber \\ &\left. \hspace{1em}
+40\zeta(2)\zeta(3) +28\zeta(6) -8\zeta(3)^2\right) \label{eq10029} \\
\sum_{k=1}^\infty \frac{H(k)}{k^{4}(k+1)^{2}(k+2)^{3}} \sumend &= \frac{1}{32}\left( -41 -41\zeta(2) +102\zeta(3) -17\zeta(4) +12\zeta(5)  \right. \nonumber \\ &\left. \hspace{1em}
-4\zeta(2)\zeta(3)\right) \label{eq10030} \\
\sum_{k=1}^\infty \frac{H(k)}{k^{3}(k+1)^{3}(k+2)^{3}} \sumend &= \frac{-1}{16}\left(  48 -35\zeta(2) +7\zeta(3) +\zeta(4)\right) \label{eq10031} \\
\sum_{k=1}^\infty \frac{H(k)}{k^{2}(k+1)^{4}(k+2)^{3}} \sumend &= \frac{1}{16}\left( -111 +13\zeta(2) +76\zeta(3) -3\zeta(4) +32\zeta(5)  \right. \nonumber \\ &\left. \hspace{1em}
-16\zeta(2)\zeta(3)\right) \label{eq10032}
\end{align}
 
\begin{align}
\sum_{k=1}^\infty \frac{H(k)}{k(k+1)^{5}(k+2)^{3}} \sumend &= \frac{-1}{8}\left(  127 -27\zeta(2) -74\zeta(3) +7\zeta(4) -32\zeta(5) +16\zeta(2)\zeta(3)  \right. \nonumber \\ &\left. \hspace{1em}
+6\zeta(6) -4\zeta(3)^2\right) \label{eq10033} \\
\sum_{k=1}^\infty \frac{H(k)}{(k+1)^{6}(k+2)^{3}} \sumend &= \frac{-1}{4}\left(  144 -28\zeta(2) -88\zeta(3) +9\zeta(4) -48\zeta(5) +24\zeta(2)\zeta(3)  \right. \nonumber \\ &\left. \hspace{1em}
+9\zeta(6) -6\zeta(3)^2 -12\zeta(7) +4\zeta(2)\zeta(5) +4\zeta(3)\zeta(4)\right) \label{eq10034} \\
\sum_{k=1}^\infty \frac{H(k)}{k^{5}(k+2)^{4}} \sumend &= \frac{1}{256}\left(  187 -23\zeta(2) -138\zeta(3) +37\zeta(4) -112\zeta(5) +40\zeta(2)\zeta(3)  \right. \nonumber \\ &\left. \hspace{1em}
+28\zeta(6) -8\zeta(3)^2\right) \label{eq10035} \\
\sum_{k=1}^\infty \frac{H(k)}{k^{4}(k+1)(k+2)^{4}} \sumend &= \frac{-1}{32}\left( -64 +35\zeta(2) -3\zeta(3) +11\zeta(4) -2\zeta(5)\right) \label{eq10036} \\
\sum_{k=1}^\infty \frac{H(k)}{k^{3}(k+1)^{2}(k+2)^{4}} \sumend &= \frac{1}{32}\left(  169 -29\zeta(2) -96\zeta(3) -5\zeta(4) -8\zeta(5)  \right. \nonumber \\ &\left. \hspace{1em}
+4\zeta(2)\zeta(3)\right) \label{eq10037} \\
\sum_{k=1}^\infty \frac{H(k)}{k^{2}(k+1)^{3}(k+2)^{4}} \sumend &= \frac{-1}{16}\left( -217 +64\zeta(2) +89\zeta(3) +4\zeta(4) +8\zeta(5)  \right. \nonumber \\ &\left. \hspace{1em}
-4\zeta(2)\zeta(3)\right) \label{eq10038} \\
\sum_{k=1}^\infty \frac{H(k)}{k(k+1)^{4}(k+2)^{4}} \sumend &= \frac{1}{16}\left(  545 -141\zeta(2) -254\zeta(3) -5\zeta(4) -48\zeta(5)  \right. \nonumber \\ &\left. \hspace{1em}
+24\zeta(2)\zeta(3)\right) \label{eq10039} \\
\sum_{k=1}^\infty \frac{H(k)}{(k+1)^{5}(k+2)^{4}} \sumend &= \frac{-1}{4}\left( -336 +84\zeta(2) +164\zeta(3) -\zeta(4) +40\zeta(5) -20\zeta(2)\zeta(3)  \right. \nonumber \\ &\left. \hspace{1em}
-3\zeta(6) +2\zeta(3)^2\right) \label{eq10040} \\
\sum_{k=1}^\infty \frac{H(k)}{k^{4}(k+2)^{5}} \sumend &= \frac{1}{256}\left( -443 +93\zeta(2) +188\zeta(3) +33\zeta(4) +104\zeta(5) -40\zeta(2)\zeta(3)  \right. \nonumber \\ &\left. \hspace{1em}
+12\zeta(6) -8\zeta(3)^2\right) \label{eq10041} \\
\sum_{k=1}^\infty \frac{H(k)}{k^{3}(k+1)(k+2)^{5}} \sumend &= \frac{-1}{128}\left(  699 -233\zeta(2) -176\zeta(3) -77\zeta(4) -96\zeta(5)  \right. \nonumber \\ &\left. \hspace{1em}
+40\zeta(2)\zeta(3) -12\zeta(6) +8\zeta(3)^2\right) \label{eq10042} \\
\sum_{k=1}^\infty \frac{H(k)}{k^{2}(k+1)^{2}(k+2)^{5}} \sumend &= \frac{1}{64}\left( -1037 +291\zeta(2) +368\zeta(3) +87\zeta(4) +112\zeta(5)  \right. \nonumber \\ &\left. \hspace{1em}
-48\zeta(2)\zeta(3) +12\zeta(6) -8\zeta(3)^2\right) \label{eq10043}
\end{align}
 
\begin{align}
\sum_{k=1}^\infty \frac{H(k)}{k(k+1)^{3}(k+2)^{5}} \sumend &= \frac{1}{32}\left( -1471 +419\zeta(2) +546\zeta(3) +95\zeta(4) +128\zeta(5)  \right. \nonumber \\ &\left. \hspace{1em}
-56\zeta(2)\zeta(3) +12\zeta(6) -8\zeta(3)^2\right) \label{eq10044} \\
\sum_{k=1}^\infty \frac{H(k)}{(k+1)^{4}(k+2)^{5}} \sumend &= \frac{-1}{4}\left(  504 -140\zeta(2) -200\zeta(3) -25\zeta(4) -44\zeta(5)  \right. \nonumber \\ &\left. \hspace{1em}
+20\zeta(2)\zeta(3) -3\zeta(6) +2\zeta(3)^2\right) \label{eq10045} \\
\sum_{k=1}^\infty \frac{H(k)}{k^{3}(k+2)^{6}} \sumend &= \frac{1}{256}\left(  825 -177\zeta(2) -222\zeta(3) -133\zeta(4) -176\zeta(5) +48\zeta(2)\zeta(3)  \right. \nonumber \\ &\left. \hspace{1em}
-68\zeta(6) +24\zeta(3)^2 -96\zeta(7) +32\zeta(2)\zeta(5) +32\zeta(3)\zeta(4)\right) \label{eq10046} \\
\sum_{k=1}^\infty \frac{H(k)}{k^{2}(k+1)(k+2)^{6}} \sumend &= \frac{-1}{64}\left( -762 +205\zeta(2) +199\zeta(3) +105\zeta(4) +136\zeta(5)  \right. \nonumber \\ &\left. \hspace{1em}
-44\zeta(2)\zeta(3) +40\zeta(6) -16\zeta(3)^2 +48\zeta(7) -16\zeta(2)\zeta(5)  \right. \nonumber \\ &\left. \hspace{1em}
-16\zeta(3)\zeta(4)\right) \label{eq10047} \\
\sum_{k=1}^\infty \frac{H(k)}{k(k+1)^{2}(k+2)^{6}} \sumend &= \frac{-1}{64}\left( -2561 +701\zeta(2) +766\zeta(3) +297\zeta(4) +384\zeta(5)  \right. \nonumber \\ &\left. \hspace{1em}
-136\zeta(2)\zeta(3) +92\zeta(6) -40\zeta(3)^2 +96\zeta(7) -32\zeta(2)\zeta(5)  \right. \nonumber \\ &\left. \hspace{1em}
-32\zeta(3)\zeta(4)\right) \label{eq10048} \\
\sum_{k=1}^\infty \frac{H(k)}{(k+1)^{3}(k+2)^{6}} \sumend &= \frac{-1}{4}\left( -504 +140\zeta(2) +164\zeta(3) +49\zeta(4) +64\zeta(5)  \right. \nonumber \\ &\left. \hspace{1em}
-24\zeta(2)\zeta(3) +13\zeta(6) -6\zeta(3)^2 +12\zeta(7) -4\zeta(2)\zeta(5) -4\zeta(3)\zeta(4)\right) \label{eq10049} \\
\sum_{k=1}^\infty \frac{H(k)}{k^{2}(k+2)^{7}} \sumend &= \frac{1}{256}\left( -1291 +233\zeta(2) +244\zeta(3) +213\zeta(4) +240\zeta(5)  \right. \nonumber \\ &\left. \hspace{1em}
-32\zeta(2)\zeta(3) +164\zeta(6) -24\zeta(3)^2 +256\zeta(7) -64\zeta(2)\zeta(5) -64\zeta(3)\zeta(4)  \right. \nonumber \\ &\left. \hspace{1em}
+80\zeta(8) -64\zeta(3)\zeta(5)\right) \label{eq10050} \\
\sum_{k=1}^\infty \frac{H(k)}{k(k+1)(k+2)^{7}} \sumend &= \frac{-1}{128}\left(  2815 -643\zeta(2) -642\zeta(3) -423\zeta(4) -512\zeta(5)  \right. \nonumber \\ &\left. \hspace{1em}
+120\zeta(2)\zeta(3) -244\zeta(6) +56\zeta(3)^2 -352\zeta(7) +96\zeta(2)\zeta(5) +96\zeta(3)\zeta(4)  \right. \nonumber \\ &\left. \hspace{1em}
-80\zeta(8) +64\zeta(3)\zeta(5)\right) \label{eq10051} \\
\sum_{k=1}^\infty \frac{H(k)}{(k+1)^{2}(k+2)^{7}} \sumend &= \frac{1}{4}\left( -336 +84\zeta(2) +88\zeta(3) +45\zeta(4) +56\zeta(5) -16\zeta(2)\zeta(3)  \right. \nonumber \\ &\left. \hspace{1em}
+21\zeta(6) -6\zeta(3)^2 +28\zeta(7) -8\zeta(2)\zeta(5) -8\zeta(3)\zeta(4) +5\zeta(8)  \right. \nonumber \\ &\left. \hspace{1em}
-4\zeta(3)\zeta(5)\right) \label{eq10052}
\end{align}
 
\begin{align}
\sum_{k=1}^\infty \frac{H(k)}{k(k+2)^{8}} \sumend &= \frac{-1}{256}\left( -1793 +253\zeta(2) +254\zeta(3) +249\zeta(4) +256\zeta(5) -8\zeta(2)\zeta(3)  \right. \nonumber \\ &\left. \hspace{1em}
+236\zeta(6) -8\zeta(3)^2 +288\zeta(7) -32\zeta(2)\zeta(5) -32\zeta(3)\zeta(4) +208\zeta(8)  \right. \nonumber \\ &\left. \hspace{1em}
-64\zeta(3)\zeta(5) +512\zeta(9) -128\zeta(3)\zeta(6) -128\zeta(4)\zeta(5) -128\zeta(2)\zeta(7)\right) \label{eq10053} \\
\sum_{k=1}^\infty \frac{H(k)}{(k+1)(k+2)^{8}} \sumend &= \frac{1}{4}\left(  144 -28\zeta(2) -28\zeta(3) -21\zeta(4) -24\zeta(5) +4\zeta(2)\zeta(3)  \right. \nonumber \\ &\left. \hspace{1em}
-15\zeta(6) +2\zeta(3)^2 -20\zeta(7) +4\zeta(2)\zeta(5) +4\zeta(3)\zeta(4) -9\zeta(8)  \right. \nonumber \\ &\left. \hspace{1em}
+4\zeta(3)\zeta(5) -16\zeta(9) +4\zeta(3)\zeta(6) +4\zeta(4)\zeta(5) +4\zeta(2)\zeta(7)\right) \label{eq10054} \\
\sum_{k=1}^\infty \frac{H(k)}{(k+2)^{9}} \sumend &= \frac{-1}{4}\left(  36 -4\zeta(2) -4\zeta(3) -4\zeta(4) -4\zeta(5) -4\zeta(6) -4\zeta(7) -4\zeta(8)  \right. \nonumber \\ &\left. \hspace{1em}
-4\zeta(9) -7\zeta(10) +4\zeta(3)\zeta(7) +2\zeta(5)^2\right) \label{eq10055} \\
\sum_{k=1}^\infty \frac{H(k)^{2}}{k^{8}} \sumend &= -\left( - M(2,8)\right) \label{eq10056} \\
\sum_{k=1}^\infty \frac{H(k)^{2}}{k^{7}(k+1)} \sumend &= \frac{-1}{24}\left( -72\zeta(3) +102\zeta(4) -84\zeta(5) +24\zeta(2)\zeta(3) +97\zeta(6)  \right. \nonumber \\ &\left. \hspace{1em}
-48\zeta(3)^2 -144\zeta(7) +24\zeta(2)\zeta(5) +60\zeta(3)\zeta(4) +24 M(2,6) -220\zeta(9)  \right. \nonumber \\ &\left. \hspace{1em}
+84\zeta(3)\zeta(6) +60\zeta(4)\zeta(5) +24\zeta(2)\zeta(7) -8\zeta(3)^3\right) \label{eq10057} \\
\sum_{k=1}^\infty \frac{H(k)^{2}}{k^{6}(k+1)^{2}} \sumend &= \frac{1}{8}\left( -144\zeta(3) +192\zeta(4) -112\zeta(5) +32\zeta(2)\zeta(3) +97\zeta(6)  \right. \nonumber \\ &\left. \hspace{1em}
-48\zeta(3)^2 -96\zeta(7) +16\zeta(2)\zeta(5) +40\zeta(3)\zeta(4) +8 M(2,6)\right) \label{eq10058} \\
\sum_{k=1}^\infty \frac{H(k)^{2}}{k^{5}(k+1)^{3}} \sumend &= \frac{1}{8}\left(  360\zeta(3) -450\zeta(4) +180\zeta(5) -56\zeta(2)\zeta(3) -97\zeta(6)  \right. \nonumber \\ &\left. \hspace{1em}
+48\zeta(3)^2 +48\zeta(7) -8\zeta(2)\zeta(5) -20\zeta(3)\zeta(4)\right) \label{eq10059} \\
\sum_{k=1}^\infty \frac{H(k)^{2}}{k^{4}(k+1)^{4}} \sumend &= \frac{-1}{12}\left(  720\zeta(3) -840\zeta(4) +240\zeta(5) -96\zeta(2)\zeta(3) -67\zeta(6)  \right. \nonumber \\ &\left. \hspace{1em}
+36\zeta(3)^2\right) \label{eq10060} \\
\sum_{k=1}^\infty \frac{H(k)^{2}}{k^{3}(k+1)^{5}} \sumend &= \frac{-1}{8}\left( -360\zeta(3) +390\zeta(4) -100\zeta(5) +56\zeta(2)\zeta(3) +37\zeta(6)  \right. \nonumber \\ &\left. \hspace{1em}
-24\zeta(3)^2 -8\zeta(7) +8\zeta(2)\zeta(5) -4\zeta(3)\zeta(4)\right) \label{eq10061} \\
\sum_{k=1}^\infty \frac{H(k)^{2}}{k^{2}(k+1)^{6}} \sumend &= \frac{-1}{8}\left(  144\zeta(3) -144\zeta(4) +48\zeta(5) -32\zeta(2)\zeta(3) -37\zeta(6)  \right. \nonumber \\ &\left. \hspace{1em}
+24\zeta(3)^2 +16\zeta(7) -16\zeta(2)\zeta(5) +8\zeta(3)\zeta(4) +28\zeta(8) -16\zeta(3)\zeta(5)  \right. \nonumber \\ &\left. \hspace{1em}
-8 M(2,6)\right) \label{eq10062}
\end{align}
 
\begin{align}
\sum_{k=1}^\infty \frac{H(k)^{2}}{k(k+1)^{7}} \sumend &= \frac{1}{24}\left(  72\zeta(3) -66\zeta(4) +36\zeta(5) -24\zeta(2)\zeta(3) -37\zeta(6)  \right. \nonumber \\ &\left. \hspace{1em}
+24\zeta(3)^2 +24\zeta(7) -24\zeta(2)\zeta(5) +12\zeta(3)\zeta(4) +84\zeta(8) -48\zeta(3)\zeta(5)  \right. \nonumber \\ &\left. \hspace{1em}
-24 M(2,6) -4\zeta(9) +36\zeta(3)\zeta(6) +12\zeta(4)\zeta(5) -24\zeta(2)\zeta(7) -8\zeta(3)^3\right) \label{eq10063} \\
\sum_{k=1}^\infty \frac{H(k)^{2}}{(k+1)^{8}} \sumend &= \frac{1}{2}\left( -9\zeta(10) +4\zeta(3)\zeta(7) +2\zeta(5)^2 +2 M(2,8)\right) \label{eq10064} \\
\sum_{k=1}^\infty \frac{H(k)^{2}}{k^{7}(k+2)} \sumend &= \frac{-1}{768}\left( -6 -6\zeta(2) -18\zeta(3) +51\zeta(4) -84\zeta(5) +24\zeta(2)\zeta(3)  \right. \nonumber \\ &\left. \hspace{1em}
+194\zeta(6) -96\zeta(3)^2 -576\zeta(7) +96\zeta(2)\zeta(5) +240\zeta(3)\zeta(4) +192 M(2,6)  \right. \nonumber \\ &\left. \hspace{1em}
-3520\zeta(9) +1344\zeta(3)\zeta(6) +960\zeta(4)\zeta(5) +384\zeta(2)\zeta(7) -128\zeta(3)^3\right) \label{eq10065} \\
\sum_{k=1}^\infty \frac{H(k)^{2}}{k^{6}(k+1)(k+2)} \sumend &= \frac{-1}{384}\left( -6 -6\zeta(2) +1134\zeta(3) -1581\zeta(4) +1260\zeta(5)  \right. \nonumber \\ &\left. \hspace{1em}
-360\zeta(2)\zeta(3) -1358\zeta(6) +672\zeta(3)^2 +1728\zeta(7) -288\zeta(2)\zeta(5)  \right. \nonumber \\ &\left. \hspace{1em}
-720\zeta(3)\zeta(4) -192 M(2,6)\right) \label{eq10066} \\
\sum_{k=1}^\infty \frac{H(k)^{2}}{k^{5}(k+1)^{2}(k+2)} \sumend &= \frac{-1}{192}\left( -6 -6\zeta(2) -2322\zeta(3) +3027\zeta(4) -1428\zeta(5)  \right. \nonumber \\ &\left. \hspace{1em}
+408\zeta(2)\zeta(3) +970\zeta(6) -480\zeta(3)^2 -576\zeta(7) +96\zeta(2)\zeta(5)  \right. \nonumber \\ &\left. \hspace{1em}
+240\zeta(3)\zeta(4)\right) \label{eq10067} \\
\sum_{k=1}^\infty \frac{H(k)^{2}}{k^{4}(k+1)^{3}(k+2)} \sumend &= \frac{-1}{96}\left( -6 -6\zeta(2) +1998\zeta(3) -2373\zeta(4) +732\zeta(5)  \right. \nonumber \\ &\left. \hspace{1em}
-264\zeta(2)\zeta(3) -194\zeta(6) +96\zeta(3)^2\right) \label{eq10068} \\
\sum_{k=1}^\infty \frac{H(k)^{2}}{k^{3}(k+1)^{4}(k+2)} \sumend &= \frac{-1}{48}\left( -6 -6\zeta(2) -882\zeta(3) +987\zeta(4) -228\zeta(5)  \right. \nonumber \\ &\left. \hspace{1em}
+120\zeta(2)\zeta(3) +74\zeta(6) -48\zeta(3)^2\right) \label{eq10069} \\
\sum_{k=1}^\infty \frac{H(k)^{2}}{k^{2}(k+1)^{5}(k+2)} \sumend &= \frac{-1}{24}\left( -6 -6\zeta(2) +198\zeta(3) -183\zeta(4) +72\zeta(5)  \right. \nonumber \\ &\left. \hspace{1em}
-48\zeta(2)\zeta(3) -37\zeta(6) +24\zeta(3)^2 +24\zeta(7) -24\zeta(2)\zeta(5)  \right. \nonumber \\ &\left. \hspace{1em}
+12\zeta(3)\zeta(4)\right) \label{eq10070} \\
\sum_{k=1}^\infty \frac{H(k)^{2}}{k(k+1)^{6}(k+2)} \sumend &= \frac{-1}{24}\left( -12 -12\zeta(2) -36\zeta(3) +66\zeta(4) +37\zeta(6) -24\zeta(3)^2  \right. \nonumber \\ &\left. \hspace{1em}
-84\zeta(8) +48\zeta(3)\zeta(5) +24 M(2,6)\right) \label{eq10071} \\
\sum_{k=1}^\infty \frac{H(k)^{2}}{(k+1)^{7}(k+2)} \sumend &= \frac{-1}{24}\left( -24 -24\zeta(2) +66\zeta(4) +36\zeta(5) -24\zeta(2)\zeta(3) +37\zeta(6)  \right. \nonumber \\ &\left. \hspace{1em}
-24\zeta(3)^2 +24\zeta(7) -24\zeta(2)\zeta(5) +12\zeta(3)\zeta(4) -84\zeta(8) +48\zeta(3)\zeta(5)  \right. \nonumber \\ &\left. \hspace{1em}
+24 M(2,6) -4\zeta(9) +36\zeta(3)\zeta(6) +12\zeta(4)\zeta(5) -24\zeta(2)\zeta(7) -8\zeta(3)^3\right) \label{eq10072}
\end{align}
 
\begin{align}
\sum_{k=1}^\infty \frac{H(k)^{2}}{k^{6}(k+2)^{2}} \sumend &= \frac{1}{128}\left( -12 -6\zeta(2) -14\zeta(3) +48\zeta(4) -56\zeta(5) +16\zeta(2)\zeta(3)  \right. \nonumber \\ &\left. \hspace{1em}
+97\zeta(6) -48\zeta(3)^2 -192\zeta(7) +32\zeta(2)\zeta(5) +80\zeta(3)\zeta(4) +32 M(2,6)\right) \label{eq10073} \\
\sum_{k=1}^\infty \frac{H(k)^{2}}{k^{5}(k+1)(k+2)^{2}} \sumend &= \frac{1}{384}\left( -78 -42\zeta(2) +1050\zeta(3) -1293\zeta(4) +924\zeta(5)  \right. \nonumber \\ &\left. \hspace{1em}
-264\zeta(2)\zeta(3) -776\zeta(6) +384\zeta(3)^2 +576\zeta(7) -96\zeta(2)\zeta(5)  \right. \nonumber \\ &\left. \hspace{1em}
-240\zeta(3)\zeta(4)\right) \label{eq10074} \\
\sum_{k=1}^\infty \frac{H(k)^{2}}{k^{4}(k+1)^{2}(k+2)^{2}} \sumend &= \frac{1}{96}\left( -42 -24\zeta(2) -636\zeta(3) +867\zeta(4) -252\zeta(5)  \right. \nonumber \\ &\left. \hspace{1em}
+72\zeta(2)\zeta(3) +97\zeta(6) -48\zeta(3)^2\right) \label{eq10075} \\
\sum_{k=1}^\infty \frac{H(k)^{2}}{k^{3}(k+1)^{3}(k+2)^{2}} \sumend &= \frac{1}{32}\left( -30 -18\zeta(2) +242\zeta(3) -213\zeta(4) +76\zeta(5)  \right. \nonumber \\ &\left. \hspace{1em}
-40\zeta(2)\zeta(3)\right) \label{eq10076} \\
\sum_{k=1}^\infty \frac{H(k)^{2}}{k^{2}(k+1)^{4}(k+2)^{2}} \sumend &= \frac{-1}{24}\left(  48 +30\zeta(2) +78\zeta(3) -174\zeta(4) -37\zeta(6)  \right. \nonumber \\ &\left. \hspace{1em}
+24\zeta(3)^2\right) \label{eq10077} \\
\sum_{k=1}^\infty \frac{H(k)^{2}}{k(k+1)^{5}(k+2)^{2}} \sumend &= \frac{-1}{24}\left(  102 +66\zeta(2) -42\zeta(3) -165\zeta(4) -72\zeta(5)  \right. \nonumber \\ &\left. \hspace{1em}
+48\zeta(2)\zeta(3) -37\zeta(6) +24\zeta(3)^2 -24\zeta(7) +24\zeta(2)\zeta(5)  \right. \nonumber \\ &\left. \hspace{1em}
-12\zeta(3)\zeta(4)\right) \label{eq10078} \\
\sum_{k=1}^\infty \frac{H(k)^{2}}{(k+1)^{6}(k+2)^{2}} \sumend &= \frac{-1}{8}\left(  72 +48\zeta(2) -16\zeta(3) -132\zeta(4) -48\zeta(5)  \right. \nonumber \\ &\left. \hspace{1em}
+32\zeta(2)\zeta(3) -37\zeta(6) +24\zeta(3)^2 -16\zeta(7) +16\zeta(2)\zeta(5) -8\zeta(3)\zeta(4)  \right. \nonumber \\ &\left. \hspace{1em}
+28\zeta(8) -16\zeta(3)\zeta(5) -8 M(2,6)\right) \label{eq10079} \\
\sum_{k=1}^\infty \frac{H(k)^{2}}{k^{5}(k+2)^{3}} \sumend &= \frac{1}{256}\left(  138 +22\zeta(2) +26\zeta(3) -229\zeta(4) +180\zeta(5)  \right. \nonumber \\ &\left. \hspace{1em}
-56\zeta(2)\zeta(3) -194\zeta(6) +96\zeta(3)^2 +192\zeta(7) -32\zeta(2)\zeta(5)  \right. \nonumber \\ &\left. \hspace{1em}
-80\zeta(3)\zeta(4)\right) \label{eq10080} \\
\sum_{k=1}^\infty \frac{H(k)^{2}}{k^{4}(k+1)(k+2)^{3}} \sumend &= \frac{1}{192}\left(  246 +54\zeta(2) -486\zeta(3) +303\zeta(4) -192\zeta(5)  \right. \nonumber \\ &\left. \hspace{1em}
+48\zeta(2)\zeta(3) +97\zeta(6) -48\zeta(3)^2\right) \label{eq10081} \\
\sum_{k=1}^\infty \frac{H(k)^{2}}{k^{3}(k+1)^{2}(k+2)^{3}} \sumend &= \frac{1}{16}\left(  48 +13\zeta(2) +25\zeta(3) -94\zeta(4) +10\zeta(5)  \right. \nonumber \\ &\left. \hspace{1em}
-4\zeta(2)\zeta(3)\right) \label{eq10082}
\end{align}
 
\begin{align}
\sum_{k=1}^\infty \frac{H(k)^{2}}{k^{2}(k+1)^{3}(k+2)^{3}} \sumend &= \frac{-1}{32}\left( -222 -70\zeta(2) +142\zeta(3) +163\zeta(4) +36\zeta(5)  \right. \nonumber \\ &\left. \hspace{1em}
-24\zeta(2)\zeta(3)\right) \label{eq10083} \\
\sum_{k=1}^\infty \frac{H(k)^{2}}{k(k+1)^{4}(k+2)^{3}} \sumend &= \frac{1}{48}\left(  762 +270\zeta(2) -270\zeta(3) -837\zeta(4) -108\zeta(5)  \right. \nonumber \\ &\left. \hspace{1em}
+72\zeta(2)\zeta(3) -74\zeta(6) +48\zeta(3)^2\right) \label{eq10084} \\
\sum_{k=1}^\infty \frac{H(k)^{2}}{(k+1)^{5}(k+2)^{3}} \sumend &= \frac{-1}{8}\left( -288 -112\zeta(2) +104\zeta(3) +334\zeta(4) +60\zeta(5)  \right. \nonumber \\ &\left. \hspace{1em}
-40\zeta(2)\zeta(3) +37\zeta(6) -24\zeta(3)^2 +8\zeta(7) -8\zeta(2)\zeta(5) +4\zeta(3)\zeta(4)\right) \label{eq10085} \\
\sum_{k=1}^\infty \frac{H(k)^{2}}{k^{4}(k+2)^{4}} \sumend &= \frac{-1}{192}\left(  384 -18\zeta(2) -90\zeta(3) -240\zeta(4) +72\zeta(5)  \right. \nonumber \\ &\left. \hspace{1em}
-24\zeta(2)\zeta(3) -67\zeta(6) +36\zeta(3)^2\right) \label{eq10086} \\
\sum_{k=1}^\infty \frac{H(k)^{2}}{k^{3}(k+1)(k+2)^{4}} \sumend &= \frac{-1}{192}\left(  1014 +18\zeta(2) -666\zeta(3) -177\zeta(4) -48\zeta(5)  \right. \nonumber \\ &\left. \hspace{1em}
-37\zeta(6) +24\zeta(3)^2\right) \label{eq10087} \\
\sum_{k=1}^\infty \frac{H(k)^{2}}{k^{2}(k+1)^{2}(k+2)^{4}} \sumend &= \frac{-1}{96}\left(  1302 +96\zeta(2) -516\zeta(3) -741\zeta(4) +12\zeta(5)  \right. \nonumber \\ &\left. \hspace{1em}
-24\zeta(2)\zeta(3) -37\zeta(6) +24\zeta(3)^2\right) \label{eq10088} \\
\sum_{k=1}^\infty \frac{H(k)^{2}}{k(k+1)^{3}(k+2)^{4}} \sumend &= \frac{-1}{96}\left(  3270 +402\zeta(2) -1458\zeta(3) -1971\zeta(4) -84\zeta(5)  \right. \nonumber \\ &\left. \hspace{1em}
+24\zeta(2)\zeta(3) -74\zeta(6) +48\zeta(3)^2\right) \label{eq10089} \\
\sum_{k=1}^\infty \frac{H(k)^{2}}{(k+1)^{4}(k+2)^{4}} \sumend &= \frac{1}{12}\left( -1008 -168\zeta(2) +432\zeta(3) +702\zeta(4) +48\zeta(5)  \right. \nonumber \\ &\left. \hspace{1em}
-24\zeta(2)\zeta(3) +37\zeta(6) -24\zeta(3)^2\right) \label{eq10090} \\
\sum_{k=1}^\infty \frac{H(k)^{2}}{k^{3}(k+2)^{5}} \sumend &= \frac{1}{256}\left(  1398 -210\zeta(2) -486\zeta(3) -371\zeta(4) -252\zeta(5)  \right. \nonumber \\ &\left. \hspace{1em}
+104\zeta(2)\zeta(3) -122\zeta(6) +80\zeta(3)^2 +32\zeta(7) -32\zeta(2)\zeta(5)  \right. \nonumber \\ &\left. \hspace{1em}
+16\zeta(3)\zeta(4)\right) \label{eq10091} \\
\sum_{k=1}^\infty \frac{H(k)^{2}}{k^{2}(k+1)(k+2)^{5}} \sumend &= \frac{-1}{384}\left( -6222 +594\zeta(2) +2790\zeta(3) +1467\zeta(4) +852\zeta(5)  \right. \nonumber \\ &\left. \hspace{1em}
-312\zeta(2)\zeta(3) +440\zeta(6) -288\zeta(3)^2 -96\zeta(7) +96\zeta(2)\zeta(5)  \right. \nonumber \\ &\left. \hspace{1em}
-48\zeta(3)\zeta(4)\right) \label{eq10092}
\end{align}
 
\begin{align}
\sum_{k=1}^\infty \frac{H(k)^{2}}{k(k+1)^{2}(k+2)^{5}} \sumend &= \frac{-1}{192}\left( -8826 +402\zeta(2) +3822\zeta(3) +2949\zeta(4) +828\zeta(5)  \right. \nonumber \\ &\left. \hspace{1em}
-264\zeta(2)\zeta(3) +514\zeta(6) -336\zeta(3)^2 -96\zeta(7) +96\zeta(2)\zeta(5)  \right. \nonumber \\ &\left. \hspace{1em}
-48\zeta(3)\zeta(4)\right) \label{eq10093} \\
\sum_{k=1}^\infty \frac{H(k)^{2}}{(k+1)^{3}(k+2)^{5}} \sumend &= \frac{-1}{8}\left( -1008 +440\zeta(3) +410\zeta(4) +76\zeta(5) -24\zeta(2)\zeta(3)  \right. \nonumber \\ &\left. \hspace{1em}
+49\zeta(6) -32\zeta(3)^2 -8\zeta(7) +8\zeta(2)\zeta(5) -4\zeta(3)\zeta(4)\right) \label{eq10094} \\
\sum_{k=1}^\infty \frac{H(k)^{2}}{k^{2}(k+2)^{6}} \sumend &= \frac{-1}{128}\left(  1524 -282\zeta(2) -498\zeta(3) -272\zeta(4) -424\zeta(5)  \right. \nonumber \\ &\left. \hspace{1em}
+160\zeta(2)\zeta(3) -165\zeta(6) +88\zeta(3)^2 -160\zeta(7) +32\zeta(2)\zeta(5) +80\zeta(3)\zeta(4)  \right. \nonumber \\ &\left. \hspace{1em}
+112\zeta(8) -64\zeta(3)\zeta(5) -32 M(2,6)\right) \label{eq10095} \\
\sum_{k=1}^\infty \frac{H(k)^{2}}{k(k+1)(k+2)^{6}} \sumend &= \frac{1}{384}\left( -15366 +2286\zeta(2) +5778\zeta(3) +3099\zeta(4) +3396\zeta(5)  \right. \nonumber \\ &\left. \hspace{1em}
-1272\zeta(2)\zeta(3) +1430\zeta(6) -816\zeta(3)^2 +864\zeta(7) -96\zeta(2)\zeta(5)  \right. \nonumber \\ &\left. \hspace{1em}
-528\zeta(3)\zeta(4) -672\zeta(8) +384\zeta(3)\zeta(5) +192 M(2,6)\right) \label{eq10096} \\
\sum_{k=1}^\infty \frac{H(k)^{2}}{(k+1)^{2}(k+2)^{6}} \sumend &= \frac{-1}{8}\left(  1008 -112\zeta(2) -400\zeta(3) -252\zeta(4) -176\zeta(5)  \right. \nonumber \\ &\left. \hspace{1em}
+64\zeta(2)\zeta(3) -81\zeta(6) +48\zeta(3)^2 -32\zeta(7) +24\zeta(3)\zeta(4) +28\zeta(8)  \right. \nonumber \\ &\left. \hspace{1em}
-16\zeta(3)\zeta(5) -8 M(2,6)\right) \label{eq10097} \\
\sum_{k=1}^\infty \frac{H(k)^{2}}{k(k+2)^{7}} \sumend &= \frac{1}{768}\left(  16890 -3090\zeta(2) -4590\zeta(3) -2757\zeta(4) -4476\zeta(5)  \right. \nonumber \\ &\left. \hspace{1em}
+1416\zeta(2)\zeta(3) -2042\zeta(6) +720\zeta(3)^2 -3744\zeta(7) +1056\zeta(2)\zeta(5)  \right. \nonumber \\ &\left. \hspace{1em}
+1200\zeta(3)\zeta(4) -288\zeta(8) +384\zeta(3)\zeta(5) -192 M(2,6) -64\zeta(9) +576\zeta(3)\zeta(6)  \right. \nonumber \\ &\left. \hspace{1em}
+192\zeta(4)\zeta(5) -384\zeta(2)\zeta(7) -128\zeta(3)^3\right) \label{eq10098} \\
\sum_{k=1}^\infty \frac{H(k)^{2}}{(k+1)(k+2)^{7}} \sumend &= \frac{-1}{24}\left( -2016 +336\zeta(2) +648\zeta(3) +366\zeta(4) +492\zeta(5)  \right. \nonumber \\ &\left. \hspace{1em}
-168\zeta(2)\zeta(3) +217\zeta(6) -96\zeta(3)^2 +288\zeta(7) -72\zeta(2)\zeta(5) -108\zeta(3)\zeta(4)  \right. \nonumber \\ &\left. \hspace{1em}
-24\zeta(8) +24 M(2,6) +4\zeta(9) -36\zeta(3)\zeta(6) -12\zeta(4)\zeta(5) +24\zeta(2)\zeta(7)  \right. \nonumber \\ &\left. \hspace{1em}
+8\zeta(3)^3\right) \label{eq10099} \\
\sum_{k=1}^\infty \frac{H(k)^{2}}{(k+2)^{8}} \sumend &= \frac{1}{2}\left( -72 +12\zeta(2) +16\zeta(3) +11\zeta(4) +16\zeta(5) -4\zeta(2)\zeta(3)  \right. \nonumber \\ &\left. \hspace{1em}
+9\zeta(6) -2\zeta(3)^2 +16\zeta(7) -4\zeta(2)\zeta(5) -4\zeta(3)\zeta(4) +7\zeta(8) -4\zeta(3)\zeta(5)  \right. \nonumber \\ &\left. \hspace{1em}
+16\zeta(9) -4\zeta(3)\zeta(6) -4\zeta(4)\zeta(5) -4\zeta(2)\zeta(7) -9\zeta(10) +4\zeta(3)\zeta(7)  \right. \nonumber \\ &\left. \hspace{1em}
+2\zeta(5)^2 +2 M(2,8)\right) \label{eq10100}
\end{align}
 
\begin{align}
\sum_{k=1}^\infty \frac{H(k)^{3}}{k^{7}} \sumend &= \frac{-1}{160}\left(  1661\zeta(10) -1280\zeta(3)\zeta(7) -80\zeta(3)^2\zeta(4)  \right. \nonumber \\ &\left. \hspace{1em}
+560\zeta(2)\zeta(3)\zeta(5) -720\zeta(5)^2 -520 M(2,8)\right) \label{eq10101} \\
\sum_{k=1}^\infty \frac{H(k)^{3}}{k^{6}(k+1)} \sumend &= \frac{1}{96}\left( -960\zeta(4) +960\zeta(5) +96\zeta(2)\zeta(3) -558\zeta(6) +240\zeta(3)^2  \right. \nonumber \\ &\left. \hspace{1em}
+1386\zeta(7) +192\zeta(2)\zeta(5) -1224\zeta(3)\zeta(4) +595\zeta(8) +120\zeta(2)\zeta(3)^2  \right. \nonumber \\ &\left. \hspace{1em}
-576\zeta(3)\zeta(5) -264 M(2,6) +2084\zeta(9) -1164\zeta(3)\zeta(6) -1224\zeta(4)\zeta(5)  \right. \nonumber \\ &\left. \hspace{1em}
+288\zeta(2)\zeta(7) +192\zeta(3)^3\right) \label{eq10102} \\
\sum_{k=1}^\infty \frac{H(k)^{3}}{k^{5}(k+1)^{2}} \sumend &= \frac{-1}{96}\left( -4800\zeta(4) +4560\zeta(5) +480\zeta(2)\zeta(3) -1674\zeta(6)  \right. \nonumber \\ &\left. \hspace{1em}
+720\zeta(3)^2 +2772\zeta(7) +384\zeta(2)\zeta(5) -2448\zeta(3)\zeta(4) +595\zeta(8)  \right. \nonumber \\ &\left. \hspace{1em}
+120\zeta(2)\zeta(3)^2 -576\zeta(3)\zeta(5) -264 M(2,6)\right) \label{eq10103} \\
\sum_{k=1}^\infty \frac{H(k)^{3}}{k^{4}(k+1)^{3}} \sumend &= \frac{-1}{16}\left(  1600\zeta(4) -1440\zeta(5) -160\zeta(2)\zeta(3) +312\zeta(6)  \right. \nonumber \\ &\left. \hspace{1em}
-152\zeta(3)^2 -231\zeta(7) -32\zeta(2)\zeta(5) +204\zeta(3)\zeta(4)\right) \label{eq10104} \\
\sum_{k=1}^\infty \frac{H(k)^{3}}{k^{3}(k+1)^{4}} \sumend &= \frac{1}{16}\left(  1600\zeta(4) -1360\zeta(5) -160\zeta(2)\zeta(3) +192\zeta(6)  \right. \nonumber \\ &\left. \hspace{1em}
-136\zeta(3)^2 -119\zeta(7) -32\zeta(2)\zeta(5) +132\zeta(3)\zeta(4)\right) \label{eq10105} \\
\sum_{k=1}^\infty \frac{H(k)^{3}}{k^{2}(k+1)^{5}} \sumend &= \frac{-1}{96}\left(  4800\zeta(4) -3840\zeta(5) -480\zeta(2)\zeta(3) +594\zeta(6)  \right. \nonumber \\ &\left. \hspace{1em}
-576\zeta(3)^2 -1428\zeta(7) -384\zeta(2)\zeta(5) +1584\zeta(3)\zeta(4) +43\zeta(8)  \right. \nonumber \\ &\left. \hspace{1em}
+120\zeta(2)\zeta(3)^2 -288\zeta(3)\zeta(5) +24 M(2,6)\right) \label{eq10106} \\
\sum_{k=1}^\infty \frac{H(k)^{3}}{k(k+1)^{6}} \sumend &= \frac{1}{96}\left(  960\zeta(4) -720\zeta(5) -96\zeta(2)\zeta(3) +198\zeta(6) -192\zeta(3)^2  \right. \nonumber \\ &\left. \hspace{1em}
-714\zeta(7) -192\zeta(2)\zeta(5) +792\zeta(3)\zeta(4) +43\zeta(8) +120\zeta(2)\zeta(3)^2  \right. \nonumber \\ &\left. \hspace{1em}
-288\zeta(3)\zeta(5) +24 M(2,6) -788\zeta(9) +444\zeta(3)\zeta(6) +792\zeta(4)\zeta(5)  \right. \nonumber \\ &\left. \hspace{1em}
-288\zeta(2)\zeta(7) -96\zeta(3)^3\right) \label{eq10107} \\
\sum_{k=1}^\infty \frac{H(k)^{3}}{(k+1)^{7}} \sumend &= \frac{1}{160}\left( -501\zeta(10) +800\zeta(3)\zeta(7) +80\zeta(3)^2\zeta(4)  \right. \nonumber \\ &\left. \hspace{1em}
-560\zeta(2)\zeta(3)\zeta(5) +480\zeta(5)^2 +40 M(2,8)\right) \label{eq10108} \\
\sum_{k=1}^\infty \frac{H(k)^{3}}{k^{6}(k+2)} \sumend &= \frac{1}{768}\left( -12 -24\zeta(2) -48\zeta(3) -120\zeta(4) +240\zeta(5) +24\zeta(2)\zeta(3)  \right. \nonumber \\ &\left. \hspace{1em}
-279\zeta(6) +120\zeta(3)^2 +1386\zeta(7) +192\zeta(2)\zeta(5) -1224\zeta(3)\zeta(4) +1190\zeta(8)  \right. \nonumber \\ &\left. \hspace{1em}
+240\zeta(2)\zeta(3)^2 -1152\zeta(3)\zeta(5) -528 M(2,6) +8336\zeta(9) -4656\zeta(3)\zeta(6)  \right. \nonumber \\ &\left. \hspace{1em}
-4896\zeta(4)\zeta(5) +1152\zeta(2)\zeta(7) +768\zeta(3)^3\right) \label{eq10109}
\end{align}
 
\begin{align}
\sum_{k=1}^\infty \frac{H(k)^{3}}{k^{5}(k+1)(k+2)} \sumend &= \frac{1}{384}\left( -12 -24\zeta(2) -48\zeta(3) +3720\zeta(4) -3600\zeta(5)  \right. \nonumber \\ &\left. \hspace{1em}
-360\zeta(2)\zeta(3) +1953\zeta(6) -840\zeta(3)^2 -4158\zeta(7) -576\zeta(2)\zeta(5)  \right. \nonumber \\ &\left. \hspace{1em}
+3672\zeta(3)\zeta(4) -1190\zeta(8) -240\zeta(2)\zeta(3)^2 +1152\zeta(3)\zeta(5) +528 M(2,6)\right) \label{eq10110} \\
\sum_{k=1}^\infty \frac{H(k)^{3}}{k^{4}(k+1)^{2}(k+2)} \sumend &= \frac{1}{64}\left( -4 -8\zeta(2) -16\zeta(3) -1960\zeta(4) +1840\zeta(5)  \right. \nonumber \\ &\left. \hspace{1em}
+200\zeta(2)\zeta(3) -465\zeta(6) +200\zeta(3)^2 +462\zeta(7) +64\zeta(2)\zeta(5)  \right. \nonumber \\ &\left. \hspace{1em}
-408\zeta(3)\zeta(4)\right) \label{eq10111} \\
\sum_{k=1}^\infty \frac{H(k)^{3}}{k^{3}(k+1)^{3}(k+2)} \sumend &= \frac{1}{32}\left( -4 -8\zeta(2) -16\zeta(3) +1240\zeta(4) -1040\zeta(5)  \right. \nonumber \\ &\left. \hspace{1em}
-120\zeta(2)\zeta(3) +159\zeta(6) -104\zeta(3)^2\right) \label{eq10112} \\
\sum_{k=1}^\infty \frac{H(k)^{3}}{k^{2}(k+1)^{4}(k+2)} \sumend &= \frac{-1}{16}\left(  4 +8\zeta(2) +16\zeta(3) +360\zeta(4) -320\zeta(5)  \right. \nonumber \\ &\left. \hspace{1em}
-40\zeta(2)\zeta(3) +33\zeta(6) -32\zeta(3)^2 -119\zeta(7) -32\zeta(2)\zeta(5)  \right. \nonumber \\ &\left. \hspace{1em}
+132\zeta(3)\zeta(4)\right) \label{eq10113} \\
\sum_{k=1}^\infty \frac{H(k)^{3}}{k(k+1)^{5}(k+2)} \sumend &= \frac{1}{96}\left( -48 -96\zeta(2) -192\zeta(3) +480\zeta(4) +198\zeta(6) -192\zeta(3)^2  \right. \nonumber \\ &\left. \hspace{1em}
+43\zeta(8) +120\zeta(2)\zeta(3)^2 -288\zeta(3)\zeta(5) +24 M(2,6)\right) \label{eq10114} \\
\sum_{k=1}^\infty \frac{H(k)^{3}}{(k+1)^{6}(k+2)} \sumend &= \frac{-1}{96}\left(  96 +192\zeta(2) +384\zeta(3) -720\zeta(5) -96\zeta(2)\zeta(3)  \right. \nonumber \\ &\left. \hspace{1em}
-198\zeta(6) +192\zeta(3)^2 -714\zeta(7) -192\zeta(2)\zeta(5) +792\zeta(3)\zeta(4) -43\zeta(8)  \right. \nonumber \\ &\left. \hspace{1em}
-120\zeta(2)\zeta(3)^2 +288\zeta(3)\zeta(5) -24 M(2,6) -788\zeta(9) +444\zeta(3)\zeta(6)  \right. \nonumber \\ &\left. \hspace{1em}
+792\zeta(4)\zeta(5) -288\zeta(2)\zeta(7) -96\zeta(3)^3\right) \label{eq10115} \\
\sum_{k=1}^\infty \frac{H(k)^{3}}{k^{5}(k+2)^{2}} \sumend &= \frac{-1}{768}\left( -156 -192\zeta(2) -264\zeta(3) -402\zeta(4) +1140\zeta(5)  \right. \nonumber \\ &\left. \hspace{1em}
+120\zeta(2)\zeta(3) -837\zeta(6) +360\zeta(3)^2 +2772\zeta(7) +384\zeta(2)\zeta(5)  \right. \nonumber \\ &\left. \hspace{1em}
-2448\zeta(3)\zeta(4) +1190\zeta(8) +240\zeta(2)\zeta(3)^2 -1152\zeta(3)\zeta(5) -528 M(2,6)\right) \label{eq10116} \\
\sum_{k=1}^\infty \frac{H(k)^{3}}{k^{4}(k+1)(k+2)^{2}} \sumend &= \frac{-1}{64}\left( -28 -36\zeta(2) -52\zeta(3) +553\zeta(4) -410\zeta(5)  \right. \nonumber \\ &\left. \hspace{1em}
-40\zeta(2)\zeta(3) +186\zeta(6) -80\zeta(3)^2 -231\zeta(7) -32\zeta(2)\zeta(5)  \right. \nonumber \\ &\left. \hspace{1em}
+204\zeta(3)\zeta(4)\right) \label{eq10117} \\
\sum_{k=1}^\infty \frac{H(k)^{3}}{k^{3}(k+1)^{2}(k+2)^{2}} \sumend &= \frac{-1}{64}\left( -60 -80\zeta(2) -120\zeta(3) -854\zeta(4) +1020\zeta(5)  \right. \nonumber \\ &\left. \hspace{1em}
+120\zeta(2)\zeta(3) -93\zeta(6) +40\zeta(3)^2\right) \label{eq10118}
\end{align}
 
\begin{align}
\sum_{k=1}^\infty \frac{H(k)^{3}}{k^{2}(k+1)^{3}(k+2)^{2}} \sumend &= \frac{-1}{16}\left( -32 -44\zeta(2) -68\zeta(3) +193\zeta(4) -10\zeta(5)  \right. \nonumber \\ &\left. \hspace{1em}
+33\zeta(6) -32\zeta(3)^2\right) \label{eq10119} \\
\sum_{k=1}^\infty \frac{H(k)^{3}}{k(k+1)^{4}(k+2)^{2}} \sumend &= \frac{1}{16}\left(  68 +96\zeta(2) +152\zeta(3) -26\zeta(4) -300\zeta(5)  \right. \nonumber \\ &\left. \hspace{1em}
-40\zeta(2)\zeta(3) -33\zeta(6) +32\zeta(3)^2 -119\zeta(7) -32\zeta(2)\zeta(5)  \right. \nonumber \\ &\left. \hspace{1em}
+132\zeta(3)\zeta(4)\right) \label{eq10120} \\
\sum_{k=1}^\infty \frac{H(k)^{3}}{(k+1)^{5}(k+2)^{2}} \sumend &= \frac{-1}{96}\left( -864 -1248\zeta(2) -2016\zeta(3) +792\zeta(4) +3600\zeta(5)  \right. \nonumber \\ &\left. \hspace{1em}
+480\zeta(2)\zeta(3) +594\zeta(6) -576\zeta(3)^2 +1428\zeta(7) +384\zeta(2)\zeta(5)  \right. \nonumber \\ &\left. \hspace{1em}
-1584\zeta(3)\zeta(4) +43\zeta(8) +120\zeta(2)\zeta(3)^2 -288\zeta(3)\zeta(5) +24 M(2,6)\right) \label{eq10121} \\
\sum_{k=1}^\infty \frac{H(k)^{3}}{k^{4}(k+2)^{3}} \sumend &= \frac{-1}{128}\left(  164 +112\zeta(2) +72\zeta(3) -4\zeta(4) -324\zeta(5)  \right. \nonumber \\ &\left. \hspace{1em}
-64\zeta(2)\zeta(3) +156\zeta(6) -76\zeta(3)^2 -231\zeta(7) -32\zeta(2)\zeta(5)  \right. \nonumber \\ &\left. \hspace{1em}
+204\zeta(3)\zeta(4)\right) \label{eq10122} \\
\sum_{k=1}^\infty \frac{H(k)^{3}}{k^{3}(k+1)(k+2)^{3}} \sumend &= \frac{-1}{64}\left(  192 +148\zeta(2) +124\zeta(3) -557\zeta(4) +86\zeta(5)  \right. \nonumber \\ &\left. \hspace{1em}
-24\zeta(2)\zeta(3) -30\zeta(6) +4\zeta(3)^2\right) \label{eq10123} \\
\sum_{k=1}^\infty \frac{H(k)^{3}}{k^{2}(k+1)^{2}(k+2)^{3}} \sumend &= \frac{1}{64}\left( -444 -376\zeta(2) -368\zeta(3) +260\zeta(4) +848\zeta(5)  \right. \nonumber \\ &\left. \hspace{1em}
+168\zeta(2)\zeta(3) -33\zeta(6) +32\zeta(3)^2\right) \label{eq10124} \\
\sum_{k=1}^\infty \frac{H(k)^{3}}{k(k+1)^{3}(k+2)^{3}} \sumend &= \frac{-1}{32}\left(  508 +464\zeta(2) +504\zeta(3) -646\zeta(4) -828\zeta(5)  \right. \nonumber \\ &\left. \hspace{1em}
-168\zeta(2)\zeta(3) -33\zeta(6) +32\zeta(3)^2\right) \label{eq10125} \\
\sum_{k=1}^\infty \frac{H(k)^{3}}{(k+1)^{4}(k+2)^{3}} \sumend &= \frac{1}{16}\left( -576 -560\zeta(2) -656\zeta(3) +672\zeta(4) +1128\zeta(5)  \right. \nonumber \\ &\left. \hspace{1em}
+208\zeta(2)\zeta(3) +66\zeta(6) -64\zeta(3)^2 +119\zeta(7) +32\zeta(2)\zeta(5)  \right. \nonumber \\ &\left. \hspace{1em}
-132\zeta(3)\zeta(4)\right) \label{eq10126} \\
\sum_{k=1}^\infty \frac{H(k)^{3}}{k^{3}(k+2)^{4}} \sumend &= \frac{1}{128}\left(  676 +216\zeta(2) -96\zeta(3) -386\zeta(4) -256\zeta(5)  \right. \nonumber \\ &\left. \hspace{1em}
-112\zeta(2)\zeta(3) +22\zeta(6) -20\zeta(3)^2 -119\zeta(7) -32\zeta(2)\zeta(5)  \right. \nonumber \\ &\left. \hspace{1em}
+132\zeta(3)\zeta(4)\right) \label{eq10127}
\end{align}
 
\begin{align}
\sum_{k=1}^\infty \frac{H(k)^{3}}{k^{2}(k+1)(k+2)^{4}} \sumend &= \frac{-1}{64}\left( -868 -364\zeta(2) -28\zeta(3) +943\zeta(4) +170\zeta(5)  \right. \nonumber \\ &\left. \hspace{1em}
+136\zeta(2)\zeta(3) +8\zeta(6) +16\zeta(3)^2 +119\zeta(7) +32\zeta(2)\zeta(5)  \right. \nonumber \\ &\left. \hspace{1em}
-132\zeta(3)\zeta(4)\right) \label{eq10128} \\
\sum_{k=1}^\infty \frac{H(k)^{3}}{k(k+1)^{2}(k+2)^{4}} \sumend &= \frac{-1}{64}\left( -2180 -1104\zeta(2) -424\zeta(3) +2146\zeta(4) +1188\zeta(5)  \right. \nonumber \\ &\left. \hspace{1em}
+440\zeta(2)\zeta(3) -17\zeta(6) +64\zeta(3)^2 +238\zeta(7) +64\zeta(2)\zeta(5)  \right. \nonumber \\ &\left. \hspace{1em}
-264\zeta(3)\zeta(4)\right) \label{eq10129} \\
\sum_{k=1}^\infty \frac{H(k)^{3}}{(k+1)^{3}(k+2)^{4}} \sumend &= \frac{1}{16}\left(  1344 +784\zeta(2) +464\zeta(3) -1396\zeta(4) -1008\zeta(5)  \right. \nonumber \\ &\left. \hspace{1em}
-304\zeta(2)\zeta(3) -8\zeta(6) -16\zeta(3)^2 -119\zeta(7) -32\zeta(2)\zeta(5)  \right. \nonumber \\ &\left. \hspace{1em}
+132\zeta(3)\zeta(4)\right) \label{eq10130} \\
\sum_{k=1}^\infty \frac{H(k)^{3}}{k^{2}(k+2)^{5}} \sumend &= \frac{1}{768}\left( -12444 -1224\zeta(2) +4320\zeta(3) +6144\zeta(4) +2232\zeta(5)  \right. \nonumber \\ &\left. \hspace{1em}
-24\zeta(2)\zeta(3) +1911\zeta(6) -1152\zeta(3)^2 +852\zeta(7) +960\zeta(2)\zeta(5)  \right. \nonumber \\ &\left. \hspace{1em}
-1872\zeta(3)\zeta(4) -86\zeta(8) -240\zeta(2)\zeta(3)^2 +576\zeta(3)\zeta(5) -48 M(2,6)\right) \label{eq10131} \\
\sum_{k=1}^\infty \frac{H(k)^{3}}{k(k+1)(k+2)^{5}} \sumend &= \frac{-1}{384}\left(  17652 +3408\zeta(2) -4152\zeta(3) -11802\zeta(4) -3252\zeta(5)  \right. \nonumber \\ &\left. \hspace{1em}
-792\zeta(2)\zeta(3) -1959\zeta(6) +1056\zeta(3)^2 -1566\zeta(7) -1152\zeta(2)\zeta(5)  \right. \nonumber \\ &\left. \hspace{1em}
+2664\zeta(3)\zeta(4) +86\zeta(8) +240\zeta(2)\zeta(3)^2 -576\zeta(3)\zeta(5) +48 M(2,6)\right) \label{eq10132} \\
\sum_{k=1}^\infty \frac{H(k)^{3}}{(k+1)^{2}(k+2)^{5}} \sumend &= \frac{1}{96}\left( -12096 -3360\zeta(2) +1440\zeta(3) +9120\zeta(4) +3408\zeta(5)  \right. \nonumber \\ &\left. \hspace{1em}
+1056\zeta(2)\zeta(3) +954\zeta(6) -432\zeta(3)^2 +1140\zeta(7) +672\zeta(2)\zeta(5)  \right. \nonumber \\ &\left. \hspace{1em}
-1728\zeta(3)\zeta(4) -43\zeta(8) -120\zeta(2)\zeta(3)^2 +288\zeta(3)\zeta(5) -24 M(2,6)\right) \label{eq10133} \\
\sum_{k=1}^\infty \frac{H(k)^{3}}{k(k+2)^{6}} \sumend &= \frac{-1}{768}\left( -30732 +720\zeta(2) +11592\zeta(3) +10758\zeta(4) +8076\zeta(5)  \right. \nonumber \\ &\left. \hspace{1em}
-2712\zeta(2)\zeta(3) +5553\zeta(6) -3360\zeta(3)^2 +2442\zeta(7) +768\zeta(2)\zeta(5)  \right. \nonumber \\ &\left. \hspace{1em}
-2808\zeta(3)\zeta(4) -4118\zeta(8) -240\zeta(2)\zeta(3)^2 +2880\zeta(3)\zeta(5) +1104 M(2,6)  \right. \nonumber \\ &\left. \hspace{1em}
+3152\zeta(9) -1776\zeta(3)\zeta(6) -3168\zeta(4)\zeta(5) +1152\zeta(2)\zeta(7) +384\zeta(3)^3\right) \label{eq10134} \\
\sum_{k=1}^\infty \frac{H(k)^{3}}{(k+1)(k+2)^{6}} \sumend &= \frac{1}{96}\left(  12096 +672\zeta(2) -3936\zeta(3) -5640\zeta(4) -2832\zeta(5)  \right. \nonumber \\ &\left. \hspace{1em}
+480\zeta(2)\zeta(3) -1878\zeta(6) +1104\zeta(3)^2 -1002\zeta(7) -480\zeta(2)\zeta(5)  \right. \nonumber \\ &\left. \hspace{1em}
+1368\zeta(3)\zeta(4) +1051\zeta(8) +120\zeta(2)\zeta(3)^2 -864\zeta(3)\zeta(5) -264 M(2,6)  \right. \nonumber \\ &\left. \hspace{1em}
-788\zeta(9) +444\zeta(3)\zeta(6) +792\zeta(4)\zeta(5) -288\zeta(2)\zeta(7) -96\zeta(3)^3\right) \label{eq10135}
\end{align}
 
\begin{align}
\sum_{k=1}^\infty \frac{H(k)^{3}}{(k+2)^{7}} \sumend &= \frac{1}{160}\left( -13440 +1120\zeta(2) +4640\zeta(3) +3520\zeta(4) +4080\zeta(5)  \right. \nonumber \\ &\left. \hspace{1em}
-1440\zeta(2)\zeta(3) +2300\zeta(6) -1200\zeta(3)^2 +2560\zeta(7) -480\zeta(2)\zeta(5)  \right. \nonumber \\ &\left. \hspace{1em}
-1200\zeta(3)\zeta(4) -1080\zeta(8) +480\zeta(3)\zeta(5) +480 M(2,6) +80\zeta(9) -720\zeta(3)\zeta(6)  \right. \nonumber \\ &\left. \hspace{1em}
-240\zeta(4)\zeta(5) +480\zeta(2)\zeta(7) +160\zeta(3)^3 -501\zeta(10) +800\zeta(3)\zeta(7)  \right. \nonumber \\ &\left. \hspace{1em}
+80\zeta(3)^2\zeta(4) -560\zeta(2)\zeta(3)\zeta(5) +480\zeta(5)^2 +40 M(2,8)\right) \label{eq10136} \\
\sum_{k=1}^\infty \frac{H(k)^{4}}{k^{6}} \sumend &= \frac{-1}{640}\left(  68823\zeta(10) -60000\zeta(3)\zeta(7) -1000\zeta(3)^2\zeta(4)  \right. \nonumber \\ &\left. \hspace{1em}
+21680\zeta(2)\zeta(3)\zeta(5) -23560\zeta(5)^2 -12120 M(2,8) -1280\zeta(2) M(2,6)\right) \label{eq10137} \\
\sum_{k=1}^\infty \frac{H(k)^{4}}{k^{5}(k+1)} \sumend &= \frac{1}{144}\left(  4320\zeta(5) +864\zeta(2)\zeta(3) -5874\zeta(6) -432\zeta(3)^2  \right. \nonumber \\ &\left. \hspace{1em}
+3330\zeta(7) +720\zeta(2)\zeta(5) -3096\zeta(3)\zeta(4) +14833\zeta(8) +4032\zeta(2)\zeta(3)^2  \right. \nonumber \\ &\left. \hspace{1em}
-16704\zeta(3)\zeta(5) -3744 M(2,6) +5232\zeta(9) -3348\zeta(3)\zeta(6) -3096\zeta(4)\zeta(5)  \right. \nonumber \\ &\left. \hspace{1em}
+1008\zeta(2)\zeta(7) +480\zeta(3)^3\right) \label{eq10138} \\
\sum_{k=1}^\infty \frac{H(k)^{4}}{k^{4}(k+1)^{2}} \sumend &= \frac{-1}{144}\left(  17280\zeta(5) +3456\zeta(2)\zeta(3) -22776\zeta(6) -1728\zeta(3)^2  \right. \nonumber \\ &\left. \hspace{1em}
+6660\zeta(7) +1440\zeta(2)\zeta(5) -6192\zeta(3)\zeta(4) +14833\zeta(8) +4032\zeta(2)\zeta(3)^2  \right. \nonumber \\ &\left. \hspace{1em}
-16704\zeta(3)\zeta(5) -3744 M(2,6)\right) \label{eq10139} \\
\sum_{k=1}^\infty \frac{H(k)^{4}}{k^{3}(k+1)^{3}} \sumend &= \frac{1}{4}\left(  720\zeta(5) +144\zeta(2)\zeta(3) -919\zeta(6) -72\zeta(3)^2 +147\zeta(7)  \right. \nonumber \\ &\left. \hspace{1em}
+40\zeta(2)\zeta(5) -160\zeta(3)\zeta(4)\right) \label{eq10140} \\
\sum_{k=1}^\infty \frac{H(k)^{4}}{k^{2}(k+1)^{4}} \sumend &= \frac{-1}{144}\left(  17280\zeta(5) +3456\zeta(2)\zeta(3) -21336\zeta(6) -1728\zeta(3)^2  \right. \nonumber \\ &\left. \hspace{1em}
+3924\zeta(7) +1440\zeta(2)\zeta(5) -5328\zeta(3)\zeta(4) +12415\zeta(8) +3312\zeta(2)\zeta(3)^2  \right. \nonumber \\ &\left. \hspace{1em}
-13824\zeta(3)\zeta(5) -3024 M(2,6)\right) \label{eq10141} \\
\sum_{k=1}^\infty \frac{H(k)^{4}}{k(k+1)^{5}} \sumend &= \frac{-1}{144}\left( -4320\zeta(5) -864\zeta(2)\zeta(3) +5154\zeta(6) +432\zeta(3)^2  \right. \nonumber \\ &\left. \hspace{1em}
-1962\zeta(7) -720\zeta(2)\zeta(5) +2664\zeta(3)\zeta(4) -12415\zeta(8) -3312\zeta(2)\zeta(3)^2  \right. \nonumber \\ &\left. \hspace{1em}
+13824\zeta(3)\zeta(5) +3024 M(2,6) -2088\zeta(9) +1188\zeta(3)\zeta(6) +2664\zeta(4)\zeta(5)  \right. \nonumber \\ &\left. \hspace{1em}
-1008\zeta(2)\zeta(7) -384\zeta(3)^3\right) \label{eq10142} \\
\sum_{k=1}^\infty \frac{H(k)^{4}}{(k+1)^{6}} \sumend &= \frac{-1}{640}\left(  48647\zeta(10) -42080\zeta(3)\zeta(7) +280\zeta(3)^2\zeta(4)  \right. \nonumber \\ &\left. \hspace{1em}
+12720\zeta(2)\zeta(3)\zeta(5) -13320\zeta(5)^2 -7640 M(2,8) -1280\zeta(2) M(2,6)\right) \label{eq10143}
\end{align}
 
\begin{align}
\sum_{k=1}^\infty \frac{H(k)^{4}}{k^{5}(k+2)} \sumend &= \frac{1}{1152}\left(  36 +108\zeta(2) +396\zeta(3) +666\zeta(4) +1080\zeta(5)  \right. \nonumber \\ &\left. \hspace{1em}
+216\zeta(2)\zeta(3) -2937\zeta(6) -216\zeta(3)^2 +3330\zeta(7) +720\zeta(2)\zeta(5)  \right. \nonumber \\ &\left. \hspace{1em}
-3096\zeta(3)\zeta(4) +29666\zeta(8) +8064\zeta(2)\zeta(3)^2 -33408\zeta(3)\zeta(5) -7488 M(2,6)  \right. \nonumber \\ &\left. \hspace{1em}
+20928\zeta(9) -13392\zeta(3)\zeta(6) -12384\zeta(4)\zeta(5) +4032\zeta(2)\zeta(7)  \right. \nonumber \\ &\left. \hspace{1em}
+1920\zeta(3)^3\right) \label{eq10144} \\
\sum_{k=1}^\infty \frac{H(k)^{4}}{k^{4}(k+1)(k+2)} \sumend &= \frac{1}{576}\left(  36 +108\zeta(2) +396\zeta(3) +666\zeta(4) -16200\zeta(5)  \right. \nonumber \\ &\left. \hspace{1em}
-3240\zeta(2)\zeta(3) +20559\zeta(6) +1512\zeta(3)^2 -9990\zeta(7) -2160\zeta(2)\zeta(5)  \right. \nonumber \\ &\left. \hspace{1em}
+9288\zeta(3)\zeta(4) -29666\zeta(8) -8064\zeta(2)\zeta(3)^2 +33408\zeta(3)\zeta(5)  \right. \nonumber \\ &\left. \hspace{1em}
+7488 M(2,6)\right) \label{eq10145} \\
\sum_{k=1}^\infty \frac{H(k)^{4}}{k^{3}(k+1)^{2}(k+2)} \sumend &= \frac{1}{32}\left(  4 +12\zeta(2) +44\zeta(3) +74\zeta(4) +2040\zeta(5)  \right. \nonumber \\ &\left. \hspace{1em}
+408\zeta(2)\zeta(3) -2777\zeta(6) -216\zeta(3)^2 +370\zeta(7) +80\zeta(2)\zeta(5)  \right. \nonumber \\ &\left. \hspace{1em}
-344\zeta(3)\zeta(4)\right) \label{eq10146} \\
\sum_{k=1}^\infty \frac{H(k)^{4}}{k^{2}(k+1)^{3}(k+2)} \sumend &= \frac{-1}{16}\left( -4 -12\zeta(2) -44\zeta(3) -74\zeta(4) +840\zeta(5)  \right. \nonumber \\ &\left. \hspace{1em}
+168\zeta(2)\zeta(3) -899\zeta(6) -72\zeta(3)^2 +218\zeta(7) +80\zeta(2)\zeta(5)  \right. \nonumber \\ &\left. \hspace{1em}
-296\zeta(3)\zeta(4)\right) \label{eq10147} \\
\sum_{k=1}^\infty \frac{H(k)^{4}}{k(k+1)^{4}(k+2)} \sumend &= \frac{-1}{144}\left( -72 -216\zeta(2) -792\zeta(3) -1332\zeta(4) -2160\zeta(5)  \right. \nonumber \\ &\left. \hspace{1em}
-432\zeta(2)\zeta(3) +5154\zeta(6) +432\zeta(3)^2 -12415\zeta(8) -3312\zeta(2)\zeta(3)^2  \right. \nonumber \\ &\left. \hspace{1em}
+13824\zeta(3)\zeta(5) +3024 M(2,6)\right) \label{eq10148} \\
\sum_{k=1}^\infty \frac{H(k)^{4}}{(k+1)^{5}(k+2)} \sumend &= \frac{1}{144}\left(  144 +432\zeta(2) +1584\zeta(3) +2664\zeta(4) -5154\zeta(6)  \right. \nonumber \\ &\left. \hspace{1em}
-432\zeta(3)^2 -1962\zeta(7) -720\zeta(2)\zeta(5) +2664\zeta(3)\zeta(4) +12415\zeta(8)  \right. \nonumber \\ &\left. \hspace{1em}
+3312\zeta(2)\zeta(3)^2 -13824\zeta(3)\zeta(5) -3024 M(2,6) -2088\zeta(9) +1188\zeta(3)\zeta(6)  \right. \nonumber \\ &\left. \hspace{1em}
+2664\zeta(4)\zeta(5) -1008\zeta(2)\zeta(7) -384\zeta(3)^3\right) \label{eq10149} \\
\sum_{k=1}^\infty \frac{H(k)^{4}}{k^{4}(k+2)^{2}} \sumend &= \frac{-1}{576}\left(  252 +504\zeta(2) +1440\zeta(3) +1404\zeta(4) +1080\zeta(5)  \right. \nonumber \\ &\left. \hspace{1em}
+288\zeta(2)\zeta(3) -5694\zeta(6) -432\zeta(3)^2 +3330\zeta(7) +720\zeta(2)\zeta(5)  \right. \nonumber \\ &\left. \hspace{1em}
-3096\zeta(3)\zeta(4) +14833\zeta(8) +4032\zeta(2)\zeta(3)^2 -16704\zeta(3)\zeta(5)  \right. \nonumber \\ &\left. \hspace{1em}
-3744 M(2,6)\right) \label{eq10150}
\end{align}
 
\begin{align}
\sum_{k=1}^\infty \frac{H(k)^{4}}{k^{3}(k+1)(k+2)^{2}} \sumend &= \frac{-1}{64}\left(  60 +124\zeta(2) +364\zeta(3) +386\zeta(4) -1560\zeta(5)  \right. \nonumber \\ &\left. \hspace{1em}
-296\zeta(2)\zeta(3) +1019\zeta(6) +72\zeta(3)^2 -370\zeta(7) -80\zeta(2)\zeta(5)  \right. \nonumber \\ &\left. \hspace{1em}
+344\zeta(3)\zeta(4)\right) \label{eq10151} \\
\sum_{k=1}^\infty \frac{H(k)^{4}}{k^{2}(k+1)^{2}(k+2)^{2}} \sumend &= \frac{-1}{16}\left(  32 +68\zeta(2) +204\zeta(3) +230\zeta(4) +240\zeta(5)  \right. \nonumber \\ &\left. \hspace{1em}
+56\zeta(2)\zeta(3) -879\zeta(6) -72\zeta(3)^2\right) \label{eq10152} \\
\sum_{k=1}^\infty \frac{H(k)^{4}}{k(k+1)^{3}(k+2)^{2}} \sumend &= \frac{-1}{16}\left(  68 +148\zeta(2) +452\zeta(3) +534\zeta(4) -360\zeta(5)  \right. \nonumber \\ &\left. \hspace{1em}
-56\zeta(2)\zeta(3) -859\zeta(6) -72\zeta(3)^2 -218\zeta(7) -80\zeta(2)\zeta(5)  \right. \nonumber \\ &\left. \hspace{1em}
+296\zeta(3)\zeta(4)\right) \label{eq10153} \\
\sum_{k=1}^\infty \frac{H(k)^{4}}{(k+1)^{4}(k+2)^{2}} \sumend &= \frac{-1}{144}\left(  1296 +2880\zeta(2) +8928\zeta(3) +10944\zeta(4) -4320\zeta(5)  \right. \nonumber \\ &\left. \hspace{1em}
-576\zeta(2)\zeta(3) -20616\zeta(6) -1728\zeta(3)^2 -3924\zeta(7) -1440\zeta(2)\zeta(5)  \right. \nonumber \\ &\left. \hspace{1em}
+5328\zeta(3)\zeta(4) +12415\zeta(8) +3312\zeta(2)\zeta(3)^2 -13824\zeta(3)\zeta(5)  \right. \nonumber \\ &\left. \hspace{1em}
-3024 M(2,6)\right) \label{eq10154} \\
\sum_{k=1}^\infty \frac{H(k)^{4}}{k^{3}(k+2)^{3}} \sumend &= \frac{1}{64}\left(  192 +244\zeta(2) +508\zeta(3) +122\zeta(4) -168\zeta(5)  \right. \nonumber \\ &\left. \hspace{1em}
-56\zeta(2)\zeta(3) -853\zeta(6) -136\zeta(3)^2 +294\zeta(7) +80\zeta(2)\zeta(5)  \right. \nonumber \\ &\left. \hspace{1em}
-320\zeta(3)\zeta(4)\right) \label{eq10155} \\
\sum_{k=1}^\infty \frac{H(k)^{4}}{k^{2}(k+1)(k+2)^{3}} \sumend &= \frac{1}{64}\left(  444 +612\zeta(2) +1380\zeta(3) +630\zeta(4) -1896\zeta(5)  \right. \nonumber \\ &\left. \hspace{1em}
-408\zeta(2)\zeta(3) -687\zeta(6) -200\zeta(3)^2 +218\zeta(7) +80\zeta(2)\zeta(5)  \right. \nonumber \\ &\left. \hspace{1em}
-296\zeta(3)\zeta(4)\right) \label{eq10156} \\
\sum_{k=1}^\infty \frac{H(k)^{4}}{k(k+1)^{2}(k+2)^{3}} \sumend &= \frac{-1}{32}\left( -508 -748\zeta(2) -1788\zeta(3) -1090\zeta(4) +1416\zeta(5)  \right. \nonumber \\ &\left. \hspace{1em}
+296\zeta(2)\zeta(3) +2445\zeta(6) +344\zeta(3)^2 -218\zeta(7) -80\zeta(2)\zeta(5)  \right. \nonumber \\ &\left. \hspace{1em}
+296\zeta(3)\zeta(4)\right) \label{eq10157} \\
\sum_{k=1}^\infty \frac{H(k)^{4}}{(k+1)^{3}(k+2)^{3}} \sumend &= \frac{-1}{2}\left( -72 -112\zeta(2) -280\zeta(3) -203\zeta(4) +222\zeta(5)  \right. \nonumber \\ &\left. \hspace{1em}
+44\zeta(2)\zeta(3) +413\zeta(6) +52\zeta(3)^2\right) \label{eq10158} \\
\sum_{k=1}^\infty \frac{H(k)^{4}}{k^{2}(k+2)^{4}} \sumend &= \frac{1}{576}\left( -7812 -5976\zeta(2) -7920\zeta(3) +5724\zeta(4) +6984\zeta(5)  \right. \nonumber \\ &\left. \hspace{1em}
+3168\zeta(2)\zeta(3) +4290\zeta(6) +1872\zeta(3)^2 +2322\zeta(7) +432\zeta(2)\zeta(5)  \right. \nonumber \\ &\left. \hspace{1em}
-2088\zeta(3)\zeta(4) -12415\zeta(8) -3312\zeta(2)\zeta(3)^2 +13824\zeta(3)\zeta(5)  \right. \nonumber \\ &\left. \hspace{1em}
+3024 M(2,6)\right) \label{eq10159}
\end{align}
 
\begin{align}
\sum_{k=1}^\infty \frac{H(k)^{4}}{k(k+1)(k+2)^{4}} \sumend &= \frac{1}{576}\left( -19620 -17460\zeta(2) -28260\zeta(3) +5778\zeta(4) +31032\zeta(5)  \right. \nonumber \\ &\left. \hspace{1em}
+10008\zeta(2)\zeta(3) +14763\zeta(6) +5544\zeta(3)^2 +2682\zeta(7) +144\zeta(2)\zeta(5)  \right. \nonumber \\ &\left. \hspace{1em}
-1512\zeta(3)\zeta(4) -24830\zeta(8) -6624\zeta(2)\zeta(3)^2 +27648\zeta(3)\zeta(5)  \right. \nonumber \\ &\left. \hspace{1em}
+6048 M(2,6)\right) \label{eq10160} \\
\sum_{k=1}^\infty \frac{H(k)^{4}}{(k+1)^{2}(k+2)^{4}} \sumend &= \frac{1}{144}\left( -12096 -12096\zeta(2) -22176\zeta(3) -2016\zeta(4) +21888\zeta(5)  \right. \nonumber \\ &\left. \hspace{1em}
+6336\zeta(2)\zeta(3) +18384\zeta(6) +4320\zeta(3)^2 +360\zeta(7) -288\zeta(2)\zeta(5)  \right. \nonumber \\ &\left. \hspace{1em}
+576\zeta(3)\zeta(4) -12415\zeta(8) -3312\zeta(2)\zeta(3)^2 +13824\zeta(3)\zeta(5) +3024 M(2,6)\right) \label{eq10161} \\
\sum_{k=1}^\infty \frac{H(k)^{4}}{k(k+2)^{5}} \sumend &= \frac{-1}{1152}\left( -52956 -22860\zeta(2) -12060\zeta(3) +44622\zeta(4) +25992\zeta(5)  \right. \nonumber \\ &\left. \hspace{1em}
+10728\zeta(2)\zeta(3) +9309\zeta(6) -1512\zeta(3)^2 +20286\zeta(7) +9648\zeta(2)\zeta(5)  \right. \nonumber \\ &\left. \hspace{1em}
-27576\zeta(3)\zeta(4) -25862\zeta(8) -9504\zeta(2)\zeta(3)^2 +34560\zeta(3)\zeta(5) +5472 M(2,6)  \right. \nonumber \\ &\left. \hspace{1em}
-8352\zeta(9) +4752\zeta(3)\zeta(6) +10656\zeta(4)\zeta(5) -4032\zeta(2)\zeta(7)  \right. \nonumber \\ &\left. \hspace{1em}
-1536\zeta(3)^3\right) \label{eq10162} \\
\sum_{k=1}^\infty \frac{H(k)^{4}}{(k+1)(k+2)^{5}} \sumend &= \frac{-1}{144}\left( -18144 -10080\zeta(2) -10080\zeta(3) +12600\zeta(4) +14256\zeta(5)  \right. \nonumber \\ &\left. \hspace{1em}
+5184\zeta(2)\zeta(3) +6018\zeta(6) +1008\zeta(3)^2 +5742\zeta(7) +2448\zeta(2)\zeta(5)  \right. \nonumber \\ &\left. \hspace{1em}
-7272\zeta(3)\zeta(4) -12673\zeta(8) -4032\zeta(2)\zeta(3)^2 +15552\zeta(3)\zeta(5) +2880 M(2,6)  \right. \nonumber \\ &\left. \hspace{1em}
-2088\zeta(9) +1188\zeta(3)\zeta(6) +2664\zeta(4)\zeta(5) -1008\zeta(2)\zeta(7) -384\zeta(3)^3\right) \label{eq10163} \\
\sum_{k=1}^\infty \frac{H(k)^{4}}{(k+2)^{6}} \sumend &= \frac{1}{1920}\left( -241920 -53760\zeta(2) +26880\zeta(3) +161280\zeta(4) +88320\zeta(5)  \right. \nonumber \\ &\left. \hspace{1em}
+7680\zeta(2)\zeta(3) +56640\zeta(6) -30720\zeta(3)^2 +57120\zeta(7) +30720\zeta(2)\zeta(5)  \right. \nonumber \\ &\left. \hspace{1em}
-82560\zeta(3)\zeta(4) -43760\zeta(8) -9600\zeta(2)\zeta(3)^2 +46080\zeta(3)\zeta(5) +9600 M(2,6)  \right. \nonumber \\ &\left. \hspace{1em}
+63040\zeta(9) -35520\zeta(3)\zeta(6) -63360\zeta(4)\zeta(5) +23040\zeta(2)\zeta(7) +7680\zeta(3)^3  \right. \nonumber \\ &\left. \hspace{1em}
-145941\zeta(10) +126240\zeta(3)\zeta(7) -840\zeta(3)^2\zeta(4) -38160\zeta(2)\zeta(3)\zeta(5)  \right. \nonumber \\ &\left. \hspace{1em}
+39960\zeta(5)^2 +22920 M(2,8) +3840\zeta(2) M(2,6)\right) \label{eq10164} \\
\sum_{k=1}^\infty \frac{H(k)^{5}}{k^{5}} \sumend &= \frac{1}{256}\left( -64433\zeta(10) +57760\zeta(3)\zeta(7) -360\zeta(3)^2\zeta(4)  \right. \nonumber \\ &\left. \hspace{1em}
-20560\zeta(2)\zeta(3)\zeta(5) +22648\zeta(5)^2 +10920 M(2,8) +1280\zeta(2) M(2,6)\right) \label{eq10165} \\
\sum_{k=1}^\infty \frac{H(k)^{5}}{k^{4}(k+1)} \sumend &= \frac{1}{288}\left( -51408\zeta(6) -6480\zeta(3)^2 +36918\zeta(7) +8208\zeta(2)\zeta(5)  \right. \nonumber \\ &\left. \hspace{1em}
+9504\zeta(3)\zeta(4) +67811\zeta(8) +19080\zeta(2)\zeta(3)^2 -78768\zeta(3)\zeta(5) -16920 M(2,6)  \right. \nonumber \\ &\left. \hspace{1em}
+37768\zeta(9) -58740\zeta(3)\zeta(6) +19008\zeta(4)\zeta(5) +9540\zeta(2)\zeta(7)  \right. \nonumber \\ &\left. \hspace{1em}
-1440\zeta(3)^3\right) \label{eq10166}
\end{align}
 
\begin{align}
\sum_{k=1}^\infty \frac{H(k)^{5}}{k^{3}(k+1)^{2}} \sumend &= \frac{-1}{288}\left( -154224\zeta(6) -19440\zeta(3)^2 +107226\zeta(7) +24624\zeta(2)\zeta(5)  \right. \nonumber \\ &\left. \hspace{1em}
+28512\zeta(3)\zeta(4) +67811\zeta(8) +19080\zeta(2)\zeta(3)^2 -78768\zeta(3)\zeta(5)  \right. \nonumber \\ &\left. \hspace{1em}
-16920 M(2,6)\right) \label{eq10167} \\
\sum_{k=1}^\infty \frac{H(k)^{5}}{k^{2}(k+1)^{3}} \sumend &= \frac{1}{288}\left( -154224\zeta(6) -19440\zeta(3)^2 +103698\zeta(7) +24624\zeta(2)\zeta(5)  \right. \nonumber \\ &\left. \hspace{1em}
+28512\zeta(3)\zeta(4) +65621\zeta(8) +17640\zeta(2)\zeta(3)^2 -72432\zeta(3)\zeta(5)  \right. \nonumber \\ &\left. \hspace{1em}
-15480 M(2,6)\right) \label{eq10168} \\
\sum_{k=1}^\infty \frac{H(k)^{5}}{k(k+1)^{4}} \sumend &= \frac{1}{288}\left(  51408\zeta(6) +6480\zeta(3)^2 -33390\zeta(7) -8208\zeta(2)\zeta(5)  \right. \nonumber \\ &\left. \hspace{1em}
-9504\zeta(3)\zeta(4) -65621\zeta(8) -17640\zeta(2)\zeta(3)^2 +72432\zeta(3)\zeta(5) +15480 M(2,6)  \right. \nonumber \\ &\left. \hspace{1em}
-28480\zeta(9) +51540\zeta(3)\zeta(6) -19008\zeta(4)\zeta(5) -9540\zeta(2)\zeta(7)  \right. \nonumber \\ &\left. \hspace{1em}
+1440\zeta(3)^3\right) \label{eq10169} \\
\sum_{k=1}^\infty \frac{H(k)^{5}}{(k+1)^{5}} \sumend &= \frac{1}{256}\left(  49901\zeta(10) -43040\zeta(3)\zeta(7) -1080\zeta(3)^2\zeta(4)  \right. \nonumber \\ &\left. \hspace{1em}
+13840\zeta(2)\zeta(3)\zeta(5) -13592\zeta(5)^2 -7560 M(2,8) -1280\zeta(2) M(2,6)\right) \label{eq10170} \\
\sum_{k=1}^\infty \frac{H(k)^{5}}{k^{4}(k+2)} \sumend &= \frac{-1}{1152}\left(  72 +288\zeta(2) +1512\zeta(3) +4518\zeta(4) +5112\zeta(5)  \right. \nonumber \\ &\left. \hspace{1em}
+1080\zeta(2)\zeta(3) +12852\zeta(6) +1620\zeta(3)^2 -18459\zeta(7) -4104\zeta(2)\zeta(5)  \right. \nonumber \\ &\left. \hspace{1em}
-4752\zeta(3)\zeta(4) -67811\zeta(8) -19080\zeta(2)\zeta(3)^2 +78768\zeta(3)\zeta(5) +16920 M(2,6)  \right. \nonumber \\ &\left. \hspace{1em}
-75536\zeta(9) +117480\zeta(3)\zeta(6) -38016\zeta(4)\zeta(5) -19080\zeta(2)\zeta(7)  \right. \nonumber \\ &\left. \hspace{1em}
+2880\zeta(3)^3\right) \label{eq10171} \\
\sum_{k=1}^\infty \frac{H(k)^{5}}{k^{3}(k+1)(k+2)} \sumend &= \frac{-1}{576}\left(  72 +288\zeta(2) +1512\zeta(3) +4518\zeta(4) +5112\zeta(5)  \right. \nonumber \\ &\left. \hspace{1em}
+1080\zeta(2)\zeta(3) -89964\zeta(6) -11340\zeta(3)^2 +55377\zeta(7) +12312\zeta(2)\zeta(5)  \right. \nonumber \\ &\left. \hspace{1em}
+14256\zeta(3)\zeta(4) +67811\zeta(8) +19080\zeta(2)\zeta(3)^2 -78768\zeta(3)\zeta(5)  \right. \nonumber \\ &\left. \hspace{1em}
-16920 M(2,6)\right) \label{eq10172} \\
\sum_{k=1}^\infty \frac{H(k)^{5}}{k^{2}(k+1)^{2}(k+2)} \sumend &= \frac{-1}{32}\left(  8 +32\zeta(2) +168\zeta(3) +502\zeta(4) +568\zeta(5)  \right. \nonumber \\ &\left. \hspace{1em}
+120\zeta(2)\zeta(3) +7140\zeta(6) +900\zeta(3)^2 -5761\zeta(7) -1368\zeta(2)\zeta(5)  \right. \nonumber \\ &\left. \hspace{1em}
-1584\zeta(3)\zeta(4)\right) \label{eq10173} \\
\sum_{k=1}^\infty \frac{H(k)^{5}}{k(k+1)^{3}(k+2)} \sumend &= \frac{-1}{288}\left(  144 +576\zeta(2) +3024\zeta(3) +9036\zeta(4) +10224\zeta(5)  \right. \nonumber \\ &\left. \hspace{1em}
+2160\zeta(2)\zeta(3) -25704\zeta(6) -3240\zeta(3)^2 +65621\zeta(8) +17640\zeta(2)\zeta(3)^2  \right. \nonumber \\ &\left. \hspace{1em}
-72432\zeta(3)\zeta(5) -15480 M(2,6)\right) \label{eq10174}
\end{align}
 
\begin{align}
\sum_{k=1}^\infty \frac{H(k)^{5}}{(k+1)^{4}(k+2)} \sumend &= \frac{1}{288}\left( -288 -1152\zeta(2) -6048\zeta(3) -18072\zeta(4) -20448\zeta(5)  \right. \nonumber \\ &\left. \hspace{1em}
-4320\zeta(2)\zeta(3) +33390\zeta(7) +8208\zeta(2)\zeta(5) +9504\zeta(3)\zeta(4) -65621\zeta(8)  \right. \nonumber \\ &\left. \hspace{1em}
-17640\zeta(2)\zeta(3)^2 +72432\zeta(3)\zeta(5) +15480 M(2,6) +28480\zeta(9) -51540\zeta(3)\zeta(6)  \right. \nonumber \\ &\left. \hspace{1em}
+19008\zeta(4)\zeta(5) +9540\zeta(2)\zeta(7) -1440\zeta(3)^3\right) \label{eq10175} \\
\sum_{k=1}^\infty \frac{H(k)^{5}}{k^{3}(k+2)^{2}} \sumend &= \frac{-1}{1152}\left( -1080 -3024\zeta(2) -12888\zeta(3) -27666\zeta(4) -14760\zeta(5)  \right. \nonumber \\ &\left. \hspace{1em}
-3960\zeta(2)\zeta(3) -12786\zeta(6) -2700\zeta(3)^2 +53613\zeta(7) +12312\zeta(2)\zeta(5)  \right. \nonumber \\ &\left. \hspace{1em}
+14256\zeta(3)\zeta(4) +67811\zeta(8) +19080\zeta(2)\zeta(3)^2 -78768\zeta(3)\zeta(5)  \right. \nonumber \\ &\left. \hspace{1em}
-16920 M(2,6)\right) \label{eq10176} \\
\sum_{k=1}^\infty \frac{H(k)^{5}}{k^{2}(k+1)(k+2)^{2}} \sumend &= \frac{1}{96}\left(  192 +552\zeta(2) +2400\zeta(3) +5364\zeta(4) +3312\zeta(5)  \right. \nonumber \\ &\left. \hspace{1em}
+840\zeta(2)\zeta(3) -12863\zeta(6) -1440\zeta(3)^2 +294\zeta(7)\right) \label{eq10177} \\
\sum_{k=1}^\infty \frac{H(k)^{5}}{k(k+1)^{2}(k+2)^{2}} \sumend &= \frac{-1}{96}\left( -408 -1200\zeta(2) -5304\zeta(3) -12234\zeta(4) -8328\zeta(5)  \right. \nonumber \\ &\left. \hspace{1em}
-2040\zeta(2)\zeta(3) +4306\zeta(6) +180\zeta(3)^2 +16695\zeta(7) +4104\zeta(2)\zeta(5)  \right. \nonumber \\ &\left. \hspace{1em}
+4752\zeta(3)\zeta(4)\right) \label{eq10178} \\
\sum_{k=1}^\infty \frac{H(k)^{5}}{(k+1)^{3}(k+2)^{2}} \sumend &= \frac{-1}{288}\left( -2592 -7776\zeta(2) -34848\zeta(3) -82440\zeta(4) -60192\zeta(5)  \right. \nonumber \\ &\left. \hspace{1em}
-14400\zeta(2)\zeta(3) +51540\zeta(6) +4320\zeta(3)^2 +100170\zeta(7) +24624\zeta(2)\zeta(5)  \right. \nonumber \\ &\left. \hspace{1em}
+28512\zeta(3)\zeta(4) -65621\zeta(8) -17640\zeta(2)\zeta(3)^2 +72432\zeta(3)\zeta(5)  \right. \nonumber \\ &\left. \hspace{1em}
+15480 M(2,6)\right) \label{eq10179} \\
\sum_{k=1}^\infty \frac{H(k)^{5}}{k^{2}(k+2)^{3}} \sumend &= \frac{1}{1152}\left( -7992 -15264\zeta(2) -51480\zeta(3) -71874\zeta(4) +4248\zeta(5)  \right. \nonumber \\ &\left. \hspace{1em}
-360\zeta(2)\zeta(3) +58584\zeta(6) +9540\zeta(3)^2 +32229\zeta(7) +5112\zeta(2)\zeta(5)  \right. \nonumber \\ &\left. \hspace{1em}
+40896\zeta(3)\zeta(4) +65621\zeta(8) +17640\zeta(2)\zeta(3)^2 -72432\zeta(3)\zeta(5)  \right. \nonumber \\ &\left. \hspace{1em}
-15480 M(2,6)\right) \label{eq10180} \\
\sum_{k=1}^\infty \frac{H(k)^{5}}{k(k+1)(k+2)^{3}} \sumend &= \frac{-1}{576}\left(  9144 +18576\zeta(2) +65880\zeta(3) +104058\zeta(4) +15624\zeta(5)  \right. \nonumber \\ &\left. \hspace{1em}
+5400\zeta(2)\zeta(3) -135762\zeta(6) -18180\zeta(3)^2 -30465\zeta(7) -5112\zeta(2)\zeta(5)  \right. \nonumber \\ &\left. \hspace{1em}
-40896\zeta(3)\zeta(4) -65621\zeta(8) -17640\zeta(2)\zeta(3)^2 +72432\zeta(3)\zeta(5)  \right. \nonumber \\ &\left. \hspace{1em}
+15480 M(2,6)\right) \label{eq10181}
\end{align}
 
\begin{align}
\sum_{k=1}^\infty \frac{H(k)^{5}}{(k+1)^{2}(k+2)^{3}} \sumend &= \frac{1}{288}\left( -10368 -22176\zeta(2) -81792\zeta(3) -140760\zeta(4) -40608\zeta(5)  \right. \nonumber \\ &\left. \hspace{1em}
-11520\zeta(2)\zeta(3) +148680\zeta(6) +18720\zeta(3)^2 +80550\zeta(7) +17424\zeta(2)\zeta(5)  \right. \nonumber \\ &\left. \hspace{1em}
+55152\zeta(3)\zeta(4) +65621\zeta(8) +17640\zeta(2)\zeta(3)^2 -72432\zeta(3)\zeta(5)  \right. \nonumber \\ &\left. \hspace{1em}
-15480 M(2,6)\right) \label{eq10182} \\
\sum_{k=1}^\infty \frac{H(k)^{5}}{k(k+2)^{4}} \sumend &= \frac{1}{1152}\left(  39240 +49968\zeta(2) +129384\zeta(3) +94806\zeta(4) -82440\zeta(5)  \right. \nonumber \\ &\left. \hspace{1em}
-25560\zeta(2)\zeta(3) -146718\zeta(6) -36540\zeta(3)^2 -675\zeta(7) +5976\zeta(2)\zeta(5)  \right. \nonumber \\ &\left. \hspace{1em}
-37152\zeta(3)\zeta(4) +182679\zeta(8) +48600\zeta(2)\zeta(3)^2 -204048\zeta(3)\zeta(5) -45000 M(2,6)  \right. \nonumber \\ &\left. \hspace{1em}
-56960\zeta(9) +103080\zeta(3)\zeta(6) -38016\zeta(4)\zeta(5) -19080\zeta(2)\zeta(7)  \right. \nonumber \\ &\left. \hspace{1em}
+2880\zeta(3)^3\right) \label{eq10183} \\
\sum_{k=1}^\infty \frac{H(k)^{5}}{(k+1)(k+2)^{4}} \sumend &= \frac{1}{288}\left(  24192 +34272\zeta(2) +97632\zeta(3) +99432\zeta(4) -33408\zeta(5)  \right. \nonumber \\ &\left. \hspace{1em}
-10080\zeta(2)\zeta(3) -141240\zeta(6) -27360\zeta(3)^2 -15570\zeta(7) +432\zeta(2)\zeta(5)  \right. \nonumber \\ &\left. \hspace{1em}
-39024\zeta(3)\zeta(4) +58529\zeta(8) +15480\zeta(2)\zeta(3)^2 -65808\zeta(3)\zeta(5) -14760 M(2,6)  \right. \nonumber \\ &\left. \hspace{1em}
-28480\zeta(9) +51540\zeta(3)\zeta(6) -19008\zeta(4)\zeta(5) -9540\zeta(2)\zeta(7)  \right. \nonumber \\ &\left. \hspace{1em}
+1440\zeta(3)^3\right) \label{eq10184} \\
\sum_{k=1}^\infty \frac{H(k)^{5}}{(k+2)^{5}} \sumend &= \frac{1}{2304}\left( -290304 -241920\zeta(2) -460800\zeta(3) -48960\zeta(4) +414720\zeta(5)  \right. \nonumber \\ &\left. \hspace{1em}
+149760\zeta(2)\zeta(3) +384960\zeta(6) +97920\zeta(3)^2 +162720\zeta(7) +57600\zeta(2)\zeta(5)  \right. \nonumber \\ &\left. \hspace{1em}
-178560\zeta(3)\zeta(4) -1003520\zeta(8) -293760\zeta(2)\zeta(3)^2 +1175040\zeta(3)\zeta(5)  \right. \nonumber \\ &\left. \hspace{1em}
+236160 M(2,6) -167040\zeta(9) +95040\zeta(3)\zeta(6) +213120\zeta(4)\zeta(5) -80640\zeta(2)\zeta(7)  \right. \nonumber \\ &\left. \hspace{1em}
-30720\zeta(3)^3 +449109\zeta(10) -387360\zeta(3)\zeta(7) -9720\zeta(3)^2\zeta(4)  \right. \nonumber \\ &\left. \hspace{1em}
+124560\zeta(2)\zeta(3)\zeta(5) -122328\zeta(5)^2 -68040 M(2,8) -11520\zeta(2) M(2,6)\right) \label{eq10185} \\
\sum_{k=1}^\infty \frac{H(k)^{6}}{k^{4}} \sumend &= \frac{-1}{128}\left(  271367\zeta(10) -176560\zeta(3)\zeta(7) +84648\zeta(3)^2\zeta(4)  \right. \nonumber \\ &\left. \hspace{1em}
+400\zeta(2)\zeta(3)\zeta(5) -121688\zeta(5)^2 -34376 M(2,8) -15040\zeta(2) M(2,6)\right) \label{eq10186} \\
\sum_{k=1}^\infty \frac{H(k)^{6}}{k^{3}(k+1)} \sumend &= \frac{-1}{24}\left( -15456\zeta(7) -3480\zeta(2)\zeta(5) -7128\zeta(3)\zeta(4) +17529\zeta(8)  \right. \nonumber \\ &\left. \hspace{1em}
-984\zeta(2)\zeta(3)^2 +11688\zeta(3)\zeta(5) +1368 M(2,6) -7474\zeta(9) +13122\zeta(3)\zeta(6)  \right. \nonumber \\ &\left. \hspace{1em}
-6048\zeta(4)\zeta(5) -1953\zeta(2)\zeta(7) +544\zeta(3)^3\right) \label{eq10187} \\
\sum_{k=1}^\infty \frac{H(k)^{6}}{k^{2}(k+1)^{2}} \sumend &= \frac{-1}{6}\left(  7728\zeta(7) +1740\zeta(2)\zeta(5) +3564\zeta(3)\zeta(4) -8639\zeta(8)  \right. \nonumber \\ &\left. \hspace{1em}
+477\zeta(2)\zeta(3)^2 -5754\zeta(3)\zeta(5) -669 M(2,6)\right) \label{eq10188}
\end{align}
 
\begin{align}
\sum_{k=1}^\infty \frac{H(k)^{6}}{k(k+1)^{3}} \sumend &= \frac{1}{24}\left(  15456\zeta(7) +3480\zeta(2)\zeta(5) +7128\zeta(3)\zeta(4) -17027\zeta(8)  \right. \nonumber \\ &\left. \hspace{1em}
+924\zeta(2)\zeta(3)^2 -11328\zeta(3)\zeta(5) -1308 M(2,6) +6146\zeta(9) -12582\zeta(3)\zeta(6)  \right. \nonumber \\ &\left. \hspace{1em}
+5832\zeta(4)\zeta(5) +1953\zeta(2)\zeta(7) -536\zeta(3)^3\right) \label{eq10189} \\
\sum_{k=1}^\infty \frac{H(k)^{6}}{(k+1)^{4}} \sumend &= \frac{-1}{128}\left(  259945\zeta(10) -163568\zeta(3)\zeta(7) +81848\zeta(3)^2\zeta(4)  \right. \nonumber \\ &\left. \hspace{1em}
-5200\zeta(2)\zeta(3)\zeta(5) -113288\zeta(5)^2 -31576 M(2,8) -15040\zeta(2) M(2,6)\right) \label{eq10190} \\
\sum_{k=1}^\infty \frac{H(k)^{6}}{k^{3}(k+2)} \sumend &= \frac{-1}{96}\left( -12 -60\zeta(2) -408\zeta(3) -1713\zeta(4) -3426\zeta(5)  \right. \nonumber \\ &\left. \hspace{1em}
-732\zeta(2)\zeta(3) -6291\zeta(6) -804\zeta(3)^2 -7728\zeta(7) -1740\zeta(2)\zeta(5)  \right. \nonumber \\ &\left. \hspace{1em}
-3564\zeta(3)\zeta(4) +17529\zeta(8) -984\zeta(2)\zeta(3)^2 +11688\zeta(3)\zeta(5) +1368 M(2,6)  \right. \nonumber \\ &\left. \hspace{1em}
-14948\zeta(9) +26244\zeta(3)\zeta(6) -12096\zeta(4)\zeta(5) -3906\zeta(2)\zeta(7)  \right. \nonumber \\ &\left. \hspace{1em}
+1088\zeta(3)^3\right) \label{eq10191} \\
\sum_{k=1}^\infty \frac{H(k)^{6}}{k^{2}(k+1)(k+2)} \sumend &= \frac{1}{16}\left(  4 +20\zeta(2) +136\zeta(3) +571\zeta(4) +1142\zeta(5)  \right. \nonumber \\ &\left. \hspace{1em}
+244\zeta(2)\zeta(3) +2097\zeta(6) +268\zeta(3)^2 -7728\zeta(7) -1740\zeta(2)\zeta(5)  \right. \nonumber \\ &\left. \hspace{1em}
-3564\zeta(3)\zeta(4) +5843\zeta(8) -328\zeta(2)\zeta(3)^2 +3896\zeta(3)\zeta(5) +456 M(2,6)\right) \label{eq10192} \\
\sum_{k=1}^\infty \frac{H(k)^{6}}{k(k+1)^{2}(k+2)} \sumend &= \frac{1}{24}\left(  12 +60\zeta(2) +408\zeta(3) +1713\zeta(4) +3426\zeta(5)  \right. \nonumber \\ &\left. \hspace{1em}
+732\zeta(2)\zeta(3) +6291\zeta(6) +804\zeta(3)^2 +7728\zeta(7) +1740\zeta(2)\zeta(5)  \right. \nonumber \\ &\left. \hspace{1em}
+3564\zeta(3)\zeta(4) -17027\zeta(8) +924\zeta(2)\zeta(3)^2 -11328\zeta(3)\zeta(5) -1308 M(2,6)\right) \label{eq10193} \\
\sum_{k=1}^\infty \frac{H(k)^{6}}{(k+1)^{3}(k+2)} \sumend &= \frac{1}{24}\left(  24 +120\zeta(2) +816\zeta(3) +3426\zeta(4) +6852\zeta(5)  \right. \nonumber \\ &\left. \hspace{1em}
+1464\zeta(2)\zeta(3) +12582\zeta(6) +1608\zeta(3)^2 -17027\zeta(8) +924\zeta(2)\zeta(3)^2  \right. \nonumber \\ &\left. \hspace{1em}
-11328\zeta(3)\zeta(5) -1308 M(2,6) -6146\zeta(9) +12582\zeta(3)\zeta(6) -5832\zeta(4)\zeta(5)  \right. \nonumber \\ &\left. \hspace{1em}
-1953\zeta(2)\zeta(7) +536\zeta(3)^3\right) \label{eq10194} \\
\sum_{k=1}^\infty \frac{H(k)^{6}}{k^{2}(k+2)^{2}} \sumend &= \frac{1}{96}\left( -192 -696\zeta(2) -3936\zeta(3) -12714\zeta(4) -16956\zeta(5)  \right. \nonumber \\ &\left. \hspace{1em}
-3912\zeta(2)\zeta(3) -12279\zeta(6) -2136\zeta(3)^2 +1239\zeta(7) +624\zeta(2)\zeta(5)  \right. \nonumber \\ &\left. \hspace{1em}
-2376\zeta(3)\zeta(4) +34556\zeta(8) -1908\zeta(2)\zeta(3)^2 +23016\zeta(3)\zeta(5)  \right. \nonumber \\ &\left. \hspace{1em}
+2676 M(2,6)\right) \label{eq10195} \\
\sum_{k=1}^\infty \frac{H(k)^{6}}{k(k+1)(k+2)^{2}} \sumend &= \frac{1}{48}\left( -204 -756\zeta(2) -4344\zeta(3) -14427\zeta(4) -20382\zeta(5)  \right. \nonumber \\ &\left. \hspace{1em}
-4644\zeta(2)\zeta(3) -18570\zeta(6) -2940\zeta(3)^2 +24423\zeta(7) +5844\zeta(2)\zeta(5)  \right. \nonumber \\ &\left. \hspace{1em}
+8316\zeta(3)\zeta(4) +17027\zeta(8) -924\zeta(2)\zeta(3)^2 +11328\zeta(3)\zeta(5) +1308 M(2,6)\right) \label{eq10196}
\end{align}
 
\begin{align}
\sum_{k=1}^\infty \frac{H(k)^{6}}{(k+1)^{2}(k+2)^{2}} \sumend &= \frac{1}{24}\left( -216 -816\zeta(2) -4752\zeta(3) -16140\zeta(4) -23808\zeta(5)  \right. \nonumber \\ &\left. \hspace{1em}
-5376\zeta(2)\zeta(3) -24861\zeta(6) -3744\zeta(3)^2 +16695\zeta(7) +4104\zeta(2)\zeta(5)  \right. \nonumber \\ &\left. \hspace{1em}
+4752\zeta(3)\zeta(4) +34054\zeta(8) -1848\zeta(2)\zeta(3)^2 +22656\zeta(3)\zeta(5)  \right. \nonumber \\ &\left. \hspace{1em}
+2616 M(2,6)\right) \label{eq10197} \\
\sum_{k=1}^\infty \frac{H(k)^{6}}{k(k+2)^{3}} \sumend &= \frac{1}{96}\left(  1524 +3948\zeta(2) +18360\zeta(3) +44121\zeta(4) +34122\zeta(5)  \right. \nonumber \\ &\left. \hspace{1em}
+8364\zeta(2)\zeta(3) -18408\zeta(6) -1692\zeta(3)^2 -32547\zeta(7) -6972\zeta(2)\zeta(5)  \right. \nonumber \\ &\left. \hspace{1em}
-24012\zeta(3)\zeta(4) -82648\zeta(8) -16716\zeta(2)\zeta(3)^2 +61104\zeta(3)\zeta(5) +14172 M(2,6)  \right. \nonumber \\ &\left. \hspace{1em}
+12292\zeta(9) -25164\zeta(3)\zeta(6) +11664\zeta(4)\zeta(5) +3906\zeta(2)\zeta(7)  \right. \nonumber \\ &\left. \hspace{1em}
-1072\zeta(3)^3\right) \label{eq10198} \\
\sum_{k=1}^\infty \frac{H(k)^{6}}{(k+1)(k+2)^{3}} \sumend &= \frac{-1}{48}\left( -1728 -4704\zeta(2) -22704\zeta(3) -58548\zeta(4) -54504\zeta(5)  \right. \nonumber \\ &\left. \hspace{1em}
-13008\zeta(2)\zeta(3) -162\zeta(6) -1248\zeta(3)^2 +56970\zeta(7) +12816\zeta(2)\zeta(5)  \right. \nonumber \\ &\left. \hspace{1em}
+32328\zeta(3)\zeta(4) +99675\zeta(8) +15792\zeta(2)\zeta(3)^2 -49776\zeta(3)\zeta(5) -12864 M(2,6)  \right. \nonumber \\ &\left. \hspace{1em}
-12292\zeta(9) +25164\zeta(3)\zeta(6) -11664\zeta(4)\zeta(5) -3906\zeta(2)\zeta(7)  \right. \nonumber \\ &\left. \hspace{1em}
+1072\zeta(3)^3\right) \label{eq10199} \\
\sum_{k=1}^\infty \frac{H(k)^{6}}{(k+2)^{4}} \sumend &= \frac{-1}{384}\left(  32256 +59136\zeta(2) +224256\zeta(3) +386496\zeta(4) +122880\zeta(5)  \right. \nonumber \\ &\left. \hspace{1em}
+26880\zeta(2)\zeta(3) -378528\zeta(6) -66432\zeta(3)^2 -167280\zeta(7) -23424\zeta(2)\zeta(5)  \right. \nonumber \\ &\left. \hspace{1em}
-225792\zeta(3)\zeta(4) -28368\zeta(8) -8640\zeta(2)\zeta(3)^2 +26496\zeta(3)\zeta(5) +2880 M(2,6)  \right. \nonumber \\ &\left. \hspace{1em}
-227840\zeta(9) +412320\zeta(3)\zeta(6) -152064\zeta(4)\zeta(5) -76320\zeta(2)\zeta(7) +11520\zeta(3)^3  \right. \nonumber \\ &\left. \hspace{1em}
+779835\zeta(10) -490704\zeta(3)\zeta(7) +245544\zeta(3)^2\zeta(4) -15600\zeta(2)\zeta(3)\zeta(5)  \right. \nonumber \\ &\left. \hspace{1em}
-339864\zeta(5)^2 -94728 M(2,8) -45120\zeta(2) M(2,6)\right) \label{eq10200} \\
\sum_{k=1}^\infty \frac{H(k)^{7}}{k^{3}} \sumend &= \frac{1}{2560}\left( -16614991\zeta(10) +10315520\zeta(3)\zeta(7) -5879160\zeta(3)^2\zeta(4)  \right. \nonumber \\ &\left. \hspace{1em}
+705040\zeta(2)\zeta(3)\zeta(5) +7710760\zeta(5)^2 +2021880 M(2,8) +1008000\zeta(2) M(2,6)\right) \label{eq10201} \\
\sum_{k=1}^\infty \frac{H(k)^{7}}{k^{2}(k+1)} \sumend &= \frac{1}{72}\left( -479096\zeta(8) -12096\zeta(2)\zeta(3)^2 -109620\zeta(3)\zeta(5)  \right. \nonumber \\ &\left. \hspace{1em}
+276341\zeta(9) +88665\zeta(3)\zeta(6) +143163\zeta(4)\zeta(5) +59166\zeta(2)\zeta(7)  \right. \nonumber \\ &\left. \hspace{1em}
+4032\zeta(3)^3\right) \label{eq10202} \\
\sum_{k=1}^\infty \frac{H(k)^{7}}{k(k+1)^{2}} \sumend &= \frac{1}{72}\left(  479096\zeta(8) +12096\zeta(2)\zeta(3)^2 +109620\zeta(3)\zeta(5)  \right. \nonumber \\ &\left. \hspace{1em}
-269402\zeta(9) -88665\zeta(3)\zeta(6) -141273\zeta(4)\zeta(5) -59166\zeta(2)\zeta(7)  \right. \nonumber \\ &\left. \hspace{1em}
-4032\zeta(3)^3\right) \label{eq10203}
\end{align}
 
\begin{align}
\sum_{k=1}^\infty \frac{H(k)^{7}}{(k+1)^{3}} \sumend &= \frac{-1}{2560}\left( -16597239\zeta(10) +9974400\zeta(3)\zeta(7) -5800760\zeta(3)^2\zeta(4)  \right. \nonumber \\ &\left. \hspace{1em}
+834960\zeta(2)\zeta(3)\zeta(5) +7473640\zeta(5)^2 +1956920 M(2,8) +1008000\zeta(2) M(2,6)\right) \label{eq10204} \\
\sum_{k=1}^\infty \frac{H(k)^{7}}{k^{2}(k+2)} \sumend &= \frac{-1}{576}\left(  144 +864\zeta(2) +7200\zeta(3) +38664\zeta(4) +108504\zeta(5)  \right. \nonumber \\ &\left. \hspace{1em}
+23184\zeta(2)\zeta(3) +352887\zeta(6) +45864\zeta(3)^2 +319554\zeta(7) +73080\zeta(2)\zeta(5)  \right. \nonumber \\ &\left. \hspace{1em}
+148932\zeta(3)\zeta(4) +958192\zeta(8) +24192\zeta(2)\zeta(3)^2 +219240\zeta(3)\zeta(5)  \right. \nonumber \\ &\left. \hspace{1em}
-1105364\zeta(9) -354660\zeta(3)\zeta(6) -572652\zeta(4)\zeta(5) -236664\zeta(2)\zeta(7)  \right. \nonumber \\ &\left. \hspace{1em}
-16128\zeta(3)^3\right) \label{eq10205} \\
\sum_{k=1}^\infty \frac{H(k)^{7}}{k(k+1)(k+2)} \sumend &= \frac{-1}{288}\left(  144 +864\zeta(2) +7200\zeta(3) +38664\zeta(4) +108504\zeta(5)  \right. \nonumber \\ &\left. \hspace{1em}
+23184\zeta(2)\zeta(3) +352887\zeta(6) +45864\zeta(3)^2 +319554\zeta(7) +73080\zeta(2)\zeta(5)  \right. \nonumber \\ &\left. \hspace{1em}
+148932\zeta(3)\zeta(4) -958192\zeta(8) -24192\zeta(2)\zeta(3)^2 -219240\zeta(3)\zeta(5)\right) \label{eq10206} \\
\sum_{k=1}^\infty \frac{H(k)^{7}}{(k+1)^{2}(k+2)} \sumend &= \frac{1}{144}\left( -144 -864\zeta(2) -7200\zeta(3) -38664\zeta(4) -108504\zeta(5)  \right. \nonumber \\ &\left. \hspace{1em}
-23184\zeta(2)\zeta(3) -352887\zeta(6) -45864\zeta(3)^2 -319554\zeta(7) -73080\zeta(2)\zeta(5)  \right. \nonumber \\ &\left. \hspace{1em}
-148932\zeta(3)\zeta(4) +538804\zeta(9) +177330\zeta(3)\zeta(6) +282546\zeta(4)\zeta(5)  \right. \nonumber \\ &\left. \hspace{1em}
+118332\zeta(2)\zeta(7) +8064\zeta(3)^3\right) \label{eq10207} \\
\sum_{k=1}^\infty \frac{H(k)^{7}}{k(k+2)^{2}} \sumend &= \frac{1}{576}\left(  2448 +10944\zeta(2) +77184\zeta(3) +331344\zeta(4) +683928\zeta(5)  \right. \nonumber \\ &\left. \hspace{1em}
+152208\zeta(2)\zeta(3) +1403655\zeta(6) +199080\zeta(3)^2 +257472\zeta(7) +46872\zeta(2)\zeta(5)  \right. \nonumber \\ &\left. \hspace{1em}
+247212\zeta(3)\zeta(4) -472076\zeta(8) +101808\zeta(2)\zeta(3)^2 -732312\zeta(3)\zeta(5)  \right. \nonumber \\ &\left. \hspace{1em}
-109872 M(2,6) -1077608\zeta(9) -354660\zeta(3)\zeta(6) -565092\zeta(4)\zeta(5) -236664\zeta(2)\zeta(7)  \right. \nonumber \\ &\left. \hspace{1em}
-16128\zeta(3)^3\right) \label{eq10208} \\
\sum_{k=1}^\infty \frac{H(k)^{7}}{(k+1)(k+2)^{2}} \sumend &= \frac{-1}{144}\left( -1296 -5904\zeta(2) -42192\zeta(3) -185004\zeta(4) -396216\zeta(5)  \right. \nonumber \\ &\left. \hspace{1em}
-87696\zeta(2)\zeta(3) -878271\zeta(6) -122472\zeta(3)^2 -288513\zeta(7) -59976\zeta(2)\zeta(5)  \right. \nonumber \\ &\left. \hspace{1em}
-198072\zeta(3)\zeta(4) +715134\zeta(8) -38808\zeta(2)\zeta(3)^2 +475776\zeta(3)\zeta(5) +54936 M(2,6)  \right. \nonumber \\ &\left. \hspace{1em}
+538804\zeta(9) +177330\zeta(3)\zeta(6) +282546\zeta(4)\zeta(5) +118332\zeta(2)\zeta(7)  \right. \nonumber \\ &\left. \hspace{1em}
+8064\zeta(3)^3\right) \label{eq10209} \\
\sum_{k=1}^\infty \frac{H(k)^{7}}{(k+2)^{3}} \sumend &= \frac{-1}{7680}\left(  276480 +913920\zeta(2) +5429760\zeta(3) +18389760\zeta(4)  \right. \nonumber \\ &\left. \hspace{1em}
+26899200\zeta(5) +6182400\zeta(2)\zeta(3) +28153920\zeta(6) +4381440\zeta(3)^2 -16691520\zeta(7)  \right. \nonumber \\ &\left. \hspace{1em}
-3951360\zeta(2)\zeta(5) -7674240\zeta(3)\zeta(4) -74888240\zeta(8) -7808640\zeta(2)\zeta(3)^2  \right. \nonumber \\ &\left. \hspace{1em}
+15187200\zeta(3)\zeta(5) +5738880 M(2,6) +13767040\zeta(9) -28183680\zeta(3)\zeta(6)  \right. \nonumber \\ &\left. \hspace{1em}
+13063680\zeta(4)\zeta(5) +4374720\zeta(2)\zeta(7) -1200640\zeta(3)^3 -49791717\zeta(10)  \right. \nonumber \\ &\left. \hspace{1em}
+29923200\zeta(3)\zeta(7) -17402280\zeta(3)^2\zeta(4) +2504880\zeta(2)\zeta(3)\zeta(5)  \right. \nonumber \\ &\left. \hspace{1em}
+22420920\zeta(5)^2 +5870760 M(2,8) +3024000\zeta(2) M(2,6)\right) \label{eq10210}
\end{align}
 
\begin{align}
\sum_{k=1}^\infty \frac{H(k)^{8}}{k^{2}} \sumend &= \frac{-1}{480}\left( -18741581\zeta(10) -6689520\zeta(3)\zeta(7) +524640\zeta(3)^2\zeta(4)  \right. \nonumber \\ &\left. \hspace{1em}
-1452480\zeta(2)\zeta(3)\zeta(5) -4247040\zeta(5)^2 -485280 M(2,8) -299520\zeta(2) M(2,6)\right) \label{eq10211} \\
\sum_{k=1}^\infty \frac{H(k)^{8}}{k(k+1)} \sumend &= \frac{-1}{6}\left( -166700\zeta(9) -88665\zeta(3)\zeta(6) -80400\zeta(4)\zeta(5)  \right. \nonumber \\ &\left. \hspace{1em}
-35091\zeta(2)\zeta(7) -4032\zeta(3)^3\right) \label{eq10212} \\
\sum_{k=1}^\infty \frac{H(k)^{8}}{(k+1)^{2}} \sumend &= \frac{-1}{240}\left( -9295879\zeta(10) -3314520\zeta(3)\zeta(7) +258540\zeta(3)^2\zeta(4)  \right. \nonumber \\ &\left. \hspace{1em}
-733800\zeta(2)\zeta(3)\zeta(5) -2098980\zeta(5)^2 -238860 M(2,8) -149760\zeta(2) M(2,6)\right) \label{eq10213} \\
\sum_{k=1}^\infty \frac{H(k)^{8}}{k(k+2)} \sumend &= \frac{1}{144}\left(  72 +504\zeta(2) +4968\zeta(3) +32400\zeta(4) +116280\zeta(5)  \right. \nonumber \\ &\left. \hspace{1em}
+24768\zeta(2)\zeta(3) +530346\zeta(6) +69264\zeta(3)^2 +849654\zeta(7) +193104\zeta(2)\zeta(5)  \right. \nonumber \\ &\left. \hspace{1em}
+404280\zeta(3)\zeta(4) +1906367\zeta(8) +48384\zeta(2)\zeta(3)^2 +436896\zeta(3)\zeta(5)  \right. \nonumber \\ &\left. \hspace{1em}
+2000400\zeta(9) +1063980\zeta(3)\zeta(6) +964800\zeta(4)\zeta(5) +421092\zeta(2)\zeta(7)  \right. \nonumber \\ &\left. \hspace{1em}
+48384\zeta(3)^3\right) \label{eq10214} \\
\sum_{k=1}^\infty \frac{H(k)^{8}}{(k+1)(k+2)} \sumend &= \frac{1}{72}\left(  72 +504\zeta(2) +4968\zeta(3) +32400\zeta(4) +116280\zeta(5)  \right. \nonumber \\ &\left. \hspace{1em}
+24768\zeta(2)\zeta(3) +530346\zeta(6) +69264\zeta(3)^2 +849654\zeta(7) +193104\zeta(2)\zeta(5)  \right. \nonumber \\ &\left. \hspace{1em}
+404280\zeta(3)\zeta(4) +1906367\zeta(8) +48384\zeta(2)\zeta(3)^2 +436896\zeta(3)\zeta(5)\right) \label{eq10215} \\
\sum_{k=1}^\infty \frac{H(k)^{8}}{(k+2)^{2}} \sumend &= \frac{1}{720}\left( -6480 -34560\zeta(2) -292320\zeta(3) -1564560\zeta(4) -4348800\zeta(5)  \right. \nonumber \\ &\left. \hspace{1em}
-950400\zeta(2)\zeta(3) -14106480\zeta(6) -1926720\zeta(3)^2 -12318480\zeta(7) -2712960\zeta(2)\zeta(5)  \right. \nonumber \\ &\left. \hspace{1em}
-6755040\zeta(3)\zeta(4) -4760990\zeta(8) -1260000\zeta(2)\zeta(3)^2 +5146560\zeta(3)\zeta(5)  \right. \nonumber \\ &\left. \hspace{1em}
+1098720 M(2,6) +21552160\zeta(9) +7093200\zeta(3)\zeta(6) +11301840\zeta(4)\zeta(5)  \right. \nonumber \\ &\left. \hspace{1em}
+4733280\zeta(2)\zeta(7) +322560\zeta(3)^3 +27887637\zeta(10) +9943560\zeta(3)\zeta(7)  \right. \nonumber \\ &\left. \hspace{1em}
-775620\zeta(3)^2\zeta(4) +2201400\zeta(2)\zeta(3)\zeta(5) +6296940\zeta(5)^2 +716580 M(2,8)  \right. \nonumber \\ &\left. \hspace{1em}
+449280\zeta(2) M(2,6)\right) \label{eq10216}
\end{align}

\newpage
Formulas for order $r = m + n + p + q = 11$:
\begin{align}
\sum_{k=1}^\infty \frac{H(k)}{k^{10}} \sumend &= \left(  6\zeta(11) -\zeta(2)\zeta(9) -\zeta(3)\zeta(8) -\zeta(4)\zeta(7) -\zeta(5)\zeta(6)\right) \label{eq11001} \\
\sum_{k=1}^\infty \frac{H(k)}{k^{9}(k+1)} \sumend &= \frac{1}{4}\left(  4\zeta(2) -8\zeta(3) +5\zeta(4) -12\zeta(5) +4\zeta(2)\zeta(3) +7\zeta(6)  \right. \nonumber \\ &\left. \hspace{1em}
-2\zeta(3)^2 -16\zeta(7) +4\zeta(2)\zeta(5) +4\zeta(3)\zeta(4) +9\zeta(8) -4\zeta(3)\zeta(5)  \right. \nonumber \\ &\left. \hspace{1em}
-20\zeta(9) +4\zeta(3)\zeta(6) +4\zeta(4)\zeta(5) +4\zeta(2)\zeta(7) +11\zeta(10) -4\zeta(3)\zeta(7)  \right. \nonumber \\ &\left. \hspace{1em}
-2\zeta(5)^2\right) \label{eq11002} \\
\sum_{k=1}^\infty \frac{H(k)}{k^{8}(k+1)^{2}} \sumend &= \frac{1}{2}\left( -16\zeta(2) +30\zeta(3) -15\zeta(4) +30\zeta(5) -10\zeta(2)\zeta(3)  \right. \nonumber \\ &\left. \hspace{1em}
-14\zeta(6) +4\zeta(3)^2 +24\zeta(7) -6\zeta(2)\zeta(5) -6\zeta(3)\zeta(4) -9\zeta(8)  \right. \nonumber \\ &\left. \hspace{1em}
+4\zeta(3)\zeta(5) +10\zeta(9) -2\zeta(3)\zeta(6) -2\zeta(4)\zeta(5) -2\zeta(2)\zeta(7)\right) \label{eq11003} \\
\sum_{k=1}^\infty \frac{H(k)}{k^{7}(k+1)^{3}} \sumend &= \frac{-1}{4}\left( -112\zeta(2) +196\zeta(3) -74\zeta(4) +120\zeta(5) -40\zeta(2)\zeta(3)  \right. \nonumber \\ &\left. \hspace{1em}
-42\zeta(6) +12\zeta(3)^2 +48\zeta(7) -12\zeta(2)\zeta(5) -12\zeta(3)\zeta(4) -9\zeta(8)  \right. \nonumber \\ &\left. \hspace{1em}
+4\zeta(3)\zeta(5)\right) \label{eq11004} \\
\sum_{k=1}^\infty \frac{H(k)}{k^{6}(k+1)^{4}} \sumend &= \frac{1}{2}\left( -112\zeta(2) +182\zeta(3) -47\zeta(4) +64\zeta(5) -22\zeta(2)\zeta(3)  \right. \nonumber \\ &\left. \hspace{1em}
-14\zeta(6) +4\zeta(3)^2 +8\zeta(7) -2\zeta(2)\zeta(5) -2\zeta(3)\zeta(4)\right) \label{eq11005} \\
\sum_{k=1}^\infty \frac{H(k)}{k^{5}(k+1)^{5}} \sumend &= -\left( -70\zeta(2) +105\zeta(3) -15\zeta(4) +25\zeta(5) -10\zeta(2)\zeta(3) -\zeta(6)\right) \label{eq11006} \\
\sum_{k=1}^\infty \frac{H(k)}{k^{4}(k+1)^{6}} \sumend &= \frac{1}{2}\left( -112\zeta(2) +154\zeta(3) -5\zeta(4) +46\zeta(5) -22\zeta(2)\zeta(3)  \right. \nonumber \\ &\left. \hspace{1em}
+6\zeta(6) -4\zeta(3)^2 +6\zeta(7) -2\zeta(2)\zeta(5) -2\zeta(3)\zeta(4)\right) \label{eq11007} \\
\sum_{k=1}^\infty \frac{H(k)}{k^{3}(k+1)^{7}} \sumend &= \frac{-1}{4}\left( -112\zeta(2) +140\zeta(3) +10\zeta(4) +80\zeta(5) -40\zeta(2)\zeta(3)  \right. \nonumber \\ &\left. \hspace{1em}
+18\zeta(6) -12\zeta(3)^2 +36\zeta(7) -12\zeta(2)\zeta(5) -12\zeta(3)\zeta(4) +5\zeta(8)  \right. \nonumber \\ &\left. \hspace{1em}
-4\zeta(3)\zeta(5)\right) \label{eq11008} \\
\sum_{k=1}^\infty \frac{H(k)}{k^{2}(k+1)^{8}} \sumend &= \frac{1}{2}\left( -16\zeta(2) +18\zeta(3) +3\zeta(4) +20\zeta(5) -10\zeta(2)\zeta(3) +6\zeta(6)  \right. \nonumber \\ &\left. \hspace{1em}
-4\zeta(3)^2 +18\zeta(7) -6\zeta(2)\zeta(5) -6\zeta(3)\zeta(4) +5\zeta(8) -4\zeta(3)\zeta(5) +8\zeta(9)  \right. \nonumber \\ &\left. \hspace{1em}
-2\zeta(3)\zeta(6) -2\zeta(4)\zeta(5) -2\zeta(2)\zeta(7)\right) \label{eq11009} \\
\sum_{k=1}^\infty \frac{H(k)}{k(k+1)^{9}} \sumend &= \frac{1}{4}\left(  4\zeta(2) -4\zeta(3) -\zeta(4) -8\zeta(5) +4\zeta(2)\zeta(3) -3\zeta(6)  \right. \nonumber \\ &\left. \hspace{1em}
+2\zeta(3)^2 -12\zeta(7) +4\zeta(2)\zeta(5) +4\zeta(3)\zeta(4) -5\zeta(8) +4\zeta(3)\zeta(5)  \right. \nonumber \\ &\left. \hspace{1em}
-16\zeta(9) +4\zeta(3)\zeta(6) +4\zeta(4)\zeta(5) +4\zeta(2)\zeta(7) -7\zeta(10) +4\zeta(3)\zeta(7)  \right. \nonumber \\ &\left. \hspace{1em}
+2\zeta(5)^2\right) \label{eq11010}
\end{align}
 
\begin{align}
\sum_{k=1}^\infty \frac{H(k)}{(k+1)^{10}} \sumend &= -\left( -5\zeta(11) +\zeta(2)\zeta(9) +\zeta(3)\zeta(8) +\zeta(4)\zeta(7) +\zeta(5)\zeta(6)\right) \label{eq11011} \\
\sum_{k=1}^\infty \frac{H(k)}{k^{9}(k+2)} \sumend &= \frac{1}{512}\left(  1 +\zeta(2) -4\zeta(3) +5\zeta(4) -24\zeta(5) +8\zeta(2)\zeta(3) +28\zeta(6)  \right. \nonumber \\ &\left. \hspace{1em}
-8\zeta(3)^2 -128\zeta(7) +32\zeta(2)\zeta(5) +32\zeta(3)\zeta(4) +144\zeta(8) -64\zeta(3)\zeta(5)  \right. \nonumber \\ &\left. \hspace{1em}
-640\zeta(9) +128\zeta(3)\zeta(6) +128\zeta(4)\zeta(5) +128\zeta(2)\zeta(7) +704\zeta(10)  \right. \nonumber \\ &\left. \hspace{1em}
-256\zeta(3)\zeta(7) -128\zeta(5)^2\right) \label{eq11012} \\
\sum_{k=1}^\infty \frac{H(k)}{k^{8}(k+1)(k+2)} \sumend &= \frac{-1}{256}\left( -1 +255\zeta(2) -508\zeta(3) +315\zeta(4) -744\zeta(5)  \right. \nonumber \\ &\left. \hspace{1em}
+248\zeta(2)\zeta(3) +420\zeta(6) -120\zeta(3)^2 -896\zeta(7) +224\zeta(2)\zeta(5) +224\zeta(3)\zeta(4)  \right. \nonumber \\ &\left. \hspace{1em}
+432\zeta(8) -192\zeta(3)\zeta(5) -640\zeta(9) +128\zeta(3)\zeta(6) +128\zeta(4)\zeta(5)  \right. \nonumber \\ &\left. \hspace{1em}
+128\zeta(2)\zeta(7)\right) \label{eq11013} \\
\sum_{k=1}^\infty \frac{H(k)}{k^{7}(k+1)^{2}(k+2)} \sumend &= \frac{-1}{128}\left( -1 -769\zeta(2) +1412\zeta(3) -645\zeta(4) +1176\zeta(5)  \right. \nonumber \\ &\left. \hspace{1em}
-392\zeta(2)\zeta(3) -476\zeta(6) +136\zeta(3)^2 +640\zeta(7) -160\zeta(2)\zeta(5) -160\zeta(3)\zeta(4)  \right. \nonumber \\ &\left. \hspace{1em}
-144\zeta(8) +64\zeta(3)\zeta(5)\right) \label{eq11014} \\
\sum_{k=1}^\infty \frac{H(k)}{k^{6}(k+1)^{3}(k+2)} \sumend &= \frac{1}{64}\left(  1 -1023\zeta(2) +1724\zeta(3) -539\zeta(4) +744\zeta(5)  \right. \nonumber \\ &\left. \hspace{1em}
-248\zeta(2)\zeta(3) -196\zeta(6) +56\zeta(3)^2 +128\zeta(7) -32\zeta(2)\zeta(5)  \right. \nonumber \\ &\left. \hspace{1em}
-32\zeta(3)\zeta(4)\right) \label{eq11015} \\
\sum_{k=1}^\infty \frac{H(k)}{k^{5}(k+1)^{4}(k+2)} \sumend &= \frac{1}{32}\left(  1 +769\zeta(2) -1188\zeta(3) +213\zeta(4) -280\zeta(5)  \right. \nonumber \\ &\left. \hspace{1em}
+104\zeta(2)\zeta(3) +28\zeta(6) -8\zeta(3)^2\right) \label{eq11016} \\
\sum_{k=1}^\infty \frac{H(k)}{k^{4}(k+1)^{5}(k+2)} \sumend &= \frac{-1}{16}\left( -1 +351\zeta(2) -492\zeta(3) +27\zeta(4) -120\zeta(5)  \right. \nonumber \\ &\left. \hspace{1em}
+56\zeta(2)\zeta(3) -12\zeta(6) +8\zeta(3)^2\right) \label{eq11017} \\
\sum_{k=1}^\infty \frac{H(k)}{k^{3}(k+1)^{6}(k+2)} \sumend &= \frac{1}{8}\left(  1 +97\zeta(2) -124\zeta(3) -7\zeta(4) -64\zeta(5) +32\zeta(2)\zeta(3)  \right. \nonumber \\ &\left. \hspace{1em}
-12\zeta(6) +8\zeta(3)^2 -24\zeta(7) +8\zeta(2)\zeta(5) +8\zeta(3)\zeta(4)\right) \label{eq11018} \\
\sum_{k=1}^\infty \frac{H(k)}{k^{2}(k+1)^{7}(k+2)} \sumend &= \frac{1}{4}\left(  1 -15\zeta(2) +16\zeta(3) +3\zeta(4) +16\zeta(5) -8\zeta(2)\zeta(3)  \right. \nonumber \\ &\left. \hspace{1em}
+6\zeta(6) -4\zeta(3)^2 +12\zeta(7) -4\zeta(2)\zeta(5) -4\zeta(3)\zeta(4) +5\zeta(8)  \right. \nonumber \\ &\left. \hspace{1em}
-4\zeta(3)\zeta(5)\right) \label{eq11019}
\end{align}
 
\begin{align}
\sum_{k=1}^\infty \frac{H(k)}{k(k+1)^{8}(k+2)} \sumend &= \frac{-1}{2}\left( -1 -\zeta(2) +2\zeta(3) +4\zeta(5) -2\zeta(2)\zeta(3) +6\zeta(7)  \right. \nonumber \\ &\left. \hspace{1em}
-2\zeta(2)\zeta(5) -2\zeta(3)\zeta(4) +8\zeta(9) -2\zeta(3)\zeta(6) -2\zeta(4)\zeta(5)  \right. \nonumber \\ &\left. \hspace{1em}
-2\zeta(2)\zeta(7)\right) \label{eq11020} \\
\sum_{k=1}^\infty \frac{H(k)}{(k+1)^{9}(k+2)} \sumend &= \frac{1}{4}\left(  4 -4\zeta(3) +\zeta(4) -8\zeta(5) +4\zeta(2)\zeta(3) +3\zeta(6) -2\zeta(3)^2  \right. \nonumber \\ &\left. \hspace{1em}
-12\zeta(7) +4\zeta(2)\zeta(5) +4\zeta(3)\zeta(4) +5\zeta(8) -4\zeta(3)\zeta(5) -16\zeta(9)  \right. \nonumber \\ &\left. \hspace{1em}
+4\zeta(3)\zeta(6) +4\zeta(4)\zeta(5) +4\zeta(2)\zeta(7) +7\zeta(10) -4\zeta(3)\zeta(7)  \right. \nonumber \\ &\left. \hspace{1em}
-2\zeta(5)^2\right) \label{eq11021} \\
\sum_{k=1}^\infty \frac{H(k)}{k^{8}(k+2)^{2}} \sumend &= \frac{-1}{256}\left(  6 +3\zeta(2) -15\zeta(3) +15\zeta(4) -60\zeta(5) +20\zeta(2)\zeta(3)  \right. \nonumber \\ &\left. \hspace{1em}
+56\zeta(6) -16\zeta(3)^2 -192\zeta(7) +48\zeta(2)\zeta(5) +48\zeta(3)\zeta(4) +144\zeta(8)  \right. \nonumber \\ &\left. \hspace{1em}
-64\zeta(3)\zeta(5) -320\zeta(9) +64\zeta(3)\zeta(6) +64\zeta(4)\zeta(5) +64\zeta(2)\zeta(7)\right) \label{eq11022} \\
\sum_{k=1}^\infty \frac{H(k)}{k^{7}(k+1)(k+2)^{2}} \sumend &= \frac{1}{256}\left( -13 +249\zeta(2) -478\zeta(3) +285\zeta(4) -624\zeta(5)  \right. \nonumber \\ &\left. \hspace{1em}
+208\zeta(2)\zeta(3) +308\zeta(6) -88\zeta(3)^2 -512\zeta(7) +128\zeta(2)\zeta(5) +128\zeta(3)\zeta(4)  \right. \nonumber \\ &\left. \hspace{1em}
+144\zeta(8) -64\zeta(3)\zeta(5)\right) \label{eq11023} \\
\sum_{k=1}^\infty \frac{H(k)}{k^{6}(k+1)^{2}(k+2)^{2}} \sumend &= \frac{-1}{64}\left(  7 +260\zeta(2) -467\zeta(3) +180\zeta(4) -276\zeta(5)  \right. \nonumber \\ &\left. \hspace{1em}
+92\zeta(2)\zeta(3) +84\zeta(6) -24\zeta(3)^2 -64\zeta(7) +16\zeta(2)\zeta(5)  \right. \nonumber \\ &\left. \hspace{1em}
+16\zeta(3)\zeta(4)\right) \label{eq11024} \\
\sum_{k=1}^\infty \frac{H(k)}{k^{5}(k+1)^{3}(k+2)^{2}} \sumend &= \frac{-1}{64}\left(  15 -503\zeta(2) +790\zeta(3) -179\zeta(4) +192\zeta(5)  \right. \nonumber \\ &\left. \hspace{1em}
-64\zeta(2)\zeta(3) -28\zeta(6) +8\zeta(3)^2\right) \label{eq11025} \\
\sum_{k=1}^\infty \frac{H(k)}{k^{4}(k+1)^{4}(k+2)^{2}} \sumend &= \frac{-1}{16}\left(  8 +133\zeta(2) -199\zeta(3) +17\zeta(4) -44\zeta(5)  \right. \nonumber \\ &\left. \hspace{1em}
+20\zeta(2)\zeta(3)\right) \label{eq11026} \\
\sum_{k=1}^\infty \frac{H(k)}{k^{3}(k+1)^{5}(k+2)^{2}} \sumend &= \frac{-1}{16}\left(  17 -85\zeta(2) +94\zeta(3) +7\zeta(4) +32\zeta(5)  \right. \nonumber \\ &\left. \hspace{1em}
-16\zeta(2)\zeta(3) +12\zeta(6) -8\zeta(3)^2\right) \label{eq11027} \\
\sum_{k=1}^\infty \frac{H(k)}{k^{2}(k+1)^{6}(k+2)^{2}} \sumend &= \frac{-1}{4}\left(  9 +6\zeta(2) -15\zeta(3) -16\zeta(5) +8\zeta(2)\zeta(3) -12\zeta(7)  \right. \nonumber \\ &\left. \hspace{1em}
+4\zeta(2)\zeta(5) +4\zeta(3)\zeta(4)\right) \label{eq11028}
\end{align}
 
\begin{align}
\sum_{k=1}^\infty \frac{H(k)}{k(k+1)^{7}(k+2)^{2}} \sumend &= \frac{-1}{4}\left(  19 -3\zeta(2) -14\zeta(3) +3\zeta(4) -16\zeta(5) +8\zeta(2)\zeta(3)  \right. \nonumber \\ &\left. \hspace{1em}
+6\zeta(6) -4\zeta(3)^2 -12\zeta(7) +4\zeta(2)\zeta(5) +4\zeta(3)\zeta(4) +5\zeta(8)  \right. \nonumber \\ &\left. \hspace{1em}
-4\zeta(3)\zeta(5)\right) \label{eq11029} \\
\sum_{k=1}^\infty \frac{H(k)}{(k+1)^{8}(k+2)^{2}} \sumend &= \frac{1}{2}\left( -20 +2\zeta(2) +16\zeta(3) -3\zeta(4) +20\zeta(5) -10\zeta(2)\zeta(3)  \right. \nonumber \\ &\left. \hspace{1em}
-6\zeta(6) +4\zeta(3)^2 +18\zeta(7) -6\zeta(2)\zeta(5) -6\zeta(3)\zeta(4) -5\zeta(8) +4\zeta(3)\zeta(5)  \right. \nonumber \\ &\left. \hspace{1em}
+8\zeta(9) -2\zeta(3)\zeta(6) -2\zeta(4)\zeta(5) -2\zeta(2)\zeta(7)\right) \label{eq11030} \\
\sum_{k=1}^\infty \frac{H(k)}{k^{7}(k+2)^{3}} \sumend &= \frac{1}{256}\left(  34 +5\zeta(2) -51\zeta(3) +37\zeta(4) -120\zeta(5) +40\zeta(2)\zeta(3)  \right. \nonumber \\ &\left. \hspace{1em}
+84\zeta(6) -24\zeta(3)^2 -192\zeta(7) +48\zeta(2)\zeta(5) +48\zeta(3)\zeta(4) +72\zeta(8)  \right. \nonumber \\ &\left. \hspace{1em}
-32\zeta(3)\zeta(5)\right) \label{eq11031} \\
\sum_{k=1}^\infty \frac{H(k)}{k^{6}(k+1)(k+2)^{3}} \sumend &= \frac{1}{256}\left(  81 -239\zeta(2) +376\zeta(3) -211\zeta(4) +384\zeta(5)  \right. \nonumber \\ &\left. \hspace{1em}
-128\zeta(2)\zeta(3) -140\zeta(6) +40\zeta(3)^2 +128\zeta(7) -32\zeta(2)\zeta(5)  \right. \nonumber \\ &\left. \hspace{1em}
-32\zeta(3)\zeta(4)\right) \label{eq11032} \\
\sum_{k=1}^\infty \frac{H(k)}{k^{5}(k+1)^{2}(k+2)^{3}} \sumend &= \frac{-1}{128}\left( -95 -281\zeta(2) +558\zeta(3) -149\zeta(4) +168\zeta(5)  \right. \nonumber \\ &\left. \hspace{1em}
-56\zeta(2)\zeta(3) -28\zeta(6) +8\zeta(3)^2\right) \label{eq11033} \\
\sum_{k=1}^\infty \frac{H(k)}{k^{4}(k+1)^{3}(k+2)^{3}} \sumend &= \frac{1}{32}\left(  55 -111\zeta(2) +116\zeta(3) -15\zeta(4) +12\zeta(5)  \right. \nonumber \\ &\left. \hspace{1em}
-4\zeta(2)\zeta(3)\right) \label{eq11034} \\
\sum_{k=1}^\infty \frac{H(k)}{k^{3}(k+1)^{4}(k+2)^{3}} \sumend &= \frac{-1}{16}\left( -63 -22\zeta(2) +83\zeta(3) -2\zeta(4) +32\zeta(5)  \right. \nonumber \\ &\left. \hspace{1em}
-16\zeta(2)\zeta(3)\right) \label{eq11035} \\
\sum_{k=1}^\infty \frac{H(k)}{k^{2}(k+1)^{5}(k+2)^{3}} \sumend &= \frac{1}{16}\left(  143 -41\zeta(2) -72\zeta(3) +11\zeta(4) -32\zeta(5)  \right. \nonumber \\ &\left. \hspace{1em}
+16\zeta(2)\zeta(3) +12\zeta(6) -8\zeta(3)^2\right) \label{eq11036} \\
\sum_{k=1}^\infty \frac{H(k)}{k(k+1)^{6}(k+2)^{3}} \sumend &= \frac{1}{8}\left(  161 -29\zeta(2) -102\zeta(3) +11\zeta(4) -64\zeta(5) +32\zeta(2)\zeta(3)  \right. \nonumber \\ &\left. \hspace{1em}
+12\zeta(6) -8\zeta(3)^2 -24\zeta(7) +8\zeta(2)\zeta(5) +8\zeta(3)\zeta(4)\right) \label{eq11037} \\
\sum_{k=1}^\infty \frac{H(k)}{(k+1)^{7}(k+2)^{3}} \sumend &= \frac{1}{4}\left(  180 -32\zeta(2) -116\zeta(3) +14\zeta(4) -80\zeta(5) +40\zeta(2)\zeta(3)  \right. \nonumber \\ &\left. \hspace{1em}
+18\zeta(6) -12\zeta(3)^2 -36\zeta(7) +12\zeta(2)\zeta(5) +12\zeta(3)\zeta(4) +5\zeta(8)  \right. \nonumber \\ &\left. \hspace{1em}
-4\zeta(3)\zeta(5)\right) \label{eq11038}
\end{align}
 
\begin{align}
\sum_{k=1}^\infty \frac{H(k)}{k^{6}(k+2)^{4}} \sumend &= \frac{1}{256}\left( -122 +9\zeta(2) +107\zeta(3) -43\zeta(4) +128\zeta(5) -44\zeta(2)\zeta(3)  \right. \nonumber \\ &\left. \hspace{1em}
-56\zeta(6) +16\zeta(3)^2 +64\zeta(7) -16\zeta(2)\zeta(5) -16\zeta(3)\zeta(4)\right) \label{eq11039} \\
\sum_{k=1}^\infty \frac{H(k)}{k^{5}(k+1)(k+2)^{4}} \sumend &= \frac{1}{256}\left( -325 +257\zeta(2) -162\zeta(3) +125\zeta(4) -128\zeta(5)  \right. \nonumber \\ &\left. \hspace{1em}
+40\zeta(2)\zeta(3) +28\zeta(6) -8\zeta(3)^2\right) \label{eq11040} \\
\sum_{k=1}^\infty \frac{H(k)}{k^{4}(k+1)^{2}(k+2)^{4}} \sumend &= \frac{-1}{32}\left(  105 +6\zeta(2) -99\zeta(3) +6\zeta(4) -10\zeta(5)  \right. \nonumber \\ &\left. \hspace{1em}
+4\zeta(2)\zeta(3)\right) \label{eq11041} \\
\sum_{k=1}^\infty \frac{H(k)}{k^{3}(k+1)^{3}(k+2)^{4}} \sumend &= \frac{1}{32}\left( -265 +99\zeta(2) +82\zeta(3) +3\zeta(4) +8\zeta(5)  \right. \nonumber \\ &\left. \hspace{1em}
-4\zeta(2)\zeta(3)\right) \label{eq11042} \\
\sum_{k=1}^\infty \frac{H(k)}{k^{2}(k+1)^{4}(k+2)^{4}} \sumend &= \frac{-1}{16}\left(  328 -77\zeta(2) -165\zeta(3) -\zeta(4) -40\zeta(5)  \right. \nonumber \\ &\left. \hspace{1em}
+20\zeta(2)\zeta(3)\right) \label{eq11043} \\
\sum_{k=1}^\infty \frac{H(k)}{k(k+1)^{5}(k+2)^{4}} \sumend &= \frac{1}{16}\left( -799 +195\zeta(2) +402\zeta(3) -9\zeta(4) +112\zeta(5)  \right. \nonumber \\ &\left. \hspace{1em}
-56\zeta(2)\zeta(3) -12\zeta(6) +8\zeta(3)^2\right) \label{eq11044} \\
\sum_{k=1}^\infty \frac{H(k)}{(k+1)^{6}(k+2)^{4}} \sumend &= \frac{-1}{2}\left(  240 -56\zeta(2) -126\zeta(3) +5\zeta(4) -44\zeta(5) +22\zeta(2)\zeta(3)  \right. \nonumber \\ &\left. \hspace{1em}
+6\zeta(6) -4\zeta(3)^2 -6\zeta(7) +2\zeta(2)\zeta(5) +2\zeta(3)\zeta(4)\right) \label{eq11045} \\
\sum_{k=1}^\infty \frac{H(k)}{k^{5}(k+2)^{5}} \sumend &= \frac{1}{256}\left(  315 -58\zeta(2) -163\zeta(3) +2\zeta(4) -108\zeta(5) +40\zeta(2)\zeta(3)  \right. \nonumber \\ &\left. \hspace{1em}
+8\zeta(6)\right) \label{eq11046} \\
\sum_{k=1}^\infty \frac{H(k)}{k^{4}(k+1)(k+2)^{5}} \sumend &= \frac{-1}{256}\left( -955 +373\zeta(2) +164\zeta(3) +121\zeta(4) +88\zeta(5)  \right. \nonumber \\ &\left. \hspace{1em}
-40\zeta(2)\zeta(3) +12\zeta(6) -8\zeta(3)^2\right) \label{eq11047} \\
\sum_{k=1}^\infty \frac{H(k)}{k^{3}(k+1)^{2}(k+2)^{5}} \sumend &= \frac{-1}{128}\left( -1375 +349\zeta(2) +560\zeta(3) +97\zeta(4) +128\zeta(5)  \right. \nonumber \\ &\left. \hspace{1em}
-56\zeta(2)\zeta(3) +12\zeta(6) -8\zeta(3)^2\right) \label{eq11048} \\
\sum_{k=1}^\infty \frac{H(k)}{k^{2}(k+1)^{3}(k+2)^{5}} \sumend &= \frac{1}{64}\left(  1905 -547\zeta(2) -724\zeta(3) -103\zeta(4) -144\zeta(5)  \right. \nonumber \\ &\left. \hspace{1em}
+64\zeta(2)\zeta(3) -12\zeta(6) +8\zeta(3)^2\right) \label{eq11049}
\end{align}
 
\begin{align}
\sum_{k=1}^\infty \frac{H(k)}{k(k+1)^{4}(k+2)^{5}} \sumend &= \frac{-1}{32}\left( -2561 +701\zeta(2) +1054\zeta(3) +105\zeta(4) +224\zeta(5)  \right. \nonumber \\ &\left. \hspace{1em}
-104\zeta(2)\zeta(3) +12\zeta(6) -8\zeta(3)^2\right) \label{eq11050} \\
\sum_{k=1}^\infty \frac{H(k)}{(k+1)^{5}(k+2)^{5}} \sumend &= -\left( -210 +56\zeta(2) +91\zeta(3) +6\zeta(4) +21\zeta(5) -10\zeta(2)\zeta(3)\right) \label{eq11051} \\
\sum_{k=1}^\infty \frac{H(k)}{k^{4}(k+2)^{6}} \sumend &= \frac{1}{256}\left( -634 +135\zeta(2) +205\zeta(3) +83\zeta(4) +140\zeta(5) -44\zeta(2)\zeta(3)  \right. \nonumber \\ &\left. \hspace{1em}
+40\zeta(6) -16\zeta(3)^2 +48\zeta(7) -16\zeta(2)\zeta(5) -16\zeta(3)\zeta(4)\right) \label{eq11052} \\
\sum_{k=1}^\infty \frac{H(k)}{k^{3}(k+1)(k+2)^{6}} \sumend &= \frac{1}{256}\left( -2223 +643\zeta(2) +574\zeta(3) +287\zeta(4) +368\zeta(5)  \right. \nonumber \\ &\left. \hspace{1em}
-128\zeta(2)\zeta(3) +92\zeta(6) -40\zeta(3)^2 +96\zeta(7) -32\zeta(2)\zeta(5)  \right. \nonumber \\ &\left. \hspace{1em}
-32\zeta(3)\zeta(4)\right) \label{eq11053} \\
\sum_{k=1}^\infty \frac{H(k)}{k^{2}(k+1)^{2}(k+2)^{6}} \sumend &= \frac{1}{64}\left( -1799 +496\zeta(2) +567\zeta(3) +192\zeta(4) +248\zeta(5)  \right. \nonumber \\ &\left. \hspace{1em}
-92\zeta(2)\zeta(3) +52\zeta(6) -24\zeta(3)^2 +48\zeta(7) -16\zeta(2)\zeta(5)  \right. \nonumber \\ &\left. \hspace{1em}
-16\zeta(3)\zeta(4)\right) \label{eq11054} \\
\sum_{k=1}^\infty \frac{H(k)}{k(k+1)^{3}(k+2)^{6}} \sumend &= \frac{-1}{64}\left(  5503 -1539\zeta(2) -1858\zeta(3) -487\zeta(4) -640\zeta(5)  \right. \nonumber \\ &\left. \hspace{1em}
+248\zeta(2)\zeta(3) -116\zeta(6) +56\zeta(3)^2 -96\zeta(7) +32\zeta(2)\zeta(5)  \right. \nonumber \\ &\left. \hspace{1em}
+32\zeta(3)\zeta(4)\right) \label{eq11055} \\
\sum_{k=1}^\infty \frac{H(k)}{(k+1)^{4}(k+2)^{6}} \sumend &= \frac{1}{2}\left( -504 +140\zeta(2) +182\zeta(3) +37\zeta(4) +54\zeta(5) -22\zeta(2)\zeta(3)  \right. \nonumber \\ &\left. \hspace{1em}
+8\zeta(6) -4\zeta(3)^2 +6\zeta(7) -2\zeta(2)\zeta(5) -2\zeta(3)\zeta(4)\right) \label{eq11056} \\
\sum_{k=1}^\infty \frac{H(k)}{k^{3}(k+2)^{7}} \sumend &= \frac{-1}{256}\left( -1058 +205\zeta(2) +233\zeta(3) +173\zeta(4) +208\zeta(5)  \right. \nonumber \\ &\left. \hspace{1em}
-40\zeta(2)\zeta(3) +116\zeta(6) -24\zeta(3)^2 +176\zeta(7) -48\zeta(2)\zeta(5) -48\zeta(3)\zeta(4)  \right. \nonumber \\ &\left. \hspace{1em}
+40\zeta(8) -32\zeta(3)\zeta(5)\right) \label{eq11057} \\
\sum_{k=1}^\infty \frac{H(k)}{k^{2}(k+1)(k+2)^{7}} \sumend &= \frac{-1}{256}\left( -4339 +1053\zeta(2) +1040\zeta(3) +633\zeta(4) +784\zeta(5)  \right. \nonumber \\ &\left. \hspace{1em}
-208\zeta(2)\zeta(3) +324\zeta(6) -88\zeta(3)^2 +448\zeta(7) -128\zeta(2)\zeta(5) -128\zeta(3)\zeta(4)  \right. \nonumber \\ &\left. \hspace{1em}
+80\zeta(8) -64\zeta(3)\zeta(5)\right) \label{eq11058} \\
\sum_{k=1}^\infty \frac{H(k)}{k(k+1)^{2}(k+2)^{7}} \sumend &= \frac{-1}{128}\left( -7937 +2045\zeta(2) +2174\zeta(3) +1017\zeta(4) +1280\zeta(5)  \right. \nonumber \\ &\left. \hspace{1em}
-392\zeta(2)\zeta(3) +428\zeta(6) -136\zeta(3)^2 +544\zeta(7) -160\zeta(2)\zeta(5) -160\zeta(3)\zeta(4)  \right. \nonumber \\ &\left. \hspace{1em}
+80\zeta(8) -64\zeta(3)\zeta(5)\right) \label{eq11059}
\end{align}
 
\begin{align}
\sum_{k=1}^\infty \frac{H(k)}{(k+1)^{3}(k+2)^{7}} \sumend &= \frac{-1}{4}\left( -840 +224\zeta(2) +252\zeta(3) +94\zeta(4) +120\zeta(5)  \right. \nonumber \\ &\left. \hspace{1em}
-40\zeta(2)\zeta(3) +34\zeta(6) -12\zeta(3)^2 +40\zeta(7) -12\zeta(2)\zeta(5) -12\zeta(3)\zeta(4)  \right. \nonumber \\ &\left. \hspace{1em}
+5\zeta(8) -4\zeta(3)\zeta(5)\right) \label{eq11060} \\
\sum_{k=1}^\infty \frac{H(k)}{k^{2}(k+2)^{8}} \sumend &= \frac{1}{256}\left( -1542 +243\zeta(2) +249\zeta(3) +231\zeta(4) +248\zeta(5)  \right. \nonumber \\ &\left. \hspace{1em}
-20\zeta(2)\zeta(3) +200\zeta(6) -16\zeta(3)^2 +272\zeta(7) -48\zeta(2)\zeta(5) -48\zeta(3)\zeta(4)  \right. \nonumber \\ &\left. \hspace{1em}
+144\zeta(8) -64\zeta(3)\zeta(5) +256\zeta(9) -64\zeta(3)\zeta(6) -64\zeta(4)\zeta(5)  \right. \nonumber \\ &\left. \hspace{1em}
-64\zeta(2)\zeta(7)\right) \label{eq11061} \\
\sum_{k=1}^\infty \frac{H(k)}{k(k+1)(k+2)^{8}} \sumend &= \frac{-1}{256}\left(  7423 -1539\zeta(2) -1538\zeta(3) -1095\zeta(4) -1280\zeta(5)  \right. \nonumber \\ &\left. \hspace{1em}
+248\zeta(2)\zeta(3) -724\zeta(6) +120\zeta(3)^2 -992\zeta(7) +224\zeta(2)\zeta(5) +224\zeta(3)\zeta(4)  \right. \nonumber \\ &\left. \hspace{1em}
-368\zeta(8) +192\zeta(3)\zeta(5) -512\zeta(9) +128\zeta(3)\zeta(6) +128\zeta(4)\zeta(5)  \right. \nonumber \\ &\left. \hspace{1em}
+128\zeta(2)\zeta(7)\right) \label{eq11062} \\
\sum_{k=1}^\infty \frac{H(k)}{(k+1)^{2}(k+2)^{8}} \sumend &= \frac{-1}{2}\left(  240 -56\zeta(2) -58\zeta(3) -33\zeta(4) -40\zeta(5) +10\zeta(2)\zeta(3)  \right. \nonumber \\ &\left. \hspace{1em}
-18\zeta(6) +4\zeta(3)^2 -24\zeta(7) +6\zeta(2)\zeta(5) +6\zeta(3)\zeta(4) -7\zeta(8)  \right. \nonumber \\ &\left. \hspace{1em}
+4\zeta(3)\zeta(5) -8\zeta(9) +2\zeta(3)\zeta(6) +2\zeta(4)\zeta(5) +2\zeta(2)\zeta(7)\right) \label{eq11063} \\
\sum_{k=1}^\infty \frac{H(k)}{k(k+2)^{9}} \sumend &= \frac{1}{512}\left(  4097 -509\zeta(2) -510\zeta(3) -505\zeta(4) -512\zeta(5) +8\zeta(2)\zeta(3)  \right. \nonumber \\ &\left. \hspace{1em}
-492\zeta(6) +8\zeta(3)^2 -544\zeta(7) +32\zeta(2)\zeta(5) +32\zeta(3)\zeta(4) -464\zeta(8)  \right. \nonumber \\ &\left. \hspace{1em}
+64\zeta(3)\zeta(5) -768\zeta(9) +128\zeta(3)\zeta(6) +128\zeta(4)\zeta(5) +128\zeta(2)\zeta(7)  \right. \nonumber \\ &\left. \hspace{1em}
-448\zeta(10) +256\zeta(3)\zeta(7) +128\zeta(5)^2\right) \label{eq11064} \\
\sum_{k=1}^\infty \frac{H(k)}{(k+1)(k+2)^{9}} \sumend &= \frac{-1}{4}\left( -180 +32\zeta(2) +32\zeta(3) +25\zeta(4) +28\zeta(5) -4\zeta(2)\zeta(3)  \right. \nonumber \\ &\left. \hspace{1em}
+19\zeta(6) -2\zeta(3)^2 +24\zeta(7) -4\zeta(2)\zeta(5) -4\zeta(3)\zeta(4) +13\zeta(8)  \right. \nonumber \\ &\left. \hspace{1em}
-4\zeta(3)\zeta(5) +20\zeta(9) -4\zeta(3)\zeta(6) -4\zeta(4)\zeta(5) -4\zeta(2)\zeta(7) +7\zeta(10)  \right. \nonumber \\ &\left. \hspace{1em}
-4\zeta(3)\zeta(7) -2\zeta(5)^2\right) \label{eq11065} \\
\sum_{k=1}^\infty \frac{H(k)}{(k+2)^{10}} \sumend &= \left( -10 +\zeta(2) +\zeta(3) +\zeta(4) +\zeta(5) +\zeta(6) +\zeta(7) +\zeta(8) +\zeta(9)  \right. \nonumber \\ &\left. \hspace{1em}
+\zeta(10) +5\zeta(11) -\zeta(2)\zeta(9) -\zeta(3)\zeta(8) -\zeta(4)\zeta(7) -\zeta(5)\zeta(6)\right) \label{eq11066} \\
\sum_{k=1}^\infty \frac{H(k)^{2}}{k^{9}} \sumend &= \frac{1}{2}\left(  26\zeta(11) -2\zeta(2)\zeta(9) -9\zeta(3)\zeta(8) -5\zeta(4)\zeta(7)  \right. \nonumber \\ &\left. \hspace{1em}
-7\zeta(5)\zeta(6) +2\zeta(3)^2\zeta(5)\right) \label{eq11067}
\end{align}
 
\begin{align}
\sum_{k=1}^\infty \frac{H(k)^{2}}{k^{8}(k+1)} \sumend &= \frac{-1}{24}\left(  72\zeta(3) -102\zeta(4) +84\zeta(5) -24\zeta(2)\zeta(3) -97\zeta(6)  \right. \nonumber \\ &\left. \hspace{1em}
+48\zeta(3)^2 +144\zeta(7) -24\zeta(2)\zeta(5) -60\zeta(3)\zeta(4) -24 M(2,6) +220\zeta(9)  \right. \nonumber \\ &\left. \hspace{1em}
-84\zeta(3)\zeta(6) -60\zeta(4)\zeta(5) -24\zeta(2)\zeta(7) +8\zeta(3)^3 -24 M(2,8)\right) \label{eq11068} \\
\sum_{k=1}^\infty \frac{H(k)^{2}}{k^{7}(k+1)^{2}} \sumend &= \frac{1}{12}\left(  252\zeta(3) -339\zeta(4) +210\zeta(5) -60\zeta(2)\zeta(3) -194\zeta(6)  \right. \nonumber \\ &\left. \hspace{1em}
+96\zeta(3)^2 +216\zeta(7) -36\zeta(2)\zeta(5) -90\zeta(3)\zeta(4) -24 M(2,6) +110\zeta(9)  \right. \nonumber \\ &\left. \hspace{1em}
-42\zeta(3)\zeta(6) -30\zeta(4)\zeta(5) -12\zeta(2)\zeta(7) +4\zeta(3)^3\right) \label{eq11069} \\
\sum_{k=1}^\infty \frac{H(k)^{2}}{k^{6}(k+1)^{3}} \sumend &= \frac{-1}{4}\left(  252\zeta(3) -321\zeta(4) +146\zeta(5) -44\zeta(2)\zeta(3) -97\zeta(6)  \right. \nonumber \\ &\left. \hspace{1em}
+48\zeta(3)^2 +72\zeta(7) -12\zeta(2)\zeta(5) -30\zeta(3)\zeta(4) -4 M(2,6)\right) \label{eq11070} \\
\sum_{k=1}^\infty \frac{H(k)^{2}}{k^{5}(k+1)^{4}} \sumend &= \frac{1}{24}\left(  2520\zeta(3) -3030\zeta(4) +1020\zeta(5) -360\zeta(2)\zeta(3)  \right. \nonumber \\ &\left. \hspace{1em}
-425\zeta(6) +216\zeta(3)^2 +144\zeta(7) -24\zeta(2)\zeta(5) -60\zeta(3)\zeta(4)\right) \label{eq11071} \\
\sum_{k=1}^\infty \frac{H(k)^{2}}{k^{4}(k+1)^{5}} \sumend &= \frac{1}{24}\left( -2520\zeta(3) +2850\zeta(4) -780\zeta(5) +360\zeta(2)\zeta(3)  \right. \nonumber \\ &\left. \hspace{1em}
+245\zeta(6) -144\zeta(3)^2 -24\zeta(7) +24\zeta(2)\zeta(5) -12\zeta(3)\zeta(4)\right) \label{eq11072} \\
\sum_{k=1}^\infty \frac{H(k)^{2}}{k^{3}(k+1)^{6}} \sumend &= \frac{-1}{4}\left( -252\zeta(3) +267\zeta(4) -74\zeta(5) +44\zeta(2)\zeta(3) +37\zeta(6)  \right. \nonumber \\ &\left. \hspace{1em}
-24\zeta(3)^2 -12\zeta(7) +12\zeta(2)\zeta(5) -6\zeta(3)\zeta(4) -14\zeta(8) +8\zeta(3)\zeta(5)  \right. \nonumber \\ &\left. \hspace{1em}
+4 M(2,6)\right) \label{eq11073} \\
\sum_{k=1}^\infty \frac{H(k)^{2}}{k^{2}(k+1)^{7}} \sumend &= \frac{-1}{12}\left(  252\zeta(3) -249\zeta(4) +90\zeta(5) -60\zeta(2)\zeta(3) -74\zeta(6)  \right. \nonumber \\ &\left. \hspace{1em}
+48\zeta(3)^2 +36\zeta(7) -36\zeta(2)\zeta(5) +18\zeta(3)\zeta(4) +84\zeta(8) -48\zeta(3)\zeta(5)  \right. \nonumber \\ &\left. \hspace{1em}
-24 M(2,6) -2\zeta(9) +18\zeta(3)\zeta(6) +6\zeta(4)\zeta(5) -12\zeta(2)\zeta(7) -4\zeta(3)^3\right) \label{eq11074} \\
\sum_{k=1}^\infty \frac{H(k)^{2}}{k(k+1)^{8}} \sumend &= \frac{1}{24}\left(  72\zeta(3) -66\zeta(4) +36\zeta(5) -24\zeta(2)\zeta(3) -37\zeta(6)  \right. \nonumber \\ &\left. \hspace{1em}
+24\zeta(3)^2 +24\zeta(7) -24\zeta(2)\zeta(5) +12\zeta(3)\zeta(4) +84\zeta(8) -48\zeta(3)\zeta(5)  \right. \nonumber \\ &\left. \hspace{1em}
-24 M(2,6) -4\zeta(9) +36\zeta(3)\zeta(6) +12\zeta(4)\zeta(5) -24\zeta(2)\zeta(7) -8\zeta(3)^3  \right. \nonumber \\ &\left. \hspace{1em}
+108\zeta(10) -48\zeta(3)\zeta(7) -24\zeta(5)^2 -24 M(2,8)\right) \label{eq11075} \\
\sum_{k=1}^\infty \frac{H(k)^{2}}{(k+1)^{9}} \sumend &= \frac{1}{2}\left(  4\zeta(11) +2\zeta(2)\zeta(9) -5\zeta(3)\zeta(8) -\zeta(4)\zeta(7)  \right. \nonumber \\ &\left. \hspace{1em}
-3\zeta(5)\zeta(6) +2\zeta(3)^2\zeta(5)\right) \label{eq11076}
\end{align}
 
\begin{align}
\sum_{k=1}^\infty \frac{H(k)^{2}}{k^{8}(k+2)} \sumend &= \frac{-1}{1536}\left(  6 +6\zeta(2) +18\zeta(3) -51\zeta(4) +84\zeta(5) -24\zeta(2)\zeta(3)  \right. \nonumber \\ &\left. \hspace{1em}
-194\zeta(6) +96\zeta(3)^2 +576\zeta(7) -96\zeta(2)\zeta(5) -240\zeta(3)\zeta(4) -192 M(2,6)  \right. \nonumber \\ &\left. \hspace{1em}
+3520\zeta(9) -1344\zeta(3)\zeta(6) -960\zeta(4)\zeta(5) -384\zeta(2)\zeta(7) +128\zeta(3)^3  \right. \nonumber \\ &\left. \hspace{1em}
-768 M(2,8)\right) \label{eq11077} \\
\sum_{k=1}^\infty \frac{H(k)^{2}}{k^{7}(k+1)(k+2)} \sumend &= \frac{1}{768}\left( -6 -6\zeta(2) +2286\zeta(3) -3213\zeta(4) +2604\zeta(5)  \right. \nonumber \\ &\left. \hspace{1em}
-744\zeta(2)\zeta(3) -2910\zeta(6) +1440\zeta(3)^2 +4032\zeta(7) -672\zeta(2)\zeta(5)  \right. \nonumber \\ &\left. \hspace{1em}
-1680\zeta(3)\zeta(4) -576 M(2,6) +3520\zeta(9) -1344\zeta(3)\zeta(6) -960\zeta(4)\zeta(5)  \right. \nonumber \\ &\left. \hspace{1em}
-384\zeta(2)\zeta(7) +128\zeta(3)^3\right) \label{eq11078} \\
\sum_{k=1}^\infty \frac{H(k)^{2}}{k^{6}(k+1)^{2}(k+2)} \sumend &= \frac{1}{384}\left( -6 -6\zeta(2) -5778\zeta(3) +7635\zeta(4) -4116\zeta(5)  \right. \nonumber \\ &\left. \hspace{1em}
+1176\zeta(2)\zeta(3) +3298\zeta(6) -1632\zeta(3)^2 -2880\zeta(7) +480\zeta(2)\zeta(5)  \right. \nonumber \\ &\left. \hspace{1em}
+1200\zeta(3)\zeta(4) +192 M(2,6)\right) \label{eq11079} \\
\sum_{k=1}^\infty \frac{H(k)^{2}}{k^{5}(k+1)^{3}(k+2)} \sumend &= \frac{-1}{192}\left(  6 +6\zeta(2) -6318\zeta(3) +7773\zeta(4) -2892\zeta(5)  \right. \nonumber \\ &\left. \hspace{1em}
+936\zeta(2)\zeta(3) +1358\zeta(6) -672\zeta(3)^2 -576\zeta(7) +96\zeta(2)\zeta(5)  \right. \nonumber \\ &\left. \hspace{1em}
+240\zeta(3)\zeta(4)\right) \label{eq11080} \\
\sum_{k=1}^\infty \frac{H(k)^{2}}{k^{4}(k+1)^{4}(k+2)} \sumend &= \frac{1}{32}\left( -2 -2\zeta(2) -1254\zeta(3) +1449\zeta(4) -396\zeta(5)  \right. \nonumber \\ &\left. \hspace{1em}
+168\zeta(2)\zeta(3) +114\zeta(6) -64\zeta(3)^2\right) \label{eq11081} \\
\sum_{k=1}^\infty \frac{H(k)^{2}}{k^{3}(k+1)^{5}(k+2)} \sumend &= \frac{1}{48}\left( -6 -6\zeta(2) +1278\zeta(3) -1353\zeta(4) +372\zeta(5)  \right. \nonumber \\ &\left. \hspace{1em}
-216\zeta(2)\zeta(3) -148\zeta(6) +96\zeta(3)^2 +48\zeta(7) -48\zeta(2)\zeta(5)  \right. \nonumber \\ &\left. \hspace{1em}
+24\zeta(3)\zeta(4)\right) \label{eq11082} \\
\sum_{k=1}^\infty \frac{H(k)^{2}}{k^{2}(k+1)^{6}(k+2)} \sumend &= \frac{1}{24}\left( -6 -6\zeta(2) -234\zeta(3) +249\zeta(4) -72\zeta(5)  \right. \nonumber \\ &\left. \hspace{1em}
+48\zeta(2)\zeta(3) +74\zeta(6) -48\zeta(3)^2 -24\zeta(7) +24\zeta(2)\zeta(5) -12\zeta(3)\zeta(4)  \right. \nonumber \\ &\left. \hspace{1em}
-84\zeta(8) +48\zeta(3)\zeta(5) +24 M(2,6)\right) \label{eq11083} \\
\sum_{k=1}^\infty \frac{H(k)^{2}}{k(k+1)^{7}(k+2)} \sumend &= \frac{1}{6}\left( -3 -3\zeta(2) +9\zeta(3) +9\zeta(5) -6\zeta(2)\zeta(3) +6\zeta(7)  \right. \nonumber \\ &\left. \hspace{1em}
-6\zeta(2)\zeta(5) +3\zeta(3)\zeta(4) -\zeta(9) +9\zeta(3)\zeta(6) +3\zeta(4)\zeta(5)  \right. \nonumber \\ &\left. \hspace{1em}
-6\zeta(2)\zeta(7) -2\zeta(3)^3\right) \label{eq11084}
\end{align}
 
\begin{align}
\sum_{k=1}^\infty \frac{H(k)^{2}}{(k+1)^{8}(k+2)} \sumend &= \frac{1}{24}\left( -24 -24\zeta(2) +66\zeta(4) +36\zeta(5) -24\zeta(2)\zeta(3) +37\zeta(6)  \right. \nonumber \\ &\left. \hspace{1em}
-24\zeta(3)^2 +24\zeta(7) -24\zeta(2)\zeta(5) +12\zeta(3)\zeta(4) -84\zeta(8) +48\zeta(3)\zeta(5)  \right. \nonumber \\ &\left. \hspace{1em}
+24 M(2,6) -4\zeta(9) +36\zeta(3)\zeta(6) +12\zeta(4)\zeta(5) -24\zeta(2)\zeta(7) -8\zeta(3)^3  \right. \nonumber \\ &\left. \hspace{1em}
-108\zeta(10) +48\zeta(3)\zeta(7) +24\zeta(5)^2 +24 M(2,8)\right) \label{eq11085} \\
\sum_{k=1}^\infty \frac{H(k)^{2}}{k^{7}(k+2)^{2}} \sumend &= \frac{1}{1536}\left(  78 +42\zeta(2) +102\zeta(3) -339\zeta(4) +420\zeta(5)  \right. \nonumber \\ &\left. \hspace{1em}
-120\zeta(2)\zeta(3) -776\zeta(6) +384\zeta(3)^2 +1728\zeta(7) -288\zeta(2)\zeta(5)  \right. \nonumber \\ &\left. \hspace{1em}
-720\zeta(3)\zeta(4) -384 M(2,6) +3520\zeta(9) -1344\zeta(3)\zeta(6) -960\zeta(4)\zeta(5)  \right. \nonumber \\ &\left. \hspace{1em}
-384\zeta(2)\zeta(7) +128\zeta(3)^3\right) \label{eq11086} \\
\sum_{k=1}^\infty \frac{H(k)^{2}}{k^{6}(k+1)(k+2)^{2}} \sumend &= \frac{-1}{384}\left( -42 -24\zeta(2) +1092\zeta(3) -1437\zeta(4) +1092\zeta(5)  \right. \nonumber \\ &\left. \hspace{1em}
-312\zeta(2)\zeta(3) -1067\zeta(6) +528\zeta(3)^2 +1152\zeta(7) -192\zeta(2)\zeta(5)  \right. \nonumber \\ &\left. \hspace{1em}
-480\zeta(3)\zeta(4) -96 M(2,6)\right) \label{eq11087} \\
\sum_{k=1}^\infty \frac{H(k)^{2}}{k^{5}(k+1)^{2}(k+2)^{2}} \sumend &= \frac{-1}{128}\left( -30 -18\zeta(2) -1198\zeta(3) +1587\zeta(4) -644\zeta(5)  \right. \nonumber \\ &\left. \hspace{1em}
+184\zeta(2)\zeta(3) +388\zeta(6) -192\zeta(3)^2 -192\zeta(7) +32\zeta(2)\zeta(5)  \right. \nonumber \\ &\left. \hspace{1em}
+80\zeta(3)\zeta(4)\right) \label{eq11088} \\
\sum_{k=1}^\infty \frac{H(k)^{2}}{k^{4}(k+1)^{3}(k+2)^{2}} \sumend &= \frac{-1}{96}\left( -48 -30\zeta(2) +1362\zeta(3) -1506\zeta(4) +480\zeta(5)  \right. \nonumber \\ &\left. \hspace{1em}
-192\zeta(2)\zeta(3) -97\zeta(6) +48\zeta(3)^2\right) \label{eq11089} \\
\sum_{k=1}^\infty \frac{H(k)^{2}}{k^{3}(k+1)^{4}(k+2)^{2}} \sumend &= \frac{1}{96}\left(  102 +66\zeta(2) +1038\zeta(3) -1335\zeta(4) +228\zeta(5)  \right. \nonumber \\ &\left. \hspace{1em}
-120\zeta(2)\zeta(3) -148\zeta(6) +96\zeta(3)^2\right) \label{eq11090} \\
\sum_{k=1}^\infty \frac{H(k)^{2}}{k^{2}(k+1)^{5}(k+2)^{2}} \sumend &= \frac{1}{8}\left(  18 +12\zeta(2) -40\zeta(3) +3\zeta(4) -24\zeta(5)  \right. \nonumber \\ &\left. \hspace{1em}
+16\zeta(2)\zeta(3) -8\zeta(7) +8\zeta(2)\zeta(5) -4\zeta(3)\zeta(4)\right) \label{eq11091} \\
\sum_{k=1}^\infty \frac{H(k)^{2}}{k(k+1)^{6}(k+2)^{2}} \sumend &= \frac{1}{24}\left(  114 +78\zeta(2) -6\zeta(3) -231\zeta(4) -72\zeta(5)  \right. \nonumber \\ &\left. \hspace{1em}
+48\zeta(2)\zeta(3) -74\zeta(6) +48\zeta(3)^2 -24\zeta(7) +24\zeta(2)\zeta(5) -12\zeta(3)\zeta(4)  \right. \nonumber \\ &\left. \hspace{1em}
+84\zeta(8) -48\zeta(3)\zeta(5) -24 M(2,6)\right) \label{eq11092} \\
\sum_{k=1}^\infty \frac{H(k)^{2}}{(k+1)^{7}(k+2)^{2}} \sumend &= \frac{1}{12}\left(  120 +84\zeta(2) -24\zeta(3) -231\zeta(4) -90\zeta(5)  \right. \nonumber \\ &\left. \hspace{1em}
+60\zeta(2)\zeta(3) -74\zeta(6) +48\zeta(3)^2 -36\zeta(7) +36\zeta(2)\zeta(5) -18\zeta(3)\zeta(4)  \right. \nonumber \\ &\left. \hspace{1em}
+84\zeta(8) -48\zeta(3)\zeta(5) -24 M(2,6) +2\zeta(9) -18\zeta(3)\zeta(6) -6\zeta(4)\zeta(5)  \right. \nonumber \\ &\left. \hspace{1em}
+12\zeta(2)\zeta(7) +4\zeta(3)^3\right) \label{eq11093}
\end{align}
 
\begin{align}
\sum_{k=1}^\infty \frac{H(k)^{2}}{k^{6}(k+2)^{3}} \sumend &= \frac{-1}{512}\left(  162 +34\zeta(2) +54\zeta(3) -325\zeta(4) +292\zeta(5)  \right. \nonumber \\ &\left. \hspace{1em}
-88\zeta(2)\zeta(3) -388\zeta(6) +192\zeta(3)^2 +576\zeta(7) -96\zeta(2)\zeta(5) -240\zeta(3)\zeta(4)  \right. \nonumber \\ &\left. \hspace{1em}
-64 M(2,6)\right) \label{eq11094} \\
\sum_{k=1}^\infty \frac{H(k)^{2}}{k^{5}(k+1)(k+2)^{3}} \sumend &= \frac{-1}{768}\left(  570 +150\zeta(2) -2022\zeta(3) +1899\zeta(4) -1308\zeta(5)  \right. \nonumber \\ &\left. \hspace{1em}
+360\zeta(2)\zeta(3) +970\zeta(6) -480\zeta(3)^2 -576\zeta(7) +96\zeta(2)\zeta(5)  \right. \nonumber \\ &\left. \hspace{1em}
+240\zeta(3)\zeta(4)\right) \label{eq11095} \\
\sum_{k=1}^\infty \frac{H(k)^{2}}{k^{4}(k+1)^{2}(k+2)^{3}} \sumend &= \frac{1}{192}\left( -330 -102\zeta(2) -786\zeta(3) +1431\zeta(4) -312\zeta(5)  \right. \nonumber \\ &\left. \hspace{1em}
+96\zeta(2)\zeta(3) +97\zeta(6) -48\zeta(3)^2\right) \label{eq11096} \\
\sum_{k=1}^\infty \frac{H(k)^{2}}{k^{3}(k+1)^{3}(k+2)^{3}} \sumend &= \frac{1}{32}\left( -126 -44\zeta(2) +192\zeta(3) -25\zeta(4) +56\zeta(5)  \right. \nonumber \\ &\left. \hspace{1em}
-32\zeta(2)\zeta(3)\right) \label{eq11097} \\
\sum_{k=1}^\infty \frac{H(k)^{2}}{k^{2}(k+1)^{4}(k+2)^{3}} \sumend &= \frac{-1}{96}\left(  858 +330\zeta(2) -114\zeta(3) -1185\zeta(4) -108\zeta(5)  \right. \nonumber \\ &\left. \hspace{1em}
+72\zeta(2)\zeta(3) -148\zeta(6) +96\zeta(3)^2\right) \label{eq11098} \\
\sum_{k=1}^\infty \frac{H(k)^{2}}{k(k+1)^{5}(k+2)^{3}} \sumend &= \frac{1}{48}\left( -966 -402\zeta(2) +354\zeta(3) +1167\zeta(4) +252\zeta(5)  \right. \nonumber \\ &\left. \hspace{1em}
-168\zeta(2)\zeta(3) +148\zeta(6) -96\zeta(3)^2 +48\zeta(7) -48\zeta(2)\zeta(5)  \right. \nonumber \\ &\left. \hspace{1em}
+24\zeta(3)\zeta(4)\right) \label{eq11099} \\
\sum_{k=1}^\infty \frac{H(k)^{2}}{(k+1)^{6}(k+2)^{3}} \sumend &= \frac{1}{4}\left( -180 -80\zeta(2) +60\zeta(3) +233\zeta(4) +54\zeta(5)  \right. \nonumber \\ &\left. \hspace{1em}
-36\zeta(2)\zeta(3) +37\zeta(6) -24\zeta(3)^2 +12\zeta(7) -12\zeta(2)\zeta(5) +6\zeta(3)\zeta(4)  \right. \nonumber \\ &\left. \hspace{1em}
-14\zeta(8) +8\zeta(3)\zeta(5) +4 M(2,6)\right) \label{eq11100} \\
\sum_{k=1}^\infty \frac{H(k)^{2}}{k^{5}(k+2)^{4}} \sumend &= \frac{1}{1536}\left(  1950 -6\zeta(2) -282\zeta(3) -1647\zeta(4) +828\zeta(5)  \right. \nonumber \\ &\left. \hspace{1em}
-264\zeta(2)\zeta(3) -850\zeta(6) +432\zeta(3)^2 +576\zeta(7) -96\zeta(2)\zeta(5)  \right. \nonumber \\ &\left. \hspace{1em}
-240\zeta(3)\zeta(4)\right) \label{eq11101} \\
\sum_{k=1}^\infty \frac{H(k)^{2}}{k^{4}(k+1)(k+2)^{4}} \sumend &= \frac{1}{64}\left(  210 +12\zeta(2) -192\zeta(3) +21\zeta(4) -40\zeta(5)  \right. \nonumber \\ &\left. \hspace{1em}
+8\zeta(2)\zeta(3) +10\zeta(6) -4\zeta(3)^2\right) \label{eq11102}
\end{align}
 
\begin{align}
\sum_{k=1}^\infty \frac{H(k)^{2}}{k^{3}(k+1)^{2}(k+2)^{4}} \sumend &= \frac{-1}{192}\left( -1590 -174\zeta(2) +366\zeta(3) +1305\zeta(4) -72\zeta(5)  \right. \nonumber \\ &\left. \hspace{1em}
+48\zeta(2)\zeta(3) +37\zeta(6) -24\zeta(3)^2\right) \label{eq11103} \\
\sum_{k=1}^\infty \frac{H(k)^{2}}{k^{2}(k+1)^{3}(k+2)^{4}} \sumend &= \frac{1}{96}\left(  1968 +306\zeta(2) -942\zeta(3) -1230\zeta(4) -96\zeta(5)  \right. \nonumber \\ &\left. \hspace{1em}
+48\zeta(2)\zeta(3) -37\zeta(6) +24\zeta(3)^2\right) \label{eq11104} \\
\sum_{k=1}^\infty \frac{H(k)^{2}}{k(k+1)^{4}(k+2)^{4}} \sumend &= \frac{-1}{32}\left( -1598 -314\zeta(2) +666\zeta(3) +1215\zeta(4) +100\zeta(5)  \right. \nonumber \\ &\left. \hspace{1em}
-56\zeta(2)\zeta(3) +74\zeta(6) -48\zeta(3)^2\right) \label{eq11105} \\
\sum_{k=1}^\infty \frac{H(k)^{2}}{(k+1)^{5}(k+2)^{4}} \sumend &= \frac{-1}{24}\left( -2880 -672\zeta(2) +1176\zeta(3) +2406\zeta(4) +276\zeta(5)  \right. \nonumber \\ &\left. \hspace{1em}
-168\zeta(2)\zeta(3) +185\zeta(6) -120\zeta(3)^2 +24\zeta(7) -24\zeta(2)\zeta(5)  \right. \nonumber \\ &\left. \hspace{1em}
+12\zeta(3)\zeta(4)\right) \label{eq11106} \\
\sum_{k=1}^\infty \frac{H(k)^{2}}{k^{4}(k+2)^{5}} \sumend &= \frac{-1}{1536}\left(  5730 -702\zeta(2) -1818\zeta(3) -2073\zeta(4) -468\zeta(5)  \right. \nonumber \\ &\left. \hspace{1em}
+216\zeta(2)\zeta(3) -634\zeta(6) +384\zeta(3)^2 +96\zeta(7) -96\zeta(2)\zeta(5)  \right. \nonumber \\ &\left. \hspace{1em}
+48\zeta(3)\zeta(4)\right) \label{eq11107} \\
\sum_{k=1}^\infty \frac{H(k)^{2}}{k^{3}(k+1)(k+2)^{5}} \sumend &= \frac{1}{768}\left( -8250 +558\zeta(2) +4122\zeta(3) +1821\zeta(4) +948\zeta(5)  \right. \nonumber \\ &\left. \hspace{1em}
-312\zeta(2)\zeta(3) +514\zeta(6) -336\zeta(3)^2 -96\zeta(7) +96\zeta(2)\zeta(5)  \right. \nonumber \\ &\left. \hspace{1em}
-48\zeta(3)\zeta(4)\right) \label{eq11108} \\
\sum_{k=1}^\infty \frac{H(k)^{2}}{k^{2}(k+1)^{2}(k+2)^{5}} \sumend &= \frac{-1}{128}\left(  3810 -70\zeta(2) -1618\zeta(3) -1477\zeta(4) -268\zeta(5)  \right. \nonumber \\ &\left. \hspace{1em}
+72\zeta(2)\zeta(3) -196\zeta(6) +128\zeta(3)^2 +32\zeta(7) -32\zeta(2)\zeta(5)  \right. \nonumber \\ &\left. \hspace{1em}
+16\zeta(3)\zeta(4)\right) \label{eq11109} \\
\sum_{k=1}^\infty \frac{H(k)^{2}}{k(k+1)^{3}(k+2)^{5}} \sumend &= \frac{-1}{192}\left(  15366 +402\zeta(2) -6738\zeta(3) -6891\zeta(4) -996\zeta(5)  \right. \nonumber \\ &\left. \hspace{1em}
+312\zeta(2)\zeta(3) -662\zeta(6) +432\zeta(3)^2 +96\zeta(7) -96\zeta(2)\zeta(5)  \right. \nonumber \\ &\left. \hspace{1em}
+48\zeta(3)\zeta(4)\right) \label{eq11110} \\
\sum_{k=1}^\infty \frac{H(k)^{2}}{(k+1)^{4}(k+2)^{5}} \sumend &= \frac{-1}{24}\left(  5040 +336\zeta(2) -2184\zeta(3) -2634\zeta(4) -324\zeta(5)  \right. \nonumber \\ &\left. \hspace{1em}
+120\zeta(2)\zeta(3) -221\zeta(6) +144\zeta(3)^2 +24\zeta(7) -24\zeta(2)\zeta(5)  \right. \nonumber \\ &\left. \hspace{1em}
+12\zeta(3)\zeta(4)\right) \label{eq11111}
\end{align}
 
\begin{align}
\sum_{k=1}^\infty \frac{H(k)^{2}}{k^{3}(k+2)^{6}} \sumend &= \frac{-1}{512}\left( -4446 +774\zeta(2) +1482\zeta(3) +915\zeta(4) +1100\zeta(5)  \right. \nonumber \\ &\left. \hspace{1em}
-424\zeta(2)\zeta(3) +452\zeta(6) -256\zeta(3)^2 +288\zeta(7) -32\zeta(2)\zeta(5) -176\zeta(3)\zeta(4)  \right. \nonumber \\ &\left. \hspace{1em}
-224\zeta(8) +128\zeta(3)\zeta(5) +64 M(2,6)\right) \label{eq11112} \\
\sum_{k=1}^\infty \frac{H(k)^{2}}{k^{2}(k+1)(k+2)^{6}} \sumend &= \frac{-1}{384}\left( -10794 +1440\zeta(2) +4284\zeta(3) +2283\zeta(4) +2124\zeta(5)  \right. \nonumber \\ &\left. \hspace{1em}
-792\zeta(2)\zeta(3) +935\zeta(6) -552\zeta(3)^2 +384\zeta(7) -288\zeta(3)\zeta(4) -336\zeta(8)  \right. \nonumber \\ &\left. \hspace{1em}
+192\zeta(3)\zeta(5) +96 M(2,6)\right) \label{eq11113} \\
\sum_{k=1}^\infty \frac{H(k)^{2}}{k(k+1)^{2}(k+2)^{6}} \sumend &= \frac{-1}{384}\left( -33018 +3090\zeta(2) +13422\zeta(3) +8997\zeta(4) +5052\zeta(5)  \right. \nonumber \\ &\left. \hspace{1em}
-1800\zeta(2)\zeta(3) +2458\zeta(6) -1488\zeta(3)^2 +672\zeta(7) +96\zeta(2)\zeta(5)  \right. \nonumber \\ &\left. \hspace{1em}
-624\zeta(3)\zeta(4) -672\zeta(8) +384\zeta(3)\zeta(5) +192 M(2,6)\right) \label{eq11114} \\
\sum_{k=1}^\infty \frac{H(k)^{2}}{(k+1)^{3}(k+2)^{6}} \sumend &= \frac{-1}{4}\left( -1008 +56\zeta(2) +420\zeta(3) +331\zeta(4) +126\zeta(5)  \right. \nonumber \\ &\left. \hspace{1em}
-44\zeta(2)\zeta(3) +65\zeta(6) -40\zeta(3)^2 +12\zeta(7) +4\zeta(2)\zeta(5) -14\zeta(3)\zeta(4)  \right. \nonumber \\ &\left. \hspace{1em}
-14\zeta(8) +8\zeta(3)\zeta(5) +4 M(2,6)\right) \label{eq11115} \\
\sum_{k=1}^\infty \frac{H(k)^{2}}{k^{2}(k+2)^{7}} \sumend &= \frac{1}{1536}\left( -26034 +4782\zeta(2) +7578\zeta(3) +4389\zeta(4) +7020\zeta(5)  \right. \nonumber \\ &\left. \hspace{1em}
-2376\zeta(2)\zeta(3) +3032\zeta(6) -1248\zeta(3)^2 +4704\zeta(7) -1248\zeta(2)\zeta(5)  \right. \nonumber \\ &\left. \hspace{1em}
-1680\zeta(3)\zeta(4) -384\zeta(8) +384 M(2,6) +64\zeta(9) -576\zeta(3)\zeta(6) -192\zeta(4)\zeta(5)  \right. \nonumber \\ &\left. \hspace{1em}
+384\zeta(2)\zeta(7) +128\zeta(3)^3\right) \label{eq11116} \\
\sum_{k=1}^\infty \frac{H(k)^{2}}{k(k+1)(k+2)^{7}} \sumend &= \frac{1}{768}\left( -47622 +7662\zeta(2) +16146\zeta(3) +8955\zeta(4) +11268\zeta(5)  \right. \nonumber \\ &\left. \hspace{1em}
-3960\zeta(2)\zeta(3) +4902\zeta(6) -2352\zeta(3)^2 +5472\zeta(7) -1248\zeta(2)\zeta(5)  \right. \nonumber \\ &\left. \hspace{1em}
-2256\zeta(3)\zeta(4) -1056\zeta(8) +384\zeta(3)\zeta(5) +576 M(2,6) +64\zeta(9) -576\zeta(3)\zeta(6)  \right. \nonumber \\ &\left. \hspace{1em}
-192\zeta(4)\zeta(5) +384\zeta(2)\zeta(7) +128\zeta(3)^3\right) \label{eq11117} \\
\sum_{k=1}^\infty \frac{H(k)^{2}}{(k+1)^{2}(k+2)^{7}} \sumend &= \frac{1}{12}\left( -2520 +336\zeta(2) +924\zeta(3) +561\zeta(4) +510\zeta(5)  \right. \nonumber \\ &\left. \hspace{1em}
-180\zeta(2)\zeta(3) +230\zeta(6) -120\zeta(3)^2 +192\zeta(7) -36\zeta(2)\zeta(5) -90\zeta(3)\zeta(4)  \right. \nonumber \\ &\left. \hspace{1em}
-54\zeta(8) +24\zeta(3)\zeta(5) +24 M(2,6) +2\zeta(9) -18\zeta(3)\zeta(6) -6\zeta(4)\zeta(5)  \right. \nonumber \\ &\left. \hspace{1em}
+12\zeta(2)\zeta(7) +4\zeta(3)^3\right) \label{eq11118} \\
\sum_{k=1}^\infty \frac{H(k)^{2}}{k(k+2)^{8}} \sumend &= \frac{-1}{1536}\left( -44538 +7698\zeta(2) +10734\zeta(3) +6981\zeta(4) +10620\zeta(5)  \right. \nonumber \\ &\left. \hspace{1em}
-2952\zeta(2)\zeta(3) +5498\zeta(6) -1488\zeta(3)^2 +9888\zeta(7) -2592\zeta(2)\zeta(5)  \right. \nonumber \\ &\left. \hspace{1em}
-2736\zeta(3)\zeta(4) +2976\zeta(8) -1920\zeta(3)\zeta(5) +192 M(2,6) +6208\zeta(9)  \right. \nonumber \\ &\left. \hspace{1em}
-2112\zeta(3)\zeta(6) -1728\zeta(4)\zeta(5) -1152\zeta(2)\zeta(7) +128\zeta(3)^3 -3456\zeta(10)  \right. \nonumber \\ &\left. \hspace{1em}
+1536\zeta(3)\zeta(7) +768\zeta(5)^2 +768 M(2,8)\right) \label{eq11119}
\end{align}
 
\begin{align}
\sum_{k=1}^\infty \frac{H(k)^{2}}{(k+1)(k+2)^{8}} \sumend &= \frac{-1}{24}\left( -2880 +480\zeta(2) +840\zeta(3) +498\zeta(4) +684\zeta(5)  \right. \nonumber \\ &\left. \hspace{1em}
-216\zeta(2)\zeta(3) +325\zeta(6) -120\zeta(3)^2 +480\zeta(7) -120\zeta(2)\zeta(5) -156\zeta(3)\zeta(4)  \right. \nonumber \\ &\left. \hspace{1em}
+60\zeta(8) -48\zeta(3)\zeta(5) +24 M(2,6) +196\zeta(9) -84\zeta(3)\zeta(6) -60\zeta(4)\zeta(5)  \right. \nonumber \\ &\left. \hspace{1em}
-24\zeta(2)\zeta(7) +8\zeta(3)^3 -108\zeta(10) +48\zeta(3)\zeta(7) +24\zeta(5)^2 +24 M(2,8)\right) \label{eq11120} \\
\sum_{k=1}^\infty \frac{H(k)^{2}}{(k+2)^{9}} \sumend &= \frac{-1}{2}\left(  90 -14\zeta(2) -18\zeta(3) -13\zeta(4) -18\zeta(5) +4\zeta(2)\zeta(3)  \right. \nonumber \\ &\left. \hspace{1em}
-11\zeta(6) +2\zeta(3)^2 -18\zeta(7) +4\zeta(2)\zeta(5) +4\zeta(3)\zeta(4) -9\zeta(8)  \right. \nonumber \\ &\left. \hspace{1em}
+4\zeta(3)\zeta(5) -18\zeta(9) +4\zeta(3)\zeta(6) +4\zeta(4)\zeta(5) +4\zeta(2)\zeta(7) -7\zeta(10)  \right. \nonumber \\ &\left. \hspace{1em}
+4\zeta(3)\zeta(7) +2\zeta(5)^2 -4\zeta(11) -2\zeta(2)\zeta(9) +5\zeta(3)\zeta(8) +\zeta(4)\zeta(7)  \right. \nonumber \\ &\left. \hspace{1em}
+3\zeta(5)\zeta(6) -2\zeta(3)^2\zeta(5)\right) \label{eq11121} \\
\sum_{k=1}^\infty \frac{H(k)^{3}}{k^{8}} \sumend &= -\left( - M(3,8)\right) \label{eq11122} \\
\sum_{k=1}^\infty \frac{H(k)^{3}}{k^{7}(k+1)} \sumend &= \frac{-1}{480}\left( -4800\zeta(4) +4800\zeta(5) +480\zeta(2)\zeta(3) -2790\zeta(6)  \right. \nonumber \\ &\left. \hspace{1em}
+1200\zeta(3)^2 +6930\zeta(7) +960\zeta(2)\zeta(5) -6120\zeta(3)\zeta(4) +2975\zeta(8)  \right. \nonumber \\ &\left. \hspace{1em}
+600\zeta(2)\zeta(3)^2 -2880\zeta(3)\zeta(5) -1320 M(2,6) +10420\zeta(9) -5820\zeta(3)\zeta(6)  \right. \nonumber \\ &\left. \hspace{1em}
-6120\zeta(4)\zeta(5) +1440\zeta(2)\zeta(7) +960\zeta(3)^3 +4983\zeta(10) -3840\zeta(3)\zeta(7)  \right. \nonumber \\ &\left. \hspace{1em}
-240\zeta(3)^2\zeta(4) +1680\zeta(2)\zeta(3)\zeta(5) -2160\zeta(5)^2 -1560 M(2,8)\right) \label{eq11123} \\
\sum_{k=1}^\infty \frac{H(k)^{3}}{k^{6}(k+1)^{2}} \sumend &= \frac{1}{48}\left( -2880\zeta(4) +2760\zeta(5) +288\zeta(2)\zeta(3) -1116\zeta(6)  \right. \nonumber \\ &\left. \hspace{1em}
+480\zeta(3)^2 +2079\zeta(7) +288\zeta(2)\zeta(5) -1836\zeta(3)\zeta(4) +595\zeta(8)  \right. \nonumber \\ &\left. \hspace{1em}
+120\zeta(2)\zeta(3)^2 -576\zeta(3)\zeta(5) -264 M(2,6) +1042\zeta(9) -582\zeta(3)\zeta(6)  \right. \nonumber \\ &\left. \hspace{1em}
-612\zeta(4)\zeta(5) +144\zeta(2)\zeta(7) +96\zeta(3)^3\right) \label{eq11124} \\
\sum_{k=1}^\infty \frac{H(k)^{3}}{k^{5}(k+1)^{3}} \sumend &= \frac{-1}{96}\left( -14400\zeta(4) +13200\zeta(5) +1440\zeta(2)\zeta(3) -3546\zeta(6)  \right. \nonumber \\ &\left. \hspace{1em}
+1632\zeta(3)^2 +4158\zeta(7) +576\zeta(2)\zeta(5) -3672\zeta(3)\zeta(4) +595\zeta(8)  \right. \nonumber \\ &\left. \hspace{1em}
+120\zeta(2)\zeta(3)^2 -576\zeta(3)\zeta(5) -264 M(2,6)\right) \label{eq11125} \\
\sum_{k=1}^\infty \frac{H(k)^{3}}{k^{4}(k+1)^{4}} \sumend &= \frac{1}{8}\left( -1600\zeta(4) +1400\zeta(5) +160\zeta(2)\zeta(3) -252\zeta(6)  \right. \nonumber \\ &\left. \hspace{1em}
+144\zeta(3)^2 +175\zeta(7) +32\zeta(2)\zeta(5) -168\zeta(3)\zeta(4)\right) \label{eq11126} \\
\sum_{k=1}^\infty \frac{H(k)^{3}}{k^{3}(k+1)^{5}} \sumend &= \frac{1}{96}\left(  14400\zeta(4) -12000\zeta(5) -1440\zeta(2)\zeta(3) +1746\zeta(6)  \right. \nonumber \\ &\left. \hspace{1em}
-1392\zeta(3)^2 -2142\zeta(7) -576\zeta(2)\zeta(5) +2376\zeta(3)\zeta(4) +43\zeta(8)  \right. \nonumber \\ &\left. \hspace{1em}
+120\zeta(2)\zeta(3)^2 -288\zeta(3)\zeta(5) +24 M(2,6)\right) \label{eq11127}
\end{align}
 
\begin{align}
\sum_{k=1}^\infty \frac{H(k)^{3}}{k^{2}(k+1)^{6}} \sumend &= \frac{-1}{48}\left(  2880\zeta(4) -2280\zeta(5) -288\zeta(2)\zeta(3) +396\zeta(6)  \right. \nonumber \\ &\left. \hspace{1em}
-384\zeta(3)^2 -1071\zeta(7) -288\zeta(2)\zeta(5) +1188\zeta(3)\zeta(4) +43\zeta(8)  \right. \nonumber \\ &\left. \hspace{1em}
+120\zeta(2)\zeta(3)^2 -288\zeta(3)\zeta(5) +24 M(2,6) -394\zeta(9) +222\zeta(3)\zeta(6)  \right. \nonumber \\ &\left. \hspace{1em}
+396\zeta(4)\zeta(5) -144\zeta(2)\zeta(7) -48\zeta(3)^3\right) \label{eq11128} \\
\sum_{k=1}^\infty \frac{H(k)^{3}}{k(k+1)^{7}} \sumend &= \frac{-1}{480}\left( -4800\zeta(4) +3600\zeta(5) +480\zeta(2)\zeta(3) -990\zeta(6)  \right. \nonumber \\ &\left. \hspace{1em}
+960\zeta(3)^2 +3570\zeta(7) +960\zeta(2)\zeta(5) -3960\zeta(3)\zeta(4) -215\zeta(8)  \right. \nonumber \\ &\left. \hspace{1em}
-600\zeta(2)\zeta(3)^2 +1440\zeta(3)\zeta(5) -120 M(2,6) +3940\zeta(9) -2220\zeta(3)\zeta(6)  \right. \nonumber \\ &\left. \hspace{1em}
-3960\zeta(4)\zeta(5) +1440\zeta(2)\zeta(7) +480\zeta(3)^3 -1503\zeta(10) +2400\zeta(3)\zeta(7)  \right. \nonumber \\ &\left. \hspace{1em}
+240\zeta(3)^2\zeta(4) -1680\zeta(2)\zeta(3)\zeta(5) +1440\zeta(5)^2 +120 M(2,8)\right) \label{eq11129} \\
\sum_{k=1}^\infty \frac{H(k)^{3}}{(k+1)^{8}} \sumend &= \frac{1}{2}\left( -44\zeta(11) +21\zeta(3)\zeta(8) +9\zeta(4)\zeta(7) +15\zeta(5)\zeta(6)  \right. \nonumber \\ &\left. \hspace{1em}
-6\zeta(3)^2\zeta(5) +2 M(3,8)\right) \label{eq11130} \\
\sum_{k=1}^\infty \frac{H(k)^{3}}{k^{7}(k+2)} \sumend &= \frac{-1}{7680}\left( -60 -120\zeta(2) -240\zeta(3) -600\zeta(4) +1200\zeta(5)  \right. \nonumber \\ &\left. \hspace{1em}
+120\zeta(2)\zeta(3) -1395\zeta(6) +600\zeta(3)^2 +6930\zeta(7) +960\zeta(2)\zeta(5)  \right. \nonumber \\ &\left. \hspace{1em}
-6120\zeta(3)\zeta(4) +5950\zeta(8) +1200\zeta(2)\zeta(3)^2 -5760\zeta(3)\zeta(5) -2640 M(2,6)  \right. \nonumber \\ &\left. \hspace{1em}
+41680\zeta(9) -23280\zeta(3)\zeta(6) -24480\zeta(4)\zeta(5) +5760\zeta(2)\zeta(7) +3840\zeta(3)^3  \right. \nonumber \\ &\left. \hspace{1em}
+39864\zeta(10) -30720\zeta(3)\zeta(7) -1920\zeta(3)^2\zeta(4) +13440\zeta(2)\zeta(3)\zeta(5)  \right. \nonumber \\ &\left. \hspace{1em}
-17280\zeta(5)^2 -12480 M(2,8)\right) \label{eq11131} \\
\sum_{k=1}^\infty \frac{H(k)^{3}}{k^{6}(k+1)(k+2)} \sumend &= \frac{1}{768}\left(  12 +24\zeta(2) +48\zeta(3) -7560\zeta(4) +7440\zeta(5)  \right. \nonumber \\ &\left. \hspace{1em}
+744\zeta(2)\zeta(3) -4185\zeta(6) +1800\zeta(3)^2 +9702\zeta(7) +1344\zeta(2)\zeta(5)  \right. \nonumber \\ &\left. \hspace{1em}
-8568\zeta(3)\zeta(4) +3570\zeta(8) +720\zeta(2)\zeta(3)^2 -3456\zeta(3)\zeta(5) -1584 M(2,6)  \right. \nonumber \\ &\left. \hspace{1em}
+8336\zeta(9) -4656\zeta(3)\zeta(6) -4896\zeta(4)\zeta(5) +1152\zeta(2)\zeta(7) +768\zeta(3)^3\right) \label{eq11132} \\
\sum_{k=1}^\infty \frac{H(k)^{3}}{k^{5}(k+1)^{2}(k+2)} \sumend &= \frac{-1}{384}\left( -12 -24\zeta(2) -48\zeta(3) -15480\zeta(4) +14640\zeta(5)  \right. \nonumber \\ &\left. \hspace{1em}
+1560\zeta(2)\zeta(3) -4743\zeta(6) +2040\zeta(3)^2 +6930\zeta(7) +960\zeta(2)\zeta(5)  \right. \nonumber \\ &\left. \hspace{1em}
-6120\zeta(3)\zeta(4) +1190\zeta(8) +240\zeta(2)\zeta(3)^2 -1152\zeta(3)\zeta(5) -528 M(2,6)\right) \label{eq11133} \\
\sum_{k=1}^\infty \frac{H(k)^{3}}{k^{4}(k+1)^{3}(k+2)} \sumend &= \frac{-1}{64}\left( -4 -8\zeta(2) -16\zeta(3) +4440\zeta(4) -3920\zeta(5)  \right. \nonumber \\ &\left. \hspace{1em}
-440\zeta(2)\zeta(3) +783\zeta(6) -408\zeta(3)^2 -462\zeta(7) -64\zeta(2)\zeta(5)  \right. \nonumber \\ &\left. \hspace{1em}
+408\zeta(3)\zeta(4)\right) \label{eq11134}
\end{align}
 
\begin{align}
\sum_{k=1}^\infty \frac{H(k)^{3}}{k^{3}(k+1)^{4}(k+2)} \sumend &= \frac{1}{32}\left(  4 +8\zeta(2) +16\zeta(3) +1960\zeta(4) -1680\zeta(5)  \right. \nonumber \\ &\left. \hspace{1em}
-200\zeta(2)\zeta(3) +225\zeta(6) -168\zeta(3)^2 -238\zeta(7) -64\zeta(2)\zeta(5)  \right. \nonumber \\ &\left. \hspace{1em}
+264\zeta(3)\zeta(4)\right) \label{eq11135} \\
\sum_{k=1}^\infty \frac{H(k)^{3}}{k^{2}(k+1)^{5}(k+2)} \sumend &= \frac{-1}{96}\left( -24 -48\zeta(2) -96\zeta(3) +2640\zeta(4) -1920\zeta(5)  \right. \nonumber \\ &\left. \hspace{1em}
-240\zeta(2)\zeta(3) +396\zeta(6) -384\zeta(3)^2 -714\zeta(7) -192\zeta(2)\zeta(5) +792\zeta(3)\zeta(4)  \right. \nonumber \\ &\left. \hspace{1em}
+43\zeta(8) +120\zeta(2)\zeta(3)^2 -288\zeta(3)\zeta(5) +24 M(2,6)\right) \label{eq11136} \\
\sum_{k=1}^\infty \frac{H(k)^{3}}{k(k+1)^{6}(k+2)} \sumend &= \frac{1}{48}\left(  24 +48\zeta(2) +96\zeta(3) +240\zeta(4) -360\zeta(5)  \right. \nonumber \\ &\left. \hspace{1em}
-48\zeta(2)\zeta(3) -357\zeta(7) -96\zeta(2)\zeta(5) +396\zeta(3)\zeta(4) -394\zeta(9)  \right. \nonumber \\ &\left. \hspace{1em}
+222\zeta(3)\zeta(6) +396\zeta(4)\zeta(5) -144\zeta(2)\zeta(7) -48\zeta(3)^3\right) \label{eq11137} \\
\sum_{k=1}^\infty \frac{H(k)^{3}}{(k+1)^{7}(k+2)} \sumend &= \frac{-1}{480}\left( -480 -960\zeta(2) -1920\zeta(3) +3600\zeta(5) +480\zeta(2)\zeta(3)  \right. \nonumber \\ &\left. \hspace{1em}
+990\zeta(6) -960\zeta(3)^2 +3570\zeta(7) +960\zeta(2)\zeta(5) -3960\zeta(3)\zeta(4) +215\zeta(8)  \right. \nonumber \\ &\left. \hspace{1em}
+600\zeta(2)\zeta(3)^2 -1440\zeta(3)\zeta(5) +120 M(2,6) +3940\zeta(9) -2220\zeta(3)\zeta(6)  \right. \nonumber \\ &\left. \hspace{1em}
-3960\zeta(4)\zeta(5) +1440\zeta(2)\zeta(7) +480\zeta(3)^3 +1503\zeta(10) -2400\zeta(3)\zeta(7)  \right. \nonumber \\ &\left. \hspace{1em}
-240\zeta(3)^2\zeta(4) +1680\zeta(2)\zeta(3)\zeta(5) -1440\zeta(5)^2 -120 M(2,8)\right) \label{eq11138} \\
\sum_{k=1}^\infty \frac{H(k)^{3}}{k^{6}(k+2)^{2}} \sumend &= \frac{-1}{768}\left(  84 +108\zeta(2) +156\zeta(3) +261\zeta(4) -690\zeta(5)  \right. \nonumber \\ &\left. \hspace{1em}
-72\zeta(2)\zeta(3) +558\zeta(6) -240\zeta(3)^2 -2079\zeta(7) -288\zeta(2)\zeta(5)  \right. \nonumber \\ &\left. \hspace{1em}
+1836\zeta(3)\zeta(4) -1190\zeta(8) -240\zeta(2)\zeta(3)^2 +1152\zeta(3)\zeta(5) +528 M(2,6)  \right. \nonumber \\ &\left. \hspace{1em}
-4168\zeta(9) +2328\zeta(3)\zeta(6) +2448\zeta(4)\zeta(5) -576\zeta(2)\zeta(7) -384\zeta(3)^3\right) \label{eq11139} \\
\sum_{k=1}^\infty \frac{H(k)^{3}}{k^{5}(k+1)(k+2)^{2}} \sumend &= \frac{-1}{768}\left(  180 +240\zeta(2) +360\zeta(3) -7038\zeta(4) +6060\zeta(5)  \right. \nonumber \\ &\left. \hspace{1em}
+600\zeta(2)\zeta(3) -3069\zeta(6) +1320\zeta(3)^2 +5544\zeta(7) +768\zeta(2)\zeta(5)  \right. \nonumber \\ &\left. \hspace{1em}
-4896\zeta(3)\zeta(4) +1190\zeta(8) +240\zeta(2)\zeta(3)^2 -1152\zeta(3)\zeta(5) -528 M(2,6)\right) \label{eq11140} \\
\sum_{k=1}^\infty \frac{H(k)^{3}}{k^{4}(k+1)^{2}(k+2)^{2}} \sumend &= \frac{-1}{64}\left(  32 +44\zeta(2) +68\zeta(3) +1407\zeta(4) -1430\zeta(5)  \right. \nonumber \\ &\left. \hspace{1em}
-160\zeta(2)\zeta(3) +279\zeta(6) -120\zeta(3)^2 -231\zeta(7) -32\zeta(2)\zeta(5)  \right. \nonumber \\ &\left. \hspace{1em}
+204\zeta(3)\zeta(4)\right) \label{eq11141} \\
\sum_{k=1}^\infty \frac{H(k)^{3}}{k^{3}(k+1)^{3}(k+2)^{2}} \sumend &= \frac{1}{64}\left( -68 -96\zeta(2) -152\zeta(3) +1626\zeta(4) -1060\zeta(5)  \right. \nonumber \\ &\left. \hspace{1em}
-120\zeta(2)\zeta(3) +225\zeta(6) -168\zeta(3)^2\right) \label{eq11142}
\end{align}
 
\begin{align}
\sum_{k=1}^\infty \frac{H(k)^{3}}{k^{2}(k+1)^{4}(k+2)^{2}} \sumend &= \frac{1}{16}\left( -36 -52\zeta(2) -84\zeta(3) -167\zeta(4) +310\zeta(5)  \right. \nonumber \\ &\left. \hspace{1em}
+40\zeta(2)\zeta(3) +119\zeta(7) +32\zeta(2)\zeta(5) -132\zeta(3)\zeta(4)\right) \label{eq11143} \\
\sum_{k=1}^\infty \frac{H(k)^{3}}{k(k+1)^{5}(k+2)^{2}} \sumend &= \frac{1}{96}\left( -456 -672\zeta(2) -1104\zeta(3) +636\zeta(4) +1800\zeta(5)  \right. \nonumber \\ &\left. \hspace{1em}
+240\zeta(2)\zeta(3) +396\zeta(6) -384\zeta(3)^2 +714\zeta(7) +192\zeta(2)\zeta(5) -792\zeta(3)\zeta(4)  \right. \nonumber \\ &\left. \hspace{1em}
+43\zeta(8) +120\zeta(2)\zeta(3)^2 -288\zeta(3)\zeta(5) +24 M(2,6)\right) \label{eq11144} \\
\sum_{k=1}^\infty \frac{H(k)^{3}}{(k+1)^{6}(k+2)^{2}} \sumend &= \frac{1}{48}\left( -480 -720\zeta(2) -1200\zeta(3) +396\zeta(4) +2160\zeta(5)  \right. \nonumber \\ &\left. \hspace{1em}
+288\zeta(2)\zeta(3) +396\zeta(6) -384\zeta(3)^2 +1071\zeta(7) +288\zeta(2)\zeta(5)  \right. \nonumber \\ &\left. \hspace{1em}
-1188\zeta(3)\zeta(4) +43\zeta(8) +120\zeta(2)\zeta(3)^2 -288\zeta(3)\zeta(5) +24 M(2,6) +394\zeta(9)  \right. \nonumber \\ &\left. \hspace{1em}
-222\zeta(3)\zeta(6) -396\zeta(4)\zeta(5) +144\zeta(2)\zeta(7) +48\zeta(3)^3\right) \label{eq11145} \\
\sum_{k=1}^\infty \frac{H(k)^{3}}{k^{5}(k+2)^{3}} \sumend &= \frac{1}{1536}\left(  1140 +864\zeta(2) +696\zeta(3) +378\zeta(4) -3084\zeta(5)  \right. \nonumber \\ &\left. \hspace{1em}
-504\zeta(2)\zeta(3) +1773\zeta(6) -816\zeta(3)^2 -4158\zeta(7) -576\zeta(2)\zeta(5)  \right. \nonumber \\ &\left. \hspace{1em}
+3672\zeta(3)\zeta(4) -1190\zeta(8) -240\zeta(2)\zeta(3)^2 +1152\zeta(3)\zeta(5) +528 M(2,6)\right) \label{eq11146} \\
\sum_{k=1}^\infty \frac{H(k)^{3}}{k^{4}(k+1)(k+2)^{3}} \sumend &= \frac{1}{128}\left(  220 +184\zeta(2) +176\zeta(3) -1110\zeta(4) +496\zeta(5)  \right. \nonumber \\ &\left. \hspace{1em}
+16\zeta(2)\zeta(3) -216\zeta(6) +84\zeta(3)^2 +231\zeta(7) +32\zeta(2)\zeta(5)  \right. \nonumber \\ &\left. \hspace{1em}
-204\zeta(3)\zeta(4)\right) \label{eq11147} \\
\sum_{k=1}^\infty \frac{H(k)^{3}}{k^{3}(k+1)^{2}(k+2)^{3}} \sumend &= \frac{-1}{64}\left( -252 -228\zeta(2) -244\zeta(3) -297\zeta(4) +934\zeta(5)  \right. \nonumber \\ &\left. \hspace{1em}
+144\zeta(2)\zeta(3) -63\zeta(6) +36\zeta(3)^2\right) \label{eq11148} \\
\sum_{k=1}^\infty \frac{H(k)^{3}}{k^{2}(k+1)^{3}(k+2)^{3}} \sumend &= \frac{1}{64}\left(  572 +552\zeta(2) +640\zeta(3) -1032\zeta(4) -808\zeta(5)  \right. \nonumber \\ &\left. \hspace{1em}
-168\zeta(2)\zeta(3) -99\zeta(6) +96\zeta(3)^2\right) \label{eq11149} \\
\sum_{k=1}^\infty \frac{H(k)^{3}}{k(k+1)^{4}(k+2)^{3}} \sumend &= \frac{-1}{32}\left( -644 -656\zeta(2) -808\zeta(3) +698\zeta(4) +1428\zeta(5)  \right. \nonumber \\ &\left. \hspace{1em}
+248\zeta(2)\zeta(3) +99\zeta(6) -96\zeta(3)^2 +238\zeta(7) +64\zeta(2)\zeta(5)  \right. \nonumber \\ &\left. \hspace{1em}
-264\zeta(3)\zeta(4)\right) \label{eq11150} \\
\sum_{k=1}^\infty \frac{H(k)^{3}}{(k+1)^{5}(k+2)^{3}} \sumend &= \frac{1}{96}\left(  4320 +4608\zeta(2) +5952\zeta(3) -4824\zeta(4) -10368\zeta(5)  \right. \nonumber \\ &\left. \hspace{1em}
-1728\zeta(2)\zeta(3) -990\zeta(6) +960\zeta(3)^2 -2142\zeta(7) -576\zeta(2)\zeta(5)  \right. \nonumber \\ &\left. \hspace{1em}
+2376\zeta(3)\zeta(4) -43\zeta(8) -120\zeta(2)\zeta(3)^2 +288\zeta(3)\zeta(5) -24 M(2,6)\right) \label{eq11151}
\end{align}
 
\begin{align}
\sum_{k=1}^\infty \frac{H(k)^{3}}{k^{4}(k+2)^{4}} \sumend &= \frac{1}{128}\left( -420 -164\zeta(2) +12\zeta(3) +195\zeta(4) +290\zeta(5)  \right. \nonumber \\ &\left. \hspace{1em}
+88\zeta(2)\zeta(3) -89\zeta(6) +48\zeta(3)^2 +175\zeta(7) +32\zeta(2)\zeta(5)  \right. \nonumber \\ &\left. \hspace{1em}
-168\zeta(3)\zeta(4)\right) \label{eq11152} \\
\sum_{k=1}^\infty \frac{H(k)^{3}}{k^{3}(k+1)(k+2)^{4}} \sumend &= \frac{-1}{128}\left(  1060 +512\zeta(2) +152\zeta(3) -1500\zeta(4) -84\zeta(5)  \right. \nonumber \\ &\left. \hspace{1em}
-160\zeta(2)\zeta(3) -38\zeta(6) -12\zeta(3)^2 -119\zeta(7) -32\zeta(2)\zeta(5)  \right. \nonumber \\ &\left. \hspace{1em}
+132\zeta(3)\zeta(4)\right) \label{eq11153} \\
\sum_{k=1}^\infty \frac{H(k)^{3}}{k^{2}(k+1)^{2}(k+2)^{4}} \sumend &= \frac{1}{64}\left( -1312 -740\zeta(2) -396\zeta(3) +1203\zeta(4) +1018\zeta(5)  \right. \nonumber \\ &\left. \hspace{1em}
+304\zeta(2)\zeta(3) -25\zeta(6) +48\zeta(3)^2 +119\zeta(7) +32\zeta(2)\zeta(5)  \right. \nonumber \\ &\left. \hspace{1em}
-132\zeta(3)\zeta(4)\right) \label{eq11154} \\
\sum_{k=1}^\infty \frac{H(k)^{3}}{k(k+1)^{3}(k+2)^{4}} \sumend &= \frac{1}{64}\left( -3196 -2032\zeta(2) -1432\zeta(3) +3438\zeta(4) +2844\zeta(5)  \right. \nonumber \\ &\left. \hspace{1em}
+776\zeta(2)\zeta(3) +49\zeta(6) +238\zeta(7) +64\zeta(2)\zeta(5) -264\zeta(3)\zeta(4)\right) \label{eq11155} \\
\sum_{k=1}^\infty \frac{H(k)^{3}}{(k+1)^{4}(k+2)^{4}} \sumend &= \frac{1}{8}\left( -960 -672\zeta(2) -560\zeta(3) +1034\zeta(4) +1068\zeta(5)  \right. \nonumber \\ &\left. \hspace{1em}
+256\zeta(2)\zeta(3) +37\zeta(6) -24\zeta(3)^2 +119\zeta(7) +32\zeta(2)\zeta(5)  \right. \nonumber \\ &\left. \hspace{1em}
-132\zeta(3)\zeta(4)\right) \label{eq11156} \\
\sum_{k=1}^\infty \frac{H(k)^{3}}{k^{3}(k+2)^{5}} \sumend &= \frac{-1}{1536}\left( -16500 -2520\zeta(2) +4896\zeta(3) +8460\zeta(4) +3768\zeta(5)  \right. \nonumber \\ &\left. \hspace{1em}
+648\zeta(2)\zeta(3) +1779\zeta(6) -1032\zeta(3)^2 +1566\zeta(7) +1152\zeta(2)\zeta(5)  \right. \nonumber \\ &\left. \hspace{1em}
-2664\zeta(3)\zeta(4) -86\zeta(8) -240\zeta(2)\zeta(3)^2 +576\zeta(3)\zeta(5) -48 M(2,6)\right) \label{eq11157} \\
\sum_{k=1}^\infty \frac{H(k)^{3}}{k^{2}(k+1)(k+2)^{5}} \sumend &= \frac{1}{768}\left(  22860 +5592\zeta(2) -3984\zeta(3) -17460\zeta(4) -4272\zeta(5)  \right. \nonumber \\ &\left. \hspace{1em}
-1608\zeta(2)\zeta(3) -2007\zeta(6) +960\zeta(3)^2 -2280\zeta(7) -1344\zeta(2)\zeta(5)  \right. \nonumber \\ &\left. \hspace{1em}
+3456\zeta(3)\zeta(4) +86\zeta(8) +240\zeta(2)\zeta(3)^2 -576\zeta(3)\zeta(5) +48 M(2,6)\right) \label{eq11158} \\
\sum_{k=1}^\infty \frac{H(k)^{3}}{k(k+1)^{2}(k+2)^{5}} \sumend &= \frac{-1}{384}\left( -30732 -10032\zeta(2) +1608\zeta(3) +24678\zeta(4) +10380\zeta(5)  \right. \nonumber \\ &\left. \hspace{1em}
+3432\zeta(2)\zeta(3) +1857\zeta(6) -672\zeta(3)^2 +2994\zeta(7) +1536\zeta(2)\zeta(5)  \right. \nonumber \\ &\left. \hspace{1em}
-4248\zeta(3)\zeta(4) -86\zeta(8) -240\zeta(2)\zeta(3)^2 +576\zeta(3)\zeta(5) -48 M(2,6)\right) \label{eq11159} \\
\sum_{k=1}^\infty \frac{H(k)^{3}}{(k+1)^{3}(k+2)^{5}} \sumend &= \frac{-1}{96}\left( -20160 -8064\zeta(2) -1344\zeta(3) +17496\zeta(4) +9456\zeta(5)  \right. \nonumber \\ &\left. \hspace{1em}
+2880\zeta(2)\zeta(3) +1002\zeta(6) -336\zeta(3)^2 +1854\zeta(7) +864\zeta(2)\zeta(5)  \right. \nonumber \\ &\left. \hspace{1em}
-2520\zeta(3)\zeta(4) -43\zeta(8) -120\zeta(2)\zeta(3)^2 +288\zeta(3)\zeta(5) -24 M(2,6)\right) \label{eq11160}
\end{align}
 
\begin{align}
\sum_{k=1}^\infty \frac{H(k)^{3}}{k^{2}(k+2)^{6}} \sumend &= \frac{1}{768}\left( -21588 -252\zeta(2) +7956\zeta(3) +8451\zeta(4) +5154\zeta(5)  \right. \nonumber \\ &\left. \hspace{1em}
-1368\zeta(2)\zeta(3) +3732\zeta(6) -2256\zeta(3)^2 +1647\zeta(7) +864\zeta(2)\zeta(5)  \right. \nonumber \\ &\left. \hspace{1em}
-2340\zeta(3)\zeta(4) -2102\zeta(8) -240\zeta(2)\zeta(3)^2 +1728\zeta(3)\zeta(5) +528 M(2,6)  \right. \nonumber \\ &\left. \hspace{1em}
+1576\zeta(9) -888\zeta(3)\zeta(6) -1584\zeta(4)\zeta(5) +576\zeta(2)\zeta(7) +192\zeta(3)^3\right) \label{eq11161} \\
\sum_{k=1}^\infty \frac{H(k)^{3}}{k(k+1)(k+2)^{6}} \sumend &= \frac{-1}{768}\left(  66036 +6096\zeta(2) -19896\zeta(3) -34362\zeta(4) -14580\zeta(5)  \right. \nonumber \\ &\left. \hspace{1em}
+1128\zeta(2)\zeta(3) -9471\zeta(6) +5472\zeta(3)^2 -5574\zeta(7) -3072\zeta(2)\zeta(5)  \right. \nonumber \\ &\left. \hspace{1em}
+8136\zeta(3)\zeta(4) +4290\zeta(8) +720\zeta(2)\zeta(3)^2 -4032\zeta(3)\zeta(5) -1008 M(2,6)  \right. \nonumber \\ &\left. \hspace{1em}
-3152\zeta(9) +1776\zeta(3)\zeta(6) +3168\zeta(4)\zeta(5) -1152\zeta(2)\zeta(7) -384\zeta(3)^3\right) \label{eq11162} \\
\sum_{k=1}^\infty \frac{H(k)^{3}}{(k+1)^{2}(k+2)^{6}} \sumend &= \frac{-1}{48}\left(  12096 +2016\zeta(2) -2688\zeta(3) -7380\zeta(4) -3120\zeta(5)  \right. \nonumber \\ &\left. \hspace{1em}
-288\zeta(2)\zeta(3) -1416\zeta(6) +768\zeta(3)^2 -1071\zeta(7) -576\zeta(2)\zeta(5)  \right. \nonumber \\ &\left. \hspace{1em}
+1548\zeta(3)\zeta(4) +547\zeta(8) +120\zeta(2)\zeta(3)^2 -576\zeta(3)\zeta(5) -120 M(2,6) -394\zeta(9)  \right. \nonumber \\ &\left. \hspace{1em}
+222\zeta(3)\zeta(6) +396\zeta(4)\zeta(5) -144\zeta(2)\zeta(7) -48\zeta(3)^3\right) \label{eq11163} \\
\sum_{k=1}^\infty \frac{H(k)^{3}}{k(k+2)^{7}} \sumend &= \frac{1}{7680}\left(  476220 -30480\zeta(2) -169320\zeta(3) -138270\zeta(4) -138300\zeta(5)  \right. \nonumber \\ &\left. \hspace{1em}
+48120\zeta(2)\zeta(3) -82965\zeta(6) +45600\zeta(3)^2 -73650\zeta(7) +7680\zeta(2)\zeta(5)  \right. \nonumber \\ &\left. \hspace{1em}
+42840\zeta(3)\zeta(4) +46510\zeta(8) +1200\zeta(2)\zeta(3)^2 -25920\zeta(3)\zeta(5) -17040 M(2,6)  \right. \nonumber \\ &\left. \hspace{1em}
-17680\zeta(9) +26160\zeta(3)\zeta(6) +21600\zeta(4)\zeta(5) -17280\zeta(2)\zeta(7) -5760\zeta(3)^3  \right. \nonumber \\ &\left. \hspace{1em}
+12024\zeta(10) -19200\zeta(3)\zeta(7) -1920\zeta(3)^2\zeta(4) +13440\zeta(2)\zeta(3)\zeta(5)  \right. \nonumber \\ &\left. \hspace{1em}
-11520\zeta(5)^2 -960 M(2,8)\right) \label{eq11164} \\
\sum_{k=1}^\infty \frac{H(k)^{3}}{(k+1)(k+2)^{7}} \sumend &= \frac{1}{480}\left(  100800 -33600\zeta(3) -38760\zeta(4) -26400\zeta(5)  \right. \nonumber \\ &\left. \hspace{1em}
+6720\zeta(2)\zeta(3) -16290\zeta(6) +9120\zeta(3)^2 -12690\zeta(7) -960\zeta(2)\zeta(5)  \right. \nonumber \\ &\left. \hspace{1em}
+10440\zeta(3)\zeta(4) +8495\zeta(8) +600\zeta(2)\zeta(3)^2 -5760\zeta(3)\zeta(5) -2760 M(2,6)  \right. \nonumber \\ &\left. \hspace{1em}
-4180\zeta(9) +4380\zeta(3)\zeta(6) +4680\zeta(4)\zeta(5) -2880\zeta(2)\zeta(7) -960\zeta(3)^3  \right. \nonumber \\ &\left. \hspace{1em}
+1503\zeta(10) -2400\zeta(3)\zeta(7) -240\zeta(3)^2\zeta(4) +1680\zeta(2)\zeta(3)\zeta(5)  \right. \nonumber \\ &\left. \hspace{1em}
-1440\zeta(5)^2 -120 M(2,8)\right) \label{eq11165} \\
\sum_{k=1}^\infty \frac{H(k)^{3}}{(k+2)^{8}} \sumend &= \frac{1}{8}\left( -960 +96\zeta(2) +304\zeta(3) +222\zeta(4) +284\zeta(5) -96\zeta(2)\zeta(3)  \right. \nonumber \\ &\left. \hspace{1em}
+157\zeta(6) -72\zeta(3)^2 +216\zeta(7) -48\zeta(2)\zeta(5) -84\zeta(3)\zeta(4) -16\zeta(8) +24 M(2,6)  \right. \nonumber \\ &\left. \hspace{1em}
+100\zeta(9) -60\zeta(3)\zeta(6) -36\zeta(4)\zeta(5) +8\zeta(3)^3 -108\zeta(10) +48\zeta(3)\zeta(7)  \right. \nonumber \\ &\left. \hspace{1em}
+24\zeta(5)^2 +24 M(2,8) -176\zeta(11) +84\zeta(3)\zeta(8) +36\zeta(4)\zeta(7) +60\zeta(5)\zeta(6)  \right. \nonumber \\ &\left. \hspace{1em}
-24\zeta(3)^2\zeta(5) +8 M(3,8)\right) \label{eq11166}
\end{align}
 
\begin{align}
\sum_{k=1}^\infty \frac{H(k)^{4}}{k^{7}} \sumend &= \frac{-1}{48}\left(  2877\zeta(11) +272\zeta(2)\zeta(9) -1190\zeta(3)\zeta(8) -1212\zeta(4)\zeta(7)  \right. \nonumber \\ &\left. \hspace{1em}
-1018\zeta(5)\zeta(6) -80\zeta(2)\zeta(3)^3 +576\zeta(3)^2\zeta(5) -176 M(3,8)\right) \label{eq11167} \\
\sum_{k=1}^\infty \frac{H(k)^{4}}{k^{6}(k+1)} \sumend &= \frac{1}{5760}\left( -172800\zeta(5) -34560\zeta(2)\zeta(3) +234960\zeta(6) +17280\zeta(3)^2  \right. \nonumber \\ &\left. \hspace{1em}
-133200\zeta(7) -28800\zeta(2)\zeta(5) +123840\zeta(3)\zeta(4) -593320\zeta(8)  \right. \nonumber \\ &\left. \hspace{1em}
-161280\zeta(2)\zeta(3)^2 +668160\zeta(3)\zeta(5) +149760 M(2,6) -209280\zeta(9)  \right. \nonumber \\ &\left. \hspace{1em}
+133920\zeta(3)\zeta(6) +123840\zeta(4)\zeta(5) -40320\zeta(2)\zeta(7) -19200\zeta(3)^3  \right. \nonumber \\ &\left. \hspace{1em}
-619407\zeta(10) +540000\zeta(3)\zeta(7) +9000\zeta(3)^2\zeta(4) -195120\zeta(2)\zeta(3)\zeta(5)  \right. \nonumber \\ &\left. \hspace{1em}
+212040\zeta(5)^2 +109080 M(2,8) +11520\zeta(2) M(2,6)\right) \label{eq11168} \\
\sum_{k=1}^\infty \frac{H(k)^{4}}{k^{5}(k+1)^{2}} \sumend &= \frac{1}{72}\left(  10800\zeta(5) +2160\zeta(2)\zeta(3) -14325\zeta(6) -1080\zeta(3)^2  \right. \nonumber \\ &\left. \hspace{1em}
+4995\zeta(7) +1080\zeta(2)\zeta(5) -4644\zeta(3)\zeta(4) +14833\zeta(8) +4032\zeta(2)\zeta(3)^2  \right. \nonumber \\ &\left. \hspace{1em}
-16704\zeta(3)\zeta(5) -3744 M(2,6) +2616\zeta(9) -1674\zeta(3)\zeta(6) -1548\zeta(4)\zeta(5)  \right. \nonumber \\ &\left. \hspace{1em}
+504\zeta(2)\zeta(7) +240\zeta(3)^3\right) \label{eq11169} \\
\sum_{k=1}^\infty \frac{H(k)^{4}}{k^{4}(k+1)^{3}} \sumend &= \frac{1}{144}\left( -43200\zeta(5) -8640\zeta(2)\zeta(3) +55860\zeta(6) +4320\zeta(3)^2  \right. \nonumber \\ &\left. \hspace{1em}
-11952\zeta(7) -2880\zeta(2)\zeta(5) +11952\zeta(3)\zeta(4) -14833\zeta(8) -4032\zeta(2)\zeta(3)^2  \right. \nonumber \\ &\left. \hspace{1em}
+16704\zeta(3)\zeta(5) +3744 M(2,6)\right) \label{eq11170} \\
\sum_{k=1}^\infty \frac{H(k)^{4}}{k^{3}(k+1)^{4}} \sumend &= \frac{-1}{144}\left( -43200\zeta(5) -8640\zeta(2)\zeta(3) +54420\zeta(6) +4320\zeta(3)^2  \right. \nonumber \\ &\left. \hspace{1em}
-9216\zeta(7) -2880\zeta(2)\zeta(5) +11088\zeta(3)\zeta(4) -12415\zeta(8) -3312\zeta(2)\zeta(3)^2  \right. \nonumber \\ &\left. \hspace{1em}
+13824\zeta(3)\zeta(5) +3024 M(2,6)\right) \label{eq11171} \\
\sum_{k=1}^\infty \frac{H(k)^{4}}{k^{2}(k+1)^{5}} \sumend &= \frac{1}{72}\left( -10800\zeta(5) -2160\zeta(2)\zeta(3) +13245\zeta(6) +1080\zeta(3)^2  \right. \nonumber \\ &\left. \hspace{1em}
-2943\zeta(7) -1080\zeta(2)\zeta(5) +3996\zeta(3)\zeta(4) -12415\zeta(8) -3312\zeta(2)\zeta(3)^2  \right. \nonumber \\ &\left. \hspace{1em}
+13824\zeta(3)\zeta(5) +3024 M(2,6) -1044\zeta(9) +594\zeta(3)\zeta(6) +1332\zeta(4)\zeta(5)  \right. \nonumber \\ &\left. \hspace{1em}
-504\zeta(2)\zeta(7) -192\zeta(3)^3\right) \label{eq11172} \\
\sum_{k=1}^\infty \frac{H(k)^{4}}{k(k+1)^{6}} \sumend &= \frac{-1}{5760}\left( -172800\zeta(5) -34560\zeta(2)\zeta(3) +206160\zeta(6) +17280\zeta(3)^2  \right. \nonumber \\ &\left. \hspace{1em}
-78480\zeta(7) -28800\zeta(2)\zeta(5) +106560\zeta(3)\zeta(4) -496600\zeta(8) -132480\zeta(2)\zeta(3)^2  \right. \nonumber \\ &\left. \hspace{1em}
+552960\zeta(3)\zeta(5) +120960 M(2,6) -83520\zeta(9) +47520\zeta(3)\zeta(6) +106560\zeta(4)\zeta(5)  \right. \nonumber \\ &\left. \hspace{1em}
-40320\zeta(2)\zeta(7) -15360\zeta(3)^3 -437823\zeta(10) +378720\zeta(3)\zeta(7)  \right. \nonumber \\ &\left. \hspace{1em}
-2520\zeta(3)^2\zeta(4) -114480\zeta(2)\zeta(3)\zeta(5) +119880\zeta(5)^2 +68760 M(2,8)  \right. \nonumber \\ &\left. \hspace{1em}
+11520\zeta(2) M(2,6)\right) \label{eq11173}
\end{align}
 
\begin{align}
\sum_{k=1}^\infty \frac{H(k)^{4}}{(k+1)^{7}} \sumend &= \frac{1}{48}\left( -237\zeta(11) -368\zeta(2)\zeta(9) +86\zeta(3)\zeta(8) +684\zeta(4)\zeta(7)  \right. \nonumber \\ &\left. \hspace{1em}
+202\zeta(5)\zeta(6) +80\zeta(2)\zeta(3)^3 -288\zeta(3)^2\zeta(5) -16 M(3,8)\right) \label{eq11174} \\
\sum_{k=1}^\infty \frac{H(k)^{4}}{k^{6}(k+2)} \sumend &= \frac{1}{11520}\left( -180 -540\zeta(2) -1980\zeta(3) -3330\zeta(4) -5400\zeta(5)  \right. \nonumber \\ &\left. \hspace{1em}
-1080\zeta(2)\zeta(3) +14685\zeta(6) +1080\zeta(3)^2 -16650\zeta(7) -3600\zeta(2)\zeta(5)  \right. \nonumber \\ &\left. \hspace{1em}
+15480\zeta(3)\zeta(4) -148330\zeta(8) -40320\zeta(2)\zeta(3)^2 +167040\zeta(3)\zeta(5) +37440 M(2,6)  \right. \nonumber \\ &\left. \hspace{1em}
-104640\zeta(9) +66960\zeta(3)\zeta(6) +61920\zeta(4)\zeta(5) -20160\zeta(2)\zeta(7) -9600\zeta(3)^3  \right. \nonumber \\ &\left. \hspace{1em}
-619407\zeta(10) +540000\zeta(3)\zeta(7) +9000\zeta(3)^2\zeta(4) -195120\zeta(2)\zeta(3)\zeta(5)  \right. \nonumber \\ &\left. \hspace{1em}
+212040\zeta(5)^2 +109080 M(2,8) +11520\zeta(2) M(2,6)\right) \label{eq11175} \\
\sum_{k=1}^\infty \frac{H(k)^{4}}{k^{5}(k+1)(k+2)} \sumend &= \frac{1}{384}\left( -12 -36\zeta(2) -132\zeta(3) -222\zeta(4) +11160\zeta(5)  \right. \nonumber \\ &\left. \hspace{1em}
+2232\zeta(2)\zeta(3) -14685\zeta(6) -1080\zeta(3)^2 +7770\zeta(7) +1680\zeta(2)\zeta(5)  \right. \nonumber \\ &\left. \hspace{1em}
-7224\zeta(3)\zeta(4) +29666\zeta(8) +8064\zeta(2)\zeta(3)^2 -33408\zeta(3)\zeta(5) -7488 M(2,6)  \right. \nonumber \\ &\left. \hspace{1em}
+6976\zeta(9) -4464\zeta(3)\zeta(6) -4128\zeta(4)\zeta(5) +1344\zeta(2)\zeta(7) +640\zeta(3)^3\right) \label{eq11176} \\
\sum_{k=1}^\infty \frac{H(k)^{4}}{k^{4}(k+1)^{2}(k+2)} \sumend &= \frac{-1}{576}\left(  36 +108\zeta(2) +396\zeta(3) +666\zeta(4) +52920\zeta(5)  \right. \nonumber \\ &\left. \hspace{1em}
+10584\zeta(2)\zeta(3) -70545\zeta(6) -5400\zeta(3)^2 +16650\zeta(7) +3600\zeta(2)\zeta(5)  \right. \nonumber \\ &\left. \hspace{1em}
-15480\zeta(3)\zeta(4) +29666\zeta(8) +8064\zeta(2)\zeta(3)^2 -33408\zeta(3)\zeta(5)  \right. \nonumber \\ &\left. \hspace{1em}
-7488 M(2,6)\right) \label{eq11177} \\
\sum_{k=1}^\infty \frac{H(k)^{4}}{k^{3}(k+1)^{3}(k+2)} \sumend &= \frac{-1}{32}\left(  4 +12\zeta(2) +44\zeta(3) +74\zeta(4) -3720\zeta(5)  \right. \nonumber \\ &\left. \hspace{1em}
-744\zeta(2)\zeta(3) +4575\zeta(6) +360\zeta(3)^2 -806\zeta(7) -240\zeta(2)\zeta(5)  \right. \nonumber \\ &\left. \hspace{1em}
+936\zeta(3)\zeta(4)\right) \label{eq11178} \\
\sum_{k=1}^\infty \frac{H(k)^{4}}{k^{2}(k+1)^{4}(k+2)} \sumend &= \frac{1}{144}\left( -36 -108\zeta(2) -396\zeta(3) -666\zeta(4) -9720\zeta(5)  \right. \nonumber \\ &\left. \hspace{1em}
-1944\zeta(2)\zeta(3) +13245\zeta(6) +1080\zeta(3)^2 -1962\zeta(7) -720\zeta(2)\zeta(5)  \right. \nonumber \\ &\left. \hspace{1em}
+2664\zeta(3)\zeta(4) -12415\zeta(8) -3312\zeta(2)\zeta(3)^2 +13824\zeta(3)\zeta(5)  \right. \nonumber \\ &\left. \hspace{1em}
+3024 M(2,6)\right) \label{eq11179} \\
\sum_{k=1}^\infty \frac{H(k)^{4}}{k(k+1)^{5}(k+2)} \sumend &= \frac{1}{24}\left( -12 -36\zeta(2) -132\zeta(3) -222\zeta(4) +360\zeta(5)  \right. \nonumber \\ &\left. \hspace{1em}
+72\zeta(2)\zeta(3) +327\zeta(7) +120\zeta(2)\zeta(5) -444\zeta(3)\zeta(4) +348\zeta(9)  \right. \nonumber \\ &\left. \hspace{1em}
-198\zeta(3)\zeta(6) -444\zeta(4)\zeta(5) +168\zeta(2)\zeta(7) +64\zeta(3)^3\right) \label{eq11180}
\end{align}
 
\begin{align}
\sum_{k=1}^\infty \frac{H(k)^{4}}{(k+1)^{6}(k+2)} \sumend &= \frac{1}{5760}\left( -5760 -17280\zeta(2) -63360\zeta(3) -106560\zeta(4) +206160\zeta(6)  \right. \nonumber \\ &\left. \hspace{1em}
+17280\zeta(3)^2 +78480\zeta(7) +28800\zeta(2)\zeta(5) -106560\zeta(3)\zeta(4) -496600\zeta(8)  \right. \nonumber \\ &\left. \hspace{1em}
-132480\zeta(2)\zeta(3)^2 +552960\zeta(3)\zeta(5) +120960 M(2,6) +83520\zeta(9) -47520\zeta(3)\zeta(6)  \right. \nonumber \\ &\left. \hspace{1em}
-106560\zeta(4)\zeta(5) +40320\zeta(2)\zeta(7) +15360\zeta(3)^3 -437823\zeta(10)  \right. \nonumber \\ &\left. \hspace{1em}
+378720\zeta(3)\zeta(7) -2520\zeta(3)^2\zeta(4) -114480\zeta(2)\zeta(3)\zeta(5) +119880\zeta(5)^2  \right. \nonumber \\ &\left. \hspace{1em}
+68760 M(2,8) +11520\zeta(2) M(2,6)\right) \label{eq11181} \\
\sum_{k=1}^\infty \frac{H(k)^{4}}{k^{5}(k+2)^{2}} \sumend &= \frac{1}{2304}\left(  540 +1116\zeta(2) +3276\zeta(3) +3474\zeta(4) +3240\zeta(5)  \right. \nonumber \\ &\left. \hspace{1em}
+792\zeta(2)\zeta(3) -14325\zeta(6) -1080\zeta(3)^2 +9990\zeta(7) +2160\zeta(2)\zeta(5)  \right. \nonumber \\ &\left. \hspace{1em}
-9288\zeta(3)\zeta(4) +59332\zeta(8) +16128\zeta(2)\zeta(3)^2 -66816\zeta(3)\zeta(5) -14976 M(2,6)  \right. \nonumber \\ &\left. \hspace{1em}
+20928\zeta(9) -13392\zeta(3)\zeta(6) -12384\zeta(4)\zeta(5) +4032\zeta(2)\zeta(7)  \right. \nonumber \\ &\left. \hspace{1em}
+1920\zeta(3)^3\right) \label{eq11182} \\
\sum_{k=1}^\infty \frac{H(k)^{4}}{k^{4}(k+1)(k+2)^{2}} \sumend &= \frac{-1}{576}\left( -288 -612\zeta(2) -1836\zeta(3) -2070\zeta(4) +15120\zeta(5)  \right. \nonumber \\ &\left. \hspace{1em}
+2952\zeta(2)\zeta(3) -14865\zeta(6) -1080\zeta(3)^2 +6660\zeta(7) +1440\zeta(2)\zeta(5)  \right. \nonumber \\ &\left. \hspace{1em}
-6192\zeta(3)\zeta(4) +14833\zeta(8) +4032\zeta(2)\zeta(3)^2 -16704\zeta(3)\zeta(5)  \right. \nonumber \\ &\left. \hspace{1em}
-3744 M(2,6)\right) \label{eq11183} \\
\sum_{k=1}^\infty \frac{H(k)^{4}}{k^{3}(k+1)^{2}(k+2)^{2}} \sumend &= \frac{-1}{64}\left( -68 -148\zeta(2) -452\zeta(3) -534\zeta(4) -2520\zeta(5)  \right. \nonumber \\ &\left. \hspace{1em}
-520\zeta(2)\zeta(3) +4535\zeta(6) +360\zeta(3)^2 -370\zeta(7) -80\zeta(2)\zeta(5)  \right. \nonumber \\ &\left. \hspace{1em}
+344\zeta(3)\zeta(4)\right) \label{eq11184} \\
\sum_{k=1}^\infty \frac{H(k)^{4}}{k^{2}(k+1)^{3}(k+2)^{2}} \sumend &= \frac{-1}{8}\left( -18 -40\zeta(2) -124\zeta(3) -152\zeta(4) +300\zeta(5)  \right. \nonumber \\ &\left. \hspace{1em}
+56\zeta(2)\zeta(3) -10\zeta(6) +109\zeta(7) +40\zeta(2)\zeta(5) -148\zeta(3)\zeta(4)\right) \label{eq11185} \\
\sum_{k=1}^\infty \frac{H(k)^{4}}{k(k+1)^{4}(k+2)^{2}} \sumend &= \frac{1}{144}\left(  684 +1548\zeta(2) +4860\zeta(3) +6138\zeta(4) -1080\zeta(5)  \right. \nonumber \\ &\left. \hspace{1em}
-72\zeta(2)\zeta(3) -12885\zeta(6) -1080\zeta(3)^2 -1962\zeta(7) -720\zeta(2)\zeta(5)  \right. \nonumber \\ &\left. \hspace{1em}
+2664\zeta(3)\zeta(4) +12415\zeta(8) +3312\zeta(2)\zeta(3)^2 -13824\zeta(3)\zeta(5)  \right. \nonumber \\ &\left. \hspace{1em}
-3024 M(2,6)\right) \label{eq11186} \\
\sum_{k=1}^\infty \frac{H(k)^{4}}{(k+1)^{5}(k+2)^{2}} \sumend &= \frac{-1}{72}\left( -720 -1656\zeta(2) -5256\zeta(3) -6804\zeta(4) +2160\zeta(5)  \right. \nonumber \\ &\left. \hspace{1em}
+288\zeta(2)\zeta(3) +12885\zeta(6) +1080\zeta(3)^2 +2943\zeta(7) +1080\zeta(2)\zeta(5)  \right. \nonumber \\ &\left. \hspace{1em}
-3996\zeta(3)\zeta(4) -12415\zeta(8) -3312\zeta(2)\zeta(3)^2 +13824\zeta(3)\zeta(5) +3024 M(2,6)  \right. \nonumber \\ &\left. \hspace{1em}
+1044\zeta(9) -594\zeta(3)\zeta(6) -1332\zeta(4)\zeta(5) +504\zeta(2)\zeta(7) +192\zeta(3)^3\right) \label{eq11187}
\end{align}
 
\begin{align}
\sum_{k=1}^\infty \frac{H(k)^{4}}{k^{4}(k+2)^{3}} \sumend &= \frac{-1}{1152}\left(  1980 +2700\zeta(2) +6012\zeta(3) +2502\zeta(4) -432\zeta(5)  \right. \nonumber \\ &\left. \hspace{1em}
-216\zeta(2)\zeta(3) -13371\zeta(6) -1656\zeta(3)^2 +5976\zeta(7) +1440\zeta(2)\zeta(5)  \right. \nonumber \\ &\left. \hspace{1em}
-5976\zeta(3)\zeta(4) +14833\zeta(8) +4032\zeta(2)\zeta(3)^2 -16704\zeta(3)\zeta(5)  \right. \nonumber \\ &\left. \hspace{1em}
-3744 M(2,6)\right) \label{eq11188} \\
\sum_{k=1}^\infty \frac{H(k)^{4}}{k^{3}(k+1)(k+2)^{3}} \sumend &= \frac{-1}{32}\left(  126 +184\zeta(2) +436\zeta(3) +254\zeta(4) -864\zeta(5)  \right. \nonumber \\ &\left. \hspace{1em}
-176\zeta(2)\zeta(3) +83\zeta(6) -32\zeta(3)^2 -38\zeta(7) +12\zeta(3)\zeta(4)\right) \label{eq11189} \\
\sum_{k=1}^\infty \frac{H(k)^{4}}{k^{2}(k+1)^{2}(k+2)^{3}} \sumend &= \frac{1}{64}\left( -572 -884\zeta(2) -2196\zeta(3) -1550\zeta(4) +936\zeta(5)  \right. \nonumber \\ &\left. \hspace{1em}
+184\zeta(2)\zeta(3) +4203\zeta(6) +488\zeta(3)^2 -218\zeta(7) -80\zeta(2)\zeta(5)  \right. \nonumber \\ &\left. \hspace{1em}
+296\zeta(3)\zeta(4)\right) \label{eq11190} \\
\sum_{k=1}^\infty \frac{H(k)^{4}}{k(k+1)^{3}(k+2)^{3}} \sumend &= \frac{-1}{32}\left(  644 +1044\zeta(2) +2692\zeta(3) +2158\zeta(4) -2136\zeta(5)  \right. \nonumber \\ &\left. \hspace{1em}
-408\zeta(2)\zeta(3) -4163\zeta(6) -488\zeta(3)^2 -218\zeta(7) -80\zeta(2)\zeta(5)  \right. \nonumber \\ &\left. \hspace{1em}
+296\zeta(3)\zeta(4)\right) \label{eq11191} \\
\sum_{k=1}^\infty \frac{H(k)^{4}}{(k+1)^{4}(k+2)^{3}} \sumend &= \frac{1}{144}\left( -6480 -10944\zeta(2) -29088\zeta(3) -25560\zeta(4) +20304\zeta(5)  \right. \nonumber \\ &\left. \hspace{1em}
+3744\zeta(2)\zeta(3) +50352\zeta(6) +5472\zeta(3)^2 +3924\zeta(7) +1440\zeta(2)\zeta(5)  \right. \nonumber \\ &\left. \hspace{1em}
-5328\zeta(3)\zeta(4) -12415\zeta(8) -3312\zeta(2)\zeta(3)^2 +13824\zeta(3)\zeta(5)  \right. \nonumber \\ &\left. \hspace{1em}
+3024 M(2,6)\right) \label{eq11192} \\
\sum_{k=1}^\infty \frac{H(k)^{4}}{k^{3}(k+2)^{4}} \sumend &= \frac{1}{1152}\left(  9540 +8172\zeta(2) +12492\zeta(3) -4626\zeta(4) -8496\zeta(5)  \right. \nonumber \\ &\left. \hspace{1em}
-3672\zeta(2)\zeta(3) -11967\zeta(6) -3096\zeta(3)^2 +324\zeta(7) +288\zeta(2)\zeta(5)  \right. \nonumber \\ &\left. \hspace{1em}
-792\zeta(3)\zeta(4) +12415\zeta(8) +3312\zeta(2)\zeta(3)^2 -13824\zeta(3)\zeta(5) -3024 M(2,6)\right) \label{eq11193} \\
\sum_{k=1}^\infty \frac{H(k)^{4}}{k^{2}(k+1)(k+2)^{4}} \sumend &= \frac{1}{576}\left(  11808 +11484\zeta(2) +20340\zeta(3) -54\zeta(4) -24048\zeta(5)  \right. \nonumber \\ &\left. \hspace{1em}
-6840\zeta(2)\zeta(3) -10473\zeta(6) -3672\zeta(3)^2 -360\zeta(7) +288\zeta(2)\zeta(5)  \right. \nonumber \\ &\left. \hspace{1em}
-576\zeta(3)\zeta(4) +12415\zeta(8) +3312\zeta(2)\zeta(3)^2 -13824\zeta(3)\zeta(5) -3024 M(2,6)\right) \label{eq11194} \\
\sum_{k=1}^\infty \frac{H(k)^{4}}{k(k+1)^{2}(k+2)^{4}} \sumend &= \frac{1}{576}\left(  28764 +30924\zeta(2) +60444\zeta(3) +13842\zeta(4) -56520\zeta(5)  \right. \nonumber \\ &\left. \hspace{1em}
-15336\zeta(2)\zeta(3) -58773\zeta(6) -11736\zeta(3)^2 +1242\zeta(7) +1296\zeta(2)\zeta(5)  \right. \nonumber \\ &\left. \hspace{1em}
-3816\zeta(3)\zeta(4) +24830\zeta(8) +6624\zeta(2)\zeta(3)^2 -27648\zeta(3)\zeta(5)  \right. \nonumber \\ &\left. \hspace{1em}
-6048 M(2,6)\right) \label{eq11195}
\end{align}
 
\begin{align}
\sum_{k=1}^\infty \frac{H(k)^{4}}{(k+1)^{3}(k+2)^{4}} \sumend &= \frac{-1}{144}\left( -17280 -20160\zeta(2) -42336\zeta(3) -16632\zeta(4) +37872\zeta(5)  \right. \nonumber \\ &\left. \hspace{1em}
+9504\zeta(2)\zeta(3) +48120\zeta(6) +8064\zeta(3)^2 +360\zeta(7) -288\zeta(2)\zeta(5)  \right. \nonumber \\ &\left. \hspace{1em}
+576\zeta(3)\zeta(4) -12415\zeta(8) -3312\zeta(2)\zeta(3)^2 +13824\zeta(3)\zeta(5) +3024 M(2,6)\right) \label{eq11196} \\
\sum_{k=1}^\infty \frac{H(k)^{4}}{k^{2}(k+2)^{5}} \sumend &= \frac{1}{2304}\left( -68580 -34812\zeta(2) -27900\zeta(3) +56070\zeta(4) +39960\zeta(5)  \right. \nonumber \\ &\left. \hspace{1em}
+17064\zeta(2)\zeta(3) +17889\zeta(6) +2232\zeta(3)^2 +24930\zeta(7) +10512\zeta(2)\zeta(5)  \right. \nonumber \\ &\left. \hspace{1em}
-31752\zeta(3)\zeta(4) -50692\zeta(8) -16128\zeta(2)\zeta(3)^2 +62208\zeta(3)\zeta(5) +11520 M(2,6)  \right. \nonumber \\ &\left. \hspace{1em}
-8352\zeta(9) +4752\zeta(3)\zeta(6) +10656\zeta(4)\zeta(5) -4032\zeta(2)\zeta(7)  \right. \nonumber \\ &\left. \hspace{1em}
-1536\zeta(3)^3\right) \label{eq11197} \\
\sum_{k=1}^\infty \frac{H(k)^{4}}{k(k+1)(k+2)^{5}} \sumend &= \frac{1}{384}\left( -30732 -19260\zeta(2) -22860\zeta(3) +18726\zeta(4) +29352\zeta(5)  \right. \nonumber \\ &\left. \hspace{1em}
+10248\zeta(2)\zeta(3) +12945\zeta(6) +3192\zeta(3)^2 +8550\zeta(7) +3312\zeta(2)\zeta(5)  \right. \nonumber \\ &\left. \hspace{1em}
-10200\zeta(3)\zeta(4) -25174\zeta(8) -7584\zeta(2)\zeta(3)^2 +29952\zeta(3)\zeta(5) +5856 M(2,6)  \right. \nonumber \\ &\left. \hspace{1em}
-2784\zeta(9) +1584\zeta(3)\zeta(6) +3552\zeta(4)\zeta(5) -1344\zeta(2)\zeta(7) -512\zeta(3)^3\right) \label{eq11198} \\
\sum_{k=1}^\infty \frac{H(k)^{4}}{(k+1)^{2}(k+2)^{5}} \sumend &= \frac{1}{72}\left( -15120 -11088\zeta(2) -16128\zeta(3) +5292\zeta(4) +18072\zeta(5)  \right. \nonumber \\ &\left. \hspace{1em}
+5760\zeta(2)\zeta(3) +12201\zeta(6) +2664\zeta(3)^2 +3051\zeta(7) +1080\zeta(2)\zeta(5)  \right. \nonumber \\ &\left. \hspace{1em}
-3348\zeta(3)\zeta(4) -12544\zeta(8) -3672\zeta(2)\zeta(3)^2 +14688\zeta(3)\zeta(5) +2952 M(2,6)  \right. \nonumber \\ &\left. \hspace{1em}
-1044\zeta(9) +594\zeta(3)\zeta(6) +1332\zeta(4)\zeta(5) -504\zeta(2)\zeta(7) -192\zeta(3)^3\right) \label{eq11199} \\
\sum_{k=1}^\infty \frac{H(k)^{4}}{k(k+2)^{6}} \sumend &= \frac{-1}{11520}\left( -990540 -275580\zeta(2) +20340\zeta(3) +706950\zeta(4) +394920\zeta(5)  \right. \nonumber \\ &\left. \hspace{1em}
+76680\zeta(2)\zeta(3) +216465\zeta(6) -99720\zeta(3)^2 +272790\zeta(7) +140400\zeta(2)\zeta(5)  \right. \nonumber \\ &\left. \hspace{1em}
-385560\zeta(3)\zeta(4) -260590\zeta(8) -76320\zeta(2)\zeta(3)^2 +311040\zeta(3)\zeta(5) +56160 M(2,6)  \right. \nonumber \\ &\left. \hspace{1em}
+147360\zeta(9) -82800\zeta(3)\zeta(6) -136800\zeta(4)\zeta(5) +48960\zeta(2)\zeta(7) +15360\zeta(3)^3  \right. \nonumber \\ &\left. \hspace{1em}
-437823\zeta(10) +378720\zeta(3)\zeta(7) -2520\zeta(3)^2\zeta(4) -114480\zeta(2)\zeta(3)\zeta(5)  \right. \nonumber \\ &\left. \hspace{1em}
+119880\zeta(5)^2 +68760 M(2,8) +11520\zeta(2) M(2,6)\right) \label{eq11200} \\
\sum_{k=1}^\infty \frac{H(k)^{4}}{(k+1)(k+2)^{6}} \sumend &= \frac{-1}{5760}\left( -1451520 -564480\zeta(2) -322560\zeta(3) +987840\zeta(4)  \right. \nonumber \\ &\left. \hspace{1em}
+835200\zeta(5) +230400\zeta(2)\zeta(3) +410640\zeta(6) -51840\zeta(3)^2 +401040\zeta(7)  \right. \nonumber \\ &\left. \hspace{1em}
+190080\zeta(2)\zeta(5) -538560\zeta(3)\zeta(4) -638200\zeta(8) -190080\zeta(2)\zeta(3)^2  \right. \nonumber \\ &\left. \hspace{1em}
+760320\zeta(3)\zeta(5) +144000 M(2,6) +105600\zeta(9) -59040\zeta(3)\zeta(6) -83520\zeta(4)\zeta(5)  \right. \nonumber \\ &\left. \hspace{1em}
+28800\zeta(2)\zeta(7) +7680\zeta(3)^3 -437823\zeta(10) +378720\zeta(3)\zeta(7) -2520\zeta(3)^2\zeta(4)  \right. \nonumber \\ &\left. \hspace{1em}
-114480\zeta(2)\zeta(3)\zeta(5) +119880\zeta(5)^2 +68760 M(2,8) +11520\zeta(2) M(2,6)\right) \label{eq11201}
\end{align}
 
\begin{align}
\sum_{k=1}^\infty \frac{H(k)^{4}}{(k+2)^{7}} \sumend &= \frac{-1}{240}\left(  50400 +6720\zeta(2) -10080\zeta(3) -27720\zeta(4) -18000\zeta(5)  \right. \nonumber \\ &\left. \hspace{1em}
+1440\zeta(2)\zeta(3) -12180\zeta(6) +6720\zeta(3)^2 -11700\zeta(7) -3360\zeta(2)\zeta(5)  \right. \nonumber \\ &\left. \hspace{1em}
+12960\zeta(3)\zeta(4) +9310\zeta(8) +1200\zeta(2)\zeta(3)^2 -7680\zeta(3)\zeta(5) -2640 M(2,6)  \right. \nonumber \\ &\left. \hspace{1em}
-8120\zeta(9) +6600\zeta(3)\zeta(6) +8640\zeta(4)\zeta(5) -4320\zeta(2)\zeta(7) -1440\zeta(3)^3  \right. \nonumber \\ &\left. \hspace{1em}
+3006\zeta(10) -4800\zeta(3)\zeta(7) -480\zeta(3)^2\zeta(4) +3360\zeta(2)\zeta(3)\zeta(5)  \right. \nonumber \\ &\left. \hspace{1em}
-2880\zeta(5)^2 -240 M(2,8) +1185\zeta(11) +1840\zeta(2)\zeta(9) -430\zeta(3)\zeta(8)  \right. \nonumber \\ &\left. \hspace{1em}
-3420\zeta(4)\zeta(7) -1010\zeta(5)\zeta(6) -400\zeta(2)\zeta(3)^3 +1440\zeta(3)^2\zeta(5)  \right. \nonumber \\ &\left. \hspace{1em}
+80 M(3,8)\right) \label{eq11202} \\
\sum_{k=1}^\infty \frac{H(k)^{5}}{k^{6}} \sumend &= \frac{-1}{576}\left(  781671\zeta(11) +88016\zeta(2)\zeta(9) -296660\zeta(3)\zeta(8)  \right. \nonumber \\ &\left. \hspace{1em}
-411984\zeta(4)\zeta(7) -220080\zeta(5)\zeta(6) -21120\zeta(2)\zeta(3)^3 +141120\zeta(3)^2\zeta(5)  \right. \nonumber \\ &\left. \hspace{1em}
-8640\zeta(3) M(2,6) -27840 M(3,8)\right) \label{eq11203} \\
\sum_{k=1}^\infty \frac{H(k)^{5}}{k^{5}(k+1)} \sumend &= \frac{-1}{2304}\left( -411264\zeta(6) -51840\zeta(3)^2 +295344\zeta(7) +65664\zeta(2)\zeta(5)  \right. \nonumber \\ &\left. \hspace{1em}
+76032\zeta(3)\zeta(4) +542488\zeta(8) +152640\zeta(2)\zeta(3)^2 -630144\zeta(3)\zeta(5) -135360 M(2,6)  \right. \nonumber \\ &\left. \hspace{1em}
+302144\zeta(9) -469920\zeta(3)\zeta(6) +152064\zeta(4)\zeta(5) +76320\zeta(2)\zeta(7) -11520\zeta(3)^3  \right. \nonumber \\ &\left. \hspace{1em}
+579897\zeta(10) -519840\zeta(3)\zeta(7) +3240\zeta(3)^2\zeta(4) +185040\zeta(2)\zeta(3)\zeta(5)  \right. \nonumber \\ &\left. \hspace{1em}
-203832\zeta(5)^2 -98280 M(2,8) -11520\zeta(2) M(2,6)\right) \label{eq11204} \\
\sum_{k=1}^\infty \frac{H(k)^{5}}{k^{4}(k+1)^{2}} \sumend &= \frac{1}{144}\left( -102816\zeta(6) -12960\zeta(3)^2 +72072\zeta(7) +16416\zeta(2)\zeta(5)  \right. \nonumber \\ &\left. \hspace{1em}
+19008\zeta(3)\zeta(4) +67811\zeta(8) +19080\zeta(2)\zeta(3)^2 -78768\zeta(3)\zeta(5) -16920 M(2,6)  \right. \nonumber \\ &\left. \hspace{1em}
+18884\zeta(9) -29370\zeta(3)\zeta(6) +9504\zeta(4)\zeta(5) +4770\zeta(2)\zeta(7)  \right. \nonumber \\ &\left. \hspace{1em}
-720\zeta(3)^3\right) \label{eq11205} \\
\sum_{k=1}^\infty \frac{H(k)^{5}}{k^{3}(k+1)^{3}} \sumend &= \frac{1}{72}\left(  77112\zeta(6) +9720\zeta(3)^2 -52731\zeta(7) -12312\zeta(2)\zeta(5)  \right. \nonumber \\ &\left. \hspace{1em}
-14256\zeta(3)\zeta(4) -33358\zeta(8) -9180\zeta(2)\zeta(3)^2 +37800\zeta(3)\zeta(5)  \right. \nonumber \\ &\left. \hspace{1em}
+8100 M(2,6)\right) \label{eq11206} \\
\sum_{k=1}^\infty \frac{H(k)^{5}}{k^{2}(k+1)^{4}} \sumend &= \frac{-1}{144}\left(  102816\zeta(6) +12960\zeta(3)^2 -68544\zeta(7) -16416\zeta(2)\zeta(5)  \right. \nonumber \\ &\left. \hspace{1em}
-19008\zeta(3)\zeta(4) -65621\zeta(8) -17640\zeta(2)\zeta(3)^2 +72432\zeta(3)\zeta(5) +15480 M(2,6)  \right. \nonumber \\ &\left. \hspace{1em}
-14240\zeta(9) +25770\zeta(3)\zeta(6) -9504\zeta(4)\zeta(5) -4770\zeta(2)\zeta(7)  \right. \nonumber \\ &\left. \hspace{1em}
+720\zeta(3)^3\right) \label{eq11207}
\end{align}
 
\begin{align}
\sum_{k=1}^\infty \frac{H(k)^{5}}{k(k+1)^{5}} \sumend &= \frac{-1}{2304}\left( -411264\zeta(6) -51840\zeta(3)^2 +267120\zeta(7) +65664\zeta(2)\zeta(5)  \right. \nonumber \\ &\left. \hspace{1em}
+76032\zeta(3)\zeta(4) +524968\zeta(8) +141120\zeta(2)\zeta(3)^2 -579456\zeta(3)\zeta(5) -123840 M(2,6)  \right. \nonumber \\ &\left. \hspace{1em}
+227840\zeta(9) -412320\zeta(3)\zeta(6) +152064\zeta(4)\zeta(5) +76320\zeta(2)\zeta(7) -11520\zeta(3)^3  \right. \nonumber \\ &\left. \hspace{1em}
+449109\zeta(10) -387360\zeta(3)\zeta(7) -9720\zeta(3)^2\zeta(4) +124560\zeta(2)\zeta(3)\zeta(5)  \right. \nonumber \\ &\left. \hspace{1em}
-122328\zeta(5)^2 -68040 M(2,8) -11520\zeta(2) M(2,6)\right) \label{eq11208} \\
\sum_{k=1}^\infty \frac{H(k)^{5}}{(k+1)^{6}} \sumend &= \frac{-1}{576}\left(  667227\zeta(11) +68816\zeta(2)\zeta(9) -248300\zeta(3)\zeta(8)  \right. \nonumber \\ &\left. \hspace{1em}
-350784\zeta(4)\zeta(7) -176280\zeta(5)\zeta(6) -16320\zeta(2)\zeta(3)^3 +112320\zeta(3)^2\zeta(5)  \right. \nonumber \\ &\left. \hspace{1em}
-8640\zeta(3) M(2,6) -23040 M(3,8)\right) \label{eq11209} \\
\sum_{k=1}^\infty \frac{H(k)^{5}}{k^{5}(k+2)} \sumend &= \frac{-1}{4608}\left( -144 -576\zeta(2) -3024\zeta(3) -9036\zeta(4) -10224\zeta(5)  \right. \nonumber \\ &\left. \hspace{1em}
-2160\zeta(2)\zeta(3) -25704\zeta(6) -3240\zeta(3)^2 +36918\zeta(7) +8208\zeta(2)\zeta(5)  \right. \nonumber \\ &\left. \hspace{1em}
+9504\zeta(3)\zeta(4) +135622\zeta(8) +38160\zeta(2)\zeta(3)^2 -157536\zeta(3)\zeta(5) -33840 M(2,6)  \right. \nonumber \\ &\left. \hspace{1em}
+151072\zeta(9) -234960\zeta(3)\zeta(6) +76032\zeta(4)\zeta(5) +38160\zeta(2)\zeta(7) -5760\zeta(3)^3  \right. \nonumber \\ &\left. \hspace{1em}
+579897\zeta(10) -519840\zeta(3)\zeta(7) +3240\zeta(3)^2\zeta(4) +185040\zeta(2)\zeta(3)\zeta(5)  \right. \nonumber \\ &\left. \hspace{1em}
-203832\zeta(5)^2 -98280 M(2,8) -11520\zeta(2) M(2,6)\right) \label{eq11210} \\
\sum_{k=1}^\infty \frac{H(k)^{5}}{k^{4}(k+1)(k+2)} \sumend &= \frac{1}{1152}\left(  72 +288\zeta(2) +1512\zeta(3) +4518\zeta(4) +5112\zeta(5)  \right. \nonumber \\ &\left. \hspace{1em}
+1080\zeta(2)\zeta(3) -192780\zeta(6) -24300\zeta(3)^2 +129213\zeta(7) +28728\zeta(2)\zeta(5)  \right. \nonumber \\ &\left. \hspace{1em}
+33264\zeta(3)\zeta(4) +203433\zeta(8) +57240\zeta(2)\zeta(3)^2 -236304\zeta(3)\zeta(5) -50760 M(2,6)  \right. \nonumber \\ &\left. \hspace{1em}
+75536\zeta(9) -117480\zeta(3)\zeta(6) +38016\zeta(4)\zeta(5) +19080\zeta(2)\zeta(7)  \right. \nonumber \\ &\left. \hspace{1em}
-2880\zeta(3)^3\right) \label{eq11211} \\
\sum_{k=1}^\infty \frac{H(k)^{5}}{k^{3}(k+1)^{2}(k+2)} \sumend &= \frac{1}{576}\left(  72 +288\zeta(2) +1512\zeta(3) +4518\zeta(4) +5112\zeta(5)  \right. \nonumber \\ &\left. \hspace{1em}
+1080\zeta(2)\zeta(3) +218484\zeta(6) +27540\zeta(3)^2 -159075\zeta(7) -36936\zeta(2)\zeta(5)  \right. \nonumber \\ &\left. \hspace{1em}
-42768\zeta(3)\zeta(4) -67811\zeta(8) -19080\zeta(2)\zeta(3)^2 +78768\zeta(3)\zeta(5)  \right. \nonumber \\ &\left. \hspace{1em}
+16920 M(2,6)\right) \label{eq11212} \\
\sum_{k=1}^\infty \frac{H(k)^{5}}{k^{2}(k+1)^{3}(k+2)} \sumend &= \frac{-1}{288}\left( -72 -288\zeta(2) -1512\zeta(3) -4518\zeta(4) -5112\zeta(5)  \right. \nonumber \\ &\left. \hspace{1em}
-1080\zeta(2)\zeta(3) +89964\zeta(6) +11340\zeta(3)^2 -51849\zeta(7) -12312\zeta(2)\zeta(5)  \right. \nonumber \\ &\left. \hspace{1em}
-14256\zeta(3)\zeta(4) -65621\zeta(8) -17640\zeta(2)\zeta(3)^2 +72432\zeta(3)\zeta(5)  \right. \nonumber \\ &\left. \hspace{1em}
+15480 M(2,6)\right) \label{eq11213}
\end{align}
 
\begin{align}
\sum_{k=1}^\infty \frac{H(k)^{5}}{k(k+1)^{4}(k+2)} \sumend &= \frac{-1}{144}\left( -72 -288\zeta(2) -1512\zeta(3) -4518\zeta(4) -5112\zeta(5)  \right. \nonumber \\ &\left. \hspace{1em}
-1080\zeta(2)\zeta(3) -12852\zeta(6) -1620\zeta(3)^2 +16695\zeta(7) +4104\zeta(2)\zeta(5)  \right. \nonumber \\ &\left. \hspace{1em}
+4752\zeta(3)\zeta(4) +14240\zeta(9) -25770\zeta(3)\zeta(6) +9504\zeta(4)\zeta(5) +4770\zeta(2)\zeta(7)  \right. \nonumber \\ &\left. \hspace{1em}
-720\zeta(3)^3\right) \label{eq11214} \\
\sum_{k=1}^\infty \frac{H(k)^{5}}{(k+1)^{5}(k+2)} \sumend &= \frac{1}{2304}\left(  2304 +9216\zeta(2) +48384\zeta(3) +144576\zeta(4) +163584\zeta(5)  \right. \nonumber \\ &\left. \hspace{1em}
+34560\zeta(2)\zeta(3) -267120\zeta(7) -65664\zeta(2)\zeta(5) -76032\zeta(3)\zeta(4) +524968\zeta(8)  \right. \nonumber \\ &\left. \hspace{1em}
+141120\zeta(2)\zeta(3)^2 -579456\zeta(3)\zeta(5) -123840 M(2,6) -227840\zeta(9)  \right. \nonumber \\ &\left. \hspace{1em}
+412320\zeta(3)\zeta(6) -152064\zeta(4)\zeta(5) -76320\zeta(2)\zeta(7) +11520\zeta(3)^3  \right. \nonumber \\ &\left. \hspace{1em}
+449109\zeta(10) -387360\zeta(3)\zeta(7) -9720\zeta(3)^2\zeta(4) +124560\zeta(2)\zeta(3)\zeta(5)  \right. \nonumber \\ &\left. \hspace{1em}
-122328\zeta(5)^2 -68040 M(2,8) -11520\zeta(2) M(2,6)\right) \label{eq11215} \\
\sum_{k=1}^\infty \frac{H(k)^{5}}{k^{4}(k+2)^{2}} \sumend &= \frac{1}{1152}\left( -576 -1656\zeta(2) -7200\zeta(3) -16092\zeta(4) -9936\zeta(5)  \right. \nonumber \\ &\left. \hspace{1em}
-2520\zeta(2)\zeta(3) -12819\zeta(6) -2160\zeta(3)^2 +36036\zeta(7) +8208\zeta(2)\zeta(5)  \right. \nonumber \\ &\left. \hspace{1em}
+9504\zeta(3)\zeta(4) +67811\zeta(8) +19080\zeta(2)\zeta(3)^2 -78768\zeta(3)\zeta(5) -16920 M(2,6)  \right. \nonumber \\ &\left. \hspace{1em}
+37768\zeta(9) -58740\zeta(3)\zeta(6) +19008\zeta(4)\zeta(5) +9540\zeta(2)\zeta(7)  \right. \nonumber \\ &\left. \hspace{1em}
-1440\zeta(3)^3\right) \label{eq11216} \\
\sum_{k=1}^\infty \frac{H(k)^{5}}{k^{3}(k+1)(k+2)^{2}} \sumend &= \frac{-1}{1152}\left(  1224 +3600\zeta(2) +15912\zeta(3) +36702\zeta(4) +24984\zeta(5)  \right. \nonumber \\ &\left. \hspace{1em}
+6120\zeta(2)\zeta(3) -167142\zeta(6) -19980\zeta(3)^2 +57141\zeta(7) +12312\zeta(2)\zeta(5)  \right. \nonumber \\ &\left. \hspace{1em}
+14256\zeta(3)\zeta(4) +67811\zeta(8) +19080\zeta(2)\zeta(3)^2 -78768\zeta(3)\zeta(5)  \right. \nonumber \\ &\left. \hspace{1em}
-16920 M(2,6)\right) \label{eq11217} \\
\sum_{k=1}^\infty \frac{H(k)^{5}}{k^{2}(k+1)^{2}(k+2)^{2}} \sumend &= \frac{-1}{96}\left(  216 +648\zeta(2) +2904\zeta(3) +6870\zeta(4) +5016\zeta(5)  \right. \nonumber \\ &\left. \hspace{1em}
+1200\zeta(2)\zeta(3) +8557\zeta(6) +1260\zeta(3)^2 -16989\zeta(7) -4104\zeta(2)\zeta(5)  \right. \nonumber \\ &\left. \hspace{1em}
-4752\zeta(3)\zeta(4)\right) \label{eq11218} \\
\sum_{k=1}^\infty \frac{H(k)^{5}}{k(k+1)^{3}(k+2)^{2}} \sumend &= \frac{-1}{288}\left(  1368 +4176\zeta(2) +18936\zeta(3) +45738\zeta(4) +35208\zeta(5)  \right. \nonumber \\ &\left. \hspace{1em}
+8280\zeta(2)\zeta(3) -38622\zeta(6) -3780\zeta(3)^2 -50085\zeta(7) -12312\zeta(2)\zeta(5)  \right. \nonumber \\ &\left. \hspace{1em}
-14256\zeta(3)\zeta(4) +65621\zeta(8) +17640\zeta(2)\zeta(3)^2 -72432\zeta(3)\zeta(5)  \right. \nonumber \\ &\left. \hspace{1em}
-15480 M(2,6)\right) \label{eq11219}
\end{align}
 
\begin{align}
\sum_{k=1}^\infty \frac{H(k)^{5}}{(k+1)^{4}(k+2)^{2}} \sumend &= \frac{1}{144}\left( -1440 -4464\zeta(2) -20448\zeta(3) -50256\zeta(4) -40320\zeta(5)  \right. \nonumber \\ &\left. \hspace{1em}
-9360\zeta(2)\zeta(3) +25770\zeta(6) +2160\zeta(3)^2 +66780\zeta(7) +16416\zeta(2)\zeta(5)  \right. \nonumber \\ &\left. \hspace{1em}
+19008\zeta(3)\zeta(4) -65621\zeta(8) -17640\zeta(2)\zeta(3)^2 +72432\zeta(3)\zeta(5) +15480 M(2,6)  \right. \nonumber \\ &\left. \hspace{1em}
+14240\zeta(9) -25770\zeta(3)\zeta(6) +9504\zeta(4)\zeta(5) +4770\zeta(2)\zeta(7)  \right. \nonumber \\ &\left. \hspace{1em}
-720\zeta(3)^3\right) \label{eq11220} \\
\sum_{k=1}^\infty \frac{H(k)^{5}}{k^{3}(k+2)^{3}} \sumend &= \frac{-1}{1152}\left( -4536 -9144\zeta(2) -32184\zeta(3) -49770\zeta(4) -5256\zeta(5)  \right. \nonumber \\ &\left. \hspace{1em}
-2160\zeta(2)\zeta(3) +22899\zeta(6) +3420\zeta(3)^2 +42921\zeta(7) +8712\zeta(2)\zeta(5)  \right. \nonumber \\ &\left. \hspace{1em}
+27576\zeta(3)\zeta(4) +66716\zeta(8) +18360\zeta(2)\zeta(3)^2 -75600\zeta(3)\zeta(5)  \right. \nonumber \\ &\left. \hspace{1em}
-16200 M(2,6)\right) \label{eq11221} \\
\sum_{k=1}^\infty \frac{H(k)^{5}}{k^{2}(k+1)(k+2)^{3}} \sumend &= \frac{1}{1152}\left(  10296 +21888\zeta(2) +80280\zeta(3) +136242\zeta(4)  \right. \nonumber \\ &\left. \hspace{1em}
+35496\zeta(5) +10440\zeta(2)\zeta(3) -212940\zeta(6) -26820\zeta(3)^2 -28701\zeta(7)  \right. \nonumber \\ &\left. \hspace{1em}
-5112\zeta(2)\zeta(5) -40896\zeta(3)\zeta(4) -65621\zeta(8) -17640\zeta(2)\zeta(3)^2  \right. \nonumber \\ &\left. \hspace{1em}
+72432\zeta(3)\zeta(5) +15480 M(2,6)\right) \label{eq11222} \\
\sum_{k=1}^\infty \frac{H(k)^{5}}{k(k+1)^{2}(k+2)^{3}} \sumend &= \frac{1}{576}\left(  11592 +25776\zeta(2) +97704\zeta(3) +177462\zeta(4) +65592\zeta(5)  \right. \nonumber \\ &\left. \hspace{1em}
+17640\zeta(2)\zeta(3) -161598\zeta(6) -19260\zeta(3)^2 -130635\zeta(7) -29736\zeta(2)\zeta(5)  \right. \nonumber \\ &\left. \hspace{1em}
-69408\zeta(3)\zeta(4) -65621\zeta(8) -17640\zeta(2)\zeta(3)^2 +72432\zeta(3)\zeta(5)  \right. \nonumber \\ &\left. \hspace{1em}
+15480 M(2,6)\right) \label{eq11223} \\
\sum_{k=1}^\infty \frac{H(k)^{5}}{(k+1)^{3}(k+2)^{3}} \sumend &= \frac{-1}{24}\left( -1080 -2496\zeta(2) -9720\zeta(3) -18600\zeta(4) -8400\zeta(5)  \right. \nonumber \\ &\left. \hspace{1em}
-2160\zeta(2)\zeta(3) +16685\zeta(6) +1920\zeta(3)^2 +15060\zeta(7) +3504\zeta(2)\zeta(5)  \right. \nonumber \\ &\left. \hspace{1em}
+6972\zeta(3)\zeta(4)\right) \label{eq11224} \\
\sum_{k=1}^\infty \frac{H(k)^{5}}{k^{2}(k+2)^{4}} \sumend &= \frac{-1}{1152}\left(  23616 +32616\zeta(2) +90432\zeta(3) +83340\zeta(4) -43344\zeta(5)  \right. \nonumber \\ &\left. \hspace{1em}
-12600\zeta(2)\zeta(3) -102651\zeta(6) -23040\zeta(3)^2 -16452\zeta(7) +432\zeta(2)\zeta(5)  \right. \nonumber \\ &\left. \hspace{1em}
-39024\zeta(3)\zeta(4) +58529\zeta(8) +15480\zeta(2)\zeta(3)^2 -65808\zeta(3)\zeta(5) -14760 M(2,6)  \right. \nonumber \\ &\left. \hspace{1em}
-28480\zeta(9) +51540\zeta(3)\zeta(6) -19008\zeta(4)\zeta(5) -9540\zeta(2)\zeta(7)  \right. \nonumber \\ &\left. \hspace{1em}
+1440\zeta(3)^3\right) \label{eq11225} \\
\sum_{k=1}^\infty \frac{H(k)^{5}}{k(k+1)(k+2)^{4}} \sumend &= \frac{1}{1152}\left( -57528 -87120\zeta(2) -261144\zeta(3) -302922\zeta(4) +51192\zeta(5)  \right. \nonumber \\ &\left. \hspace{1em}
+14760\zeta(2)\zeta(3) +418242\zeta(6) +72900\zeta(3)^2 +61605\zeta(7) +4248\zeta(2)\zeta(5)  \right. \nonumber \\ &\left. \hspace{1em}
+118944\zeta(3)\zeta(4) -51437\zeta(8) -13320\zeta(2)\zeta(3)^2 +59184\zeta(3)\zeta(5) +14040 M(2,6)  \right. \nonumber \\ &\left. \hspace{1em}
+56960\zeta(9) -103080\zeta(3)\zeta(6) +38016\zeta(4)\zeta(5) +19080\zeta(2)\zeta(7)  \right. \nonumber \\ &\left. \hspace{1em}
-2880\zeta(3)^3\right) \label{eq11226}
\end{align}
 
\begin{align}
\sum_{k=1}^\infty \frac{H(k)^{5}}{(k+1)^{2}(k+2)^{4}} \sumend &= \frac{1}{72}\left( -8640 -14112\zeta(2) -44856\zeta(3) -60048\zeta(4) -1800\zeta(5)  \right. \nonumber \\ &\left. \hspace{1em}
-360\zeta(2)\zeta(3) +72480\zeta(6) +11520\zeta(3)^2 +24030\zeta(7) +4248\zeta(2)\zeta(5)  \right. \nonumber \\ &\left. \hspace{1em}
+23544\zeta(3)\zeta(4) +1773\zeta(8) +540\zeta(2)\zeta(3)^2 -1656\zeta(3)\zeta(5) -180 M(2,6)  \right. \nonumber \\ &\left. \hspace{1em}
+7120\zeta(9) -12885\zeta(3)\zeta(6) +4752\zeta(4)\zeta(5) +2385\zeta(2)\zeta(7) -360\zeta(3)^3\right) \label{eq11227} \\
\sum_{k=1}^\infty \frac{H(k)^{5}}{k(k+2)^{5}} \sumend &= \frac{-1}{4608}\left( -368784 -341856\zeta(2) -719568\zeta(3) -238572\zeta(4) +579600\zeta(5)  \right. \nonumber \\ &\left. \hspace{1em}
+200880\zeta(2)\zeta(3) +678396\zeta(6) +171000\zeta(3)^2 +164070\zeta(7) +45648\zeta(2)\zeta(5)  \right. \nonumber \\ &\left. \hspace{1em}
-104256\zeta(3)\zeta(4) -1368878\zeta(8) -390960\zeta(2)\zeta(3)^2 +1583136\zeta(3)\zeta(5)  \right. \nonumber \\ &\left. \hspace{1em}
+326160 M(2,6) -53120\zeta(9) -111120\zeta(3)\zeta(6) +289152\zeta(4)\zeta(5) -42480\zeta(2)\zeta(7)  \right. \nonumber \\ &\left. \hspace{1em}
-36480\zeta(3)^3 +449109\zeta(10) -387360\zeta(3)\zeta(7) -9720\zeta(3)^2\zeta(4)  \right. \nonumber \\ &\left. \hspace{1em}
+124560\zeta(2)\zeta(3)\zeta(5) -122328\zeta(5)^2 -68040 M(2,8) -11520\zeta(2) M(2,6)\right) \label{eq11228} \\
\sum_{k=1}^\infty \frac{H(k)^{5}}{(k+1)(k+2)^{5}} \sumend &= \frac{-1}{2304}\left( -483840 -516096\zeta(2) -1241856\zeta(3) -844416\zeta(4)  \right. \nonumber \\ &\left. \hspace{1em}
+681984\zeta(5) +230400\zeta(2)\zeta(3) +1514880\zeta(6) +316800\zeta(3)^2 +287280\zeta(7)  \right. \nonumber \\ &\left. \hspace{1em}
+54144\zeta(2)\zeta(5) +133632\zeta(3)\zeta(4) -1471752\zeta(8) -417600\zeta(2)\zeta(3)^2  \right. \nonumber \\ &\left. \hspace{1em}
+1701504\zeta(3)\zeta(5) +354240 M(2,6) +60800\zeta(9) -317280\zeta(3)\zeta(6) +365184\zeta(4)\zeta(5)  \right. \nonumber \\ &\left. \hspace{1em}
-4320\zeta(2)\zeta(7) -42240\zeta(3)^3 +449109\zeta(10) -387360\zeta(3)\zeta(7) -9720\zeta(3)^2\zeta(4)  \right. \nonumber \\ &\left. \hspace{1em}
+124560\zeta(2)\zeta(3)\zeta(5) -122328\zeta(5)^2 -68040 M(2,8) -11520\zeta(2) M(2,6)\right) \label{eq11229} \\
\sum_{k=1}^\infty \frac{H(k)^{5}}{(k+2)^{6}} \sumend &= \frac{-1}{1152}\left(  290304 +177408\zeta(2) +266112\zeta(3) -87552\zeta(4) -298368\zeta(5)  \right. \nonumber \\ &\left. \hspace{1em}
-97920\zeta(2)\zeta(3) -236112\zeta(6) -28800\zeta(3)^2 -161280\zeta(7) -69120\zeta(2)\zeta(5)  \right. \nonumber \\ &\left. \hspace{1em}
+201600\zeta(3)\zeta(4) +547240\zeta(8) +161280\zeta(2)\zeta(3)^2 -645120\zeta(3)\zeta(5)  \right. \nonumber \\ &\left. \hspace{1em}
-126720 M(2,6) -11040\zeta(9) +5760\zeta(3)\zeta(6) -11520\zeta(4)\zeta(5) +5760\zeta(2)\zeta(7)  \right. \nonumber \\ &\left. \hspace{1em}
+3840\zeta(3)^3 +437823\zeta(10) -378720\zeta(3)\zeta(7) +2520\zeta(3)^2\zeta(4)  \right. \nonumber \\ &\left. \hspace{1em}
+114480\zeta(2)\zeta(3)\zeta(5) -119880\zeta(5)^2 -68760 M(2,8) -11520\zeta(2) M(2,6) +1334454\zeta(11)  \right. \nonumber \\ &\left. \hspace{1em}
+137632\zeta(2)\zeta(9) -496600\zeta(3)\zeta(8) -701568\zeta(4)\zeta(7) -352560\zeta(5)\zeta(6)  \right. \nonumber \\ &\left. \hspace{1em}
-32640\zeta(2)\zeta(3)^3 +224640\zeta(3)^2\zeta(5) -17280\zeta(3) M(2,6) -46080 M(3,8)\right) \label{eq11230} \\
\sum_{k=1}^\infty \frac{H(k)^{6}}{k^{5}} \sumend &= \frac{1}{192}\left( -734643\zeta(11) -83472\zeta(2)\zeta(9) +271244\zeta(3)\zeta(8)  \right. \nonumber \\ &\left. \hspace{1em}
+395088\zeta(4)\zeta(7) +205424\zeta(5)\zeta(6) +19360\zeta(2)\zeta(3)^3 -130176\zeta(3)^2\zeta(5)  \right. \nonumber \\ &\left. \hspace{1em}
+9120\zeta(3) M(2,6) +25600 M(3,8)\right) \label{eq11231}
\end{align}
 
\begin{align}
\sum_{k=1}^\infty \frac{H(k)^{6}}{k^{4}(k+1)} \sumend &= \frac{-1}{384}\left(  247296\zeta(7) +55680\zeta(2)\zeta(5) +114048\zeta(3)\zeta(4)  \right. \nonumber \\ &\left. \hspace{1em}
-280464\zeta(8) +15744\zeta(2)\zeta(3)^2 -187008\zeta(3)\zeta(5) -21888 M(2,6) +119584\zeta(9)  \right. \nonumber \\ &\left. \hspace{1em}
-209952\zeta(3)\zeta(6) +96768\zeta(4)\zeta(5) +31248\zeta(2)\zeta(7) -8704\zeta(3)^3 +814101\zeta(10)  \right. \nonumber \\ &\left. \hspace{1em}
-529680\zeta(3)\zeta(7) +253944\zeta(3)^2\zeta(4) +1200\zeta(2)\zeta(3)\zeta(5) -365064\zeta(5)^2  \right. \nonumber \\ &\left. \hspace{1em}
-103128 M(2,8) -45120\zeta(2) M(2,6)\right) \label{eq11232} \\
\sum_{k=1}^\infty \frac{H(k)^{6}}{k^{3}(k+1)^{2}} \sumend &= \frac{1}{24}\left(  46368\zeta(7) +10440\zeta(2)\zeta(5) +21384\zeta(3)\zeta(4)  \right. \nonumber \\ &\left. \hspace{1em}
-52085\zeta(8) +2892\zeta(2)\zeta(3)^2 -34704\zeta(3)\zeta(5) -4044 M(2,6) +7474\zeta(9)  \right. \nonumber \\ &\left. \hspace{1em}
-13122\zeta(3)\zeta(6) +6048\zeta(4)\zeta(5) +1953\zeta(2)\zeta(7) -544\zeta(3)^3\right) \label{eq11233} \\
\sum_{k=1}^\infty \frac{H(k)^{6}}{k^{2}(k+1)^{3}} \sumend &= \frac{1}{24}\left( -46368\zeta(7) -10440\zeta(2)\zeta(5) -21384\zeta(3)\zeta(4)  \right. \nonumber \\ &\left. \hspace{1em}
+51583\zeta(8) -2832\zeta(2)\zeta(3)^2 +34344\zeta(3)\zeta(5) +3984 M(2,6) -6146\zeta(9)  \right. \nonumber \\ &\left. \hspace{1em}
+12582\zeta(3)\zeta(6) -5832\zeta(4)\zeta(5) -1953\zeta(2)\zeta(7) +536\zeta(3)^3\right) \label{eq11234} \\
\sum_{k=1}^\infty \frac{H(k)^{6}}{k(k+1)^{4}} \sumend &= \frac{-1}{384}\left( -247296\zeta(7) -55680\zeta(2)\zeta(5) -114048\zeta(3)\zeta(4)  \right. \nonumber \\ &\left. \hspace{1em}
+272432\zeta(8) -14784\zeta(2)\zeta(3)^2 +181248\zeta(3)\zeta(5) +20928 M(2,6) -98336\zeta(9)  \right. \nonumber \\ &\left. \hspace{1em}
+201312\zeta(3)\zeta(6) -93312\zeta(4)\zeta(5) -31248\zeta(2)\zeta(7) +8576\zeta(3)^3 -779835\zeta(10)  \right. \nonumber \\ &\left. \hspace{1em}
+490704\zeta(3)\zeta(7) -245544\zeta(3)^2\zeta(4) +15600\zeta(2)\zeta(3)\zeta(5) +339864\zeta(5)^2  \right. \nonumber \\ &\left. \hspace{1em}
+94728 M(2,8) +45120\zeta(2) M(2,6)\right) \label{eq11235} \\
\sum_{k=1}^\infty \frac{H(k)^{6}}{(k+1)^{5}} \sumend &= \frac{1}{192}\left(  686799\zeta(11) +74512\zeta(2)\zeta(9) -262484\zeta(3)\zeta(8)  \right. \nonumber \\ &\left. \hspace{1em}
-362208\zeta(4)\zeta(7) -182584\zeta(5)\zeta(6) -18080\zeta(2)\zeta(3)^3 +120384\zeta(3)^2\zeta(5)  \right. \nonumber \\ &\left. \hspace{1em}
-8160\zeta(3) M(2,6) -23360 M(3,8)\right) \label{eq11236} \\
\sum_{k=1}^\infty \frac{H(k)^{6}}{k^{4}(k+2)} \sumend &= \frac{1}{768}\left( -48 -240\zeta(2) -1632\zeta(3) -6852\zeta(4) -13704\zeta(5)  \right. \nonumber \\ &\left. \hspace{1em}
-2928\zeta(2)\zeta(3) -25164\zeta(6) -3216\zeta(3)^2 -30912\zeta(7) -6960\zeta(2)\zeta(5)  \right. \nonumber \\ &\left. \hspace{1em}
-14256\zeta(3)\zeta(4) +70116\zeta(8) -3936\zeta(2)\zeta(3)^2 +46752\zeta(3)\zeta(5) +5472 M(2,6)  \right. \nonumber \\ &\left. \hspace{1em}
-59792\zeta(9) +104976\zeta(3)\zeta(6) -48384\zeta(4)\zeta(5) -15624\zeta(2)\zeta(7) +4352\zeta(3)^3  \right. \nonumber \\ &\left. \hspace{1em}
-814101\zeta(10) +529680\zeta(3)\zeta(7) -253944\zeta(3)^2\zeta(4) -1200\zeta(2)\zeta(3)\zeta(5)  \right. \nonumber \\ &\left. \hspace{1em}
+365064\zeta(5)^2 +103128 M(2,8) +45120\zeta(2) M(2,6)\right) \label{eq11237} \\
\sum_{k=1}^\infty \frac{H(k)^{6}}{k^{3}(k+1)(k+2)} \sumend &= \frac{1}{96}\left( -12 -60\zeta(2) -408\zeta(3) -1713\zeta(4) -3426\zeta(5)  \right. \nonumber \\ &\left. \hspace{1em}
-732\zeta(2)\zeta(3) -6291\zeta(6) -804\zeta(3)^2 +54096\zeta(7) +12180\zeta(2)\zeta(5)  \right. \nonumber \\ &\left. \hspace{1em}
+24948\zeta(3)\zeta(4) -52587\zeta(8) +2952\zeta(2)\zeta(3)^2 -35064\zeta(3)\zeta(5) -4104 M(2,6)  \right. \nonumber \\ &\left. \hspace{1em}
+14948\zeta(9) -26244\zeta(3)\zeta(6) +12096\zeta(4)\zeta(5) +3906\zeta(2)\zeta(7)  \right. \nonumber \\ &\left. \hspace{1em}
-1088\zeta(3)^3\right) \label{eq11238}
\end{align}
 
\begin{align}
\sum_{k=1}^\infty \frac{H(k)^{6}}{k^{2}(k+1)^{2}(k+2)} \sumend &= \frac{1}{48}\left( -12 -60\zeta(2) -408\zeta(3) -1713\zeta(4) -3426\zeta(5)  \right. \nonumber \\ &\left. \hspace{1em}
-732\zeta(2)\zeta(3) -6291\zeta(6) -804\zeta(3)^2 -38640\zeta(7) -8700\zeta(2)\zeta(5)  \right. \nonumber \\ &\left. \hspace{1em}
-17820\zeta(3)\zeta(4) +51583\zeta(8) -2832\zeta(2)\zeta(3)^2 +34344\zeta(3)\zeta(5)  \right. \nonumber \\ &\left. \hspace{1em}
+3984 M(2,6)\right) \label{eq11239} \\
\sum_{k=1}^\infty \frac{H(k)^{6}}{k(k+1)^{3}(k+2)} \sumend &= \frac{1}{24}\left( -12 -60\zeta(2) -408\zeta(3) -1713\zeta(4) -3426\zeta(5)  \right. \nonumber \\ &\left. \hspace{1em}
-732\zeta(2)\zeta(3) -6291\zeta(6) -804\zeta(3)^2 +7728\zeta(7) +1740\zeta(2)\zeta(5)  \right. \nonumber \\ &\left. \hspace{1em}
+3564\zeta(3)\zeta(4) +6146\zeta(9) -12582\zeta(3)\zeta(6) +5832\zeta(4)\zeta(5) +1953\zeta(2)\zeta(7)  \right. \nonumber \\ &\left. \hspace{1em}
-536\zeta(3)^3\right) \label{eq11240} \\
\sum_{k=1}^\infty \frac{H(k)^{6}}{(k+1)^{4}(k+2)} \sumend &= \frac{1}{384}\left( -384 -1920\zeta(2) -13056\zeta(3) -54816\zeta(4) -109632\zeta(5)  \right. \nonumber \\ &\left. \hspace{1em}
-23424\zeta(2)\zeta(3) -201312\zeta(6) -25728\zeta(3)^2 +272432\zeta(8) -14784\zeta(2)\zeta(3)^2  \right. \nonumber \\ &\left. \hspace{1em}
+181248\zeta(3)\zeta(5) +20928 M(2,6) +98336\zeta(9) -201312\zeta(3)\zeta(6) +93312\zeta(4)\zeta(5)  \right. \nonumber \\ &\left. \hspace{1em}
+31248\zeta(2)\zeta(7) -8576\zeta(3)^3 -779835\zeta(10) +490704\zeta(3)\zeta(7)  \right. \nonumber \\ &\left. \hspace{1em}
-245544\zeta(3)^2\zeta(4) +15600\zeta(2)\zeta(3)\zeta(5) +339864\zeta(5)^2 +94728 M(2,8)  \right. \nonumber \\ &\left. \hspace{1em}
+45120\zeta(2) M(2,6)\right) \label{eq11241} \\
\sum_{k=1}^\infty \frac{H(k)^{6}}{k^{3}(k+2)^{2}} \sumend &= \frac{1}{192}\left(  204 +756\zeta(2) +4344\zeta(3) +14427\zeta(4) +20382\zeta(5)  \right. \nonumber \\ &\left. \hspace{1em}
+4644\zeta(2)\zeta(3) +18570\zeta(6) +2940\zeta(3)^2 +6489\zeta(7) +1116\zeta(2)\zeta(5)  \right. \nonumber \\ &\left. \hspace{1em}
+5940\zeta(3)\zeta(4) -52085\zeta(8) +2892\zeta(2)\zeta(3)^2 -34704\zeta(3)\zeta(5) -4044 M(2,6)  \right. \nonumber \\ &\left. \hspace{1em}
+14948\zeta(9) -26244\zeta(3)\zeta(6) +12096\zeta(4)\zeta(5) +3906\zeta(2)\zeta(7)  \right. \nonumber \\ &\left. \hspace{1em}
-1088\zeta(3)^3\right) \label{eq11242} \\
\sum_{k=1}^\infty \frac{H(k)^{6}}{k^{2}(k+1)(k+2)^{2}} \sumend &= \frac{1}{96}\left(  216 +816\zeta(2) +4752\zeta(3) +16140\zeta(4) +23808\zeta(5)  \right. \nonumber \\ &\left. \hspace{1em}
+5376\zeta(2)\zeta(3) +24861\zeta(6) +3744\zeta(3)^2 -47607\zeta(7) -11064\zeta(2)\zeta(5)  \right. \nonumber \\ &\left. \hspace{1em}
-19008\zeta(3)\zeta(4) +502\zeta(8) -60\zeta(2)\zeta(3)^2 +360\zeta(3)\zeta(5) +60 M(2,6)\right) \label{eq11243} \\
\sum_{k=1}^\infty \frac{H(k)^{6}}{k(k+1)^{2}(k+2)^{2}} \sumend &= \frac{-1}{16}\left( -76 -292\zeta(2) -1720\zeta(3) -5951\zeta(4) -9078\zeta(5)  \right. \nonumber \\ &\left. \hspace{1em}
-2036\zeta(2)\zeta(3) -10384\zeta(6) -1516\zeta(3)^2 +2989\zeta(7) +788\zeta(2)\zeta(5)  \right. \nonumber \\ &\left. \hspace{1em}
+396\zeta(3)\zeta(4) +17027\zeta(8) -924\zeta(2)\zeta(3)^2 +11328\zeta(3)\zeta(5) +1308 M(2,6)\right) \label{eq11244} \\
\sum_{k=1}^\infty \frac{H(k)^{6}}{(k+1)^{3}(k+2)^{2}} \sumend &= \frac{1}{24}\left(  240 +936\zeta(2) +5568\zeta(3) +19566\zeta(4) +30660\zeta(5)  \right. \nonumber \\ &\left. \hspace{1em}
+6840\zeta(2)\zeta(3) +37443\zeta(6) +5352\zeta(3)^2 -16695\zeta(7) -4104\zeta(2)\zeta(5)  \right. \nonumber \\ &\left. \hspace{1em}
-4752\zeta(3)\zeta(4) -51081\zeta(8) +2772\zeta(2)\zeta(3)^2 -33984\zeta(3)\zeta(5) -3924 M(2,6)  \right. \nonumber \\ &\left. \hspace{1em}
-6146\zeta(9) +12582\zeta(3)\zeta(6) -5832\zeta(4)\zeta(5) -1953\zeta(2)\zeta(7) +536\zeta(3)^3\right) \label{eq11245}
\end{align}
 
\begin{align}
\sum_{k=1}^\infty \frac{H(k)^{6}}{k^{2}(k+2)^{3}} \sumend &= \frac{1}{192}\left( -1716 -4644\zeta(2) -22296\zeta(3) -56835\zeta(4) -51078\zeta(5)  \right. \nonumber \\ &\left. \hspace{1em}
-12276\zeta(2)\zeta(3) +6129\zeta(6) -444\zeta(3)^2 +33786\zeta(7) +7596\zeta(2)\zeta(5)  \right. \nonumber \\ &\left. \hspace{1em}
+21636\zeta(3)\zeta(4) +117204\zeta(8) +14808\zeta(2)\zeta(3)^2 -38088\zeta(3)\zeta(5) -11496 M(2,6)  \right. \nonumber \\ &\left. \hspace{1em}
-12292\zeta(9) +25164\zeta(3)\zeta(6) -11664\zeta(4)\zeta(5) -3906\zeta(2)\zeta(7)  \right. \nonumber \\ &\left. \hspace{1em}
+1072\zeta(3)^3\right) \label{eq11246} \\
\sum_{k=1}^\infty \frac{H(k)^{6}}{k(k+1)(k+2)^{3}} \sumend &= \frac{-1}{96}\left(  1932 +5460\zeta(2) +27048\zeta(3) +72975\zeta(4) +74886\zeta(5)  \right. \nonumber \\ &\left. \hspace{1em}
+17652\zeta(2)\zeta(3) +18732\zeta(6) +4188\zeta(3)^2 -81393\zeta(7) -18660\zeta(2)\zeta(5)  \right. \nonumber \\ &\left. \hspace{1em}
-40644\zeta(3)\zeta(4) -116702\zeta(8) -14868\zeta(2)\zeta(3)^2 +38448\zeta(3)\zeta(5) +11556 M(2,6)  \right. \nonumber \\ &\left. \hspace{1em}
+12292\zeta(9) -25164\zeta(3)\zeta(6) +11664\zeta(4)\zeta(5) +3906\zeta(2)\zeta(7)  \right. \nonumber \\ &\left. \hspace{1em}
-1072\zeta(3)^3\right) \label{eq11247} \\
\sum_{k=1}^\infty \frac{H(k)^{6}}{(k+1)^{2}(k+2)^{3}} \sumend &= \frac{1}{48}\left( -2160 -6336\zeta(2) -32208\zeta(3) -90828\zeta(4) -102120\zeta(5)  \right. \nonumber \\ &\left. \hspace{1em}
-23760\zeta(2)\zeta(3) -49884\zeta(6) -8736\zeta(3)^2 +90360\zeta(7) +21024\zeta(2)\zeta(5)  \right. \nonumber \\ &\left. \hspace{1em}
+41832\zeta(3)\zeta(4) +167783\zeta(8) +12096\zeta(2)\zeta(3)^2 -4464\zeta(3)\zeta(5) -7632 M(2,6)  \right. \nonumber \\ &\left. \hspace{1em}
-12292\zeta(9) +25164\zeta(3)\zeta(6) -11664\zeta(4)\zeta(5) -3906\zeta(2)\zeta(7)  \right. \nonumber \\ &\left. \hspace{1em}
+1072\zeta(3)^3\right) \label{eq11248} \\
\sum_{k=1}^\infty \frac{H(k)^{6}}{k(k+2)^{4}} \sumend &= \frac{-1}{768}\left( -38352 -74928\zeta(2) -297696\zeta(3) -562980\zeta(4) -259368\zeta(5)  \right. \nonumber \\ &\left. \hspace{1em}
-60336\zeta(2)\zeta(3) +452160\zeta(6) +73200\zeta(3)^2 +297468\zeta(7) +51312\zeta(2)\zeta(5)  \right. \nonumber \\ &\left. \hspace{1em}
+321840\zeta(3)\zeta(4) +358960\zeta(8) +75504\zeta(2)\zeta(3)^2 -270912\zeta(3)\zeta(5) -59568 M(2,6)  \right. \nonumber \\ &\left. \hspace{1em}
+178672\zeta(9) -311664\zeta(3)\zeta(6) +105408\zeta(4)\zeta(5) +60696\zeta(2)\zeta(7) -7232\zeta(3)^3  \right. \nonumber \\ &\left. \hspace{1em}
-779835\zeta(10) +490704\zeta(3)\zeta(7) -245544\zeta(3)^2\zeta(4) +15600\zeta(2)\zeta(3)\zeta(5)  \right. \nonumber \\ &\left. \hspace{1em}
+339864\zeta(5)^2 +94728 M(2,8) +45120\zeta(2) M(2,6)\right) \label{eq11249} \\
\sum_{k=1}^\infty \frac{H(k)^{6}}{(k+1)(k+2)^{4}} \sumend &= \frac{1}{384}\left(  46080 +96768\zeta(2) +405888\zeta(3) +854880\zeta(4) +558912\zeta(5)  \right. \nonumber \\ &\left. \hspace{1em}
+130944\zeta(2)\zeta(3) -377232\zeta(6) -56448\zeta(3)^2 -623040\zeta(7) -125952\zeta(2)\zeta(5)  \right. \nonumber \\ &\left. \hspace{1em}
-484416\zeta(3)\zeta(4) -825768\zeta(8) -134976\zeta(2)\zeta(3)^2 +424704\zeta(3)\zeta(5)  \right. \nonumber \\ &\left. \hspace{1em}
+105792 M(2,6) -129504\zeta(9) +211008\zeta(3)\zeta(6) -58752\zeta(4)\zeta(5) -45072\zeta(2)\zeta(7)  \right. \nonumber \\ &\left. \hspace{1em}
+2944\zeta(3)^3 +779835\zeta(10) -490704\zeta(3)\zeta(7) +245544\zeta(3)^2\zeta(4)  \right. \nonumber \\ &\left. \hspace{1em}
-15600\zeta(2)\zeta(3)\zeta(5) -339864\zeta(5)^2 -94728 M(2,8) -45120\zeta(2) M(2,6)\right) \label{eq11250}
\end{align}
 
\begin{align}
\sum_{k=1}^\infty \frac{H(k)^{6}}{(k+2)^{5}} \sumend &= \frac{-1}{384}\left(  80640 +112896\zeta(2) +368256\zeta(3) +494496\zeta(4) +35136\zeta(5)  \right. \nonumber \\ &\left. \hspace{1em}
-5760\zeta(2)\zeta(3) -565536\zeta(6) -111360\zeta(3)^2 -197280\zeta(7) -31104\zeta(2)\zeta(5)  \right. \nonumber \\ &\left. \hspace{1em}
-202752\zeta(3)\zeta(4) +471672\zeta(8) +133440\zeta(2)\zeta(3)^2 -549504\zeta(3)\zeta(5)  \right. \nonumber \\ &\left. \hspace{1em}
-116160 M(2,6) -144320\zeta(9) +364800\zeta(3)\zeta(6) -258624\zeta(4)\zeta(5) -36000\zeta(2)\zeta(7)  \right. \nonumber \\ &\left. \hspace{1em}
+26880\zeta(3)^3 -449109\zeta(10) +387360\zeta(3)\zeta(7) +9720\zeta(3)^2\zeta(4)  \right. \nonumber \\ &\left. \hspace{1em}
-124560\zeta(2)\zeta(3)\zeta(5) +122328\zeta(5)^2 +68040 M(2,8) +11520\zeta(2) M(2,6) -1373598\zeta(11)  \right. \nonumber \\ &\left. \hspace{1em}
-149024\zeta(2)\zeta(9) +524968\zeta(3)\zeta(8) +724416\zeta(4)\zeta(7) +365168\zeta(5)\zeta(6)  \right. \nonumber \\ &\left. \hspace{1em}
+36160\zeta(2)\zeta(3)^3 -240768\zeta(3)^2\zeta(5) +16320\zeta(3) M(2,6) +46720 M(3,8)\right) \label{eq11251} \\
\sum_{k=1}^\infty \frac{H(k)^{7}}{k^{4}} \sumend &= \frac{1}{1152}\left( -16370805\zeta(11) -1684144\zeta(2)\zeta(9) -5889744\zeta(3)\zeta(8)  \right. \nonumber \\ &\left. \hspace{1em}
+10724760\zeta(4)\zeta(7) +10480104\zeta(5)\zeta(6) -844032\zeta(2)\zeta(3)^3  \right. \nonumber \\ &\left. \hspace{1em}
+2330496\zeta(3)^2\zeta(5) +1431360\zeta(3) M(2,6) +630336 M(3,8)\right) \label{eq11252} \\
\sum_{k=1}^\infty \frac{H(k)^{7}}{k^{3}(k+1)} \sumend &= \frac{-1}{23040}\left( -153310720\zeta(8) -3870720\zeta(2)\zeta(3)^2 -35078400\zeta(3)\zeta(5)  \right. \nonumber \\ &\left. \hspace{1em}
+88429120\zeta(9) +28372800\zeta(3)\zeta(6) +45812160\zeta(4)\zeta(5) +18933120\zeta(2)\zeta(7)  \right. \nonumber \\ &\left. \hspace{1em}
+1290240\zeta(3)^3 +149534919\zeta(10) -92839680\zeta(3)\zeta(7) +52912440\zeta(3)^2\zeta(4)  \right. \nonumber \\ &\left. \hspace{1em}
-6345360\zeta(2)\zeta(3)\zeta(5) -69396840\zeta(5)^2 -18196920 M(2,8) -9072000\zeta(2) M(2,6)\right) \label{eq11253} \\
\sum_{k=1}^\infty \frac{H(k)^{7}}{k^{2}(k+1)^{2}} \sumend &= \frac{-1}{72}\left(  958192\zeta(8) +24192\zeta(2)\zeta(3)^2 +219240\zeta(3)\zeta(5)  \right. \nonumber \\ &\left. \hspace{1em}
-545743\zeta(9) -177330\zeta(3)\zeta(6) -284436\zeta(4)\zeta(5) -118332\zeta(2)\zeta(7)  \right. \nonumber \\ &\left. \hspace{1em}
-8064\zeta(3)^3\right) \label{eq11254} \\
\sum_{k=1}^\infty \frac{H(k)^{7}}{k(k+1)^{3}} \sumend &= \frac{1}{23040}\left(  153310720\zeta(8) +3870720\zeta(2)\zeta(3)^2 +35078400\zeta(3)\zeta(5)  \right. \nonumber \\ &\left. \hspace{1em}
-86208640\zeta(9) -28372800\zeta(3)\zeta(6) -45207360\zeta(4)\zeta(5) -18933120\zeta(2)\zeta(7)  \right. \nonumber \\ &\left. \hspace{1em}
-1290240\zeta(3)^3 -149375151\zeta(10) +89769600\zeta(3)\zeta(7) -52206840\zeta(3)^2\zeta(4)  \right. \nonumber \\ &\left. \hspace{1em}
+7514640\zeta(2)\zeta(3)\zeta(5) +67262760\zeta(5)^2 +17612280 M(2,8) +9072000\zeta(2) M(2,6)\right) \label{eq11255} \\
\sum_{k=1}^\infty \frac{H(k)^{7}}{(k+1)^{4}} \sumend &= \frac{1}{1152}\left( -16196565\zeta(11) -1630384\zeta(2)\zeta(9) -5721072\zeta(3)\zeta(8)  \right. \nonumber \\ &\left. \hspace{1em}
+10468728\zeta(4)\zeta(7) +10317144\zeta(5)\zeta(6) -837312\zeta(2)\zeta(3)^3  \right. \nonumber \\ &\left. \hspace{1em}
+2330496\zeta(3)^2\zeta(5) +1411200\zeta(3) M(2,6) +616896 M(3,8)\right) \label{eq11256} \\
\sum_{k=1}^\infty \frac{H(k)^{7}}{k^{3}(k+2)} \sumend &= \frac{-1}{46080}\left( -5760 -34560\zeta(2) -288000\zeta(3) -1546560\zeta(4) -4340160\zeta(5)  \right. \nonumber \\ &\left. \hspace{1em}
-927360\zeta(2)\zeta(3) -14115480\zeta(6) -1834560\zeta(3)^2 -12782160\zeta(7) -2923200\zeta(2)\zeta(5)  \right. \nonumber \\ &\left. \hspace{1em}
-5957280\zeta(3)\zeta(4) -38327680\zeta(8) -967680\zeta(2)\zeta(3)^2 -8769600\zeta(3)\zeta(5)  \right. \nonumber \\ &\left. \hspace{1em}
+44214560\zeta(9) +14186400\zeta(3)\zeta(6) +22906080\zeta(4)\zeta(5) +9466560\zeta(2)\zeta(7)  \right. \nonumber \\ &\left. \hspace{1em}
+645120\zeta(3)^3 +149534919\zeta(10) -92839680\zeta(3)\zeta(7) +52912440\zeta(3)^2\zeta(4)  \right. \nonumber \\ &\left. \hspace{1em}
-6345360\zeta(2)\zeta(3)\zeta(5) -69396840\zeta(5)^2 -18196920 M(2,8) -9072000\zeta(2) M(2,6)\right) \label{eq11257}
\end{align}
 
\begin{align}
\sum_{k=1}^\infty \frac{H(k)^{7}}{k^{2}(k+1)(k+2)} \sumend &= \frac{-1}{576}\left( -144 -864\zeta(2) -7200\zeta(3) -38664\zeta(4) -108504\zeta(5)  \right. \nonumber \\ &\left. \hspace{1em}
-23184\zeta(2)\zeta(3) -352887\zeta(6) -45864\zeta(3)^2 -319554\zeta(7) -73080\zeta(2)\zeta(5)  \right. \nonumber \\ &\left. \hspace{1em}
-148932\zeta(3)\zeta(4) +2874576\zeta(8) +72576\zeta(2)\zeta(3)^2 +657720\zeta(3)\zeta(5)  \right. \nonumber \\ &\left. \hspace{1em}
-1105364\zeta(9) -354660\zeta(3)\zeta(6) -572652\zeta(4)\zeta(5) -236664\zeta(2)\zeta(7)  \right. \nonumber \\ &\left. \hspace{1em}
-16128\zeta(3)^3\right) \label{eq11258} \\
\sum_{k=1}^\infty \frac{H(k)^{7}}{k(k+1)^{2}(k+2)} \sumend &= \frac{1}{288}\left(  144 +864\zeta(2) +7200\zeta(3) +38664\zeta(4) +108504\zeta(5)  \right. \nonumber \\ &\left. \hspace{1em}
+23184\zeta(2)\zeta(3) +352887\zeta(6) +45864\zeta(3)^2 +319554\zeta(7) +73080\zeta(2)\zeta(5)  \right. \nonumber \\ &\left. \hspace{1em}
+148932\zeta(3)\zeta(4) +958192\zeta(8) +24192\zeta(2)\zeta(3)^2 +219240\zeta(3)\zeta(5)  \right. \nonumber \\ &\left. \hspace{1em}
-1077608\zeta(9) -354660\zeta(3)\zeta(6) -565092\zeta(4)\zeta(5) -236664\zeta(2)\zeta(7)  \right. \nonumber \\ &\left. \hspace{1em}
-16128\zeta(3)^3\right) \label{eq11259} \\
\sum_{k=1}^\infty \frac{H(k)^{7}}{(k+1)^{3}(k+2)} \sumend &= \frac{-1}{23040}\left( -23040 -138240\zeta(2) -1152000\zeta(3) -6186240\zeta(4)  \right. \nonumber \\ &\left. \hspace{1em}
-17360640\zeta(5) -3709440\zeta(2)\zeta(3) -56461920\zeta(6) -7338240\zeta(3)^2 -51128640\zeta(7)  \right. \nonumber \\ &\left. \hspace{1em}
-11692800\zeta(2)\zeta(5) -23829120\zeta(3)\zeta(4) +86208640\zeta(9) +28372800\zeta(3)\zeta(6)  \right. \nonumber \\ &\left. \hspace{1em}
+45207360\zeta(4)\zeta(5) +18933120\zeta(2)\zeta(7) +1290240\zeta(3)^3 -149375151\zeta(10)  \right. \nonumber \\ &\left. \hspace{1em}
+89769600\zeta(3)\zeta(7) -52206840\zeta(3)^2\zeta(4) +7514640\zeta(2)\zeta(3)\zeta(5)  \right. \nonumber \\ &\left. \hspace{1em}
+67262760\zeta(5)^2 +17612280 M(2,8) +9072000\zeta(2) M(2,6)\right) \label{eq11260} \\
\sum_{k=1}^\infty \frac{H(k)^{7}}{k^{2}(k+2)^{2}} \sumend &= \frac{1}{576}\left( -1296 -5904\zeta(2) -42192\zeta(3) -185004\zeta(4) -396216\zeta(5)  \right. \nonumber \\ &\left. \hspace{1em}
-87696\zeta(2)\zeta(3) -878271\zeta(6) -122472\zeta(3)^2 -288513\zeta(7) -59976\zeta(2)\zeta(5)  \right. \nonumber \\ &\left. \hspace{1em}
-198072\zeta(3)\zeta(4) -243058\zeta(8) -63000\zeta(2)\zeta(3)^2 +256536\zeta(3)\zeta(5) +54936 M(2,6)  \right. \nonumber \\ &\left. \hspace{1em}
+1091486\zeta(9) +354660\zeta(3)\zeta(6) +568872\zeta(4)\zeta(5) +236664\zeta(2)\zeta(7)  \right. \nonumber \\ &\left. \hspace{1em}
+16128\zeta(3)^3\right) \label{eq11261} \\
\sum_{k=1}^\infty \frac{H(k)^{7}}{k(k+1)(k+2)^{2}} \sumend &= \frac{1}{576}\left( -2736 -12672\zeta(2) -91584\zeta(3) -408672\zeta(4) -900936\zeta(5)  \right. \nonumber \\ &\left. \hspace{1em}
-198576\zeta(2)\zeta(3) -2109429\zeta(6) -290808\zeta(3)^2 -896580\zeta(7) -193032\zeta(2)\zeta(5)  \right. \nonumber \\ &\left. \hspace{1em}
-545076\zeta(3)\zeta(4) +2388460\zeta(8) -53424\zeta(2)\zeta(3)^2 +1170792\zeta(3)\zeta(5)  \right. \nonumber \\ &\left. \hspace{1em}
+109872 M(2,6) +1077608\zeta(9) +354660\zeta(3)\zeta(6) +565092\zeta(4)\zeta(5) +236664\zeta(2)\zeta(7)  \right. \nonumber \\ &\left. \hspace{1em}
+16128\zeta(3)^3\right) \label{eq11262} \\
\sum_{k=1}^\infty \frac{H(k)^{7}}{(k+1)^{2}(k+2)^{2}} \sumend &= \frac{1}{144}\left( -1440 -6768\zeta(2) -49392\zeta(3) -223668\zeta(4) -504720\zeta(5)  \right. \nonumber \\ &\left. \hspace{1em}
-110880\zeta(2)\zeta(3) -1231158\zeta(6) -168336\zeta(3)^2 -608067\zeta(7) -133056\zeta(2)\zeta(5)  \right. \nonumber \\ &\left. \hspace{1em}
-347004\zeta(3)\zeta(4) +715134\zeta(8) -38808\zeta(2)\zeta(3)^2 +475776\zeta(3)\zeta(5) +54936 M(2,6)  \right. \nonumber \\ &\left. \hspace{1em}
+1077608\zeta(9) +354660\zeta(3)\zeta(6) +565092\zeta(4)\zeta(5) +236664\zeta(2)\zeta(7)  \right. \nonumber \\ &\left. \hspace{1em}
+16128\zeta(3)^3\right) \label{eq11263}
\end{align}
 
\begin{align}
\sum_{k=1}^\infty \frac{H(k)^{7}}{k(k+2)^{3}} \sumend &= \frac{1}{46080}\left(  927360 +3179520\zeta(2) +19376640\zeta(3) +68423040\zeta(4)  \right. \nonumber \\ &\left. \hspace{1em}
+108054720\zeta(5) +24635520\zeta(2)\zeta(3) +140607960\zeta(6) +21107520\zeta(3)^2 -39775680\zeta(7)  \right. \nonumber \\ &\left. \hspace{1em}
-9979200\zeta(2)\zeta(5) -13134240\zeta(3)\zeta(4) -243547760\zeta(8) -19353600\zeta(2)\zeta(3)^2  \right. \nonumber \\ &\left. \hspace{1em}
+16269120\zeta(3)\zeta(5) +12821760 M(2,6) -1803200\zeta(9) -98737440\zeta(3)\zeta(6)  \right. \nonumber \\ &\left. \hspace{1em}
+16587360\zeta(4)\zeta(5) +3657600\zeta(2)\zeta(7) -4247040\zeta(3)^3 -149375151\zeta(10)  \right. \nonumber \\ &\left. \hspace{1em}
+89769600\zeta(3)\zeta(7) -52206840\zeta(3)^2\zeta(4) +7514640\zeta(2)\zeta(3)\zeta(5)  \right. \nonumber \\ &\left. \hspace{1em}
+67262760\zeta(5)^2 +17612280 M(2,8) +9072000\zeta(2) M(2,6)\right) \label{eq11264} \\
\sum_{k=1}^\infty \frac{H(k)^{7}}{(k+1)(k+2)^{3}} \sumend &= \frac{1}{23040}\left(  1036800 +3686400\zeta(2) +23040000\zeta(3) +84769920\zeta(4)  \right. \nonumber \\ &\left. \hspace{1em}
+144092160\zeta(5) +32578560\zeta(2)\zeta(3) +224985120\zeta(6) +32739840\zeta(3)^2 -3912480\zeta(7)  \right. \nonumber \\ &\left. \hspace{1em}
-2257920\zeta(2)\zeta(5) +8668800\zeta(3)\zeta(4) -339086160\zeta(8) -17216640\zeta(2)\zeta(3)^2  \right. \nonumber \\ &\left. \hspace{1em}
-30562560\zeta(3)\zeta(5) +8426880 M(2,6) -44907520\zeta(9) -112923840\zeta(3)\zeta(6)  \right. \nonumber \\ &\left. \hspace{1em}
-6016320\zeta(4)\zeta(5) -5808960\zeta(2)\zeta(7) -4892160\zeta(3)^3 -149375151\zeta(10)  \right. \nonumber \\ &\left. \hspace{1em}
+89769600\zeta(3)\zeta(7) -52206840\zeta(3)^2\zeta(4) +7514640\zeta(2)\zeta(3)\zeta(5)  \right. \nonumber \\ &\left. \hspace{1em}
+67262760\zeta(5)^2 +17612280 M(2,8) +9072000\zeta(2) M(2,6)\right) \label{eq11265} \\
\sum_{k=1}^\infty \frac{H(k)^{7}}{(k+2)^{4}} \sumend &= \frac{1}{1152}\left( -138240 -354816\zeta(2) -1838592\zeta(3) -5175072\zeta(4) -5821632\zeta(5)  \right. \nonumber \\ &\left. \hspace{1em}
-1338624\zeta(2)\zeta(3) -2929008\zeta(6) -443520\zeta(3)^2 +4509648\zeta(7) +959616\zeta(2)\zeta(5)  \right. \nonumber \\ &\left. \hspace{1em}
+3108672\zeta(3)\zeta(4) +13269200\zeta(8) +1725696\zeta(2)\zeta(3)^2 -4491648\zeta(3)\zeta(5)  \right. \nonumber \\ &\left. \hspace{1em}
-1314432 M(2,6) +327264\zeta(9) -101808\zeta(3)\zeta(6) -362880\zeta(4)\zeta(5) +145152\zeta(2)\zeta(7)  \right. \nonumber \\ &\left. \hspace{1em}
+59136\zeta(3)^3 -16376535\zeta(10) +10304784\zeta(3)\zeta(7) -5156424\zeta(3)^2\zeta(4)  \right. \nonumber \\ &\left. \hspace{1em}
+327600\zeta(2)\zeta(3)\zeta(5) +7137144\zeta(5)^2 +1989288 M(2,8) +947520\zeta(2) M(2,6)  \right. \nonumber \\ &\left. \hspace{1em}
-16196565\zeta(11) -1630384\zeta(2)\zeta(9) -5721072\zeta(3)\zeta(8) +10468728\zeta(4)\zeta(7)  \right. \nonumber \\ &\left. \hspace{1em}
+10317144\zeta(5)\zeta(6) -837312\zeta(2)\zeta(3)^3 +2330496\zeta(3)^2\zeta(5) +1411200\zeta(3) M(2,6)  \right. \nonumber \\ &\left. \hspace{1em}
+616896 M(3,8)\right) \label{eq11266} \\
\sum_{k=1}^\infty \frac{H(k)^{8}}{k^{3}} \sumend &= \frac{1}{72}\left( -2824380\zeta(11) -277304\zeta(2)\zeta(9) -1926401\zeta(3)\zeta(8)  \right. \nonumber \\ &\left. \hspace{1em}
+1998972\zeta(4)\zeta(7) +2270310\zeta(5)\zeta(6) -243648\zeta(2)\zeta(3)^3 +803808\zeta(3)^2\zeta(5)  \right. \nonumber \\ &\left. \hspace{1em}
+341280\zeta(3) M(2,6) +113760 M(3,8)\right) \label{eq11267} \\
\sum_{k=1}^\infty \frac{H(k)^{8}}{k^{2}(k+1)} \sumend &= \frac{1}{480}\left( -13336000\zeta(9) -7093200\zeta(3)\zeta(6) -6432000\zeta(4)\zeta(5)  \right. \nonumber \\ &\left. \hspace{1em}
-2807280\zeta(2)\zeta(7) -322560\zeta(3)^3 +18741581\zeta(10) +6689520\zeta(3)\zeta(7)  \right. \nonumber \\ &\left. \hspace{1em}
-524640\zeta(3)^2\zeta(4) +1452480\zeta(2)\zeta(3)\zeta(5) +4247040\zeta(5)^2 +485280 M(2,8)  \right. \nonumber \\ &\left. \hspace{1em}
+299520\zeta(2) M(2,6)\right) \label{eq11268}
\end{align}
 
\begin{align}
\sum_{k=1}^\infty \frac{H(k)^{8}}{k(k+1)^{2}} \sumend &= \frac{-1}{240}\left( -6668000\zeta(9) -3546600\zeta(3)\zeta(6) -3216000\zeta(4)\zeta(5)  \right. \nonumber \\ &\left. \hspace{1em}
-1403640\zeta(2)\zeta(7) -161280\zeta(3)^3 +9295879\zeta(10) +3314520\zeta(3)\zeta(7)  \right. \nonumber \\ &\left. \hspace{1em}
-258540\zeta(3)^2\zeta(4) +733800\zeta(2)\zeta(3)\zeta(5) +2098980\zeta(5)^2 +238860 M(2,8)  \right. \nonumber \\ &\left. \hspace{1em}
+149760\zeta(2) M(2,6)\right) \label{eq11269} \\
\sum_{k=1}^\infty \frac{H(k)^{8}}{(k+1)^{3}} \sumend &= \frac{-1}{72}\left( -2839707\zeta(11) -274424\zeta(2)\zeta(9) -1906367\zeta(3)\zeta(8)  \right. \nonumber \\ &\left. \hspace{1em}
+1976076\zeta(4)\zeta(7) +2252940\zeta(5)\zeta(6) -242208\zeta(2)\zeta(3)^3 +798912\zeta(3)^2\zeta(5)  \right. \nonumber \\ &\left. \hspace{1em}
+339120\zeta(3) M(2,6) +113040 M(3,8)\right) \label{eq11270} \\
\sum_{k=1}^\infty \frac{H(k)^{8}}{k^{2}(k+2)} \sumend &= \frac{1}{2880}\left( -720 -5040\zeta(2) -49680\zeta(3) -324000\zeta(4) -1162800\zeta(5)  \right. \nonumber \\ &\left. \hspace{1em}
-247680\zeta(2)\zeta(3) -5303460\zeta(6) -692640\zeta(3)^2 -8496540\zeta(7) -1931040\zeta(2)\zeta(5)  \right. \nonumber \\ &\left. \hspace{1em}
-4042800\zeta(3)\zeta(4) -19063670\zeta(8) -483840\zeta(2)\zeta(3)^2 -4368960\zeta(3)\zeta(5)  \right. \nonumber \\ &\left. \hspace{1em}
-20004000\zeta(9) -10639800\zeta(3)\zeta(6) -9648000\zeta(4)\zeta(5) -4210920\zeta(2)\zeta(7)  \right. \nonumber \\ &\left. \hspace{1em}
-483840\zeta(3)^3 +56224743\zeta(10) +20068560\zeta(3)\zeta(7) -1573920\zeta(3)^2\zeta(4)  \right. \nonumber \\ &\left. \hspace{1em}
+4357440\zeta(2)\zeta(3)\zeta(5) +12741120\zeta(5)^2 +1455840 M(2,8) +898560\zeta(2) M(2,6)\right) \label{eq11271} \\
\sum_{k=1}^\infty \frac{H(k)^{8}}{k(k+1)(k+2)} \sumend &= \frac{-1}{144}\left(  72 +504\zeta(2) +4968\zeta(3) +32400\zeta(4) +116280\zeta(5)  \right. \nonumber \\ &\left. \hspace{1em}
+24768\zeta(2)\zeta(3) +530346\zeta(6) +69264\zeta(3)^2 +849654\zeta(7) +193104\zeta(2)\zeta(5)  \right. \nonumber \\ &\left. \hspace{1em}
+404280\zeta(3)\zeta(4) +1906367\zeta(8) +48384\zeta(2)\zeta(3)^2 +436896\zeta(3)\zeta(5)  \right. \nonumber \\ &\left. \hspace{1em}
-2000400\zeta(9) -1063980\zeta(3)\zeta(6) -964800\zeta(4)\zeta(5) -421092\zeta(2)\zeta(7)  \right. \nonumber \\ &\left. \hspace{1em}
-48384\zeta(3)^3\right) \label{eq11272} \\
\sum_{k=1}^\infty \frac{H(k)^{8}}{(k+1)^{2}(k+2)} \sumend &= \frac{1}{720}\left( -720 -5040\zeta(2) -49680\zeta(3) -324000\zeta(4) -1162800\zeta(5)  \right. \nonumber \\ &\left. \hspace{1em}
-247680\zeta(2)\zeta(3) -5303460\zeta(6) -692640\zeta(3)^2 -8496540\zeta(7) -1931040\zeta(2)\zeta(5)  \right. \nonumber \\ &\left. \hspace{1em}
-4042800\zeta(3)\zeta(4) -19063670\zeta(8) -483840\zeta(2)\zeta(3)^2 -4368960\zeta(3)\zeta(5)  \right. \nonumber \\ &\left. \hspace{1em}
+27887637\zeta(10) +9943560\zeta(3)\zeta(7) -775620\zeta(3)^2\zeta(4) +2201400\zeta(2)\zeta(3)\zeta(5)  \right. \nonumber \\ &\left. \hspace{1em}
+6296940\zeta(5)^2 +716580 M(2,8) +449280\zeta(2) M(2,6)\right) \label{eq11273} \\
\sum_{k=1}^\infty \frac{H(k)^{8}}{k(k+2)^{2}} \sumend &= \frac{1}{1440}\left(  6840 +37080\zeta(2) +317160\zeta(3) +1726560\zeta(4) +4930200\zeta(5)  \right. \nonumber \\ &\left. \hspace{1em}
+1074240\zeta(2)\zeta(3) +16758210\zeta(6) +2273040\zeta(3)^2 +16566750\zeta(7)  \right. \nonumber \\ &\left. \hspace{1em}
+3678480\zeta(2)\zeta(5) +8776440\zeta(3)\zeta(4) +14292825\zeta(8) +1501920\zeta(2)\zeta(3)^2  \right. \nonumber \\ &\left. \hspace{1em}
-2962080\zeta(3)\zeta(5) -1098720 M(2,6) -11550160\zeta(9) -1773300\zeta(3)\zeta(6)  \right. \nonumber \\ &\left. \hspace{1em}
-6477840\zeta(4)\zeta(5) -2627820\zeta(2)\zeta(7) -80640\zeta(3)^3 -27887637\zeta(10)  \right. \nonumber \\ &\left. \hspace{1em}
-9943560\zeta(3)\zeta(7) +775620\zeta(3)^2\zeta(4) -2201400\zeta(2)\zeta(3)\zeta(5) -6296940\zeta(5)^2  \right. \nonumber \\ &\left. \hspace{1em}
-716580 M(2,8) -449280\zeta(2) M(2,6)\right) \label{eq11274}
\end{align}
 
\begin{align}
\sum_{k=1}^\infty \frac{H(k)^{8}}{(k+1)(k+2)^{2}} \sumend &= \frac{-1}{720}\left( -7200 -39600\zeta(2) -342000\zeta(3) -1888560\zeta(4) -5511600\zeta(5)  \right. \nonumber \\ &\left. \hspace{1em}
-1198080\zeta(2)\zeta(3) -19409940\zeta(6) -2619360\zeta(3)^2 -20815020\zeta(7)  \right. \nonumber \\ &\left. \hspace{1em}
-4644000\zeta(2)\zeta(5) -10797840\zeta(3)\zeta(4) -23824660\zeta(8) -1743840\zeta(2)\zeta(3)^2  \right. \nonumber \\ &\left. \hspace{1em}
+777600\zeta(3)\zeta(5) +1098720 M(2,6) +21552160\zeta(9) +7093200\zeta(3)\zeta(6)  \right. \nonumber \\ &\left. \hspace{1em}
+11301840\zeta(4)\zeta(5) +4733280\zeta(2)\zeta(7) +322560\zeta(3)^3 +27887637\zeta(10)  \right. \nonumber \\ &\left. \hspace{1em}
+9943560\zeta(3)\zeta(7) -775620\zeta(3)^2\zeta(4) +2201400\zeta(2)\zeta(3)\zeta(5) +6296940\zeta(5)^2  \right. \nonumber \\ &\left. \hspace{1em}
+716580 M(2,8) +449280\zeta(2) M(2,6)\right) \label{eq11275} \\
\sum_{k=1}^\infty \frac{H(k)^{8}}{(k+2)^{3}} \sumend &= \frac{1}{2880}\left( -129600 -541440\zeta(2) -4008960\zeta(3) -18103680\zeta(4) -40584960\zeta(5)  \right. \nonumber \\ &\left. \hspace{1em}
-9020160\zeta(2)\zeta(3) -98753280\zeta(6) -13881600\zeta(3)^2 -48610080\zeta(7)  \right. \nonumber \\ &\left. \hspace{1em}
-10391040\zeta(2)\zeta(5) -28817280\zeta(3)\zeta(4) +74914520\zeta(8) +1733760\zeta(2)\zeta(3)^2  \right. \nonumber \\ &\left. \hspace{1em}
+18086400\zeta(3)\zeta(5) +120960 M(2,6) +65558080\zeta(9) +70648320\zeta(3)\zeta(6)  \right. \nonumber \\ &\left. \hspace{1em}
+25611840\zeta(4)\zeta(5) +12371040\zeta(2)\zeta(7) +3091200\zeta(3)^3 +149375151\zeta(10)  \right. \nonumber \\ &\left. \hspace{1em}
-89769600\zeta(3)\zeta(7) +52206840\zeta(3)^2\zeta(4) -7514640\zeta(2)\zeta(3)\zeta(5)  \right. \nonumber \\ &\left. \hspace{1em}
-67262760\zeta(5)^2 -17612280 M(2,8) -9072000\zeta(2) M(2,6) +113588280\zeta(11)  \right. \nonumber \\ &\left. \hspace{1em}
+10976960\zeta(2)\zeta(9) +76254680\zeta(3)\zeta(8) -79043040\zeta(4)\zeta(7) -90117600\zeta(5)\zeta(6)  \right. \nonumber \\ &\left. \hspace{1em}
+9688320\zeta(2)\zeta(3)^3 -31956480\zeta(3)^2\zeta(5) -13564800\zeta(3) M(2,6) -4521600 M(3,8)\right) \label{eq11276} \\
\sum_{k=1}^\infty \frac{H(k)^{9}}{k^{2}} \sumend &= \frac{-1}{64}\left( -7739347\zeta(11) -2048432\zeta(2)\zeta(9) -5357920\zeta(3)\zeta(8)  \right. \nonumber \\ &\left. \hspace{1em}
-8811792\zeta(4)\zeta(7) -10526056\zeta(5)\zeta(6) +294208\zeta(2)\zeta(3)^3 -2064192\zeta(3)^2\zeta(5)  \right. \nonumber \\ &\left. \hspace{1em}
-540096\zeta(3) M(2,6) -199936 M(3,8)\right) \label{eq11277} \\
\sum_{k=1}^\infty \frac{H(k)^{9}}{k(k+1)} \sumend &= \frac{1}{40}\left(  17039209\zeta(10) +3158190\zeta(3)\zeta(7) +704820\zeta(3)^2\zeta(4)  \right. \nonumber \\ &\left. \hspace{1em}
+928080\zeta(2)\zeta(3)\zeta(5) +1767000\zeta(5)^2 +37320\zeta(2) M(2,6)\right) \label{eq11278} \\
\sum_{k=1}^\infty \frac{H(k)^{9}}{(k+1)^{2}} \sumend &= \frac{1}{64}\left(  7676163\zeta(11) +2050992\zeta(2)\zeta(9) +5357920\zeta(3)\zeta(8)  \right. \nonumber \\ &\left. \hspace{1em}
+8776416\zeta(4)\zeta(7) +10489496\zeta(5)\zeta(6) -292928\zeta(2)\zeta(3)^3 +2058432\zeta(3)^2\zeta(5)  \right. \nonumber \\ &\left. \hspace{1em}
+538176\zeta(3) M(2,6) +199296 M(3,8)\right) \label{eq11279} \\
\sum_{k=1}^\infty \frac{H(k)^{9}}{k(k+2)} \sumend &= \frac{1}{160}\left(  80 +640\zeta(2) +7280\zeta(3) +55780\zeta(4) +243200\zeta(5)  \right. \nonumber \\ &\left. \hspace{1em}
+51600\zeta(2)\zeta(3) +1416700\zeta(6) +185160\zeta(3)^2 +3186110\zeta(7) +721920\zeta(2)\zeta(5)  \right. \nonumber \\ &\left. \hspace{1em}
+1525680\zeta(3)\zeta(4) +12880535\zeta(8) +365400\zeta(2)\zeta(3)^2 +2739360\zeta(3)\zeta(5)  \right. \nonumber \\ &\left. \hspace{1em}
-37320 M(2,6) +9964440\zeta(9) +5312700\zeta(3)\zeta(6) +4812120\zeta(4)\zeta(5)  \right. \nonumber \\ &\left. \hspace{1em}
+2105460\zeta(2)\zeta(7) +241760\zeta(3)^3 +34078418\zeta(10) +6316380\zeta(3)\zeta(7)  \right. \nonumber \\ &\left. \hspace{1em}
+1409640\zeta(3)^2\zeta(4) +1856160\zeta(2)\zeta(3)\zeta(5) +3534000\zeta(5)^2  \right. \nonumber \\ &\left. \hspace{1em}
+74640\zeta(2) M(2,6)\right) \label{eq11280}
\end{align}
 
\begin{align}
\sum_{k=1}^\infty \frac{H(k)^{9}}{(k+1)(k+2)} \sumend &= \frac{-1}{16}\left( -16 -128\zeta(2) -1456\zeta(3) -11156\zeta(4) -48640\zeta(5)  \right. \nonumber \\ &\left. \hspace{1em}
-10320\zeta(2)\zeta(3) -283340\zeta(6) -37032\zeta(3)^2 -637222\zeta(7) -144384\zeta(2)\zeta(5)  \right. \nonumber \\ &\left. \hspace{1em}
-305136\zeta(3)\zeta(4) -2576107\zeta(8) -73080\zeta(2)\zeta(3)^2 -547872\zeta(3)\zeta(5) +7464 M(2,6)  \right. \nonumber \\ &\left. \hspace{1em}
-1992888\zeta(9) -1062540\zeta(3)\zeta(6) -962424\zeta(4)\zeta(5) -421092\zeta(2)\zeta(7)  \right. \nonumber \\ &\left. \hspace{1em}
-48352\zeta(3)^3\right) \label{eq11281} \\
\sum_{k=1}^\infty \frac{H(k)^{9}}{(k+2)^{2}} \sumend &= \frac{-1}{320}\left(  3200 +20160\zeta(2) +200960\zeta(3) +1307040\zeta(4) +4662400\zeta(5)  \right. \nonumber \\ &\left. \hspace{1em}
+1005120\zeta(2)\zeta(3) +21224080\zeta(6) +2841600\zeta(3)^2 +33558720\zeta(7)  \right. \nonumber \\ &\left. \hspace{1em}
+7513920\zeta(2)\zeta(5) +16977600\zeta(3)\zeta(4) +83974040\zeta(8) +3958080\zeta(2)\zeta(3)^2  \right. \nonumber \\ &\left. \hspace{1em}
+9227520\zeta(3)\zeta(5) -1763520 M(2,6) -3246560\zeta(9) +7064400\zeta(3)\zeta(6)  \right. \nonumber \\ &\left. \hspace{1em}
-3355200\zeta(4)\zeta(5) -1044720\zeta(2)\zeta(7) +321920\zeta(3)^3 -111550548\zeta(10)  \right. \nonumber \\ &\left. \hspace{1em}
-39774240\zeta(3)\zeta(7) +3102480\zeta(3)^2\zeta(4) -8805600\zeta(2)\zeta(3)\zeta(5)  \right. \nonumber \\ &\left. \hspace{1em}
-25187760\zeta(5)^2 -2866320 M(2,8) -1797120\zeta(2) M(2,6) -38380815\zeta(11)  \right. \nonumber \\ &\left. \hspace{1em}
-10254960\zeta(2)\zeta(9) -26789600\zeta(3)\zeta(8) -43882080\zeta(4)\zeta(7) -52447480\zeta(5)\zeta(6)  \right. \nonumber \\ &\left. \hspace{1em}
+1464640\zeta(2)\zeta(3)^3 -10292160\zeta(3)^2\zeta(5) -2690880\zeta(3) M(2,6) -996480 M(3,8)\right) \label{eq11282}
\end{align}

\newpage
Formulas for order $r = m + n = 12$:
\begin{align}
\sum_{k=1}^\infty \frac{H(k)}{k^{11}} \sumend &= \frac{-1}{4}\left( -13\zeta(12) +4\zeta(3)\zeta(9) +4\zeta(5)\zeta(7)\right) \label{eq12001} \\
\sum_{k=1}^\infty \frac{H(k)}{k^{10}(k+1)} \sumend &= \frac{-1}{4}\left(  4\zeta(2) -8\zeta(3) +5\zeta(4) -12\zeta(5) +4\zeta(2)\zeta(3) +7\zeta(6)  \right. \nonumber \\ &\left. \hspace{1em}
-2\zeta(3)^2 -16\zeta(7) +4\zeta(2)\zeta(5) +4\zeta(3)\zeta(4) +9\zeta(8) -4\zeta(3)\zeta(5)  \right. \nonumber \\ &\left. \hspace{1em}
-20\zeta(9) +4\zeta(3)\zeta(6) +4\zeta(4)\zeta(5) +4\zeta(2)\zeta(7) +11\zeta(10) -4\zeta(3)\zeta(7)  \right. \nonumber \\ &\left. \hspace{1em}
-2\zeta(5)^2 -24\zeta(11) +4\zeta(2)\zeta(9) +4\zeta(3)\zeta(8) +4\zeta(4)\zeta(7)  \right. \nonumber \\ &\left. \hspace{1em}
+4\zeta(5)\zeta(6)\right) \label{eq12002} \\
\sum_{k=1}^\infty \frac{H(k)}{k^{9}(k+1)^{2}} \sumend &= \frac{-1}{4}\left( -36\zeta(2) +68\zeta(3) -35\zeta(4) +72\zeta(5) -24\zeta(2)\zeta(3)  \right. \nonumber \\ &\left. \hspace{1em}
-35\zeta(6) +10\zeta(3)^2 +64\zeta(7) -16\zeta(2)\zeta(5) -16\zeta(3)\zeta(4) -27\zeta(8)  \right. \nonumber \\ &\left. \hspace{1em}
+12\zeta(3)\zeta(5) +40\zeta(9) -8\zeta(3)\zeta(6) -8\zeta(4)\zeta(5) -8\zeta(2)\zeta(7) -11\zeta(10)  \right. \nonumber \\ &\left. \hspace{1em}
+4\zeta(3)\zeta(7) +2\zeta(5)^2\right) \label{eq12003} \\
\sum_{k=1}^\infty \frac{H(k)}{k^{8}(k+1)^{3}} \sumend &= \frac{1}{4}\left( -144\zeta(2) +256\zeta(3) -104\zeta(4) +180\zeta(5) -60\zeta(2)\zeta(3)  \right. \nonumber \\ &\left. \hspace{1em}
-70\zeta(6) +20\zeta(3)^2 +96\zeta(7) -24\zeta(2)\zeta(5) -24\zeta(3)\zeta(4) -27\zeta(8)  \right. \nonumber \\ &\left. \hspace{1em}
+12\zeta(3)\zeta(5) +20\zeta(9) -4\zeta(3)\zeta(6) -4\zeta(4)\zeta(5) -4\zeta(2)\zeta(7)\right) \label{eq12004} \\
\sum_{k=1}^\infty \frac{H(k)}{k^{7}(k+1)^{4}} \sumend &= \frac{-1}{4}\left( -336\zeta(2) +560\zeta(3) -168\zeta(4) +248\zeta(5) -84\zeta(2)\zeta(3)  \right. \nonumber \\ &\left. \hspace{1em}
-70\zeta(6) +20\zeta(3)^2 +64\zeta(7) -16\zeta(2)\zeta(5) -16\zeta(3)\zeta(4) -9\zeta(8)  \right. \nonumber \\ &\left. \hspace{1em}
+4\zeta(3)\zeta(5)\right) \label{eq12005} \\
\sum_{k=1}^\infty \frac{H(k)}{k^{6}(k+1)^{5}} \sumend &= \frac{1}{2}\left( -252\zeta(2) +392\zeta(3) -77\zeta(4) +114\zeta(5) -42\zeta(2)\zeta(3)  \right. \nonumber \\ &\left. \hspace{1em}
-16\zeta(6) +4\zeta(3)^2 +8\zeta(7) -2\zeta(2)\zeta(5) -2\zeta(3)\zeta(4)\right) \label{eq12006} \\
\sum_{k=1}^\infty \frac{H(k)}{k^{5}(k+1)^{6}} \sumend &= \frac{-1}{2}\left( -252\zeta(2) +364\zeta(3) -35\zeta(4) +96\zeta(5) -42\zeta(2)\zeta(3)  \right. \nonumber \\ &\left. \hspace{1em}
+4\zeta(6) -4\zeta(3)^2 +6\zeta(7) -2\zeta(2)\zeta(5) -2\zeta(3)\zeta(4)\right) \label{eq12007} \\
\sum_{k=1}^\infty \frac{H(k)}{k^{4}(k+1)^{7}} \sumend &= \frac{-1}{4}\left(  336\zeta(2) -448\zeta(3) -172\zeta(5) +84\zeta(2)\zeta(3) -30\zeta(6)  \right. \nonumber \\ &\left. \hspace{1em}
+20\zeta(3)^2 -48\zeta(7) +16\zeta(2)\zeta(5) +16\zeta(3)\zeta(4) -5\zeta(8) +4\zeta(3)\zeta(5)\right) \label{eq12008} \\
\sum_{k=1}^\infty \frac{H(k)}{k^{3}(k+1)^{8}} \sumend &= \frac{-1}{4}\left( -144\zeta(2) +176\zeta(3) +16\zeta(4) +120\zeta(5) -60\zeta(2)\zeta(3)  \right. \nonumber \\ &\left. \hspace{1em}
+30\zeta(6) -20\zeta(3)^2 +72\zeta(7) -24\zeta(2)\zeta(5) -24\zeta(3)\zeta(4) +15\zeta(8)  \right. \nonumber \\ &\left. \hspace{1em}
-12\zeta(3)\zeta(5) +16\zeta(9) -4\zeta(3)\zeta(6) -4\zeta(4)\zeta(5) -4\zeta(2)\zeta(7)\right) \label{eq12009}
\end{align}
 
\begin{align}
\sum_{k=1}^\infty \frac{H(k)}{k^{2}(k+1)^{9}} \sumend &= \frac{-1}{4}\left(  36\zeta(2) -40\zeta(3) -7\zeta(4) -48\zeta(5) +24\zeta(2)\zeta(3)  \right. \nonumber \\ &\left. \hspace{1em}
-15\zeta(6) +10\zeta(3)^2 -48\zeta(7) +16\zeta(2)\zeta(5) +16\zeta(3)\zeta(4) -15\zeta(8)  \right. \nonumber \\ &\left. \hspace{1em}
+12\zeta(3)\zeta(5) -32\zeta(9) +8\zeta(3)\zeta(6) +8\zeta(4)\zeta(5) +8\zeta(2)\zeta(7) -7\zeta(10)  \right. \nonumber \\ &\left. \hspace{1em}
+4\zeta(3)\zeta(7) +2\zeta(5)^2\right) \label{eq12010} \\
\sum_{k=1}^\infty \frac{H(k)}{k(k+1)^{10}} \sumend &= \frac{1}{4}\left(  4\zeta(2) -4\zeta(3) -\zeta(4) -8\zeta(5) +4\zeta(2)\zeta(3) -3\zeta(6)  \right. \nonumber \\ &\left. \hspace{1em}
+2\zeta(3)^2 -12\zeta(7) +4\zeta(2)\zeta(5) +4\zeta(3)\zeta(4) -5\zeta(8) +4\zeta(3)\zeta(5)  \right. \nonumber \\ &\left. \hspace{1em}
-16\zeta(9) +4\zeta(3)\zeta(6) +4\zeta(4)\zeta(5) +4\zeta(2)\zeta(7) -7\zeta(10) +4\zeta(3)\zeta(7)  \right. \nonumber \\ &\left. \hspace{1em}
+2\zeta(5)^2 -20\zeta(11) +4\zeta(2)\zeta(9) +4\zeta(3)\zeta(8) +4\zeta(4)\zeta(7)  \right. \nonumber \\ &\left. \hspace{1em}
+4\zeta(5)\zeta(6)\right) \label{eq12011} \\
\sum_{k=1}^\infty \frac{H(k)}{(k+1)^{11}} \sumend &= \frac{-1}{4}\left( -9\zeta(12) +4\zeta(3)\zeta(9) +4\zeta(5)\zeta(7)\right) \label{eq12012} \\
\sum_{k=1}^\infty \frac{H(k)^{2}}{k^{10}} \sumend &= -\left( - M(2,10)\right) \label{eq12013} \\
\sum_{k=1}^\infty \frac{H(k)^{2}}{k^{9}(k+1)} \sumend &= \frac{1}{24}\left(  72\zeta(3) -102\zeta(4) +84\zeta(5) -24\zeta(2)\zeta(3) -97\zeta(6)  \right. \nonumber \\ &\left. \hspace{1em}
+48\zeta(3)^2 +144\zeta(7) -24\zeta(2)\zeta(5) -60\zeta(3)\zeta(4) -24 M(2,6) +220\zeta(9)  \right. \nonumber \\ &\left. \hspace{1em}
-84\zeta(3)\zeta(6) -60\zeta(4)\zeta(5) -24\zeta(2)\zeta(7) +8\zeta(3)^3 -24 M(2,8) +312\zeta(11)  \right. \nonumber \\ &\left. \hspace{1em}
-24\zeta(2)\zeta(9) -108\zeta(3)\zeta(8) -60\zeta(4)\zeta(7) -84\zeta(5)\zeta(6)  \right. \nonumber \\ &\left. \hspace{1em}
+24\zeta(3)^2\zeta(5)\right) \label{eq12014} \\
\sum_{k=1}^\infty \frac{H(k)^{2}}{k^{8}(k+1)^{2}} \sumend &= \frac{1}{24}\left( -576\zeta(3) +780\zeta(4) -504\zeta(5) +144\zeta(2)\zeta(3) +485\zeta(6)  \right. \nonumber \\ &\left. \hspace{1em}
-240\zeta(3)^2 -576\zeta(7) +96\zeta(2)\zeta(5) +240\zeta(3)\zeta(4) +72 M(2,6) -440\zeta(9)  \right. \nonumber \\ &\left. \hspace{1em}
+168\zeta(3)\zeta(6) +120\zeta(4)\zeta(5) +48\zeta(2)\zeta(7) -16\zeta(3)^3 +24 M(2,8)\right) \label{eq12015} \\
\sum_{k=1}^\infty \frac{H(k)^{2}}{k^{7}(k+1)^{3}} \sumend &= \frac{-1}{12}\left( -1008\zeta(3) +1302\zeta(4) -648\zeta(5) +192\zeta(2)\zeta(3)  \right. \nonumber \\ &\left. \hspace{1em}
+485\zeta(6) -240\zeta(3)^2 -432\zeta(7) +72\zeta(2)\zeta(5) +180\zeta(3)\zeta(4) +36 M(2,6)  \right. \nonumber \\ &\left. \hspace{1em}
-110\zeta(9) +42\zeta(3)\zeta(6) +30\zeta(4)\zeta(5) +12\zeta(2)\zeta(7) -4\zeta(3)^3\right) \label{eq12016} \\
\sum_{k=1}^\infty \frac{H(k)^{2}}{k^{6}(k+1)^{4}} \sumend &= \frac{1}{24}\left( -4032\zeta(3) +4956\zeta(4) -1896\zeta(5) +624\zeta(2)\zeta(3)  \right. \nonumber \\ &\left. \hspace{1em}
+1007\zeta(6) -504\zeta(3)^2 -576\zeta(7) +96\zeta(2)\zeta(5) +240\zeta(3)\zeta(4) +24 M(2,6)\right) \label{eq12017}
\end{align}
 
\begin{align}
\sum_{k=1}^\infty \frac{H(k)^{2}}{k^{5}(k+1)^{5}} \sumend &= \frac{1}{12}\left(  2520\zeta(3) -2940\zeta(4) +900\zeta(5) -360\zeta(2)\zeta(3)  \right. \nonumber \\ &\left. \hspace{1em}
-335\zeta(6) +180\zeta(3)^2 +84\zeta(7) -24\zeta(2)\zeta(5) -24\zeta(3)\zeta(4)\right) \label{eq12018} \\
\sum_{k=1}^\infty \frac{H(k)^{2}}{k^{4}(k+1)^{6}} \sumend &= \frac{1}{24}\left( -4032\zeta(3) +4452\zeta(4) -1224\zeta(5) +624\zeta(2)\zeta(3)  \right. \nonumber \\ &\left. \hspace{1em}
+467\zeta(6) -288\zeta(3)^2 -96\zeta(7) +96\zeta(2)\zeta(5) -48\zeta(3)\zeta(4) -84\zeta(8)  \right. \nonumber \\ &\left. \hspace{1em}
+48\zeta(3)\zeta(5) +24 M(2,6)\right) \label{eq12019} \\
\sum_{k=1}^\infty \frac{H(k)^{2}}{k^{3}(k+1)^{7}} \sumend &= \frac{-1}{12}\left( -1008\zeta(3) +1050\zeta(4) -312\zeta(5) +192\zeta(2)\zeta(3)  \right. \nonumber \\ &\left. \hspace{1em}
+185\zeta(6) -120\zeta(3)^2 -72\zeta(7) +72\zeta(2)\zeta(5) -36\zeta(3)\zeta(4) -126\zeta(8)  \right. \nonumber \\ &\left. \hspace{1em}
+72\zeta(3)\zeta(5) +36 M(2,6) +2\zeta(9) -18\zeta(3)\zeta(6) -6\zeta(4)\zeta(5) +12\zeta(2)\zeta(7)  \right. \nonumber \\ &\left. \hspace{1em}
+4\zeta(3)^3\right) \label{eq12020} \\
\sum_{k=1}^\infty \frac{H(k)^{2}}{k^{2}(k+1)^{8}} \sumend &= \frac{1}{24}\left( -576\zeta(3) +564\zeta(4) -216\zeta(5) +144\zeta(2)\zeta(3) +185\zeta(6)  \right. \nonumber \\ &\left. \hspace{1em}
-120\zeta(3)^2 -96\zeta(7) +96\zeta(2)\zeta(5) -48\zeta(3)\zeta(4) -252\zeta(8) +144\zeta(3)\zeta(5)  \right. \nonumber \\ &\left. \hspace{1em}
+72 M(2,6) +8\zeta(9) -72\zeta(3)\zeta(6) -24\zeta(4)\zeta(5) +48\zeta(2)\zeta(7) +16\zeta(3)^3  \right. \nonumber \\ &\left. \hspace{1em}
-108\zeta(10) +48\zeta(3)\zeta(7) +24\zeta(5)^2 +24 M(2,8)\right) \label{eq12021} \\
\sum_{k=1}^\infty \frac{H(k)^{2}}{k(k+1)^{9}} \sumend &= \frac{1}{24}\left(  72\zeta(3) -66\zeta(4) +36\zeta(5) -24\zeta(2)\zeta(3) -37\zeta(6)  \right. \nonumber \\ &\left. \hspace{1em}
+24\zeta(3)^2 +24\zeta(7) -24\zeta(2)\zeta(5) +12\zeta(3)\zeta(4) +84\zeta(8) -48\zeta(3)\zeta(5)  \right. \nonumber \\ &\left. \hspace{1em}
-24 M(2,6) -4\zeta(9) +36\zeta(3)\zeta(6) +12\zeta(4)\zeta(5) -24\zeta(2)\zeta(7) -8\zeta(3)^3  \right. \nonumber \\ &\left. \hspace{1em}
+108\zeta(10) -48\zeta(3)\zeta(7) -24\zeta(5)^2 -24 M(2,8) -48\zeta(11) -24\zeta(2)\zeta(9)  \right. \nonumber \\ &\left. \hspace{1em}
+60\zeta(3)\zeta(8) +12\zeta(4)\zeta(7) +36\zeta(5)\zeta(6) -24\zeta(3)^2\zeta(5)\right) \label{eq12022} \\
\sum_{k=1}^\infty \frac{H(k)^{2}}{(k+1)^{10}} \sumend &= \frac{-1}{2}\left(  11\zeta(12) -4\zeta(3)\zeta(9) -4\zeta(5)\zeta(7) -2 M(2,10)\right) \label{eq12023} \\
\sum_{k=1}^\infty \frac{H(k)^{3}}{k^{9}} \sumend &= \frac{-1}{22112}\left(  355355\zeta(12) -221120\zeta(3)\zeta(9) -265344\zeta(5)\zeta(7)  \right. \nonumber \\ &\left. \hspace{1em}
-33168\zeta(3)^2\zeta(6) +5528\zeta(3)^4 +49752\zeta(2)\zeta(5)^2 +99504\zeta(2)\zeta(3)\zeta(7)  \right. \nonumber \\ &\left. \hspace{1em}
-82920 M(2,10)\right) \label{eq12024} \\
\sum_{k=1}^\infty \frac{H(k)^{3}}{k^{8}(k+1)} \sumend &= \frac{-1}{480}\left(  4800\zeta(4) -4800\zeta(5) -480\zeta(2)\zeta(3) +2790\zeta(6)  \right. \nonumber \\ &\left. \hspace{1em}
-1200\zeta(3)^2 -6930\zeta(7) -960\zeta(2)\zeta(5) +6120\zeta(3)\zeta(4) -2975\zeta(8)  \right. \nonumber \\ &\left. \hspace{1em}
-600\zeta(2)\zeta(3)^2 +2880\zeta(3)\zeta(5) +1320 M(2,6) -10420\zeta(9) +5820\zeta(3)\zeta(6)  \right. \nonumber \\ &\left. \hspace{1em}
+6120\zeta(4)\zeta(5) -1440\zeta(2)\zeta(7) -960\zeta(3)^3 -4983\zeta(10) +3840\zeta(3)\zeta(7)  \right. \nonumber \\ &\left. \hspace{1em}
+240\zeta(3)^2\zeta(4) -1680\zeta(2)\zeta(3)\zeta(5) +2160\zeta(5)^2 +1560 M(2,8) -480 M(3,8)\right) \label{eq12025}
\end{align}
 
\begin{align}
\sum_{k=1}^\infty \frac{H(k)^{3}}{k^{7}(k+1)^{2}} \sumend &= \frac{1}{480}\left(  33600\zeta(4) -32400\zeta(5) -3360\zeta(2)\zeta(3) +13950\zeta(6)  \right. \nonumber \\ &\left. \hspace{1em}
-6000\zeta(3)^2 -27720\zeta(7) -3840\zeta(2)\zeta(5) +24480\zeta(3)\zeta(4) -8925\zeta(8)  \right. \nonumber \\ &\left. \hspace{1em}
-1800\zeta(2)\zeta(3)^2 +8640\zeta(3)\zeta(5) +3960 M(2,6) -20840\zeta(9) +11640\zeta(3)\zeta(6)  \right. \nonumber \\ &\left. \hspace{1em}
+12240\zeta(4)\zeta(5) -2880\zeta(2)\zeta(7) -1920\zeta(3)^3 -4983\zeta(10) +3840\zeta(3)\zeta(7)  \right. \nonumber \\ &\left. \hspace{1em}
+240\zeta(3)^2\zeta(4) -1680\zeta(2)\zeta(3)\zeta(5) +2160\zeta(5)^2 +1560 M(2,8)\right) \label{eq12026} \\
\sum_{k=1}^\infty \frac{H(k)^{3}}{k^{6}(k+1)^{3}} \sumend &= \frac{1}{96}\left( -20160\zeta(4) +18720\zeta(5) +2016\zeta(2)\zeta(3) -5778\zeta(6)  \right. \nonumber \\ &\left. \hspace{1em}
+2592\zeta(3)^2 +8316\zeta(7) +1152\zeta(2)\zeta(5) -7344\zeta(3)\zeta(4) +1785\zeta(8)  \right. \nonumber \\ &\left. \hspace{1em}
+360\zeta(2)\zeta(3)^2 -1728\zeta(3)\zeta(5) -792 M(2,6) +2084\zeta(9) -1164\zeta(3)\zeta(6)  \right. \nonumber \\ &\left. \hspace{1em}
-1224\zeta(4)\zeta(5) +288\zeta(2)\zeta(7) +192\zeta(3)^3\right) \label{eq12027} \\
\sum_{k=1}^\infty \frac{H(k)^{3}}{k^{5}(k+1)^{4}} \sumend &= \frac{-1}{96}\left( -33600\zeta(4) +30000\zeta(5) +3360\zeta(2)\zeta(3) -6570\zeta(6)  \right. \nonumber \\ &\left. \hspace{1em}
+3360\zeta(3)^2 +6258\zeta(7) +960\zeta(2)\zeta(5) -5688\zeta(3)\zeta(4) +595\zeta(8)  \right. \nonumber \\ &\left. \hspace{1em}
+120\zeta(2)\zeta(3)^2 -576\zeta(3)\zeta(5) -264 M(2,6)\right) \label{eq12028} \\
\sum_{k=1}^\infty \frac{H(k)^{3}}{k^{4}(k+1)^{5}} \sumend &= \frac{1}{96}\left( -33600\zeta(4) +28800\zeta(5) +3360\zeta(2)\zeta(3) -4770\zeta(6)  \right. \nonumber \\ &\left. \hspace{1em}
+3120\zeta(3)^2 +4242\zeta(7) +960\zeta(2)\zeta(5) -4392\zeta(3)\zeta(4) -43\zeta(8)  \right. \nonumber \\ &\left. \hspace{1em}
-120\zeta(2)\zeta(3)^2 +288\zeta(3)\zeta(5) -24 M(2,6)\right) \label{eq12029} \\
\sum_{k=1}^\infty \frac{H(k)^{3}}{k^{3}(k+1)^{6}} \sumend &= \frac{1}{96}\left(  20160\zeta(4) -16560\zeta(5) -2016\zeta(2)\zeta(3) +2538\zeta(6)  \right. \nonumber \\ &\left. \hspace{1em}
-2160\zeta(3)^2 -4284\zeta(7) -1152\zeta(2)\zeta(5) +4752\zeta(3)\zeta(4) +129\zeta(8)  \right. \nonumber \\ &\left. \hspace{1em}
+360\zeta(2)\zeta(3)^2 -864\zeta(3)\zeta(5) +72 M(2,6) -788\zeta(9) +444\zeta(3)\zeta(6)  \right. \nonumber \\ &\left. \hspace{1em}
+792\zeta(4)\zeta(5) -288\zeta(2)\zeta(7) -96\zeta(3)^3\right) \label{eq12030} \\
\sum_{k=1}^\infty \frac{H(k)^{3}}{k^{2}(k+1)^{7}} \sumend &= \frac{1}{480}\left( -33600\zeta(4) +26400\zeta(5) +3360\zeta(2)\zeta(3) -4950\zeta(6)  \right. \nonumber \\ &\left. \hspace{1em}
+4800\zeta(3)^2 +14280\zeta(7) +3840\zeta(2)\zeta(5) -15840\zeta(3)\zeta(4) -645\zeta(8)  \right. \nonumber \\ &\left. \hspace{1em}
-1800\zeta(2)\zeta(3)^2 +4320\zeta(3)\zeta(5) -360 M(2,6) +7880\zeta(9) -4440\zeta(3)\zeta(6)  \right. \nonumber \\ &\left. \hspace{1em}
-7920\zeta(4)\zeta(5) +2880\zeta(2)\zeta(7) +960\zeta(3)^3 -1503\zeta(10) +2400\zeta(3)\zeta(7)  \right. \nonumber \\ &\left. \hspace{1em}
+240\zeta(3)^2\zeta(4) -1680\zeta(2)\zeta(3)\zeta(5) +1440\zeta(5)^2 +120 M(2,8)\right) \label{eq12031} \\
\sum_{k=1}^\infty \frac{H(k)^{3}}{k(k+1)^{8}} \sumend &= \frac{1}{480}\left(  4800\zeta(4) -3600\zeta(5) -480\zeta(2)\zeta(3) +990\zeta(6) -960\zeta(3)^2  \right. \nonumber \\ &\left. \hspace{1em}
-3570\zeta(7) -960\zeta(2)\zeta(5) +3960\zeta(3)\zeta(4) +215\zeta(8) +600\zeta(2)\zeta(3)^2  \right. \nonumber \\ &\left. \hspace{1em}
-1440\zeta(3)\zeta(5) +120 M(2,6) -3940\zeta(9) +2220\zeta(3)\zeta(6) +3960\zeta(4)\zeta(5)  \right. \nonumber \\ &\left. \hspace{1em}
-1440\zeta(2)\zeta(7) -480\zeta(3)^3 +1503\zeta(10) -2400\zeta(3)\zeta(7) -240\zeta(3)^2\zeta(4)  \right. \nonumber \\ &\left. \hspace{1em}
+1680\zeta(2)\zeta(3)\zeta(5) -1440\zeta(5)^2 -120 M(2,8) +10560\zeta(11) -5040\zeta(3)\zeta(8)  \right. \nonumber \\ &\left. \hspace{1em}
-2160\zeta(4)\zeta(7) -3600\zeta(5)\zeta(6) +1440\zeta(3)^2\zeta(5) -480 M(3,8)\right) \label{eq12032}
\end{align}
 
\begin{align}
\sum_{k=1}^\infty \frac{H(k)^{3}}{(k+1)^{9}} \sumend &= \frac{1}{22112}\left( -161875\zeta(12) +154784\zeta(3)\zeta(9) +199008\zeta(5)\zeta(7)  \right. \nonumber \\ &\left. \hspace{1em}
+33168\zeta(3)^2\zeta(6) -5528\zeta(3)^4 -49752\zeta(2)\zeta(5)^2 -99504\zeta(2)\zeta(3)\zeta(7)  \right. \nonumber \\ &\left. \hspace{1em}
+16584 M(2,10)\right) \label{eq12033} \\
\sum_{k=1}^\infty \frac{H(k)^{4}}{k^{8}} \sumend &= -\left( - M(4,8)\right) \label{eq12034} \\
\sum_{k=1}^\infty \frac{H(k)^{4}}{k^{7}(k+1)} \sumend &= \frac{1}{5760}\left(  172800\zeta(5) +34560\zeta(2)\zeta(3) -234960\zeta(6) -17280\zeta(3)^2  \right. \nonumber \\ &\left. \hspace{1em}
+133200\zeta(7) +28800\zeta(2)\zeta(5) -123840\zeta(3)\zeta(4) +593320\zeta(8)  \right. \nonumber \\ &\left. \hspace{1em}
+161280\zeta(2)\zeta(3)^2 -668160\zeta(3)\zeta(5) -149760 M(2,6) +209280\zeta(9)  \right. \nonumber \\ &\left. \hspace{1em}
-133920\zeta(3)\zeta(6) -123840\zeta(4)\zeta(5) +40320\zeta(2)\zeta(7) +19200\zeta(3)^3  \right. \nonumber \\ &\left. \hspace{1em}
+619407\zeta(10) -540000\zeta(3)\zeta(7) -9000\zeta(3)^2\zeta(4) +195120\zeta(2)\zeta(3)\zeta(5)  \right. \nonumber \\ &\left. \hspace{1em}
-212040\zeta(5)^2 -109080 M(2,8) -11520\zeta(2) M(2,6) -345240\zeta(11) -32640\zeta(2)\zeta(9)  \right. \nonumber \\ &\left. \hspace{1em}
+142800\zeta(3)\zeta(8) +145440\zeta(4)\zeta(7) +122160\zeta(5)\zeta(6) +9600\zeta(2)\zeta(3)^3  \right. \nonumber \\ &\left. \hspace{1em}
-69120\zeta(3)^2\zeta(5) +21120 M(3,8)\right) \label{eq12035} \\
\sum_{k=1}^\infty \frac{H(k)^{4}}{k^{6}(k+1)^{2}} \sumend &= \frac{-1}{1920}\left(  345600\zeta(5) +69120\zeta(2)\zeta(3) -460320\zeta(6)  \right. \nonumber \\ &\left. \hspace{1em}
-34560\zeta(3)^2 +177600\zeta(7) +38400\zeta(2)\zeta(5) -165120\zeta(3)\zeta(4) +593320\zeta(8)  \right. \nonumber \\ &\left. \hspace{1em}
+161280\zeta(2)\zeta(3)^2 -668160\zeta(3)\zeta(5) -149760 M(2,6) +139520\zeta(9) -89280\zeta(3)\zeta(6)  \right. \nonumber \\ &\left. \hspace{1em}
-82560\zeta(4)\zeta(5) +26880\zeta(2)\zeta(7) +12800\zeta(3)^3 +206469\zeta(10) -180000\zeta(3)\zeta(7)  \right. \nonumber \\ &\left. \hspace{1em}
-3000\zeta(3)^2\zeta(4) +65040\zeta(2)\zeta(3)\zeta(5) -70680\zeta(5)^2 -36360 M(2,8)  \right. \nonumber \\ &\left. \hspace{1em}
-3840\zeta(2) M(2,6)\right) \label{eq12036} \\
\sum_{k=1}^\infty \frac{H(k)^{4}}{k^{5}(k+1)^{3}} \sumend &= \frac{-1}{48}\left( -21600\zeta(5) -4320\zeta(2)\zeta(3) +28170\zeta(6) +2160\zeta(3)^2  \right. \nonumber \\ &\left. \hspace{1em}
-7314\zeta(7) -1680\zeta(2)\zeta(5) +7080\zeta(3)\zeta(4) -14833\zeta(8) -4032\zeta(2)\zeta(3)^2  \right. \nonumber \\ &\left. \hspace{1em}
+16704\zeta(3)\zeta(5) +3744 M(2,6) -1744\zeta(9) +1116\zeta(3)\zeta(6) +1032\zeta(4)\zeta(5)  \right. \nonumber \\ &\left. \hspace{1em}
-336\zeta(2)\zeta(7) -160\zeta(3)^3\right) \label{eq12037} \\
\sum_{k=1}^\infty \frac{H(k)^{4}}{k^{4}(k+1)^{4}} \sumend &= \frac{1}{18}\left( -10800\zeta(5) -2160\zeta(2)\zeta(3) +13785\zeta(6) +1080\zeta(3)^2  \right. \nonumber \\ &\left. \hspace{1em}
-2646\zeta(7) -720\zeta(2)\zeta(5) +2880\zeta(3)\zeta(4) -3406\zeta(8) -918\zeta(2)\zeta(3)^2  \right. \nonumber \\ &\left. \hspace{1em}
+3816\zeta(3)\zeta(5) +846 M(2,6)\right) \label{eq12038} \\
\sum_{k=1}^\infty \frac{H(k)^{4}}{k^{3}(k+1)^{5}} \sumend &= \frac{1}{48}\left(  21600\zeta(5) +4320\zeta(2)\zeta(3) -26970\zeta(6) -2160\zeta(3)^2  \right. \nonumber \\ &\left. \hspace{1em}
+5034\zeta(7) +1680\zeta(2)\zeta(5) -6360\zeta(3)\zeta(4) +12415\zeta(8) +3312\zeta(2)\zeta(3)^2  \right. \nonumber \\ &\left. \hspace{1em}
-13824\zeta(3)\zeta(5) -3024 M(2,6) +696\zeta(9) -396\zeta(3)\zeta(6) -888\zeta(4)\zeta(5)  \right. \nonumber \\ &\left. \hspace{1em}
+336\zeta(2)\zeta(7) +128\zeta(3)^3\right) \label{eq12039}
\end{align}
 
\begin{align}
\sum_{k=1}^\infty \frac{H(k)^{4}}{k^{2}(k+1)^{6}} \sumend &= \frac{1}{1920}\left( -345600\zeta(5) -69120\zeta(2)\zeta(3) +421920\zeta(6) +34560\zeta(3)^2  \right. \nonumber \\ &\left. \hspace{1em}
-104640\zeta(7) -38400\zeta(2)\zeta(5) +142080\zeta(3)\zeta(4) -496600\zeta(8)  \right. \nonumber \\ &\left. \hspace{1em}
-132480\zeta(2)\zeta(3)^2 +552960\zeta(3)\zeta(5) +120960 M(2,6) -55680\zeta(9) +31680\zeta(3)\zeta(6)  \right. \nonumber \\ &\left. \hspace{1em}
+71040\zeta(4)\zeta(5) -26880\zeta(2)\zeta(7) -10240\zeta(3)^3 -145941\zeta(10) +126240\zeta(3)\zeta(7)  \right. \nonumber \\ &\left. \hspace{1em}
-840\zeta(3)^2\zeta(4) -38160\zeta(2)\zeta(3)\zeta(5) +39960\zeta(5)^2 +22920 M(2,8)  \right. \nonumber \\ &\left. \hspace{1em}
+3840\zeta(2) M(2,6)\right) \label{eq12040} \\
\sum_{k=1}^\infty \frac{H(k)^{4}}{k(k+1)^{7}} \sumend &= \frac{-1}{5760}\left( -172800\zeta(5) -34560\zeta(2)\zeta(3) +206160\zeta(6) +17280\zeta(3)^2  \right. \nonumber \\ &\left. \hspace{1em}
-78480\zeta(7) -28800\zeta(2)\zeta(5) +106560\zeta(3)\zeta(4) -496600\zeta(8) -132480\zeta(2)\zeta(3)^2  \right. \nonumber \\ &\left. \hspace{1em}
+552960\zeta(3)\zeta(5) +120960 M(2,6) -83520\zeta(9) +47520\zeta(3)\zeta(6) +106560\zeta(4)\zeta(5)  \right. \nonumber \\ &\left. \hspace{1em}
-40320\zeta(2)\zeta(7) -15360\zeta(3)^3 -437823\zeta(10) +378720\zeta(3)\zeta(7)  \right. \nonumber \\ &\left. \hspace{1em}
-2520\zeta(3)^2\zeta(4) -114480\zeta(2)\zeta(3)\zeta(5) +119880\zeta(5)^2 +68760 M(2,8)  \right. \nonumber \\ &\left. \hspace{1em}
+11520\zeta(2) M(2,6) -28440\zeta(11) -44160\zeta(2)\zeta(9) +10320\zeta(3)\zeta(8)  \right. \nonumber \\ &\left. \hspace{1em}
+82080\zeta(4)\zeta(7) +24240\zeta(5)\zeta(6) +9600\zeta(2)\zeta(3)^3 -34560\zeta(3)^2\zeta(5)  \right. \nonumber \\ &\left. \hspace{1em}
-1920 M(3,8)\right) \label{eq12041} \\
\sum_{k=1}^\infty \frac{H(k)^{4}}{(k+1)^{8}} \sumend &= \frac{-1}{5528}\left( -289019\zeta(12) +199008\zeta(3)\zeta(9) +243232\zeta(5)\zeta(7)  \right. \nonumber \\ &\left. \hspace{1em}
+33168\zeta(3)^2\zeta(6) -5528\zeta(3)^4 -49752\zeta(2)\zeta(5)^2 -99504\zeta(2)\zeta(3)\zeta(7)  \right. \nonumber \\ &\left. \hspace{1em}
+49752 M(2,10) -5528 M(4,8)\right) \label{eq12042} \\
\sum_{k=1}^\infty \frac{H(k)^{5}}{k^{7}} \sumend &= \frac{1}{265344}\left(  3612841\zeta(12) -884480\zeta(3)\zeta(9) -597024\zeta(5)\zeta(7)  \right. \nonumber \\ &\left. \hspace{1em}
+364848\zeta(3)^2\zeta(6) +221120\zeta(3)^4 +364848\zeta(2)\zeta(5)^2 +729696\zeta(2)\zeta(3)\zeta(7)  \right. \nonumber \\ &\left. \hspace{1em}
-3250464\zeta(4)\zeta(3)\zeta(5) -1028208 M(2,10) +663360 M(4,8)\right) \label{eq12043} \\
\sum_{k=1}^\infty \frac{H(k)^{5}}{k^{6}(k+1)} \sumend &= \frac{-1}{2304}\left(  411264\zeta(6) +51840\zeta(3)^2 -295344\zeta(7) -65664\zeta(2)\zeta(5)  \right. \nonumber \\ &\left. \hspace{1em}
-76032\zeta(3)\zeta(4) -542488\zeta(8) -152640\zeta(2)\zeta(3)^2 +630144\zeta(3)\zeta(5) +135360 M(2,6)  \right. \nonumber \\ &\left. \hspace{1em}
-302144\zeta(9) +469920\zeta(3)\zeta(6) -152064\zeta(4)\zeta(5) -76320\zeta(2)\zeta(7) +11520\zeta(3)^3  \right. \nonumber \\ &\left. \hspace{1em}
-579897\zeta(10) +519840\zeta(3)\zeta(7) -3240\zeta(3)^2\zeta(4) -185040\zeta(2)\zeta(3)\zeta(5)  \right. \nonumber \\ &\left. \hspace{1em}
+203832\zeta(5)^2 +98280 M(2,8) +11520\zeta(2) M(2,6) +3126684\zeta(11) +352064\zeta(2)\zeta(9)  \right. \nonumber \\ &\left. \hspace{1em}
-1186640\zeta(3)\zeta(8) -1647936\zeta(4)\zeta(7) -880320\zeta(5)\zeta(6) -84480\zeta(2)\zeta(3)^3  \right. \nonumber \\ &\left. \hspace{1em}
+564480\zeta(3)^2\zeta(5) -34560\zeta(3) M(2,6) -111360 M(3,8)\right) \label{eq12044}
\end{align}
 
\begin{align}
\sum_{k=1}^\infty \frac{H(k)^{5}}{k^{5}(k+1)^{2}} \sumend &= \frac{1}{2304}\left(  2056320\zeta(6) +259200\zeta(3)^2 -1448496\zeta(7)  \right. \nonumber \\ &\left. \hspace{1em}
-328320\zeta(2)\zeta(5) -380160\zeta(3)\zeta(4) -1627464\zeta(8) -457920\zeta(2)\zeta(3)^2  \right. \nonumber \\ &\left. \hspace{1em}
+1890432\zeta(3)\zeta(5) +406080 M(2,6) -604288\zeta(9) +939840\zeta(3)\zeta(6) -304128\zeta(4)\zeta(5)  \right. \nonumber \\ &\left. \hspace{1em}
-152640\zeta(2)\zeta(7) +23040\zeta(3)^3 -579897\zeta(10) +519840\zeta(3)\zeta(7)  \right. \nonumber \\ &\left. \hspace{1em}
-3240\zeta(3)^2\zeta(4) -185040\zeta(2)\zeta(3)\zeta(5) +203832\zeta(5)^2 +98280 M(2,8)  \right. \nonumber \\ &\left. \hspace{1em}
+11520\zeta(2) M(2,6)\right) \label{eq12045} \\
\sum_{k=1}^\infty \frac{H(k)^{5}}{k^{4}(k+1)^{3}} \sumend &= \frac{1}{144}\left( -257040\zeta(6) -32400\zeta(3)^2 +177534\zeta(7) +41040\zeta(2)\zeta(5)  \right. \nonumber \\ &\left. \hspace{1em}
+47520\zeta(3)\zeta(4) +134527\zeta(8) +37440\zeta(2)\zeta(3)^2 -154368\zeta(3)\zeta(5) -33120 M(2,6)  \right. \nonumber \\ &\left. \hspace{1em}
+18884\zeta(9) -29370\zeta(3)\zeta(6) +9504\zeta(4)\zeta(5) +4770\zeta(2)\zeta(7)  \right. \nonumber \\ &\left. \hspace{1em}
-720\zeta(3)^3\right) \label{eq12046} \\
\sum_{k=1}^\infty \frac{H(k)^{5}}{k^{3}(k+1)^{4}} \sumend &= \frac{1}{144}\left(  257040\zeta(6) +32400\zeta(3)^2 -174006\zeta(7) -41040\zeta(2)\zeta(5)  \right. \nonumber \\ &\left. \hspace{1em}
-47520\zeta(3)\zeta(4) -132337\zeta(8) -36000\zeta(2)\zeta(3)^2 +148032\zeta(3)\zeta(5) +31680 M(2,6)  \right. \nonumber \\ &\left. \hspace{1em}
-14240\zeta(9) +25770\zeta(3)\zeta(6) -9504\zeta(4)\zeta(5) -4770\zeta(2)\zeta(7)  \right. \nonumber \\ &\left. \hspace{1em}
+720\zeta(3)^3\right) \label{eq12047} \\
\sum_{k=1}^\infty \frac{H(k)^{5}}{k^{2}(k+1)^{5}} \sumend &= \frac{1}{2304}\left( -2056320\zeta(6) -259200\zeta(3)^2 +1363824\zeta(7)  \right. \nonumber \\ &\left. \hspace{1em}
+328320\zeta(2)\zeta(5) +380160\zeta(3)\zeta(4) +1574904\zeta(8) +423360\zeta(2)\zeta(3)^2  \right. \nonumber \\ &\left. \hspace{1em}
-1738368\zeta(3)\zeta(5) -371520 M(2,6) +455680\zeta(9) -824640\zeta(3)\zeta(6) +304128\zeta(4)\zeta(5)  \right. \nonumber \\ &\left. \hspace{1em}
+152640\zeta(2)\zeta(7) -23040\zeta(3)^3 +449109\zeta(10) -387360\zeta(3)\zeta(7)  \right. \nonumber \\ &\left. \hspace{1em}
-9720\zeta(3)^2\zeta(4) +124560\zeta(2)\zeta(3)\zeta(5) -122328\zeta(5)^2 -68040 M(2,8)  \right. \nonumber \\ &\left. \hspace{1em}
-11520\zeta(2) M(2,6)\right) \label{eq12048} \\
\sum_{k=1}^\infty \frac{H(k)^{5}}{k(k+1)^{6}} \sumend &= \frac{-1}{2304}\left( -411264\zeta(6) -51840\zeta(3)^2 +267120\zeta(7) +65664\zeta(2)\zeta(5)  \right. \nonumber \\ &\left. \hspace{1em}
+76032\zeta(3)\zeta(4) +524968\zeta(8) +141120\zeta(2)\zeta(3)^2 -579456\zeta(3)\zeta(5) -123840 M(2,6)  \right. \nonumber \\ &\left. \hspace{1em}
+227840\zeta(9) -412320\zeta(3)\zeta(6) +152064\zeta(4)\zeta(5) +76320\zeta(2)\zeta(7) -11520\zeta(3)^3  \right. \nonumber \\ &\left. \hspace{1em}
+449109\zeta(10) -387360\zeta(3)\zeta(7) -9720\zeta(3)^2\zeta(4) +124560\zeta(2)\zeta(3)\zeta(5)  \right. \nonumber \\ &\left. \hspace{1em}
-122328\zeta(5)^2 -68040 M(2,8) -11520\zeta(2) M(2,6) -2668908\zeta(11) -275264\zeta(2)\zeta(9)  \right. \nonumber \\ &\left. \hspace{1em}
+993200\zeta(3)\zeta(8) +1403136\zeta(4)\zeta(7) +705120\zeta(5)\zeta(6) +65280\zeta(2)\zeta(3)^3  \right. \nonumber \\ &\left. \hspace{1em}
-449280\zeta(3)^2\zeta(5) +34560\zeta(3) M(2,6) +92160 M(3,8)\right) \label{eq12049}
\end{align}
 
\begin{align}
\sum_{k=1}^\infty \frac{H(k)^{6}}{k^{5}(k+1)} \sumend &= \frac{1}{384}\left(  247296\zeta(7) +55680\zeta(2)\zeta(5) +114048\zeta(3)\zeta(4)  \right. \nonumber \\ &\left. \hspace{1em}
-280464\zeta(8) +15744\zeta(2)\zeta(3)^2 -187008\zeta(3)\zeta(5) -21888 M(2,6) +119584\zeta(9)  \right. \nonumber \\ &\left. \hspace{1em}
-209952\zeta(3)\zeta(6) +96768\zeta(4)\zeta(5) +31248\zeta(2)\zeta(7) -8704\zeta(3)^3 +814101\zeta(10)  \right. \nonumber \\ &\left. \hspace{1em}
-529680\zeta(3)\zeta(7) +253944\zeta(3)^2\zeta(4) +1200\zeta(2)\zeta(3)\zeta(5) -365064\zeta(5)^2  \right. \nonumber \\ &\left. \hspace{1em}
-103128 M(2,8) -45120\zeta(2) M(2,6) -1469286\zeta(11) -166944\zeta(2)\zeta(9) +542488\zeta(3)\zeta(8)  \right. \nonumber \\ &\left. \hspace{1em}
+790176\zeta(4)\zeta(7) +410848\zeta(5)\zeta(6) +38720\zeta(2)\zeta(3)^3 -260352\zeta(3)^2\zeta(5)  \right. \nonumber \\ &\left. \hspace{1em}
+18240\zeta(3) M(2,6) +51200 M(3,8)\right) \label{eq12050} \\
\sum_{k=1}^\infty \frac{H(k)^{6}}{k^{4}(k+1)^{2}} \sumend &= \frac{1}{384}\left( -989184\zeta(7) -222720\zeta(2)\zeta(5) -456192\zeta(3)\zeta(4)  \right. \nonumber \\ &\left. \hspace{1em}
+1113824\zeta(8) -62016\zeta(2)\zeta(3)^2 +742272\zeta(3)\zeta(5) +86592 M(2,6) -239168\zeta(9)  \right. \nonumber \\ &\left. \hspace{1em}
+419904\zeta(3)\zeta(6) -193536\zeta(4)\zeta(5) -62496\zeta(2)\zeta(7) +17408\zeta(3)^3  \right. \nonumber \\ &\left. \hspace{1em}
-814101\zeta(10) +529680\zeta(3)\zeta(7) -253944\zeta(3)^2\zeta(4) -1200\zeta(2)\zeta(3)\zeta(5)  \right. \nonumber \\ &\left. \hspace{1em}
+365064\zeta(5)^2 +103128 M(2,8) +45120\zeta(2) M(2,6)\right) \label{eq12051} \\
\sum_{k=1}^\infty \frac{H(k)^{6}}{k^{3}(k+1)^{3}} \sumend &= \frac{1}{4}\left(  15456\zeta(7) +3480\zeta(2)\zeta(5) +7128\zeta(3)\zeta(4) -17278\zeta(8)  \right. \nonumber \\ &\left. \hspace{1em}
+954\zeta(2)\zeta(3)^2 -11508\zeta(3)\zeta(5) -1338 M(2,6) +2270\zeta(9) -4284\zeta(3)\zeta(6)  \right. \nonumber \\ &\left. \hspace{1em}
+1980\zeta(4)\zeta(5) +651\zeta(2)\zeta(7) -180\zeta(3)^3\right) \label{eq12052} \\
\sum_{k=1}^\infty \frac{H(k)^{6}}{k^{2}(k+1)^{4}} \sumend &= \frac{1}{384}\left( -989184\zeta(7) -222720\zeta(2)\zeta(5) -456192\zeta(3)\zeta(4)  \right. \nonumber \\ &\left. \hspace{1em}
+1097760\zeta(8) -60096\zeta(2)\zeta(3)^2 +730752\zeta(3)\zeta(5) +84672 M(2,6) -196672\zeta(9)  \right. \nonumber \\ &\left. \hspace{1em}
+402624\zeta(3)\zeta(6) -186624\zeta(4)\zeta(5) -62496\zeta(2)\zeta(7) +17152\zeta(3)^3  \right. \nonumber \\ &\left. \hspace{1em}
-779835\zeta(10) +490704\zeta(3)\zeta(7) -245544\zeta(3)^2\zeta(4) +15600\zeta(2)\zeta(3)\zeta(5)  \right. \nonumber \\ &\left. \hspace{1em}
+339864\zeta(5)^2 +94728 M(2,8) +45120\zeta(2) M(2,6)\right) \label{eq12053} \\
\sum_{k=1}^\infty \frac{H(k)^{6}}{k(k+1)^{5}} \sumend &= \frac{-1}{384}\left( -247296\zeta(7) -55680\zeta(2)\zeta(5) -114048\zeta(3)\zeta(4)  \right. \nonumber \\ &\left. \hspace{1em}
+272432\zeta(8) -14784\zeta(2)\zeta(3)^2 +181248\zeta(3)\zeta(5) +20928 M(2,6) -98336\zeta(9)  \right. \nonumber \\ &\left. \hspace{1em}
+201312\zeta(3)\zeta(6) -93312\zeta(4)\zeta(5) -31248\zeta(2)\zeta(7) +8576\zeta(3)^3 -779835\zeta(10)  \right. \nonumber \\ &\left. \hspace{1em}
+490704\zeta(3)\zeta(7) -245544\zeta(3)^2\zeta(4) +15600\zeta(2)\zeta(3)\zeta(5) +339864\zeta(5)^2  \right. \nonumber \\ &\left. \hspace{1em}
+94728 M(2,8) +45120\zeta(2) M(2,6) +1373598\zeta(11) +149024\zeta(2)\zeta(9) -524968\zeta(3)\zeta(8)  \right. \nonumber \\ &\left. \hspace{1em}
-724416\zeta(4)\zeta(7) -365168\zeta(5)\zeta(6) -36160\zeta(2)\zeta(3)^3 +240768\zeta(3)^2\zeta(5)  \right. \nonumber \\ &\left. \hspace{1em}
-16320\zeta(3) M(2,6) -46720 M(3,8)\right) \label{eq12054}
\end{align}

\end{document}